\documentclass[10pt]{amsart}
\usepackage{latexsym}
\usepackage{amsmath}
\usepackage{amssymb}
\usepackage{amscd}
\usepackage{array}
\usepackage{hhline}
\usepackage{color}
\usepackage[all]{xy}
\usepackage{graphicx}
\theoremstyle{break}
\newtheorem{definition}{Definition}[section]
\newtheorem{theorem}[definition]{Theorem}
\newtheorem{proposition}[definition]{Proposition}
\newtheorem{remark}[definition]{Remark}
\newtheorem{lemma}[definition]{Lemma}
\newtheorem{corollary}[definition]{Corollary}

\newtheorem{problem}[definition]{Problem}

\def\C{{\mathbb{C}}}
\def\Q{{\mathbb{Q}}}
\def\N{{\mathbb{N}}}
\def\Z{{\mathbb{Z}}}
\def\c{{\mathcal{C}}}
\def\deg{{\mathrm{deg}}}
\def\dim{{\mathrm{dim}}}

\def\gdim{{\mathrm{gdim}}}
\def\mod{{\mathrm{mod\;}}}

\def\End{{\mathrm{End}}}
\def\END{{\mathrm{END}}}
\def\Hom{{\mathrm{Hom}}}
\def\HOM{{\mathrm{HOM}}}

\def\EXT{{\mathrm{EXT}}}
\def\ob{{\mathrm{Ob}}}
\def\mor{{\mathrm{Mor}}}

\def\H{{\mathrm{H}}}
\def\cyc{{\mathrm{Z}}}
\def\kom{{\mathrm{Kom}}}
\def\d{{\mathrm{d}}}
\def\k{{\mathcal{K}}}
\def\id{{\mathrm{Id}}}
\def\MF{{\mathrm{MF}}}
\def\HMF{{\mathrm{HMF}}}

\def\mfsim{{\hspace{.1cm}\stackrel{mf}{\sim}\hspace{.1cm}}}

\def\ostimes{{\,\otimes\hspace{-0.7em}\raisebox{-0.5ex}{${}_{{}_{S}}$}\,}}

\def\oqtimes{{\,\otimes\hspace{-0.6em}\raisebox{-0.5ex}{${}_{{}_{\Q}}$}\,}}
\newcommand{\en}[1]{\txt{\begin{picture}(0,0)(0,0)\circle{1.5}\put(-1.5,0){\makebox(0,0){\emph{$#1$}}}\end{picture}}}

\title[Quantum ($\mathfrak{sl}_n$, $\land V_n$) link invariant and matrix factorizations]{Quantum ($\mathfrak{sl}_n$,$\land V_n$) link invariant and matrix factorizations}
\author{Yasuyoshi Yonezawa}
\address{Graduate School of Mathematics, Nagoya University\\ 464-8602 Furocho, Chikusaku, Nagoya, Japan }
\email{yasuyoshi.yonezawa@math.nagoya-u.ac.jp}
\begin{document}
\maketitle
\begin{abstract}
M. Khovanov and L. Rozansky gave a categorification of the HOMFLY-PT polynomial.
This study is a generalization of the Khovanov-Rozansky homology.
We define a homology associated to the quantum $(\mathfrak{sl}_n,\land V_n)$ link invariant, where $\land V_n$ is the set of the fundamental representations of the quantum group of $\mathfrak{sl}_n$.
In the case of an oriented link diagram with $[1,k]$-colorings and $[k,1]$-colorings, we prove that its homology has invariance under colored Reidemeister moves composed of $[1,k]$-crossings and $[k,1]$-crossings.
In the case of an $[i,j]$-colored oriented link diagram, we define a normalized Poincar\'e polynomial of its homology and prove the polynomial is a link invariant.
\end{abstract}
\tableofcontents
%
%
%
%
\section{Introduction}\label{1-cat}
\begin{figure}[hbt]
\input{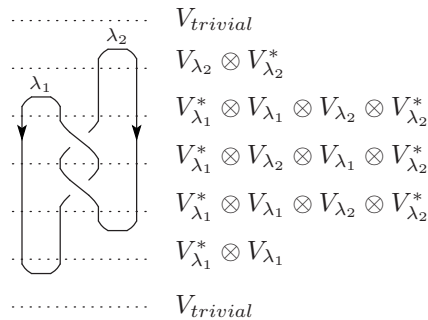}
\caption{Sliced link diagram}\label{slice-diag}
\end{figure}

\indent
M. Khovanov constructed a homological link invariant whose Euler characteristic is the Jones polynomial via a category of complexes of $\Z$-graded modules \cite{K1}.
In general, constructing a homological link invariant whose Euler characteristic is a link invariant by using objects of a category is called a categorification of the link invariant.\\
\indent
We understand the Jones polynomial to be the simplest quantum link invariant, which is obtained from the quantum group $U_q(\mathfrak{sl}_2)$ at a generic $q$ and its vector representation $V_2$.
N. Reshetikhin and V. Turaev generally constructed a link invariant associated with the quantum group $U_q(\mathfrak{g})$, where $\mathfrak{g}$ is a simple Lie algebra, and its representations, called the quantum $\mathfrak{g}$ link invariant or the Reshetikhin-Turaev $\mathfrak{g}$ link invariant \cite{RT}.\\
\indent
For a given oriented link diagram $D$, we obtain the quantum $\mathfrak{g}$ link invariant by the following procedures.
We fix a simple Lie algebra $\mathfrak{g}$ and assign one of the irreducible representations $V_{\lambda}$ of the quantum group $U_q(\mathfrak{g})$ to each component of the link diagram $D$, where $V_{\lambda}$ is the highest weight representation corresponding to a highest weight $\lambda$.
It does not matter that each assigned representation of components is different.
On a component of an oriented link diagram the marking $\lambda$ often represents assigning $V_{\lambda}$ to the component.
The marking $\lambda$ is called a coloring.
The horizontal line sweeps across the link diagram from the bottom to the top.
Then, we slice the link diagram every time a state of intersections of the horizontal line and the diagram changes, see Figure \ref{slice-diag}.
A partial sliced diagram in an interval between horizontal lines is called a tangle diagram.
A tangle diagram in an interval between neighboring horizontal lines can be considered as an intertwiner of two representations of $U_q(\mathfrak{g})$ since a representation associated to the link diagram exists on each intersection of the link diagram and a horizontal line.
Taking the composition of the intertwiners in all intervals, we obtain a Laurent polynomial of the variable $q$.
If we choose a suitable intertwiner for each tangle diagram, then the Laurent polynomial has the same evaluation for oriented diagrams transforming to each other under the Reidemeister moves.
The quantum $\mathfrak{g}$ link invariant is such an obtained Laurent polynomial.
When we consider particular representations $V_1$,...,$V_k$ of $U_q(\mathfrak{g})$ only, the quantum $\mathfrak{g}$ link invariant is called the quantum ($\mathfrak{g}$,$\mathbb{V}$) link invariant, where $\mathbb{V}=\{V_1,...,V_k\}$.
\\
\indent
We consider the following natural question.
\begin{problem}
\hspace{1cm}\\
Can we construct a homological link invariant whose Euler characteristic is a given quantum link invariant?
\end{problem}
%
%
\subsection{Khovanov-Rozansky homology}
\begin{figure}[hbt]
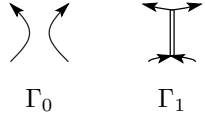

$$
\input{figure/figsmoothing1sln-p}\hspace{1cm}\input{figure/figsmoothing2sln1-p}
$$
\caption{Planar diagrams}\label{planar-single-double}
\end{figure}

\begin{figure}[hbt]
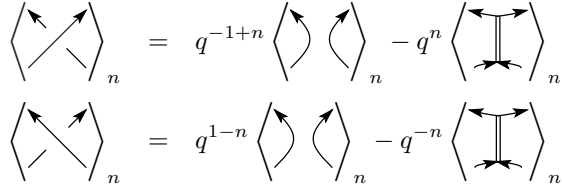

\begin{eqnarray}
\nonumber
\left<\input{figure/figplus}\right>_n&=&q^{-1+n}\left<\input{figure/figsmoothing1sln}\right>_n-q^{n}\left<\input{figure/figsmoothing2sln1}\right>_n\\[-0.1em]
\nonumber
\left<\input{figure/figminus}\right>_n&=&q^{1-n}\left<\input{figure/figsmoothing1sln}\right>_n-q^{-n}\left<\input{figure/figsmoothing2sln1}\right>_n
\end{eqnarray}
\caption{Reductions for single crossings}\label{HOMFLY-PT}
\end{figure}

\begin{figure}[htb]
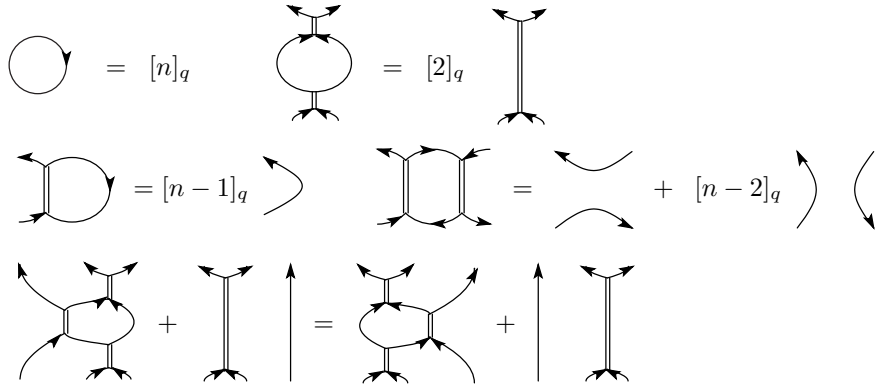

\begin{eqnarray}
\nonumber
&&\input{figure/figcircle1}\hspace{.2cm}=\hspace{.2cm}[n]_q\hspace{1cm} \input{figure/figsmoothbubble-rel}\hspace{.2cm}=\hspace{.2cm}[2]_q\hspace{.3cm}\input{figure/figsmoothbubble-res}\\[1em]
\nonumber
&&\input{figure/figsmoothbubble2}\hspace{.2cm}=[n-1]_q\hspace{.2cm}\input{figure/figsmoothing4sln1}
\hspace{1cm}\input{figure/figsmoothbubble3}\hspace{.2cm}=\hspace{.2cm}\input{figure/figsmoothing4sln}\hspace{.2cm}+\hspace{.2cm}[n-2]_q\hspace{.2cm}\input{figure/figsmoothing3sln}\\[1.5em]
\nonumber
&&\input{figure/figsmoothingreid1}
\hspace{.2cm}+\hspace{.2cm}\input{figure/figsmoothingreid2}
\hspace{.2cm}=\hspace{.2cm}\input{figure/figsmoothingreid3}
\hspace{.2cm}+\hspace{.2cm}\input{figure/figsmoothingreid4}\\
\nonumber
\end{eqnarray}
\caption{Relation of planar diagrams}\label{relation-HOMFLY-PT}
\end{figure}

\begin{figure}[hbt]
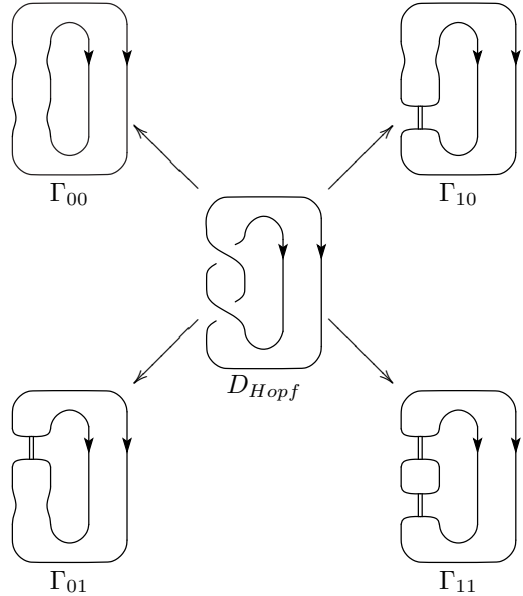

\begin{equation}
\nonumber
\xymatrix{
\input{figure/hopf-00}&&\input{figure/hopf-10}\\
&\input{figure/hopf-link}\ar[ru]\ar[rd]\ar[lu]\ar[ld]&\\
\input{figure/hopf-01}&&\input{figure/hopf-11}
}
\end{equation}\\[2em]
\caption{Planar diagrams derived from Hopf link diagram}\label{reduction-hopf}
\end{figure}

\begin{figure}[htb]
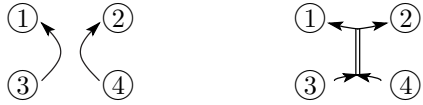

\begin{center}
\input{figure/figsmoothing1sln-mf1}
\hspace{2cm}
\input{figure/figsmoothing2sln-mf1}
\end{center}
\caption{$\Gamma_0$ and $\Gamma_1$ assigned indexes}\label{potential-KR}
\end{figure}

\begin{figure}[hbt]
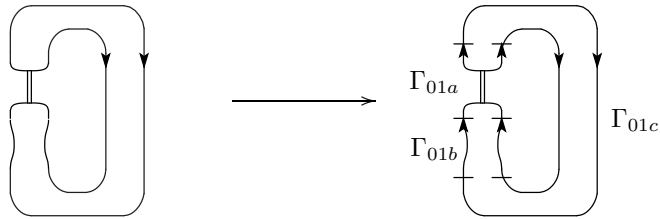

\begin{equation}
\nonumber
\xymatrix{
\input{figure/hopf-01-decomp}&\ar[rr]&&&\input{figure/hopf-01-mark}
}
\end{equation}\\[2em]
\caption{Decomposition of planar diagram $\Gamma_{10}$}\label{decomp-hopf}
\end{figure}

\indent
The HOMFLY-PT polynomial is the quantum ($\mathfrak{sl}_n$,$V_n$) link invariant, where $V_n$ is the vector representation of $U_q(\mathfrak{sl}_n)$.
In fact, M. Khovanov and L. Rozansky constructed a homological link invariant whose Euler characteristic is the HOMFLY-PT polynomial via a category of complexes of $\Z$-graded matrix factorizations, denoted by $\k^{b}(\HMF^{gr})$; see Section \ref{sec2.11}.\\
\indent
H. Murakami, T. Ohtsuki and S. Yamada gave the state model of the HOMFLY-PT polynomial by using trivalent planar diagrams, see Appendix \ref{NMOY} in the case that coloring has $1$ and $2$ only and see \cite{MOY}.
For a given oriented link diagram $D$, the HOMFLY-PT polynomial of the diagram $D$ can be calculated combinatorially only by the state model.
It is calculated by transforming to each single crossing into planar diagrams $\Gamma_0$ and $\Gamma_1$ in Figure \ref{planar-single-double} (in the case of Hopf link, see Figure \ref{reduction-hopf}), then evaluating the obtained closed planar diagrams as a Laurent polynomial using relations in Figure \ref{relation-HOMFLY-PT} and summing the Laurent polynomials by the reduction in Figure \ref{HOMFLY-PT}.
\\
\\
\\
\indent
For a $\Z$-graded polynomial ring $R$ with finite variables and a homogeneous $\Z$-graded polynomial $\omega$ called a potential, we obtain $\Z$-graded matrix factorizations with the potential $\omega$.
The matrix factorization (factorization for short ) $\overline{M}$ is a two-periodic chain composed of $\Z$-graded $R$-modules $M_0$, $M_1$ and $R$-module morphisms $d_0$, $d_1$ satisfying that $d_1d_0=\omega\id_{M_0}$ and $d_0d_1=\omega\id_{M_1}$,
\begin{equation}
\nonumber
\xymatrix{M_0\ar[rr]^{d_0}&&M_1\ar[rr]^{d_1}&&M_0}.
\end{equation}
For two factorizations $\overline{M}$ with a potential $\omega$ and $\overline{N}$ a potential $\omega'$, a tensor product $\overline{M}\boxtimes\overline{N}$ is defined; see Section \ref{sec2.8}.
Its potential is $\omega+\omega'$.\\
\indent
Khovanov and Rozansky defined a complex of matrix factorization for an oriented link diagram as follows.
First, we define a polynomial called a potential for a planar diagram.
We assign a different index $i$ to each end of the diagram, moreover, place the polynomial $x_i^{n+1}$ on the $i$-assigned end if the orientation of the edge is the direction from the diagram to the end, place $-x_i^{n+1}$ if the orientation is the opposite direction and sum these polynomials.
Remark that $n$ is associated to the quantum $(\mathfrak{sl}_n,V_n)$ link invariant.
For example, we assign indexes $1$, $2$, $3$, $4$ to $\Gamma_0$ and $\Gamma_1$ as Figure \ref{potential-KR}.
Then, the potential for $\Gamma_0$ and $\Gamma_1$ is $x_1^{n+1}+x_2^{n+1}-x_3^{n+1}-x_4^{n+1}$.
\\
\indent
For the diagrams $\Gamma_0$ and $\Gamma_1$, we define matrix factorizations $\c(\Gamma_0)_n$ and $\c(\Gamma_1)_n$ with the potential of the planar diagrams $x_1^{n+1}+x_2^{n+1}-x_3^{n+1}-x_4^{n+1}$.\\
We consider a general planar diagram formed by gluing some planar diagrams $\Gamma_0$ and $\Gamma_1$ at ends of edges with preserving an orientation.
A matrix factorization for such a planar diagram is defined by taking a tensor product of factorizations for the diagrams $\Gamma_0$, $\Gamma_1$ and making two indexes on the glued ends the same.
In particular, for a closed planar diagrams $\Gamma$, the potential of its matrix factorization $\c(\Gamma)_n$ is $0$.
Therefore, $\c(\Gamma)_n$ is a two-periodic complex of $\Q$-vector spaces.
For example, $\Gamma_{10}$ in Figure \ref{reduction-hopf} has a decomposition shown as Figure \ref{decomp-hopf}.
Then, we have a matrix factorization $\c(\Gamma_{10})_n$ as follows,
\begin{equation}
\nonumber
\c(\Gamma_{01})_n= \c(\Gamma_{01a})_n\boxtimes\c(\Gamma_{01b})_n\boxtimes\c(\Gamma_{01c})_n.
\end{equation}
\\
\indent
In the homotopy category of matrix factorizations $\HMF^{gr}$, we show that there exist isomorphisms between matrix factorizations corresponding to relations of planar diagrams in Figure \ref{relation-HOMFLY-PT}.
\\
\indent
Using the matrix factorizations, the state model of single crossing in Figure \ref{HOMFLY-PT} is represented as an object of $\k^b(\HMF^{gr})$:
\begin{eqnarray}
\nonumber
&&\xymatrix{
\c\Big( \input{figure/figplus} \Big)_n  := \Big( \ar[r]&0\ar[r] 
&{\c^{-1}\Big( \input{figure/figsmoothing2sln1-mf} \Big)_n}\ar[r]^{\chi_+}
&{\c^{0}\Big( \input{figure/figsmoothing-sln-mf} \Big)_n} \ar[r] 
&0\ar[r] & \Big), } \\[-0.5em]
\nonumber
&&
\xymatrix{\c\Big( \input{figure/figminus} \Big)_n  := \Big(\ar[r] 
&0\ar[r] 
&{\c^{0}\Big( \input{figure/figsmoothing-sln-mf} \Big)_n}\ar[r]^{\chi_-} 
&{\c^{1}\Big( \input{figure/figsmoothing2sln1-mf} \Big)_n}\ar[r]&0\ar[r] &\Big).
}
\end{eqnarray}
For an oriented link diagram $D$, we define a complex $\c(D)_n$ by exchanging every single crossing into the above complexes of matrix factorizations and taking a tensor product of the complexes for all single crossings.
The complex $\c(D)_n$ is a complex of two-periodic complexes of $\Q$-vector spaces since the matrix factorization for a closed planar diagram is a two-periodic complex of $\Q$-vector spaces.\\
\indent
For the Hopf link in Figure \ref{reduction-hopf}, we obtain the complex
\begin{equation}
\nonumber
\xymatrix{
&&0\ar@{..}[d]&&1\ar@{..}[d]&&2\ar@{..}[d]&
\\
\ar[r]&0\ar[r]
&\hspace{.5cm}{\c(\Gamma_{00})_n}\hspace{.5cm}
\ar[rr]^{\left(\begin{array}{c}{}_{\chi_-\boxtimes \id} \\ {}_{\id\boxtimes\chi_-} \end{array}\right)}
&&\hspace{.5cm}{\begin{array}{c}\c(\Gamma_{10})_n\\ \oplus\\ \c(\Gamma_{01})_n\end{array}}\hspace{.5cm}
\ar[rr]^{(\id\boxtimes\chi_-,-\chi_-\boxtimes\id)}
&&\hspace{.5cm}{\c(\Gamma_{11})_n}\ar[r]\hspace{.5cm}&0\ar[r]&
}.
\end{equation}
\indent
M. Khovanov and L. Rozansky introduced such a $\Z\oplus\Z\oplus\Z_2$-graded complex $\c(D)$ for an oriented link diagram $D$, where these gradings consist of the complex grading, the $\Z$-grading induced by a $\Z$-graded factorization, and the two-periodic grading of a factorization.
Then, they proved that the $\Z\oplus\Z\oplus\Z_2$-graded homology of the complex $\c(D)$ is a link invariant.
That is, they showed that if link diagrams $D$ and $D'$ transform to each other under the Reidemeister moves, then these $\Z\oplus\Z\oplus\Z_2$-graded homologies are isomorphic.
This homology for $D$ is called the Khovanov-Rozansky homology.
%
%
\subsection{Result of the present paper I: Matrix factorizations and colored planar diagrams}
\begin{figure}[htb]
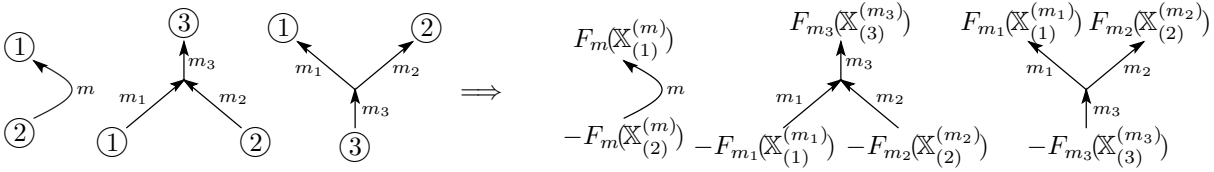

\begin{equation}
\nonumber
\input{figure/figmoy3-en} \hspace{.5cm} \input{figure/figgluing-in-3valent-en} \hspace{.5cm} \input{figure/figgluing-out-3valent-en}
\hspace{.5cm}\Longrightarrow \hspace{.5cm}
\input{figure/figmoy3-mf1} \hspace{.6cm} \input{figure/figgluing-in-3valent-mf1} \hspace{.5cm} \input{figure/figgluing-out-3valent-mf1}
\end{equation}\\[2em]
\caption{Assigned polynomials for essential planar diagrams}\label{essential-planar-function}
\end{figure}

\begin{figure}[hbt]
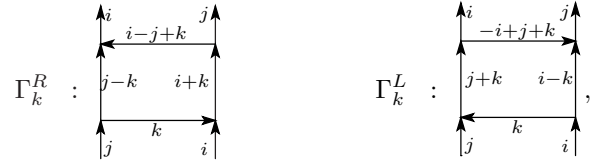

\begin{eqnarray}
\nonumber
\Gamma^R_k\hspace{.2cm}:\hspace{.2cm}\input{figure/figmoysmooth4}\hspace{1cm} \hspace{1cm} 
\Gamma^L_k\hspace{.2cm}:\hspace{.2cm}\input{figure/figmoysmooth3},\\[1em]
\nonumber
\end{eqnarray}\\[-1.5em]
\caption{Resolved planar diagrams for $[i,j]$-crossing}\label{colored-planar}
\end{figure}

\begin{figure}[htb]
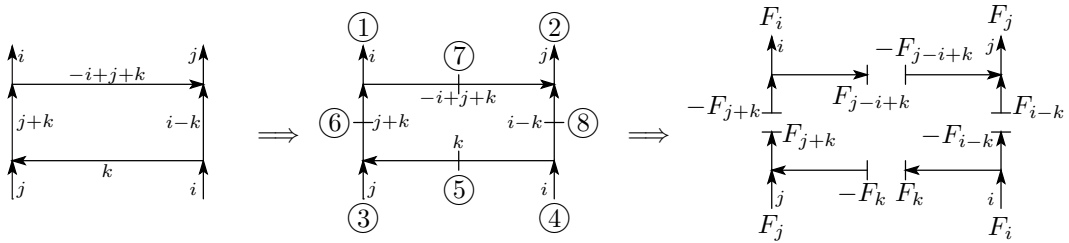

\begin{equation}
\nonumber
\input{figure/figmoysmooth-index0}\hspace{1cm}\Longrightarrow \hspace{1cm}\input{figure/figmoysmooth-index1}\hspace{1cm}\Longrightarrow \hspace{1cm}\input{figure/figmoysmooth-index2}
\end{equation}\\[2.5em]
\caption{Decomposition of $\Gamma^L_k$ and assignment of polynomial}\label{decomp-assign-planar}
\end{figure}

\indent
Khovanov and  Rozansky gave the homological link invariant whose Euler characteristic is the HOMFLY-PT polynomial by using a matrix factorization.
However, there still exist a lot of quantum link invariants which are not yet categorified.
Let $\land V_n$ be the set of the fundamental representations of $U_q(\mathfrak{sl}_n)$, that is, $\land V_n=\{V_n,\land^2 V_n, ..., \land^{n-1} V_n\}$.
This paper concerns the quantum ($\mathfrak{sl}_n$, $\land V_n$) link invariant which has not previously been categorified.
In the paper \cite{MOY}, H. Murakami, T. Ohtsuki and S. Yamada gave the state model of the quantum ($\mathfrak{sl}_n$, $\land V_n$) link invariant by using planar diagrams with colorings from $1$ to $n$
, see Appendix \ref{NMOY}.
A coloring $m$ ($1\leq m\leq n$) on an edge represents assigning $\land^{m}V_n$ to the edge.
$\land^n V_n$ is the trivial representation of $U_q(\mathfrak{sl}_n)$.
This state model is often called the MOY bracket.
This is a generalization of the state model of the HOMFLY-PT polynomial.
If we consider a single crossing with coloring $1$ and a planar diagram with colorings $1$ and $2$ only, the MOY bracket is equal to the state model of the HOMFLY-PT polynomial; see Figures \ref{HOMFLY-PT} and \ref{relation-HOMFLY-PT}.
An edge with coloring $1$ corresponds to a single edge in the state model of the HOMFLY-PT polynomial and an edge with coloring $2$ corresponds to a double edge.\\
\indent
It is a natural problem to construct a homological link invariant whose Euler characteristic is the quantum ($\mathfrak{sl}_n$, $\land V_n$) link invariant generalizing the Khovanov-Rozansky homology.
The purpose of this study is to construct such a homological link invariant by using matrix factorizations.
Unfortunately, we do not define the homological link invariant whose Euler characteristic is the quantum ($\mathfrak{sl}_n$, $\land V_n$) link invariant in the present paper.
However, we define a new polynomial link invariant which is the same with Poincar\'e polynomial of the homological link invariant, see Section \ref{sec1.3}, Section \ref{sec5} and Section \ref{sec6}.
\\
\indent
The calculation of the MOY bracket is similar to the calculation of the state model of the HOMFLY-PT polynomial.
For a colored oriented link diagram $D$, the quantum $(\mathfrak{sl}_n$, $\land V_n)$ link invariant of $D$ is obtained by transforming each $[i,j]$-colored single crossing, called $[i,j]$-crossing for short, into a linear combination of colored planar diagrams in Figure \ref{MOY-crossing}, and evaluating these closed colored planar diagrams as a Laurent polynomial by using relations in Figure \ref{MOY-relation} and then summing these Laurent polynomials.
See Appendix \ref{NMOY}.
\begin{figure}[htb]
\begin{eqnarray}
\nonumber
\left<\hspace{.3cm} \input{figure/figplus-color}\hspace{.3cm}\right>_n &=& 
\sum_{k=0}^{j} (-1)^{-k+j-i}q^{k+in-i^2+(i-j)^2+2(i-j)}
\left<\hspace{.3cm} \input{figure/figmoysmooth2}\hspace{.3cm}\right>_n \hspace{1cm} {\rm for}\,\, i \geq j\\[1em]
\nonumber
\left<\hspace{.3cm} \input{figure/figminus-color}\hspace{.3cm}\right>_n &=& 
\sum_{k=0}^{i} (-1)^{k+j-i}q^{-k-jn+j^2-(j-i)^2-2(j-i)}
\left<\hspace{.3cm} \input{figure/figmoysmooth1}\hspace{.3cm}\right>_n \hspace{1cm} {\rm for}\,\, i\leq j
\end{eqnarray}
\caption{Reduction for $[i,j]$-crossing of MOY bracket}\label{MOY-crossing}
\end{figure}
\begin{figure}[htb]
$$
\input{figure/figcircle-weight} = \left[ n \atop i \right]_q
\hspace{0.5cm}
\input{figure/fig-bubble-color}=
\left[ i_3 \atop i_1 \right]_q  \input{figure/fig-line-color}
\hspace{0.2cm}
\label{planar-relation1}
\input{figure/fig-bubble-color1}=
\left[ n - i_1 \atop i_2 \right]_q \input{figure/fig-line-color1}
$$
$$
\input{figure/fig-ass-dia1}=\input{figure/fig-ass-dia2}
\hspace{1cm}
\input{figure/fig-coass-dia1}=\input{figure/fig-coass-dia2}
$$
$$
\input{figure/figsquare1j--k--k+1j-k--k--1j}= 
\left[ i_1-1 \atop i_2 \right]_q \input{figure/figsquare1i-1} +
\left[ i_1-1 \atop i_2-1 \right]_q \input{figure/figsquare1i-1--i-1+1--1i-1}
$$
$$
\input{figure/figsquare1j--j+1--j1--j+1--1j-rev}=
\input{figure/figsquare1j-rev} +
\left[ n-j-1 \right]_q \input{figure/figsquare1j--j-1--1j-rev}
$$\\[1.5em]
\caption{Relations of colored planar diagrams}\label{MOY-relation}
\end{figure}\\
\indent
An $[i,j]$-crossing is expanded into complicated planar diagrams as shown in Figure \ref{MOY-crossing}.
However, these colored planar diagrams locally consist of the colored oriented lines and the colored oriented trivalent diagrams $\Gamma_L$, $\Gamma_{\Lambda}$ and $\Gamma_V$ shown in Figure \ref{essential-planar}, which we call essential.
We consider a colored closed planar diagram $\Gamma$ obtained by expanding a colored oriented link diagram $D$ with the reductions in Figure \ref{MOY-crossing}.
Therefore, we find that the diagram $\Gamma$ also consists of some essential planar diagrams.\\[-1.5em]
\begin{figure}[htb]
\input{figure/figmoy3} \hspace{1cm} \input{figure/figgluing-in-3valent} \hspace{1cm} \input{figure/figgluing-out-3valent}\\
\begin{equation}
\nonumber
(1\leq m\leq n, 1\leq m_1,m_2\leq n-1, m_3=m_1+m_2\leq n)
\end{equation}\\[-1em]
\caption{Essential planar diagrams}\label{essential-planar}
\end{figure}\\
%
%
\indent
A matrix factorization for a colored closed planar diagram $\Gamma$ is defined by gluing factorizations for essential diagrams.
First, we define matrix factorizations for essential planar diagrams.
The polynomial $F_m$ is the expression of the power sum $\sum_{k=1}^{m} t_{k}^{n+1}$, where $n$ is derived from $\mathfrak{sl}_n$, in the elementary symmetric functions $x_{l}=\sum_{1 \leq k_1<\ldots<k_l\leq m}t_{k_1}\ldots t_{k_l}$ ($1\leq l\leq m$):
\begin{equation}
\nonumber
F_m(x_1,x_2,\ldots,x_m) =t_1^{n+1}+t_2^{n+1}+\ldots +t_m^{n+1}.
\end{equation}
\indent
We assign different formal indexes to ends of the colored planar diagram.
We consider an $i$-assigned end of an $m$-colored edge with the orientation from the inside diagram to the outside end.
On the end we place the function $F_m(x_{1,i},...,x_{m,i})$ with $i$-assigned variables $x_{j,i}$ $(j=1,...,m)$.
We consider an $i$-assigned end of an $m$-colored edge with the opposite orientation.
On the end we place the function $-F_m(x_{1,i},...,x_{m,i})$.
We denote the sequence $(x_{i,1},...,x_{i,m})$ by $\mathbb{X}_{(i)}^{(m)}$ and the function $F_m(x_{i,1},...,x_{i,m})$ by $F_m(\mathbb{X}_{(i)}^{(m)})$ for short.
The potential for a colored planar diagram is the sum of assigned polynomials of all ends.
\\
\indent
For instance, to an essential planar diagram we assign polynomials as shown in Figure \ref{essential-planar-function}.
Therefore, potentials for the essential planar diagrams $\Gamma_L$, $\Gamma_{\Lambda}$ and $\Gamma_V$ are $F_m(\mathbb{X}_{(1)}^{(m)})-F_m(\mathbb{X}_{(2)}^{(m)})$, $-F_{m_1}(\mathbb{X}_{(1)}^{(m_1)})-F_{m_2}(\mathbb{X}_{(2)}^{(m_2)})+F_{m_3}(\mathbb{X}_{(3)}^{(m_3)})$ and $F_{m_1}(\mathbb{X}_{(1)}^{(m_1)})+F_{m_2}(\mathbb{X}_{(2)}^{(m_2)})-F_{m_3}(\mathbb{X}_{(3)}^{(m_3)})$, respectively.\\[-1em]
\\
\indent
For 
 diagrams $\Gamma_L$, $\Gamma_{\Lambda}$ and $\Gamma_V$, we define matrix factorizations $\c(\Gamma_L)_n$, $\c(\Gamma_{\Lambda})_n$ and $\c(\Gamma_V)_n$ with potentials of the essential diagrams in Section \ref{sec4.2}.
For a general colored planar diagrams $\Gamma$, we consider a decomposition of $\Gamma$ in essential planar diagrams.
A matrix factorization for the diagram $\Gamma$ is defined by the tensor product of factorizations for the essential diagrams of the decomposition.
In the homotopy category of factorizations $\HMF^{gr}$, we find that there exist isomorphisms corresponding to relations of colored planar diagrams in Figure \ref{MOY-relation}; see Section \ref{sec4.3} and \ref{sec4.4}.\\
\indent
For the diagram $\Gamma^L_k$ in Figure \ref{colored-planar}, we can define a matrix factorization $\c(\Gamma^L_k)$ using factorizations for essential diagrams.
We decompose $\Gamma^L_k$ into essential planar diagrams using markings, assign different indexes to the markings and end points; see middle in Figure \ref{decomp-assign-planar}, and then place the polynomial $\pm F_m(\mathbb{X}^{(m)}_{(i)})$ on the marking and end points; see the right-hand side of Figure \ref{decomp-assign-planar}.
The factorization $\c(\Gamma^L_k)$ is defined to be the tensor product of factorizations for all essential diagrams in the decomposition of $\Gamma^L_k$.
For the diagram $\Gamma^R_k$, the factorization $\c(\Gamma^R_k)$ can be also defined by decomposing $\Gamma^R_k$ into essential diagrams and taking the tensor product of factorizations for all essential diagrams in the decomposition.\\[-2.5em]
\\[-3em]
%
%
\subsection{Result of the present paper II: Complex of matrix factorizations and $[i,j]$-crossing}\label{sec1.3}
\begin{figure}[htb]
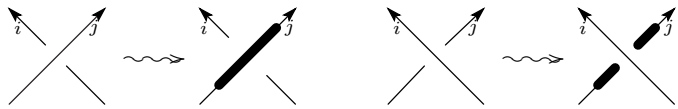

\input{figure/figplus-color}$\xymatrix{\ar@{~>}[r]&}$\input{figure/figplus-approximate-color}
\hspace{1cm}
\input{figure/figminus-color}$\xymatrix{\ar@{~>}[r]&}$\input{figure/figminus-approximate-color}
\\[1em]
\caption{Approximate diagram of $[i,j]$-crossing}\label{approximate-diag}
\end{figure}

\begin{figure}[htb]
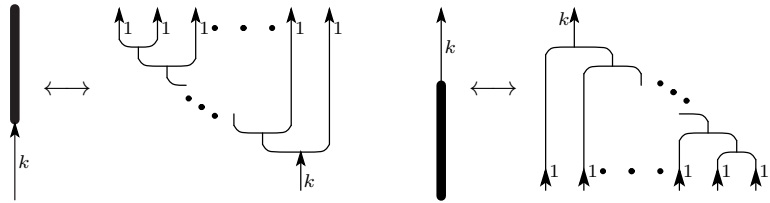

\input{figure/fig-decomp2-text} $\longleftrightarrow$ \input{figure/fig-decomp1-text} \hspace{1cm} 
\input{figure/fig-decomp4-text} $\longleftrightarrow$ \input{figure/fig-decomp3-text}\\[3em]
\caption{Wide edge and bundle of one-colored edges}\label{wide-bundle}
\end{figure}

\indent
First, we consider a complex of factorizations for an oriented link diagram with only $[1,k]$-crossings and $[k,1]$-crossings,
In the case of a $[k,1]$-crossing, the state model of the quantum $(\mathfrak{sl}_n,\land V_n)$ link invariant takes the following forms
\begin{eqnarray}
\nonumber
\left<\hspace{.1cm}\input{figure/figplus-k-1-intro}\hspace{.1cm}\right>_n&=&
(-1)^{1-k}q^{kn-1}\left<\hspace{.3cm}\input{figure/figsquare1k--k-1--k1}\hspace{.3cm}\right>_n+(-1)^{-k}q^{kn}\left<\hspace{.3cm}\input{figure/figsquare1k--k+1--k1}\hspace{.3cm}\right>_n,\\
\nonumber
\left<\hspace{.1cm}\input{figure/figminus-k-1-intro}\hspace{.1cm}\right>_n&=&
(-1)^{k-1}q^{-kn+1}\left<\hspace{.3cm}\input{figure/figsquare1k--k-1--k1}\hspace{.3cm}\right>_n+(-1)^{k}q^{-kn}\left<\hspace{.3cm}\input{figure/figsquare1k--k+1--k1}\hspace{.3cm}\right>_n.
\end{eqnarray}
\begin{figure}[htb]
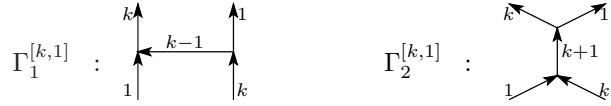

\begin{equation}
\nonumber
\Gamma^{[k,1]}_1\hspace{.2cm}:\hspace{.4cm}\input{figure/figsquare1k--k-1--k1-big}\hspace{2cm}
\Gamma^{[k,1]}_2\hspace{.2cm}:\hspace{.4cm}\input{figure/figsquare1k--k+1--k1-big}
\end{equation}
\caption{Colored planar diagrams in reduction of $[k,1]$-crossing}\label{resolved-planar}
\end{figure}\\
For $[1,k]$-crossings \input{figure/figminus-1-k-intro1} and \input{figure/figplus-1-k-intro1}, these brackets have a similar form, see Appendix \ref{NMOY}.
In Section \ref{sec5}, we define $\Z$-grading-preserving morphisms between the factorizations $\c(\Gamma^{[k,1]}_1)_n$ and $\c(\Gamma^{[k,1]}_2)_n$ for diagrams in Figure \ref{resolved-planar}:
\begin{eqnarray}
\nonumber
\chi^{[k,1]}_{+}:\c\left(\hspace{.3cm}\input{figure/figsquare1k--k+1--k1}\hspace{.3cm}\right)_n\longrightarrow\c\left(\hspace{.3cm}\input{figure/figsquare1k--k-1--k1}\hspace{.3cm}\right)_n,\\
\nonumber
\chi^{[k,1]}_{-}:\c\left(\hspace{.3cm}\input{figure/figsquare1k--k-1--k1}\hspace{.3cm}\right)_n\longrightarrow\c\left(\hspace{.3cm}\input{figure/figsquare1k--k+1--k1}\hspace{.3cm}\right)_n.
\end{eqnarray}
\indent
Using these morphisms, a complex for a single $[k,1]$-crossing is defined as an object of $\k^b(\HMF^{gr})$,
\begin{eqnarray}
\label{k-1-complex-p}
&&\xymatrix{\c\left(\input{figure/figplus-k-1-intro}\right)_n=
\ldots
\ar[r]&0
\ar[r]&\c^{-k}\left(\input{figure/figsquare1k--k+1--k1}\right)_n
\ar[r]^{\chi^{[k,1]}_{+}}&\c^{1-k}\left(\input{figure/figsquare1k--k-1--k1}\right)_n
\ar[r]&0
\ar[r]&
\ldots},\\
\label{k-1-complex-m}
&&\xymatrix{\c\left(\input{figure/figminus-k-1-intro}\right)_n=
\ldots
\ar[r]&0
\ar[r]&\c^{k-1}\left(\input{figure/figsquare1k--k-1--k1}\right)_n
\ar[r]^{\chi^{[k,1]}_{-}}&\c^{k}\left(\input{figure/figsquare1k--k+1--k1}\right)_n
\ar[r]&0
\ar[r]&
\ldots}.
\end{eqnarray}
We remark that this construction is a generalization of a complex for a $[1,2]$-crossing given by Rozansky \cite{Roz}.\\
\indent
To an oriented link diagram $D$ with $[1,k]$-crossings and $[k,1]$-crossings, we define a complex of matrix factorizations to be decomposing $D$ into single $[1,k]$-crossings and $[k,1]$-crossings and taking the tensor product of complexes for all single $[1,k]$-crossings and $[k,1]$-crossings in the decomposition.
The obtained complex is a complex of matrix factorizations with potential zero.
Then, the complex for the diagram $D$ gives rise to a $\Z\oplus\Z\oplus\Z_2$-graded homology.\\
\indent
If two oriented link diagrams with $[1,k]$-crossings and $[k,1]$-crossings transform to each other under colored Reidemeister moves which are composed of $[1,k]$-crossings and $[k,1]$-crossings only, then the associated $\Z\oplus\Z\oplus\Z_2$-graded homologies are isomorphic.
More precisely, for tangle diagrams with $[1,k]$-crossings and $[k,1]$-crossings transforming to each other under colored Reidemeister moves composed of $[1,k]$-crossings and $[k,1]$-crossings, those complexes of matrix factorizations are isomorphic in $\k^{b}(\HMF^{gr})$.
\begin{theorem}[Theorem \ref{main1} in Section \ref{sec5.2} (In the case $k=1$, Khovanov-Rozansky\cite{KR1})]
We consider tangle diagrams with $[1,k]$-crossings and $[k,1]$-crossings transforming to each other under colored Reidemeister moves composed of $[1,k]$-crossings and $[k,1]$-crossings.
Complexes of factorizations for these tangle diagrams are isomorphic in $\k^{b}(\HMF^{gr})$:
\begin{eqnarray}
\nonumber
&&\c\left(\input{figure/r1p-1}\right)_n\simeq\c\left(\input{figure/r1c-1}\right)_n\simeq\c\left(\input{figure/r1m-1}\right)_n,\\
\nonumber
&&\c\left(\input{figure/r2-1kl}\right)_n\simeq\c\left(\input{figure/r2-1kc}\right)_n\simeq\c\left(\input{figure/r2-1kr}\right)_n,\,
\c\left(\input{figure/r2-k1l}\right)_n\simeq\c\left(\,\input{figure/r2-k1c}\right)_n\simeq\c\left(\,\input{figure/r2-k1r}\right)_n,\\
\nonumber
&&\c\left(\input{figure/r2-1klrev}\right)_n\simeq\c\left(\input{figure/r2-1kcrev}\right)_n\simeq\c\left(\input{figure/r2-1krrev}\right)_n,\,
\c\left(\input{figure/r2-1klrevori}\right)_n\simeq\c\left(\input{figure/r2-1kcrevori}\right)_n\simeq\c\left(\input{figure/r2-1krrevori}\right)_n,\\
\nonumber
&&\c\left(\input{figure/r3-k11}\right)_n\simeq\c\left(\input{figure/r3-k11rev}\right)_n,
\c\left(\input{figure/r3-1k1}\right)_n\simeq\c\left(\input{figure/r3-1k1rev}\right)_n,
\c\left(\input{figure/r3-11k}\right)_n\simeq\c\left(\input{figure/r3-11krev}\right)_n.
\end{eqnarray}
\end{theorem}
\indent
For a colored oriented link diagram $D$ with $[1,k]$-crossings and $[k,1]$-crossings, we can explicitly calculate the complex $\c\left(D\right)_n$ and the $\Z\oplus\Z\oplus\Z_2$-graded link homology $\H^{i,j,k}(D)$.
We evaluate the Poincar\'e polynomial of the homology $\H^{i,j,k}(D_{Hopf})$ for an oriented Hopf link diagram $D_{Hopf}$ with a $[1,k]$-crossing and a $[k,1]$-crossing in Section \ref{sec5.5}.
\\
\indent
In the case of general $[i,j]$-crossings, it is difficult both to define boundary maps of a complex of matrix factorization for the $[i,j]$-crossing explicitly and to show that there are isomorphisms between complexes for the colored tangle diagrams that transform to each other under colored Reidemeister moves in $\k^{b}(\HMF^{gr})$.
Instead of this construction of the homological link invariant, we introduce an approximate $[i,j]$-crossing and define a complex for the approximate crossing in Figure \ref{approximate-diag}.
The wide edge of the approximate $[i,j]$-crossing represents a bundle of one-colored edges, see Figure \ref{wide-bundle}.
We arrange an $[i,j]$-crossing in the orientation from bottom to up and change a colored edge from the left-bottom to the right-top into a wide edge at an over crossing or an under crossing.
Therefore, we can define a complex for the approximate crossing using the definition of the complex for an $[i,1]$-crossing since every crossing of the approximate crossing is an $[i,1]$-crossing.
\\
\indent
We consider the homology of this complex.
The homology is not a link invariant.
However, we can obtain a link invariant as a normalized Poincar\'e polynomial of the homology.
The polynomial link invariant is a polynomial in $\Q[\,t^{\pm 1},q^{\pm 1},s]/\langle s^2-1\rangle$ and the quantum $(\mathfrak{sl}_n,\land V_n)$ link invariant by specializing $t$ to $-1$ and $s$ to $1$.
\\
\\
\indent
We consider approximate tangle diagrams transforming to each other under the Reidemeister moves composed of the approximate crossings.
We find the following isomorphisms in $\k^b(\HMF^{gr})$ as follows:
\begin{theorem}[Theorem \ref{main2} in Section \ref{sec6.2}]\label{approximate}
For approximate tangle diagrams transforming to each other under the Reidemeister moves composed of the approximate crossings, complexes of matrix factorizations for these approximate tangle diagrams are isomorphic in $\k^{b}(\HMF^{gr})$:
\begin{eqnarray}
\nonumber
&&
\c\left(\input{figure/r1p-i-wide}\right)_n
\hspace{.1cm}\simeq\hspace{.1cm}
\c\left(\input{figure/r1c-i-wide}\right)_n
\hspace{.1cm}\simeq\hspace{.1cm}
\c\left(\input{figure/r1m-i-wide}\right)_n,
\hspace{.4cm}
\c\left(\input{figure/r2-ijl-wide}\right)_n
\hspace{.1cm}\simeq\hspace{.1cm}
\c\left(\input{figure/r2-ijc-wide}\right)_n
\hspace{.1cm}\simeq\hspace{.1cm}
\c\left(\input{figure/r2-ijr-wide}\right)_n,\\
\nonumber
&&
\c\left(\input{figure/r2-ijlrev-wide}\right)_n
\hspace{.1cm}\simeq\hspace{.1cm}
\c\left(\input{figure/r2-ijcrev-wide}\right)_n,
\hspace{.4cm}
\c\left(\input{figure/r2-ijlrevori-wide}\right)_n
\hspace{.1cm}\simeq\hspace{.1cm}
\c\left(\input{figure/r2-ijcrevori-wide}\right)_n,
\hspace{.4cm}
\c\left(\input{figure/r3-ijk-wide}\right)_n
\hspace{.1cm}\simeq\hspace{.1cm}
\c\left(\input{figure/r3-ijkrev-wide}\right)_n.
\end{eqnarray}
\end{theorem}
We have not got produced isomorphisms between complexes for colored tangle diagrams transforming to each other under colored Reidemeister moves, but we hope to return to this question in a future paper.
However, for a colored oriented link diagram $D$, we consider a $\Z\oplus\Z\oplus\Z_2$-graded homology $H^{i,j,k}(D)$ through the complex for the approximate diagram of $D$.
Then, we define a polynomial $\overline{P}(D)$ to be the Poincar\'e polynomial of the homology $H^{i,j,k}(D)$
\begin{equation}
\nonumber
\overline{P}(D):=\sum_{i,j,k}t^i q^j s^k\dim_{\Q}H^{i,j,k}(D) \in \Q[t^{\pm 1},q^{\pm 1},s]/\langle s^2-1\rangle.
\end{equation}
A polynomial link invariant can be obtained by normalizing the Poincar\'e polynomial $\overline{P}(D)$ as follows.
For a colored oriented link diagram $D$, we define a function $\mathrm{Cr}_k$ $(k=1,...,n-1)$ to be
\begin{center}
$\mathrm{Cr}_k(D):=$ the number of $[\ast,k]$-crossing of $D$.
\end{center}
We then define the normalized Poincar\'e polynomial $P(D)$ to be
\begin{equation}
\nonumber
P(D):=\overline{P}(D)\prod_{k=1}^{n-1}\frac{1}{\left([k]_q!\right)^{\mathrm{Cr}_k(D)}}.
\end{equation}
By the construction and Theorem \ref{approximate}, we find that $P(D)$ is a link invariant.
\begin{theorem}[Corollary \ref{main3} in Section \ref{sec6.2}]\label{new-invariant}
Two colored oriented link diagrams $D$ and $D'$ that transform to each other under colored Reidemeister moves have the same normalized Poincar\'e polynomial,
\begin{equation}
\nonumber
P(D)=P(D').
\end{equation}
That is, we have the following equations for evaluations of colored oriented link diagrams:
\begin{eqnarray}
\nonumber
&&
P\left(\input{figure/r1p-i}\right)
\hspace{.1cm}=\hspace{.1cm}
P\left(\input{figure/r1c-i}\right)
\hspace{.1cm}=\hspace{.1cm}
P\left(\input{figure/r1m-i}\right),
\hspace{.4cm}
P\left(\input{figure/r2-ijl}\right)
\hspace{.1cm}=\hspace{.1cm}
P\left(\input{figure/r2-ijc}\right)
\hspace{.1cm}=\hspace{.1cm}
P\left(\input{figure/r2-ijr}\right),\\
\nonumber
&&
P\left(\input{figure/r2-ijlrev}\right)
\hspace{.1cm}=\hspace{.1cm}
P\left(\input{figure/r2-ijcrev}\right),
\hspace{.4cm}
P\left(\input{figure/r2-ijlrevori}\right)
\hspace{.1cm}=\hspace{.1cm}
P\left(\input{figure/r2-ijcrevori}\right),
\hspace{.4cm}
P\left(\input{figure/r3-ijk}\right)
\hspace{.1cm}=\hspace{.1cm}
P\left(\input{figure/r3-ijkrev}\right).
\end{eqnarray}
The outsides of colored tangle diagrams in each equation have the same picture.
\end{theorem}
The polynomial $P(D)$ is a refined link invariant of the quantum $(\mathfrak{sl}_n,\land V_n)$ link invariant since $P(D)$ is the quantum $(\mathfrak{sl}_n,\land V_n)$ link invariant by specializing $t$ to $-1$ and $s$ to $1$.
\subsection{Organization of paper}
\indent
The present paper consists of seven sections and two appendixes.
In Section \ref{sec2}, we recall definition of a matrix factorization, basic properties and theorems in a category of factorizations $\MF^{gr}$ and its homotopy category $\HMF^{gr}$.
For propositions given by the author, proofs are described and, for other important propositions, reference are given.
Then, we define the complex category $\kom(\HMF^{gr})$ and its homotopy category $\k^b(\HMF^{gr})$ in Section \ref{sec2.11} and \ref{sec2.12}.
The homology for a colored oriented link diagram is evaluated as an object of $\k^b(\HMF^{gr})$.
In Section \ref{sec3}, symmetric functions and its generating function are introduced.
The symmetric functions are used for defining a matrix factorization for a colored planar diagram.
Section \ref{sec4}, \ref{sec5} and \ref{sec6} are the main part of this paper and include the author's original results.
In Section \ref{sec4}, we define factorizations for essential planar diagrams using the symmetric functions.
Then, we show that there exist isomorphisms in $\HMF^{gr}$ corresponding to most MOY relations.
In Section \ref{sec5}, we define complexes of factorizations for $[1,k]$-crossing and $[k,1]$-crossing in $\ob(\k^b(\HMF^{gr}))$ and, for tangle diagrams with $[k,1]$-crossings and $[k,1]$-crossings transforming to each other under colored Reidemeister moves composed of $[k,1]$-crossings and $[k,1]$-crossings, we show isomorphisms in $\k^b(\HMF^{gr})$ between matrix factorizations of these tangle diagrams.
For the Hopf link with $[1,k]$-crossing and $[k,1]$-crossing, the Poincar\'e polynomial of the homological invariant is shown in Section \ref{sec5.5}.
In Section \ref{sec6}, we introduce a wide edge and an approximate $[i,j]$-crossing with the wide edges.
Then, we define a matrix factorization for the approximate $[i,j]$-crossing with the wide edges using a complex for $[1,k]$-crossing and show isomorphisms between complexes for approximate tangle diagrams transforming to each other under Reidemeister moves of approximate crossings.
For a colored oriented link diagram $D$, we construct a polynomial link invariant $P(D)$ as a normalized Poincar\'e polynomial of the homology of a complex for an approximate link diagram associated to the diagram $D$.
In Section \ref{sec7}, we give proofs of some properties and theorems which are skipped in order to understand this study.
In Appendix \ref{app1}, we remark that a generalization of this study to virtual link theory introduced by Kauffman.
In Appendix \ref{NMOY}, we recall the quantum ($\mathfrak{sl}_n$,$\land V_n$) link invariant by a normalizing MOY bracket\footnote{MOY bracket is a regular link invariant. That is, it satisfies invariance under colored Reidemeister moves II and III only. However, we can generally obtain the link invariant by normalizing the regular invariant.} given by H. Murakami, T. Ohtsuki and S. Yamada.
Then, we show that the normalized MOY bracket satisfies invariance under colored Reidemeister move I.
\\\\
\indent
{\bf Acknowledgement.}
The author is heartfelt thanks to Osamu Iyama, Hiroyuki Ochiai and Akihiro Tsuchiya for their appropriate advice and helpful comment of this study and deeply grateful thanks to Lev Rozansky for showing him how to construct a complex for $[1,2]$-crossing by using matrix factorizations explicitly.
He want to thank Lars Hesselholt for correcting the wording of Section \ref{1-cat} of the second version.
The author is also thanks to K. Aragane, T. Araya, H. Awata, T. Ikeda, H. Kanno, T. Kawamura, K. Iijima, H. Miyachi, H. Murakami, Y. Nohara, N. Sato, R. Takahashi, K. Wada, H. Wu, K. Yoshida and Y. Yoshino for their helpful and interesting comment.
This work is partly supported by the Grant-in-Aid for JSPS Fellows (20-2330) from Japan Society for the Promotion of Science.
%
%
%
%
\section{$\Z$-graded matrix factorization}\label{sec2}
In this section, we recall definitions, properties and theorems given by \cite{KR1}\cite{KR2}\cite{KR3}\cite{Ras}\cite{Wu1}\cite{Yone1}\cite{Yoshi}.
Many essential facts are given by Mikhail Khovanov and Lev Rozansky\cite{KR1}.
%
%
%
%
\subsection{$\Z$-graded module}\label{sec2.1}
Let $R=\Q[x_1,...,x_r]$ be a polynomial ring such that the degree $\deg(x_i)\in \Z$ is a positive integer given for each $i=1,...,r$.
Then, $R$ has a $\Z$-grading decomposition $\oplus_{i}R^i$ such that $R^iR^j\subset R^{i+j}$ and $R^0=\Q$.
A maximal ideal generated by homogeneous polynomials is unique, denoted the maximal ideal by $\mathfrak{m}$.
We consider a free $\Z$-graded $R$-module $M$ with a $\Z$-grading decomposition $\oplus_i M^i$ such that $R^jM^i\subset M^{i+j}$ for any $i\in\Z$.
Remark we allow a $\Z$-graded $R$-module to be infinite-rank.
\\
\indent
A {\bf $\Z$-grading shift} $\{ m \}$ ($m \in \Z$) is an operator up $\Z$-grading $m$ on an $R$-module.
That is, for a $\Z$-graded $R$-module $M$ with a $\Z$-grading decomposition $\bigoplus_i M^i$, we have an equality as a $\Q$-vector space
\begin{equation}
\nonumber
\left( M \{m\} \right)^i = M^{i-m}.
\end{equation}
For a Laurent polynomial $f(q)=\sum a_iq^i\in \N_{\geq 0}[q,q^{-1}]$, we define $M\{f(q)\}_{q}$ to be
\begin{equation}
\nonumber
M\{f(q)\}_{q}:=\bigoplus_{i}\left(M\{i\}\right)^{\oplus a_i}.
\end{equation}
\indent
For $R$-modules $M$ and $N$, the set of $\Z$-grading preserving morphisms from $M$ to $N$ is denoted by $\Hom_R(M,N)$.
When $N=M$ we denote it by $\End_R(M)$.
We find that a $\Z$-grading preserving morphism of $\Hom_R(M\{m\},N)$ is an $R$-module morphism with $\Z$-grading $m$ from $M$ to $N$.
We often regard $\Hom_R(M\{m\},N)$ as the set of $R$-module morphisms with $\Z$-grading $m$ from $M$ to $N$.
We consider the set of all of $R$-module morphisms from $M$ to $N$ denoted by $\HOM_{R}(M,N)$, that is
\begin{equation}
\nonumber
\HOM_{R}(M,N):=\bigoplus_{m\in \Z}\Hom_R(M\{m\},N).
\end{equation}
When $N=M$ it is denoted by $\END_R(M)$.
These sets $\HOM_{R}(M,N)$ and $\END_{R}(M)$ naturally have a $\Z$-graded $R$-module structure by
\begin{equation}
\nonumber
\xymatrix{
r_k:\Hom_{R}(M\{m\},N)\ar@{}[d]|-{\rotatebox{90}{$\in$}}\ar[r]&\Hom_{R}(M\{m+k\},N)\ar@{}[d]|-{\rotatebox{90}{$\in$}}\\
f_m\ar@{|->}[r]&r_k\,f_m
,}
\end{equation}
where $r_k\in R$ has a $\Z$-grading $k$.
\begin{proposition}\label{finite-mod}
If $M$ and $N$ are free $R$-modules of finite rank, then $\HOM_{R}(M,N)$ is a free $R$-module of finite rank.
\end{proposition}
For a $\Z$-graded $R$-module $M$, we define $M^{\ast}$ to be $\HOM_{R}(M,R)$,
\begin{equation}
\nonumber
M^{\ast}:=\HOM_{R}(M,R).
\end{equation}
The $R$-module $M^{\ast}$ is called {\bf dual} of $M$.
\begin{corollary}
If $M$ is a free $R$-module of finite rank, $M^{\ast}$ is a free $R$-module of finite rank.
\end{corollary}
\begin{proposition}\label{duality-mod}
If $M$ is a free $R$-module of finite rank, $M^{\ast\ast}\simeq M$.
\end{proposition}
\begin{proposition}\label{isom-hom}
If $M$ is a free $R$-module of finite rank, we have a $\Z$-graded $R$-module isomorphism 
\begin{equation}
\nonumber
\HOM_R(M,N)\simeq N\otimes_R M^{\ast}.
\end{equation}
\end{proposition}
%
%
%
%
\subsection{Potential and Jacobian ring}\label{sec2.2}
\indent
For a homogeneous $\Z$-graded polynomial $\omega\in R$, we define a quotient ring $R_{\omega}$ to be $R\left/I_{\omega}\right.$, where $I_{\omega}$ is the ideal generated by partial derivatives $\frac{\partial \omega}{\partial x_k}$ ($1\leq k\leq r$).
The quotient ring $R_{\omega}$ is called {\bf Jacobian ring} of $\omega$.
A homogeneous element $\omega\in \mathfrak{m}$ is a {\bf potential} of $R$ if the Jacobian ring $R_{\omega}$ is finitely dimensional as a $\Q$ vector space.
That is, the partial derivatives $\frac{\partial \omega}{\partial x_k}$ ($1\leq k\leq r$) form a regular sequence in $R$.
%
%
%
%
\subsection{$\Z$-graded matrix factorization}\label{sec2.3}
For a polynomial $\omega$ with an even homogeneous $\Z$-grading, we consider the polynomial ring $R$ generated by variables included in $\omega$.
Assume that a given polynomial $\omega$ is a potential of $R$.
We allow $\omega$ to be zero and, in such a case, we consider $\Q$ as a polynomial ring.
In this setting, we define a $\Z$-graded matrix factorization as follows \cite{Yoshi}.\\
\indent
We suppose a 4-tuple $\overline{M}=(M_0,M_1,d_{M_0},d_{M_1})$ is a two-periodic chain
\begin{equation}
\nonumber
\xymatrix{M_0\ar[r]^{d_{M_0}}&M_1\ar[r]^{d_{M_1}}&M_0},
\end{equation}
where $M_0$ and $M_1$ are $\Z$-graded free $R$-modules permitted to be infinite-rank and
$d_{M_0}:M_0\to M_1$ and $d_{M_1}:M_1\to M_0$ are homogeneous $\Z$-graded morphisms of $R$-modules.
When we fix $R$-module bases of $M_0$ and $M_1$, $d_{M_0}$ and $d_{M_1}$ can be represented as a matrix form.
The matrix forms of $d_{M_0}$ and $d_{M_1}$ are denoted by the capital letter $D_{M_0}$ and $D_{M_1}$.\\
\indent
We say that a 4-tuple $\overline{M}$ is a {\bf matrix factorization} with a potential $\omega \in \mathfrak{m}$ ({\bf factorization} for short ), if the composition of $d_{M_0}$ and $d_{M_1}$ satisfies that $d_{M_1} d_{M_0} = \omega\,\id_{M_0}$, $d_{M_0} d_{M_1} = \omega\,\id_{M_1}$ and $d_{M_1}$, $d_{M_2}$ have the $\Z$-grading $\frac{1}{2}\deg \, \omega$.\\
\indent
\indent
\begin{definition}
We define a category $\MF^{gr,all}_{R,\omega}$ of $\Z$-graded matrix factorizations as follows.
\begin{itemize}
\item
An object in $\MF^{gr,all}_{R,\omega}$ is a matrix factorization $\overline{M}=(M_0,M_1,d_{M_0},d_{M_1})$ with the potential $\omega$, where $M_0$, $M_1$ are $R$-modules and $d_{M_0}$, $d_{M_1}$ are $R$-module morphisms with the $\Z$-grading $\frac{1}{2}\deg \, \omega$.
\item
A morphism in $\MF^{gr,all}_{R,\omega}$ from
$\overline{M}=(M_0,M_1,d_{M_0},d_{M_1})$ to
$\overline{N}=(N_0,N_1,d_{N_0},d_{N_1})$ is a pair $\overline{f}=(f_0,f_1)$ of $\Z$-grading preserving morphisms of $R$-modules $f_0:M_0\to N_0$ and $f_1:M_1\to N_1$ to give a commutative diagram, that is, the morphism $\overline{f}=(f_0,f_1)$ satisfies $d_{N_0} f_0 = f_1 d_{M_0}$ and $d_{N_1} f_1 = f_0 d_{M_1}$,
\begin{equation}
\nonumber
\xymatrix{
M_0\ar[r]^{d_{M_0}}\ar[d]_{f_0}&M_1\ar[r]^{d_{M_1}}\ar[d]_{f_1}&M_0\ar[d]_{f_0}\\
N_0\ar[r]^{d_{N_0}}&N_1\ar[r]^{d_{N_1}}&N_0.
}
\end{equation}
For any matrix factorizations $\overline{M}$ and $\overline{N}$, we denote the set of $\Z$-grading preserving morphisms from $\overline{M}$ to $\overline{N}$ by $\Hom_{\MF}(\overline{M},\overline{N})$.
When $\overline{N}=\overline{M}$ we denote it by $\End_{\MF}(\overline{M})$.
\item
The composition $\overline{f}\overline{g}$ of morphisms $\overline{f}=(f_0,f_1)$ and $\overline{g}=(g_0,g_1)$ is defined by $(f_0g_0,f_1g_1)$.
\end{itemize}
\end{definition}
\indent
A matrix factorization $(M_0,M_1,d_0,d_1)$ $\in \MF^{gr,all}_{R,\omega}$ is {\bf finite} if $M_0$ and $M_1$ are free $R$-modules of finite rank.\\
\indent
Let $\MF^{gr,fin}_{R,\omega}$ be the full subcategory in $\MF^{gr,all}_{R,\omega}$ whose object is a finite factorizations.\\
\indent
A morphism $\overline{f}=(f_1,f_2): \overline{M} \to \overline{N}$ in $\MF^{gr,all}_{R,\omega}$ is {\bf null-homotopic} 
if morphisms $h_{0}:M_0 \to N_1$ and $h_{1}:M_1 \to N_0$ exist 
such that $f_0 = h_1 d_{M_0} + d_{N_1} h_0$ and $f_1 = h_0 d_{M_1} + d_{N_0} h_1$,
\begin{equation}
\nonumber
\xymatrix{
M_0\ar[r]^{d_{M_0}}\ar[d]_{f_0}&M_1\ar[r]^{d_{M_1}}\ar[d]_{f_1}\ar@{.>}[ld]_{h_1}&M_0\ar[d]_{f_0}\ar@{.>}[ld]_{h_0}\\
N_0\ar[r]_{d_{N_0}}&N_1\ar[r]_{d_{N_1}}&N_0.
}
\end{equation}
Two morphisms $\overline{f},\overline{g}: \overline{M} \to \overline{N}$ in $\MF^{gr,all}_{R,\omega}$ are {\bf homotopic} :
if $\overline{f} - \overline{g}$ is null-homotopic.
We often denote $\overline{f}\mfsim\overline{g}$ for short if $\overline{f}$ and $\overline{g}$ are homotopic.
Two matrix factorizations $\overline{M}$ and $\overline{N}$ in $\MF^{gr,all}_{R,\omega}$ are {\bf homotopy equivalent} if there exist morphisms $\overline{f}:\overline{M}\to \overline{N}$ and $\overline{g}:\overline{N}\to \overline{M}$ such that $\overline{f}\overline{g}\mfsim\id_{\overline{N}}$ and $\overline{g}\overline{f}\mfsim\id_{\overline{M}}$.
\begin{definition}
We define the homotopy category $\HMF^{gr,all}_{R,\omega}$ of $\MF^{gr,all}_{R,\omega}$ as follows;
\begin{itemize}
\item
$\HMF^{gr,all}_{R,\omega}$ has the same objects of $\MF^{gr,all}_{R,\omega}$
\item
A morphism in $\HMF^{gr,all}_{R,\omega}$ is a morphism in $\MF^{gr,all}_{R,\omega}$ modulo null-homotopic.
Denote the set of $\Z$-grading preserving morphisms from $\overline{M}$ to $\overline{N}$ in $\HMF^{gr,all}_{R,\omega}$ by $\Hom_{\HMF}(\overline{M},\overline{N})$, that is,
\begin{equation}
\nonumber
\Hom_{\HMF}(\overline{M},\overline{N})=\Hom_{\MF}(\overline{M},\overline{N})\left/\{null-homotopic\}\right..
\end{equation}
When $\overline{N}=\overline{M}$ we denote it by $\End_{\HMF}(\overline{M})$.
\item
The composition of morphisms is defined by the same way to $\MF^{gr,all}_{R,\omega}$
\end{itemize}
\end{definition}
It is obvious that homotopy equivalent is isomorphic in $\HMF^{gr,all}_{R,\omega}$.
Let $\HMF^{gr,fin}_{R,\omega}$ be the full subcategory in $\HMF^{gr,all}_{R,\omega}$ whose object is a finite factorizations in $\HMF^{gr,all}_{R,\omega}$.\\
\indent
A matrix factorization is {\bf contractible} (or called trivial in \cite{Yoshi}) if it is isomorphic in the homotopy category to the zero matrix factorization
$$
\Big(
\xymatrix{
*{0}\ar[rr]&&*{0}\ar[rr]&&*{0} 
}
\Big).
$$
Especially, $(R,R\{\frac{1}{2}\deg\,\omega\},1,\omega)$ and $(R,R\{-\frac{1}{2}\deg\,\omega\},\omega,1)$ are contractible.
In general, a contractible matrix factorization is a direct sum of the factorizations $(R,R\{\frac{1}{2}\deg\,\omega\},1,\omega)$ and $(R,R\{-\frac{1}{2}\deg\,\omega\},\omega,1)$.
A matrix factorization is {\bf essential} (or called reduced in \cite{Yoshi}) if it does not include any contractible matrix factorizations.
For a matrix factorization $\overline{M}$, we denote its essential summand by $\overline{M}_{es}$ and its contractible summand by $\overline{M}_{c}$.\\
\indent
We define a {\bf $\Z$-grading shift} $\{ m \}$ ($m \in \Z$) on $\overline{M}=(M_0,M_1,d_{M_0},d_{M_1})$ to be
\begin{eqnarray}
\nonumber
\overline{M}\{m\}&=&(M_0\{ m \},M_1\{ m\},d_{M_0},d_{M_1}).
\end{eqnarray}
For a Laurent polynomial $f(q)=\sum a_iq^i\in \N_{\geq 0}[q,q^{-1}]$, we define $\overline{M}\{f(q)\}_{q}$ to be
\begin{equation}
\nonumber
\overline{M}\{f(q)\}_{q}=\bigoplus_{i}\left(\overline{M}\{i\}\right)^{\oplus a_i}.
\end{equation}
\indent
The {\bf translation functor} $\left< 1 \right>$ changes a matrix factorization $\overline{M} = (M_{0},M_{1},d_{M_0},d_{M_1})\in \ob(\MF^{gr,all}_{R,\omega})$ into 
\begin{eqnarray}
\nonumber
\overline{M} \left< 1 \right> &=&
(
M_{1},M_{0},-d_{M_1},-d_{M_0})\in \ob(\MF^{gr,all}_{R,\omega})
.
\end{eqnarray}
The functor $\left< 2 \right> (= \left< 1 \right>^2)$ is the identity functor.
In general, we denote $\left<1\right>^k$ by $\left<k\right>$.
%
%
\subsection{Cohomology of matrix factorization}\label{sec2.4}
Let $\mathfrak{m}$ be a unique maximal ideal generated by homogeneous $\Z$-graded polynomials of $R$.
For a matrix factorization $\overline{M}=(M_0,M_1,d_0,d_1)\in \ob(\MF^{gr,all}_{R,\omega})$ (resp. $\ob(\HMF^{gr,all}_{R,\omega})$), we define a quotient $\overline{M}\left/\mathfrak{m}\overline{M}\right.$ to be a two-periodic complex of vector spaces over $R\left/\mathfrak{m}\right.\simeq \Q$,
\begin{equation}
\nonumber
\xymatrix{M_0\left/\mathfrak{m}M_0\right.\ar[r]^{d_0}&M_1\left/\mathfrak{m}M_1\right.\ar[r]^{d_1}&M_0\left/\mathfrak{m}M_0\right.}.
\end{equation}
The compositions of $d_0$ and $d_1$ satisfy that $d_1d_0=0$ and $d_0d_1=0$ since $\omega\in\mathfrak{m}$.
We denote the cohomology of the quotient complex $\overline{M}\left/\mathfrak{m}\overline{M}\right.$ by $\H(\overline{M})$ and call it the {\bf cohomology of matrix factorization}.
The cohomology of $\overline{M}$ is a $\Z\oplus\Z_2$-graded $\Q$-vector space, denoted by $\H(\overline{M})= \H^0(\overline{M})\oplus \H^1(\overline{M})$.
A morphism $\overline{f}=(f_0,f_1):\overline{M}\to\overline{N}$ naturally induces a morphism from the cohomology $\H(\overline{M})$ to $\H(\overline{N})$, denoted by $\H(\overline{f})$ or $(\H(f_0),\H(f_1))$.
\begin{proposition}\label{null-0}
A null-homotopic morphism $\overline{f}=(f_0,f_1):\overline{M}\to\overline{N}$ induces the morphism $0$ from the cohomology $\H(\overline{M})$ to $\H(\overline{N})$.
\end{proposition}
\begin{proposition}
The cohomology of a contractible matrix factorization is the zero.
\end{proposition}
\begin{definition}
\hspace{1cm}\\
\begin{itemize}
\item[(1)]
We define a {\bf $\Z$-graded dimension} of a $\Z$-graded $\Q$-vector space $V$ to be
\begin{equation}
\nonumber
\gdim(V):=\sum_{i\in\Z}q^i\dim_{\Q}V^i,
\end{equation}
where $V^i$ is the $i$-graded subspace of $V$.\\
\item[(2)]
We define a {\bf $\Z_2\oplus\Z$-graded dimension} of a matrix factorization $\overline{M}$ to be
\begin{equation}
\nonumber
\gdim\left(\overline{M}\right):=\gdim\,\H^0(\overline{M})+s\,\gdim\,\H^1(\overline{M}),
\end{equation} 
where the variable $s$ satisfies $s^2=1$.
\end{itemize}
\end{definition}
\begin{proposition}[Proposition $7$ \cite{KR1}]
The following conditions on $\overline{M}=(M_0,M_1,d_0,d_1)\in \MF_{R,\omega}^{gr,all}$ are equivalent.
\begin{itemize}
\item[(1)]$\overline{M}$ is isomorphic in $\MF_{R,\omega}^{gr,all}$ to a direct sum of $(R\{m_1\},R\{m_1+\frac{1}{2}\deg\,\omega\},1,\omega)$ and $(R\{m_2\},R\{m_2-\frac{1}{2}\deg\,\omega\},\omega,1)$, where $m_1,m_2\in\Z$.
\item[(2)]$\overline{M}$ is isomorphic in $\HMF_{R,\omega}^{gr,all}$ to the zero factorization $(\xymatrix{0\ar[r]&0\ar[r]&0})$.
\item[(3)]$\H(\overline{M})=0$.
\end{itemize}
\end{proposition}
\begin{proposition}[Proposition $8$ \cite{KR1}]
The following conditions of a morphism $(f_0,f_1) : \xymatrix{\overline{M}\ar[r]&\overline{N}}$ are equivalent.
\begin{itemize}
\item[(1)] $(f_0,f_1)$ is an isomorphism in $\HMF^{gr,all}_{R,\omega}$.
\item[(2)] $(f_0,f_1)$ induces an isomorphism $(\H(f_0),\H(f_1))$ between the cohomologies $\H(\overline{M})$ and $\H(\overline{N})$.
\end{itemize}
\end{proposition}
\begin{corollary}[Corollary 3\cite{KR1}]\label{unique-decomp}
For a matrix factorization $\overline{M}\in\ob(\MF^{gr,all}_{R,\omega})$,  a decomposition of $\overline{M}$ into into an essential factorization and a contractible factorization $\overline{M}_{es}\oplus\overline{M}_{c}$ is unique up to isomorphism.
\end{corollary}
\begin{corollary}[Corollary 4\cite{KR1}]\label{unique-finite}
A matrix factorization $\overline{M}\in\ob(\MF^{gr,all}_{R,\omega})$ with finite-dimensional cohomology is a direct sum of an essential finite factorization and a contractible factorization.
\end{corollary}
Let $\MF^{gr}_{R,\omega}$ be the full subcategory in $\MF^{gr,all}_{R,\omega}$ whose objects are matrix factorizations with finite-dimensional cohomology.
Let $\HMF^{gr}_{R,\omega}$ be the full subcategory of $\HMF^{gr,all}_{R,\omega}$ whose object is a matrix factorizations with finite-dimensional cohomology.
Since a contractible matrix factorization is homotopic to the zero factorization in the homotopy category, Corollary \ref{unique-finite} implies the following corollary.
\begin{corollary}[Corollary 5\cite{KR1}]
The homotopy categories $\HMF^{gr,fin}_{R,\omega}$ and $\HMF^{gr}_{R,\omega}$ are categorical equivalent.
\end{corollary}
%
%
\subsection{Morphism of matrix factorization}\label{sec2.5}
\indent
Let $\overline{M}$ and $\overline{N}$ be matrix factorizations with $\omega\in R$.
The set of $\Z$-grading preserving morphisms from $\overline{M}$ to $\overline{N}$ in $\MF^{gr,-}_{R,\omega}$ ($-$ is filled with "all", "fin" or the empty) is denoted by $\Hom_{\MF}(\overline{M},\overline{N})$ and the set of $\Z$-grading preserving morphisms in $\HMF^{gr,-}_{R,\omega}$ is denoted by $\Hom_{\HMF}(\overline{M},\overline{N})$.
The set of all morphisms between the factorizations in $\MF^{gr,-}_{R,\omega}$ is denoted by
\begin{eqnarray}
\nonumber
\HOM_{\MF}(\overline{M},\overline{N})&:=&\bigoplus_{m\in\Z}\Hom_{\MF}(\overline{M}\{m\},\overline{N})
\end{eqnarray}
and the set of all morphisms in $\HMF^{gr,-}_{R,\omega}$ is denoted by
\begin{eqnarray}
\nonumber
\HOM_{\HMF}(\overline{M},\overline{N})&:=&\bigoplus_{m\in\Z}\Hom_{\HMF}(\overline{M}\{m\},\overline{N}).
\end{eqnarray}
When $\overline{N}=\overline{M}$ we denote these sets by $\End_{\MF}(\overline{M})$, $\End_{\HMF}(\overline{M})$, $\END_{\MF}(\overline{M})$ and $\END_{\HMF}(\overline{M})$.\\
\indent
By definition, we have the following properties.
\begin{proposition}
\begin{align}
\nonumber
&\HOM_{\MF}(\overline{M}\{k\},\overline{N})=\HOM_{\MF}(\overline{M},\overline{N})\{-k\},&&
\HOM_{\HMF}(\overline{M}\{k\},\overline{N})=\HOM_{\HMF}(\overline{M},\overline{N})\{-k\},\\
\nonumber
&\HOM_{\MF}(\overline{M},\overline{N}\{k\})=\HOM_{\MF}(\overline{M},\overline{N})\{k\},&&
\HOM_{\HMF}(\overline{M},\overline{N}\{k\})=\HOM_{\HMF}(\overline{M},\overline{N})\{k\}.
\end{align}
\end{proposition}
\indent
The set $\HOM_{\MF}(\overline{M},\overline{N})$ (resp. $\END_{\MF}(\overline{M})$) has an $R$-module structure.
The action of $R$ is defined by $r(f_0,f_1):=(r f_0,r f_1)$ ($r\in R$).
We immediately find the following properties by propositions of $R$-modules in Section \ref{sec2.1}.
\begin{proposition}\label{preserve-finite}
If $\overline{M}$ and $\overline{N}$ are finite factorizations, $\HOM_{\MF}(\overline{M},\overline{N})$ is a free $R$-module of finite rank.
\end{proposition}
\begin{proposition}[Proposition $5$ \cite{KR1}]
For matrix factorizations $\overline{M}$ and $\overline{N}$ in $\HMF^{gr,all}_{R,\omega}$, the action of $R$ on $\HOM_{\HMF}(\overline{M},\overline{N})$ factors through the action of the Jacobi ring $R_{\omega}$.
\end{proposition}
\indent
Since the action $\partial_{x_i}\omega$ $(i=1,...,r)$ is null-homotopic, we have the following proposition.
\begin{proposition}[Proposition 5 \cite{KR1}]
The set $\HOM_{\HMF}(\overline{M},\overline{N})$ $($resp. $\END_{\HMF}(\overline{M})$ $)$ has an $R_{\omega}$-module structure.
\end{proposition}
\subsection{Duality of matrix factorization}\label{sec2.6}
For a matrix factorization $\overline{M}=(M_0,M_1,d_0,d_1)\in\ob(\MF^{gr,all}_{R,\omega})$, we define a matrix factorization $\overline{M}^{\ast}\in\ob(\MF^{gr,all}_{R,-\omega})$ to be $(M_0^{\ast},M_1^{\ast},-d_1^{\ast},d_0^{\ast})$,
\begin{equation}
\nonumber
\overline{M}^{\ast}=
(\xymatrix{
M_0^{\ast}
\ar[r]^{-d_1^{\ast}}
&
M_1^{\ast}
\ar[r]^{d_0^{\ast}}
&
M_0^{\ast}}),
\end{equation}
where $d_{0}^{\ast}f_1:=f_1 d_{0}$ and $d_{1}^{\ast}f_0:=f_0 d_{1}$ ($f_{0}\in M_{0}^{\ast}$, $f_{1}\in M_{1}^{\ast}$).
The factorization $\overline{M}^{\ast}$ is called {\bf dual} of $\overline{M}$.
The following propositions is obvious by Proposition \ref{preserve-finite} and Proposition \ref{duality-mod}.
\begin{proposition}
If $\overline{M}$ is finite, $\overline{M}^{\ast}$ is also finite.
\end{proposition}
\begin{proposition}
If $\overline{M}$ is finite, $\overline{M}^{\ast\ast}\simeq \overline{M}$.
\end{proposition}
We find that the dual of factorization preserves contractible.
\begin{proposition}
If $\overline{M}$ is contractible, $\overline{M}^{\ast}$ is also contractible.
\end{proposition}
Thus, we have a proposition for a matrix factorization with finite-dimensional cohomology.
\begin{proposition}
If $\overline{M}$ is a factorization of $\HMF^{gr}_{R,\omega}$, $\overline{M}^{\ast\ast}\simeq \overline{M}$.
\end{proposition}
\indent
The map from a factorization to the dual of the factorization can be viewed as a contravariant functor between categories of factorizations by the above propositions,
\begin{eqnarray}
\nonumber
?{}^{\ast}&:&\MF_{R,\omega}^{gr,-}\longrightarrow\MF_{R,-\omega}^{gr,-},\\
\nonumber
?{}^{\ast}&:&\HMF_{R,\omega}^{gr,-}\longrightarrow\HMF_{R,-\omega}^{gr,-},
\end{eqnarray}
where $-$ is filled with "all", "fin" or the empty.
Especially, this functor is a categorical equivalence for $\MF_{R,\omega}^{gr,fin}$, $\HMF_{R,\omega}^{gr,fin}$ and $\HMF_{R,\omega}^{gr}$.\\
\indent
These category $\MF_{R,\omega}^{gr,-}$ and $\HMF_{R,\omega}^{gr,-}$ are Krull-Schmidt categories, that is, any matrix factorization has the unique
decomposition property.
\subsection{Two-periodic complex of $R$-module morphisms and extension}\label{sec2.7}
For matrix factorizations $\overline{M}=\left(M_0,M_1,d_{M_0},d_{M_1}\right)$ in $\ob(\MF^{gr,all}_{R,\omega})$ and $\overline{N}=\left(N_0,N_1,d_{N_0},d_{N_1}\right)$ in $\ob(\MF^{gr,all}_{R,\omega'})$, we define a factorization $\HOM_{R}\left(\overline{M},\overline{N}\right)$ of $\MF^{gr,all}_{R,\omega'-\omega}$ to be
\begin{equation}
\nonumber
\left(
\left(
\begin{array}{c}
\HOM_{R}(M_0,N_0)\\
\HOM_{R}(M_1,N_1)
\end{array}
\right),
\left(
\begin{array}{c}
\begin{array}{c}
\HOM_{R}(M_0,N_1)\\
\HOM_{R}(M_1,N_0)
\end{array}
\end{array}
\right),
\left(
\begin{array}{cc}
d_{N_0}&-d_{M_0}^{\ast}\\
-d_{M_1}^{\ast}&d_{N_1}
\end{array}
\right),
\left(
\begin{array}{cc}
d_{N_1}&d_{M_0}^{\ast}\\
d_{M_1}^{\ast}&d_{N_0}
\end{array}
\right)
\right),
\end{equation}
where $d_{N_i}$ $(i\in\Z_2)$ is defined by, for $\sum_m f_m$ $(f_m\in \Hom_R(M_{i}\{m\},N_j)\{m\})$,
\begin{equation}
\nonumber
d_{N_i}(\sum_m f_m):=\sum_m d_{N_i}f_m
\end{equation}
and $d_{M_i}^{\ast}$ $(i\in\Z_2)$ is defined by, for $\sum_m g_m$ $(g_m\in\Hom_R(M_{i-1}\{m\},N_j)\{m\})$,
\begin{equation}
\nonumber
d_{M_i}^{\ast}(\sum_m g_m):=\sum_m g_m\,d_{M_i}.
\end{equation}
\begin{proposition}:
If $\overline{M}$ is a contractible matrix factorization, $\HOM_R(\overline{N},\overline{M})$ and $\HOM_R(\overline{M},\overline{N})$ are also contractible for any factorization $\overline{N}$.
\end{proposition}
Then, the $\HOM_R(\hspace{0.2cm},\hspace{0.2cm})$ is a bifunctor:
\begin{eqnarray}
\nonumber
&&\HOM_R(\hspace{0.2cm},\hspace{0.2cm}):\MF_{R,\omega}^{gr,all}\times\MF_{R,\omega'}^{gr,all}\to\MF_{R,\omega'-\omega}^{gr,all}\,,\\
\nonumber
&&\HOM_R(\hspace{0.2cm},\hspace{0.2cm}):\HMF_{R,\omega}^{gr,all}\times\HMF_{R,\omega'}^{gr,all}\to\HMF_{R,\omega'-\omega}^{gr,all}\,.
\end{eqnarray}
\begin{remark}
Denote the unit object by $\overline{R}=(\xymatrix{R\ar[r]&0\ar[r]&R})$.
For a matrix factorization $\overline{M}$ of $\MF^{gr,all}_{R,\omega}$, we have
$$
\HOM_R(\overline{M},\overline{R})=\overline{M}^{\ast}.
$$
\end{remark}
\indent
For matrix factorizations $\overline{M}$ and $\overline{N}$ with the same potential, the $\HOM_R(\overline{M},\overline{N})$ is a two periodic complex.
Moreover, $\HOM_R(\overline{M},\overline{N})$ is useful for calculus of $\dim_{\Q}\HOM_{\MF}(\overline{M},\overline{N})$ and $\dim_{\Q}\HOM_{\HMF}(\overline{M},\overline{N})$.\\
\indent
If $\overline{M}$ and $\overline{N}$ are objects of $\MF_{R,\omega}^{gr,all}$, $\HOM_R(\overline{M},\overline{N})$ is a two-periodic complex of $\Z$-graded $R$-modules.
Therefore, we can take the bigraded cohomology of $\HOM_R(\overline{M},\overline{N})$ denoted by $\EXT_R(\overline{M},\overline{N})=\EXT_R^0(\overline{M},\overline{N})\oplus \EXT_R^1(\overline{M},\overline{N})$,
where 
\begin{eqnarray}
\nonumber
\EXT_R^0(\overline{M},\overline{N})&=&
\H^0(\,\HOM_R(\overline{M},\overline{N})\,)
,\\
\nonumber
\EXT_R^1(\overline{M},\overline{N})&=&
\H^1(\,\HOM_R(\overline{M},\overline{N})\,)
.
\end{eqnarray}
We have the following isomorphisms as a $\Z$-graded $\Q$-vector space.
\begin{align}
\nonumber
&\HOM_{\MF}(\overline{M},\overline{N})\simeq\cyc^0(\HOM_R(\overline{M},\overline{N})\,)
,&&
\HOM_{\MF}(\overline{M},\overline{N}\left<1\right>)\simeq\cyc^1(\HOM_R(\overline{M},\overline{N})\,)
,&\\
\nonumber
&\HOM_{\HMF}(\overline{M},\overline{N})\simeq\EXT_R^0(\overline{M},\overline{N}),&&
\HOM_{\HMF}(\overline{M},\overline{N}\left<1\right>)\simeq\EXT_R^1(\overline{M},\overline{N}),&
\end{align}
where $\cyc^i(\HOM_R(\overline{M},\overline{N})\,)$ $(i=0,1)$ is the cycle of the two-periodic complex $\HOM_R(\overline{M},\overline{N})$.
\begin{definition}\label{gdim}
For $\overline{M}$ and $\overline{N}$ in $\MF_{R,\omega}^{gr}$, we define 
\begin{equation}
\nonumber
d(q,s):=\gdim\,\EXT_R(\overline{M},\overline{N}).
\end{equation}
We have
\begin{equation}
\nonumber
\left. \frac{1}{k!}\left(\frac{\d}{\d q}\right)^k\right|_{q=0} \hspace{-0.2cm}d(q,s)=\dim_{\Q}\Hom_{\HMF}(\overline{M}\{k\},\overline{N})+s\,\dim_{\Q}\Hom_{\HMF}(\overline{M}\{k\}\langle1\rangle,\overline{N}).
\end{equation}
\end{definition}
%
%
\subsection{Tensor product of matrix factorization}\label{sec2.8}
\indent
$\mathbb{X}=\{x_1,\ldots,x_r\}$ and $\mathbb{Y}=\{y_1,\ldots,y_s\}$ are two sets of variables.
$\mathbb{W}=\{w_1,\ldots,w_t\}$ is the common variables included in $\mathbb{X}$ and $\mathbb{Y}$.
We consider $\Z$-graded rings generated by $\mathbb{X}=\{x_1,\ldots,x_r\}$, $\mathbb{Y}=\{y_1,\ldots,y_s\}$ and $\mathbb{W}=\{w_1,\ldots,w_t\}$, denoted by $R=\Q[\mathbb{X}]$, $R'=\Q[\mathbb{Y}]$ and $S=\Q[\mathbb{W}]$.
We always take a tensor product of $R$ and $R'$ over the ring $S$ generated by the common variables of $R$ and $R'$,
\begin{equation}
\nonumber
R\ostimes R'=R\oqtimes R'/\{rs\oqtimes r'-r\oqtimes sr'\,|\,r\in R,r'\in R',s\in S\}.
\end{equation}
Even if the common variables of $R$ and $R'$ is non-empty, we simply denote $R\ostimes R'$ by the description $R\otimes R'$ without notice.
For an $R$-module $M$ and an $R'$-module $N$, we also take these tensor products over the ring $S$ generated by the set of the common variables of $R$ and $R'$,
\begin{equation}
\nonumber
M\otimes N=M\oqtimes N/\{ms\oqtimes n-m\oqtimes sn\,|\,m\in M,n\in N,s\in S\}.
\end{equation}
\indent
For $\overline{M}=(M_0,M_1,d_{M_0},d_{M_1})$ in $\MF^{gr,all}_{R,\omega}$ and $\overline{N}=(N_0,N_1,d_{N_0},d_{N_1})$ in $\MF^{gr,all}_{R',\omega'}$, 
we define a {\bf tensor product} of matrix factorizations $\overline{M}\boxtimes \overline{N}$ in $\MF^{gr,all}_{R\otimes R',\omega + \omega '}$ by
\begin{eqnarray}
\nonumber
\overline{M}\boxtimes \overline{N} &:=&
\left(
\left(\begin{array}{c}
		M_0\otimes N_0\\
		M_1\otimes N_1
	\end{array}
\right),
\left(\begin{array}{c}
		M_1\otimes N_0\\
		M_0\otimes N_1
	\end{array}
\right),
\left(
	\begin{array}{cc}
		d_{M_0}&-d_{N_1}\\
		d_{N_0}&d_{M_1}
	\end{array}
\right),
\left(
	\begin{array}{cc}
		d_{M_1}&d_{N_1}\\
		-d_{N_0}&d_{M_0}
	\end{array}
\right)\right)
.
\end{eqnarray}
%
As a bifunctor, the tensor product can be viewed
\begin{equation}
\nonumber
\boxtimes:\MF^{gr,all}_{R,\omega}\times\MF^{gr,all}_{R',\omega'}\longrightarrow\MF^{gr,all}_{R\otimes R',\omega+\omega'}.
\end{equation}
\indent
This tensor product $\boxtimes$ has commutativity, additivity and associativity.
Moreover, there exists the unit object for the tensor product.\\
\begin{proposition}\label{com-ass}
{\rm (Commutativity)}For $\overline{M}$ in $\MF^{gr,all}_{R,\omega}$ and $\overline{N}$ in $\MF^{gr,all}_{R',\omega'}$, 
there exists an isomorphism in $\MF^{gr,all}_{R\otimes R',\omega + \omega '}$
$$\overline{M}\boxtimes \overline{N} \simeq \overline{N}\boxtimes \overline{M}.$$ 
{\rm (Additivity)}For $\overline{M}_1$, $\overline{M}_2$ in $\MF^{gr,all}_{R,\omega}$ and $\overline{N}$ in $\MF^{gr,all}_{R',\omega'}$, 
there exists an isomorphism in $\MF^{gr,all}_{R\otimes R',\omega + \omega'}$ 
$$(\overline{M}_1\oplus \overline{M}_2)\boxtimes \overline{N} \simeq \overline{M}_1\boxtimes\overline{N}\oplus \overline{M}_2\boxtimes \overline{N}.$$ 
{\rm (Associativity)}For $\overline{L}$ in $\MF^{gr,all}_{R,\omega}$, $\overline{M}$ in $\MF^{gr,all}_{R',\omega'}$ and $\overline{N}$ in $\MF^{gr,all}_{R '',\omega ''}$, 
there exists an isomorphism between the factorizations in $\MF^{gr,all}_{R\otimes R'\otimes R '',\omega + \omega' +\omega ''}$ 
$$(\overline{L}\boxtimes \overline{M})\boxtimes \overline{N} \simeq \overline{L}\boxtimes (\overline{M}\boxtimes \overline{N}).$$ 
\end{proposition}
%
%
\begin{proof}
See Lemma 2.1, 2.2, 2.7\cite{Yoshi2}.\\
%
\end{proof}
\begin{remark}
As from here, $\overline{M_1}\boxtimes\overline{M_2}\boxtimes\overline{M_3}\boxtimes\ldots\boxtimes\overline{M_k}$ is defined by $(\ldots((\overline{M_1}\boxtimes\overline{M_2})\boxtimes\overline{M_3})\boxtimes\ldots)\boxtimes\overline{M_k}$:
\begin{equation}
\nonumber
\overline{M_1}\boxtimes\overline{M_2}\boxtimes\overline{M_3}\boxtimes\ldots\boxtimes\overline{M_k}
=(\ldots((\overline{M_1}\boxtimes\overline{M_2})\boxtimes\overline{M_3})\boxtimes\ldots)\boxtimes\overline{M_k}.
\end{equation}
\end{remark}
\begin{proposition}\label{identity}
The matrix factorization $(\xymatrix{R \ar[r]^0&0\ar[r]^0&R})$, denoted by $\overline{R}$, is the unit object for the tensor product to any factorization in $\MF^{gr,-}_{R,\omega}$, where $-$ is filled with "all", "fin" or the empty.
In brief, for a matrix factorization $\overline{M} \in \ob(\MF^{gr,-}_{R,\omega})$ we have
$$
\overline{M} \boxtimes \overline{R} \simeq \overline{M}.
$$
\end{proposition}
\begin{proof}
We have this isomorphism by direct calculation.
\end{proof}
We directly find the following isomorphism and identity.
\begin{proposition}
For $\overline{M} \in \ob(\MF^{gr,-}_{R,\omega})$ and $\overline{N} \in \ob(\MF^{gr,-}_{R',\omega'})$, 
there exists an isomorphism in $\MF^{gr,-}_{R\otimes R',\omega + \omega '}$
\begin{eqnarray}
\nonumber
(\overline{M}\boxtimes \overline{N})\left< 1 \right> &=& (\overline{M}\left< 1 \right> )\boxtimes \overline{N}\\[-0.1em]
\nonumber
&\simeq& \overline{M}\boxtimes (\overline{N}\left< 1 \right> ).
\end{eqnarray}
\end{proposition}
\begin{proposition}\label{functor1}
For $\overline{M} \in \ob(\MF^{gr,-}_{R,\omega})$ and $\overline{N} \in \ob(\MF^{gr,-}_{R',\omega'})$, 
there exists an equality in $\MF^{gr,-}_{R\otimes R',\omega + \omega '}$
\begin{eqnarray}
\nonumber
(\overline{M}\boxtimes \overline{N})\{ m\} &=& (\overline{M}\{ m\} )\boxtimes \overline{N}\\[-0.1em]
\nonumber
&=& \overline{M}\boxtimes (\overline{N}\{ m\} )
\end{eqnarray}
\end{proposition}
By Proposition \ref{isom-hom}, we have the following proposition.
\begin{proposition}
If $\overline{M}$ is a finite factorization of $\MF^{gr,fin}_{R,\omega}$ and $\overline{N}$ is a factorization of $\MF^{gr,all}_{R,\omega}$, then there exists an isomorphism as a two-periodic complex
\begin{equation}
\nonumber
\HOM_R(\overline{M},\overline{N})\simeq \overline{N}\boxtimes \overline{M}^{\ast}.
\end{equation}
\end{proposition}
\begin{corollary}\label{isom-ext-cohom}
If $\overline{M}$ is a finite factorization of $\MF^{gr}_{R,\omega}$ and $\overline{N}$ is a factorization of $\MF^{gr}_{R,\omega}$, then there exists an isomorphism as a $\Z\oplus\Z_2$-graded vector space over $\Q$
\begin{equation}
\nonumber
\EXT_R(\overline{M},\overline{N})\simeq \H(\overline{N}\boxtimes \overline{M}^{\ast}).
\end{equation}
\end{corollary}
It is obvious by definition of a contractible matrix factorization that we have the following proposition.
\begin{proposition}
Let $\overline{N}$ be a contractible matrix factorization.
For any matrix factorization $\overline{M}$, the tensor product $\overline{N}\boxtimes \overline{M}$ is also contractible.
\end{proposition}
This proposition means that the tensor product of matrix factorizations can be also defined in the homotopy category $\HMF^{gr,all}$.
\begin{corollary}
$\boxtimes:\HMF^{gr,all}_{R,\omega}\times\HMF^{gr,all}_{R',\omega'}\to\HMF^{gr,all}_{R\otimes R',\omega+\omega'}$ is well-defined.
\end{corollary}
\indent
We consider the special case of the tensor product of two matrix factorizations.
Let $\omega(x_{1},\ldots, x_{i})$, $\omega'(y_{1},\ldots, y_{j})$ and $\omega''(z_{1},\ldots, z_{k})$ be potentials of polynomial rings $R=\Q[x_{1},\ldots, x_{i}]$, $R'=\Q[y_{1},\ldots, y_{j}]$ and $R''=\Q[z_{1},\ldots, z_{k}]$ respectively.
One of matrix factorizations is an object of $\MF_{R\otimes R',\omega-\omega'}^{gr,all}$ denoted by $\overline{M}$.
The other is an object of $\MF_{R'\otimes R'',\omega'-\omega''}^{gr,all}$ denoted by $\overline{N}$.
Their tensor product $\overline{M} \boxtimes \overline{N}$ is an object of $\MF_{R\otimes R'\otimes R'',\omega-\omega''}^{gr,all}$ by definition.
The matrix factorization $\overline{M} \boxtimes \overline{N}$ is also an object of $\MF_{R\otimes R'',\omega-\omega''}^{gr,all}$ since the polynomial $\omega-\omega''$ is a potential of $R\otimes R'$.
Then, we can regard the tensor product as a bifunctor to $\MF^{gr,all}_{R\otimes R'',\omega-\omega''}$ through $\MF^{gr,all}_{R\otimes R'\otimes R'',\omega-\omega''}$
\begin{equation}
\nonumber
\boxtimes:\MF^{gr,all}_{R\otimes R',\omega-\omega'}\times\MF^{gr,all}_{R'\otimes R'',\omega'-\omega''}\to\MF^{gr,all}_{R\otimes R'',\omega-\omega''}.
\end{equation}
Moreover, a contractible factorization of $\MF_{R\otimes R'\otimes R'',\omega-\omega''}^{gr,all}$ is also contractible in $\MF_{R'\otimes R'',\omega-\omega''}^{gr,all}$.
Therefore, we have
\begin{equation}
\nonumber
\boxtimes:\HMF^{gr,all}_{R\otimes R',\omega-\omega'}\times\HMF^{gr,all}_{R'\otimes R'',\omega'-\omega''}\to\HMF^{gr,all}_{R\otimes R'',\omega-\omega''}.
\end{equation}
The tensor product does not preserve finiteness of a matrix factorization.
However, we have the proposition
\begin{proposition}[Proposition 13 \cite{KR1}]\label{preserve-fdc}
If $\overline{M}$ is a factorization of $\MF^{gr}_{R\otimes R',\omega-\omega'}$ and $\overline{N}$ is a factorization of $\MF^{gr}_{R'\otimes R'',\omega'-\omega''}$ (that is, the cohomologies of these factorizations have finitely dimensional), then the tensor product $\overline{M}\boxtimes\overline{N}$ is also a factorization with finite-dimensional cohomology of $\MF^{gr}_{R\otimes R'',\omega-\omega''}$.
\end{proposition}
Therefore, the tensor product preserves finite-dimensional cohomology of a matrix factorization by this proposition.
Thus, we can regard the tensor product as a bifunctor from categories of factorizations with finite-dimensional cohomology to a category of factorization with finite-dimensional cohomology:
\begin{eqnarray}
\nonumber
&\boxtimes&:\MF^{gr}_{R\otimes R',\omega-\omega'}\times\MF^{gr}_{R'\otimes R'',\omega'-\omega''}\to\MF^{gr}_{R\otimes R'',\omega-\omega''},\\
\nonumber
&\boxtimes&:\HMF^{gr}_{R\otimes R',\omega-\omega'}\times\HMF^{gr}_{R'\otimes R'',\omega'-\omega''}\to\HMF^{gr}_{R\otimes R'',\omega-\omega''}.
\end{eqnarray}
%
%
%
%
\subsection{Koszul matrix factorization}\label{sec2.9}
Let $R$ be a $\Z$-graded polynomial ring over $\Q$.
For homogeneous $\Z$-graded polynomials $a$, $b \in R$ and a $\Z$-graded $R$-module $M$, we define a matrix factorization $K(a;b)_{M}$ with the potential $a b$ by
\begin{eqnarray}
\nonumber
K(a;b)_{M}&:=& 
(M,M\{ \frac{1}{2}(\, \deg (b)-\deg (a)\, )\},a,b)\\
\nonumber
&=& 
\Big(
\xymatrix{
*{M}\ar[rr]^(.3){a}&&*{M\{ \frac{1}{2}(\, \deg (b)-\deg (a)\, )\}}\ar[rr]^(.7){b}&&*{M}
}
\Big)
,
\end{eqnarray}
where $\deg (a)$ and $\deg (b)$ are $\Z$-gradings of the homogeneous $\Z$-graded polynomials $a$, $b \in R$.
\\
\indent
In general, for sequences $\mathbf{a}={}^t(a_1, a_2, \ldots, a_k)$, $\mathbf{b}={}^t(b_1, b_2, \ldots, b_k)$ of homogeneous $\Z$-graded polynomials in $R$ and an $R$-module $M$, 
a matrix factorization $K\left( \mathbf{a} ; \mathbf{b} \right)_{M}$ with the potential $\sum_{i=1}^k a_i b_i$ is defined by 
\[
K\left( \mathbf{a} ; \mathbf{b} \right)_{M}
=
\mathop{\boxtimes}_{i=1}^k K(a_i;b_i)_{R}\boxtimes(M,0,0,0)
.
\]
This matrix factorization is called a {\bf Koszul matrix factorization} \cite{KR1}.
We easily find the following propositions by a change of bases of $R$-modules.
\begin{proposition}\label{equiv2}
Let $c$ be a non-zero element in $\Q$.
There is an isomorphism in $\MF^{gr,all}_{R,a b}$
$$
K(a;b)_M \simeq K(c a;c^{-1} b)_M.
$$
\end{proposition}
\begin{proposition}\label{functor2}
\begin{eqnarray}
\nonumber
K(a;b)_{M}\,\left< 1 \right> &=& K(-b;-a)_{M}\,\{ \frac{1}{2}(\, \deg (b) -\deg (a)\, ) \} \\[-0.1em]
\nonumber
&\simeq&K(b;a)_{M}\,\{ \frac{1}{2}(\, \deg (b) -\deg (a)\, ) \} 
\end{eqnarray}
\end{proposition}
The dual of Koszul matrix factorization can be explicitly represented as follows.
\begin{proposition}
Let $M$ be an $R$-module and let $a$ and $b$ be homogeneous $\Z$-graded polynomials of $R$.
We have
\begin{equation}
\nonumber
(K(a;b)_M)^{\ast}\simeq K(-b;a)_{M^{\ast}}.
\end{equation}
When $M=R$ we have
\begin{equation}
\nonumber
(K(a;b)_R)^{\ast}\simeq K(-b;a)_R.
\end{equation}
\end{proposition}
\begin{proposition}[Rasmussen, Section 3.3 in \cite{Ras}]\label{equiv}
Let $a_i$ and $b_i$ be homogeneous $\Z$-graded polynomials such that $\deg (a_1) +\deg (b_1) =\deg (a_2) +\deg (b_2)$ and let $\lambda_i$ {\rm ($i=1,2$)} be homogeneous $\Z$-graded polynomials such that $\deg ( \lambda_1 )=\deg (a_2) - \deg (a_1)$ and $\deg ( \lambda_2 )=-\deg (b_1) + \deg (a_2)$.\\
{\rm (1)} There is an isomorphism in $\MF^{gr,all}_{R,a_1 b_1 + a_2 b_2}$ 
$$
K\left(\left(
\begin{array}{c}
	 a_1\\
	 a_2
\end{array}
\right);
\left(
\begin{array}{c}
	 b_1\\
	 b_2
\end{array}
\right)\right)_{M}
\simeq
K\left(\left(
\begin{array}{c}
	 a_1\\
	 a_2 + \lambda_1 a_1
\end{array}
\right);
\left(
\begin{array}{c}
	 b_1 - \lambda_1 b_2\\
	 b_2
\end{array}
\right)\right)_{M}.
$$
{\rm (2)} There is an isomorphism in $\MF^{gr,all}_{R,a_1 b_1 + a_2 b_2}$ 
$$
K\left(\left(
\begin{array}{c}
	 a_1\\
	 a_2
\end{array}
\right);
\left(
\begin{array}{c}
	 b_1\\
	 b_2
\end{array}
\right)\right)_{M}
\simeq
K\left(\left(
\begin{array}{c}
	 a_1 + \lambda_2 b_2\\
	 a_2 - \lambda_2 b_1
\end{array}
\right);
\left(
\begin{array}{c}
	 b_1\\
	 b_2
\end{array}
\right)\right)_{M}.
$$
\end{proposition}
\indent
We immediately find properties for the bifunctor $\HOM_R(\hspace{.2cm},\hspace{.2cm})$ of Koszul matrix factorization by this proposition.
\begin{proposition}
Let $a$, $b$ and $c$ be non-zero homogeneous $\Z$-graded polynomials of $R$.
We have isomorphisms as a two-periodic complex.
\begin{itemize}
\item[(1)]
\begin{equation}
\nonumber
\HOM_R(K(a;b)_R,K(a;b)_R)
\simeq
K\left(\left(\begin{array}{c}0\\0\end{array}\right);\left(\begin{array}{c}a\\b\end{array}\right)\right)_{R}.
\end{equation}
\item[(2)]
\begin{equation}
\nonumber
\HOM_R(K(a;bc)_R,K(ab;c)_R)
\simeq
K\left(\left(\begin{array}{c}0\\0\end{array}\right);\left(\begin{array}{c}a\\c\end{array}\right)\right)_{R}.
\end{equation}
\item[(3)]
\begin{equation}
\nonumber
\HOM_R(K(ab;c)_R,K(a;bc)_R\{-\deg\,(b)\})
\simeq
K\left(\left(\begin{array}{c}0\\0\end{array}\right);\left(\begin{array}{c}a\\c\end{array}\right)\right)_{R}.
\end{equation}
\end{itemize}
\end{proposition}
By this proposition, we find the following isomorphisms.
\begin{corollary}
We have isomorphisms as a $\Z$-graded $\Q$-vector space.\\
\begin{itemize}
\item[(1)]
\begin{align}
\nonumber
&\HOM_{\MF}(K(a;b)_R,K(a;b)_R)\simeq R,&
&\HOM_{\MF}(K(a;b)_R,K(a;b)_R\langle 1 \rangle)\simeq 0,&\\
\nonumber
&\HOM_{\HMF}(K(a;b)_R,K(a;b)_R)\simeq R\left/\langle a,b \rangle\right.,&
&\HOM_{\HMF}(K(a;b)_R,K(a;b)_R\langle 1 \rangle)\simeq 0.&
\end{align}
\item[(2)]
\begin{align}
\nonumber
&\HOM_{\MF}(K(a;bc)_R,K(ab;c)_R)\simeq R,&
&\HOM_{\MF}(K(a;bc)_R,K(ab;c)_R\langle 1 \rangle)\simeq 0,&\\
\nonumber
&\HOM_{\HMF}(K(a;bc)_R,K(ab;c)_R)\simeq R\left/\langle a,c \rangle\right.,&
&\HOM_{\HMF}(K(a;bc)_R,K(ab;c)_R\langle 1 \rangle)\simeq 0.&
\end{align}
\item[(3)]
\begin{align}
\nonumber
&\HOM_{\MF}(K(ab;c)_R,K(a;bc)_R\{-\deg\,(b)\})\simeq R,&
&\HOM_{\MF}(K(ab;c)_R,K(a;bc)_R\{-\deg\,(b)\}\langle 1 \rangle)\simeq 0,&\\
\nonumber
&\HOM_{\HMF}(K(ab;c)_R,K(a;bc)_R\{-\deg\,(b)\})\simeq R\left/\langle a,c \rangle\right.,&
&\HOM_{\HMF}(K(ab;c)_R,K(a;bc)_R\{-\deg\,(b)\}\langle 1 \rangle)\simeq 0.&
\end{align}
\end{itemize}
\end{corollary}
We find a dimension of $\Hom_{\MF}$ and $\Hom_{\HMF}$ as a $\Q$-vector space by this corollary and Proposition \ref{gdim}.
\begin{corollary}\label{hom-dim}
We find dimension of a $\Q$-vector space of $\Z$-grading preserving morphisms between factorizations.
\begin{itemize}
\item[(1)]
\begin{align}
\nonumber
&\dim_{\Q}\, \Hom_{\MF}(K(a;b)_R,K(a;b)_R)=1,&
&\dim_{\Q}\, \Hom_{\HMF}(K(a;b)_R,K(a;b)_R)=1.&
\end{align}
\item[(2)]
\begin{align}
\nonumber
&\dim_{\Q}\, \Hom_{\MF}(K(a;bc)_R,K(ab;c)_R)=1,&
&\dim_{\Q}\, \Hom_{\HMF}(K(a;bc)_R,K(ab;c)_R)=1.&
\end{align}
\item[(3)]
\begin{align}
\nonumber
&\dim_{\Q}\, \Hom_{\MF}(K(ab;c)_R,K(a;bc)_R\{-\deg\,(b)\})=1,&
&\dim_{\Q}\, \Hom_{\HMF}(K(ab;c)_R,K(a;bc)_R\{-\deg\,(b)\})=1.&
\end{align}
\end{itemize}
\end{corollary}
\begin{theorem}[Khovanov-Rozansky, Theorem 2.1 \cite{KR3}]\label{reg-eq}
Let $a_i$, $b_i$ and $b_i'$ $(i=1,\ldots ,m)$ be homogeneous $\Z$-graded polynomials in $R$ and let $M$ be an $R$-module.
If $a_1,\ldots ,a_m$ $\in R$ form a regular sequence and
\begin{equation}
\nonumber
\sum_{i=1}^m a_i b_i=\sum_{i=1}^m a_i b_i' \,(=:\omega),
\end{equation}
there exists an isomorphism in $\MF^{gr,all}_{R,\omega}$
\begin{equation}
\nonumber
K\left(
\left(
\begin{array}{c}
a_1\\
\vdots\\
a_m
\end{array}
\right);
\left(
\begin{array}{c}
b_1\\
\vdots\\
b_m
\end{array}
\right)
\right)_M
\simeq
K\left(
\left(
\begin{array}{c}
a_1\\
\vdots\\
a_m
\end{array}
\right);
\left(
\begin{array}{c}
b_1'\\
\vdots\\
b_m'
\end{array}
\right)
\right)_M.
\end{equation}
\end{theorem}
\begin{corollary}\label{induce-sq1}
Put $R=\Q[x_1,x_2,\ldots ,x_k]$ and $R_y=R[y]\left/\left<y^l+\alpha_1y^{l-1} +\alpha_2y^{l-2} +\ldots +\alpha_l\right>\right.$, where $\alpha_i\in R$ such that $\deg (\alpha_i)=i \,\deg(y)$.\\
$(1)$Let $a_i$ be a homogeneous $\Z$-graded polynomial $\in R_y$ $(i=1,\ldots ,m)$, $b_i$ be a homogeneous $\Z$-graded polynomial $\in R$ $(i=2,\ldots ,m)$ and let $b_1$, $\beta$ be homogeneous $\Z$-graded polynomials $\in R_y$ with the property $(y+\beta)b_1\in R$.
If these polynomials hold the following conditions:
\begin{itemize}
\item[(i)]$(y+\beta)b_1$, $b_2$, $\ldots$, $b_m$ form a regular sequence in $R$,
\item[(ii)] $\displaystyle a_1b_1(y+\beta)+\sum_{i=2}^m a_ib_i$ $(=:\omega)$ $\in R$,
\end{itemize}
then there exist homogeneous $\Z$-graded polynomials $a_i' \in R$ $(i=1,\ldots ,m)$ and we have an isomorphism in $\MF_{R,\omega}^{gr,all}$
\begin{equation}
\nonumber
K\left(\left(
\begin{array}{c}
(y+\beta)a_1\\
a_2\\
\vdots\\
a_m
\end{array}
\right);
\left(
\begin{array}{c}
b_1\\
b_2\\
\vdots\\
b_m
\end{array}
\right)\right)_{R_y}
\simeq
K\left(\left(
\begin{array}{c}
(y+\beta)a_1'\\
a_2'\\
\vdots\\
a_m'
\end{array}
\right);
\left(
\begin{array}{c}
b_1\\
b_2\\
\vdots\\
b_m
\end{array}
\right)\right)_{R_y}.
\end{equation}
\\
\noindent
$(2)$Let $a_i$ be a homogeneous $\Z$-graded polynomial $\in R_y$ $(i=1,\ldots ,m)$, $b_i$ be a homogeneous $\Z$-graded polynomial $\in R$ $(i=1,\ldots ,m)$ and $\beta$ be a homogeneous $\Z$-graded polynomial $\in R$.
If these polynomials hold the following conditions:
\begin{itemize}
\item[(i)]$b_1$, $b_2$, $\ldots$, $b_m$ form a regular sequence in $R$,
\item[(ii)] $\displaystyle a_1b_1(y+\beta)+\sum_{i=2}^m a_ib_i$ $(=:\omega')$ $\in R$,
\end{itemize}
then there exist homogeneous $\Z$-graded polynomials $a_1' \in R_y$ and $a_i' \in R$ $(i=2,\ldots ,m)$ and we have an isomorphism in $\MF_{R,\omega'}^{gr,all}$
\begin{equation}
\nonumber
K\left(\left(
\begin{array}{c}
a_1\\
a_2\\
\vdots\\
a_m
\end{array}
\right);
\left(
\begin{array}{c}
b_1(y+\beta)\\
b_2\\
\vdots\\
b_m
\end{array}
\right)\right)_{R_y}
\simeq
K\left(\left(
\begin{array}{c}
a_1'\\
a_2'\\
\vdots\\
a_m'
\end{array}
\right);
\left(
\begin{array}{c}
b_1(y+\beta)\\
b_2\\
\vdots\\
b_m
\end{array}
\right)\right)_{R_y}.
\end{equation}
\end{corollary}
\begin{proof}
This corollary can be proved by using Theorem \ref{reg-eq} and the relation of the quotient $R_y$.
\end{proof}
\begin{remark}
Put $R_y=R[y]\left/\left<y^l+\alpha_1y^{l-1} +\alpha_2y^{l-2} +\ldots +\alpha_l\right>\right.$ $(\alpha_i\in R)$.
If the variable $y$ remains in a homogeneous $\Z$-graded polynomial $p$ then a matrix form of $p$ to $R_y$ is complicated as an $R$-module morphism.
However, if the variable $y$ dose not exist in a polynomial $p$ then a matrix form of $p$ is simply a diagonal map as an $R$-module morphism,
\begin{equation}
\nonumber
(\xymatrix{
R_y\ar[r]^{p}&R_y
}
)
=
(
\xymatrix{
\hbox{$
\left(
\begin{array}{r}
\beta_0 R\\
\beta_1 R\\
\vdots\\
\beta_{l-1} R
\end{array}
\right)
$}
\ar[rrrr]_{
\hbox{$
\left(
\begin{array}{cccc}
p&0&\cdots&0\\
0&p&&0\\
\vdots&&\ddots&\vdots\\
0&0&\cdots&p
\end{array}
\right)
$}
}&&&&
\hbox{
$\left(
\begin{array}{r}
\beta_0 R\\
\beta_1 R\\
\vdots\\
\beta_{l-1} R
\end{array}
\right)
$}
}
).
\end{equation}
where $\beta_0$, $\beta_1$, $\ldots$, $\beta_{l-1}$ form a basis of $R_y$ as an $R$-module.
\end{remark}
In the next section, Corollary \ref{induce-sq1} is useful for decomposing a matrix factorization into a direct sum of matrix factorizations.\\
\indent
Theorem \ref{exclude} is a generalization to multivariable of Theorem 2.2 given by Khovanov and Rozansky \cite{KR3}.
\begin{theorem}[Generalization of Khovanov-Rozansky's Theorem 2.2 \cite{KR3}]\label{exclude}
We put $R=\Q [\underline{x}]$, where $\underline{x}=(x_1,x_2,\ldots ,x_l)$.
Let $a_i$ and $b_i$ $(1\leq i\leq k)$ be homogeneous $\Z$-graded polynomials in $R [\underline{y}]$, where $\underline{y}=(y_1,y_2,\ldots ,y_m)$ and let $M$ be an $R [\underline{y}]$-module.
We suppose that sequences $\mathbf{a}={}^t( a_1 , a_2 , \ldots , a_k )$ and $\mathbf{b}={}^t( b_1 , b_2 , \ldots , b_k )$ satisfy the conditions
\begin{itemize}
\item[(i)] $\sum_{i=1}^{k}a_i b_i \,(=:\omega) \in R$,\\
\item[(ii)] There exists $j$ such that $b_j=c y_1^{m_1}y_2^{m_2}\ldots y_l^{m_l} + p$, where 
$c$ is a non-zero constant and $p \in R [\underline{y}]$ does not include the monomial $y_1^{m_1}y_2^{m_2}\ldots y_l^{m_l}$.
\end{itemize}
Then, there exists an isomorphism in $\HMF^{gr,all}_{R,\omega }$,
$$
K(\mathbf{a};\mathbf{b})_{M} \simeq K(\stackrel{j}{\check{\mathbf{a}}};\stackrel{j}{\check{\mathbf{b}}})_{M/b_j M},
$$
where $\stackrel{j}{\check{\mathbf{a}}}$ and $\stackrel{j}{\check{\mathbf{b}}}$ are the sequences omitted the $j$-th entry of $\mathbf{a}$ and $\mathbf{b}$.
\end{theorem}
\begin{proof}
This theorem is proved by a similar way to the proof of Khovanov-Rozansky's Theorem 2.2 \cite{KR3}.
\end{proof}
\begin{corollary}\label{cor2-10}
We put $R=\Q [\underline{x}]$, where $\underline{x}=(x_1,x_2,\ldots ,x_l)$.
Let $a_i$ and $b_i$ $(1\leq i\leq k)$ be homogeneous $\Z$-graded polynomials in $R [\underline{y}]$, where $\underline{y}=(y_1,y_2,\ldots ,y_m)$ and let $M$ be an $R [\underline{y}]$-module.
We suppose that sequences $\mathbf{a}={}^t( a_1 , a_2 , \ldots , a_k )$ and $\mathbf{b}={}^t( b_1 , b_2 , \ldots , b_k )$ satisfy the conditions
\begin{itemize}
\item[(i)] $\sum_{i=1}^{k}a_i b_i \,(=:\omega) \in R$,\\
\item[($\ast$)] There exists $j$ such that a homogeneous $\Z$-graded polynomial $b_j(\underline{x},\underline{y})\in R[\underline{y}]$ satisfies  $b_j(\underline{0},\underline{y})\not=0$.
\end{itemize}
Then, there exists an isomorphism in $\HMF^{gr,all}_{R,\omega }$,
$$
K(\mathbf{a};\mathbf{b})_{M} \simeq K(\stackrel{j}{\check{\mathbf{a}}};\stackrel{j}{\check{\mathbf{b}}})_{M/b_jM}.
$$
\end{corollary}
\begin{proof}
Each monomial of $b_j(\underline{0},\underline{y})$ forms $y_1^{i_1}y_2^{i_2}\ldots y_m^{i_m}$.
Then, $b_j(\underline{x},\underline{y})$ satisfies the condition (ii) of Theorem \ref{exclude}.
Thus, we obtain this corollary.
\end{proof}
\begin{corollary}\label{cor2-11}
We put $R=\Q [\underline{x}]$, where $\underline{x}=(x_1,x_2,\ldots ,x_l)$.
Let $a_i$ and $b_i$ $(1\leq i\leq k)$ be homogeneous $\Z$-graded polynomials in $R [\underline{y}]$, where $\underline{y}=(y_1,y_2,\ldots ,y_m)$ and let $M$ be an $R [\underline{y}]$-module.
We suppose that sequences $\mathbf{a}={}^t( a_1 , a_2 , \ldots , a_k )$ and $\mathbf{b}={}^t( b_1 , b_2 , \ldots , b_k )$ satisfy the conditions
\begin{itemize}
\item[(i)] $\sum_{i=1}^{k}a_i b_i \,(=:\omega) \in R$,\\
\item[(ii)] There are homogeneous $\Z$-graded polynomials $b_{j_1}(\underline{x},\underline{y}),b_{j_2}(\underline{x},\underline{y}),\ldots ,b_{j_r}(\underline{x},\underline{y})\in R[\underline{y}]$ such that the sequence $(b_{j_1}(\underline{0},\underline{y}),b_{j_2}(\underline{0},\underline{y}),\ldots ,b_{j_r}(\underline{0},\underline{y}))$ is regular in $\Q[\underline{y}]$,
\end{itemize}
then there exists an isomorphism in $\HMF^{gr,all}_{R,\omega }$,
$$
K(\mathbf{a};\mathbf{b})_{M}
\simeq K(\stackrel{j_1,j_2,\ldots ,j_r}{\check{\mathbf{a}}};\stackrel{j_1,j_2,\ldots ,j_r}{\check{\mathbf{b}}})_{M/\left< b_{j_1},b_{j_2},\ldots ,b_{j_r}\right> M}.
$$
\end{corollary}
\begin{proof}
The sequence $(b_{j_1}(\underline{0},\underline{y}),b_{j_2}(\underline{0},\underline{y}),\ldots ,b_{j_r}(\underline{0},\underline{y}))$ is regular by the assumption.
Then, by applying Corollary \ref{cor2-10} to the polynomial $b_{j_r}(\underline{x},\underline{y})$, the sequences $\stackrel{j_r}{\check{\mathbf{a}}}$ and $\stackrel{j_r}{\check{\mathbf{b}}}$ still satisfy the conditions (i) and (ii).
Thus, we can prove this corollary repeating this operation.
\end{proof}
%
%
%
%
\subsection{Complex category over a graded additive category}\label{sec2.10}
\indent
In general, for a graded additive category $\mathcal{A}$, we can define the complex category over $\mathcal{A}$ and its homotopy category.
Moreover, when $\mathcal{A}$ has tensor product structure we can define tensor product in the complex category.
\begin{definition}
Let $\mathcal{A}$ be a graded additive category.
The category of complexes bounded below and above over $\mathcal{A}$, denoted by $\kom^b(\mathcal{A})$, is defined as follow.
\begin{itemize}
\item
An object of $\kom^b(\mathcal{A})$ forms
\begin{equation}
\nonumber
X^{\bullet}:\xymatrix{\ldots\ar[r]^{{d_c}_{X^{i-2}}}&X^{i-1}\ar[r]^{{d_c}_{X^{i-1}}}&X^i\ar[r]^{{d_c}_{X^i}}&X^{i+1}\ar[r]^{{d_c}_{X^{i+1}}}&\ldots},
\end{equation}
where $X^i$ is an object of $\mathcal{A}$ for any $i$, $X^i=0$ for $i\ll0$, $i\gg0$ and the boundary map ${d_c}_{X^i}$ is a $\Z$-grading preserving morphism such that ${d_c}_{X^{i+1}}{d_c}_{X^i}=0$ in $\mor(\mathcal{A})$ for any $i$.
\item
For objects $X^{\bullet}$ and $Y^{\bullet}$, a morphism from $X^{\bullet}$ to $Y^{\bullet}$ is a collection $(\ldots,f^{i-1},f^{i},f^{i+1},\ldots)$, denoted by $f^{\bullet}$, of a $\Z$-grading preserving morphism $f^i$ of $\Hom_{\mathcal{A}}(X^i,Y^i)$ such that ${d_c}_{Y^i}f^i=f^{i+1}{d_c}_{X^i}$ for every $i$:
\begin{equation}
\nonumber
\xymatrix{
\ldots\ar[r]^{{d_c}_{X^{i-2}}}&X^{i-1}\ar[d]^{f^{i-1}}\ar[r]^{{d_c}_{X^{i-1}}}&X^i\ar[d]^{f^{i}}\ar[r]^{{d_c}_{X^i}}&X^{i+1}\ar[d]^{f^{i+1}}\ar[r]^{{d_c}_{X^{i+1}}}&\ldots\\
\ldots\ar[r]_{{d_c}_{Y^{i-2}}}&Y^{i-1}\ar[r]_{{d_c}_{Y^{i-1}}}&Y^i\ar[r]_{{d_c}_{Y^i}}&Y^{i+1}\ar[r]_{{d_c}_{Y^{i+1}}}&\ldots\hspace{1mm}.
}
\end{equation}
The set of morphisms from $X^{\bullet}$ to $Y^{\bullet}$ is denoted by $\Hom_{\mathcal{A}}(X^{\bullet},Y^{\bullet})$ and $\Hom_{\mathcal{A}}(X^{\bullet},X^{\bullet})$ is denoted by $\End_{\mathcal{A}}(X^{\bullet})$ for short.
\item
The composition of morphisms $f^{\bullet}g^{\bullet}$ is defined by $(\ldots,f^{i-1}g^{i-1},f^{i}g^{i},f^{i+1}g^{i+1},\ldots)$.
\end{itemize}
\end{definition}
\indent
We define complex null-homotopic in $\kom^b(\mathcal{A})$.
A morphism $f^{\bullet}:X^{\bullet}\to Y^{\bullet}$ is {\bf complex null-homotopic} if a collection ${h_c}^{\bullet}=(\ldots,{h_c}^{i-1},{h_c}^{i},{h_c}^{i+1},\ldots)$ of $\Z$-grading preserving morphisms ${h_c}^i:X^i\to Y^{i-1}$ exists such that $f^i={h_c}^{i+1}{d_c}_{X^i}+{d_c}_{Y^{i-1}}{h_c}^i$ in $\mor(\mathcal{A})$ for every $i$:
\begin{equation}
\nonumber
\xymatrix{
\ldots\ar[r]^{{d_c}_{X^{i-2}}}&X^{i-1}\ar[dl]|-{{h_c}^{i-1}}\ar[d]|-{f^{i-1}}\ar[r]^{{d_c}_{X^{i-1}}}&X^i\ar[dl]|-{{h_c}^{i}}\ar[d]|-{f^{i}}\ar[r]^{{d_c}_{X^i}}&X^{i+1}\ar[dl]|-{{h_c}^{i+1}}\ar[d]|-{f^{i+1}}\ar[r]^{{d_c}_{X^{i+1}}}&\ldots\\
\ldots\ar[r]_{{d_c}_{Y^{i-2}}}&Y^{i-1}\ar[r]_{{d_c}_{Y^{i-1}}}&Y^i\ar[r]_{{d_c}_{Y^i}}&Y^{i+1}\ar[r]_{{d_c}_{Y^{i+1}}}&\ldots\hspace{1mm}.
}
\end{equation}
Morphisms $f^{\bullet},g^{\bullet}:X^{\bullet}\to Y^{\bullet}$ is {\bf complex homotopic}, denoted by $f^{\bullet}\stackrel{cpx}{\sim} g^{\bullet}$, if $f^{\bullet}-g^{\bullet}$ is complex null-homotopic.
\begin{definition}
The homotopy category of $\kom^b(\mathcal{A})$, denoted by $\k^b(\mathcal{A})$, is defined as follow.
\begin{itemize}
\item
$\ob (\k^b(\mathcal{A}))=\ob (\kom^b(\mathcal{A}))$,
\item
$\mor (\k^b(\mathcal{A}))=\mor (\kom^b(\mathcal{A})/\{ complex \,\, null$-$homotopic\}$.
\item
The composition of morphisms is defined as the same in $\kom^b(\mathcal{A})$.
\end{itemize}
\end{definition}
The {\bf complex translation functor}\footnote{This definition of complex translation functor is different from the ordinary definition $(X^{\bullet}[k])^i=X^{i+k}$. This definition matches with Poincar\'e polynomial $P(D)$ of $\Z\oplus\Z\oplus\Z_2$-graded homology, see Remark \ref{rem-shift}.} $[k]$ ($k\in\Z$) changes a complex $X^{\bullet}$ into
$$
(X^{\bullet}[k])^i=X^{i-k}.
$$
\begin{definition}
We assume that a category $\mathcal{A}$ has tensor product structure,
$$
\otimes:\mathcal{A}\times\mathcal{A}\to\mathcal{A}.
$$
For complexes $X^{\bullet}$ and $Y^{\bullet}$, we define $X^{\bullet}\otimes Y^{\bullet}$ to be
\begin{eqnarray}
\nonumber
(X^{\bullet}\otimes Y^{\bullet})_k:=\bigoplus_{i+j=k}X^i\otimes Y^j,
\hspace{2mm}
{d_c}_{(X^{\bullet}\otimes Y^{\bullet})_k}=\sum_{i+j=k}({d_c}_{X^i}\otimes\id_{Y^j}+(-1)^{i}\id_{X^i}\otimes {d_c}_{Y^j}).
\end{eqnarray}
\end{definition}
%
%
%
%
\subsection{Complex category of a $\Z$-graded matrix factorizations}\label{sec2.11}
\indent
$\HMF^{gr,-}_{R,\omega}$ ($-$ is filled with "all", "fin" or the empty) is a graded additive category with  tensor product $\boxtimes$.
We consider the complex category $\kom^b(\HMF^{gr,-}_{R,\omega})$ and the homotopy category $\k^b(\HMF^{gr,-}_{R,\omega})$.
Moreover, we consider a full subcategory of $\HMF^{gr,-}_{R,\omega}$ whose objects are essential factorizations, denoted by $\HMF^{gr,-,es}_{R,\omega}$.
We also consider the complex category $\kom^b(\HMF^{gr,-,es}_{R,\omega})$ and the homotopy category $\k^b(\HMF^{gr,-,es}_{R,\omega})$.\\
\indent
We find the following proposition by Corollary \ref{unique-decomp}.
\begin{corollary}
$\HMF^{gr,-}_{R,\omega}$ and $\HMF^{gr,-,es}_{R,\omega}$ are categorical equivalent.
\end{corollary}
\indent
By Corollary \ref{unique-decomp}, we know that a matrix factorization $\overline{M}\in\ob(\MF^{gr,-}_{R,\omega})$ is a direct sum of an essential factorization and a contractible factorization $\overline{M}_{es}\oplus\overline{M}_{c}$.
Then, we also find that $\overline{M}^{\bullet}$ in $\ob(\kom^{b}(\HMF^{gr,-}_{R,\omega}))$ is describe by the following complex
\begin{equation}
\nonumber
\xymatrix{
\ldots\ar[rr]_(.45){
\left(
\begin{array}{cc}
d_{\overline{M}^{i-2}_{es}}&\ast\\
\ast&\ast
\end{array}
\right)}
&&
{\begin{array}{c}
\overline{M}^{i-1}_{es}\\
\oplus\\
\overline{M}^{i-1}_{c}
\end{array}}
\ar[rr]_{
\left(
\begin{array}{cc}
d_{\overline{M}^{i-1}_{es}}&\ast\\
\ast&\ast
\end{array}
\right)}
&&
{\begin{array}{c}
\overline{M}^{i}_{es}\\
\oplus\\
\overline{M}^{i}_{c}
\end{array}}
\ar[rr]_{
\left(
\begin{array}{cc}
d_{\overline{M}^{i}_{es}}&\ast\\
\ast&\ast
\end{array}
\right)}
&&
{\begin{array}{c}
\overline{M}^{i+1}_{es}\\
\oplus\\
\overline{M}^{i+1}_{c}
\end{array}}
\ar[rr]_(.55){
\left(
\begin{array}{cc}
d_{\overline{M}^{i+1}_{es}}&\ast\\
\ast&\ast
\end{array}
\right)}&&\ldots
}.
\end{equation}
Entries denoted by $\ast$ of boundary morphisms are null-homotopic since any morphism from a contractible factorization or to a contractible factorization is null-homotopic.
Then, this complex is isomorphic in $\kom^{b}(\HMF^{gr,-}_{R,\omega})$ to
\begin{equation}
\nonumber
\overline{M}^{\bullet}_{es}=\xymatrix{
\ldots
\ar[rr]^{d_{\overline{M}^{i-2}_{es}}}
&&
\overline{M}^{i-1}_{es}
\ar[rr]^{d_{\overline{M}^{i-1}_{es}}}
&&
\overline{M}^{i}_{es}
\ar[rr]^{d_{\overline{M}^{i}_{es}}}
&&
\overline{M}^{i+1}_{es}
\ar[rr]^{d_{\overline{M}^{i+1}_{es}}}
&&\ldots
}.
\end{equation}
For a morphism $f^{\bullet}:\overline{M}^{\bullet}\to\overline{N}^{\bullet}$, we have
\begin{equation}
\nonumber
\xymatrix{\overline{M}^{\bullet}_{es}\ar[d]^{f^{\bullet}_{es}}&
\ldots
\ar[rr]^{d_{\overline{M}^{i-2}_{es}}}
&&
\overline{M}^{i-1}_{es}\ar[d]^{f^{i-1}_{es}}
\ar[rr]^{d_{\overline{M}^{i-1}_{es}}}
&&
\overline{M}^{i}_{es}\ar[d]^{f^{i}_{es}}
\ar[rr]^{d_{\overline{M}^{i}_{es}}}
&&
\overline{M}^{i+1}_{es}\ar[d]^{f^{i+1}_{es}}
\ar[rr]^{d_{\overline{M}^{i+1}_{es}}}
&&\ldots
\\
\overline{N}^{\bullet}_{es}&\ldots
\ar[rr]^{d_{\overline{N}^{i-2}_{es}}}
&&
\overline{N}^{i-1}_{es}
\ar[rr]^{d_{\overline{N}^{i-1}_{es}}}
&&
\overline{N}^{i}_{es}
\ar[rr]^{d_{\overline{N}^{i}_{es}}}
&&
\overline{N}^{i+1}_{es}
\ar[rr]^{d_{\overline{N}^{i+1}_{es}}}
&&\ldots
}.
\end{equation}
The map from $\overline{M}^{\bullet}$ to $\overline{M}^{\bullet}_{es}$ and from $\overline{f}^{\bullet}$ to $\overline{f}^{\bullet}_{es}$ is functorial from $\kom^b(\HMF^{gr,-}_{R,\omega})$ to $\kom^b_{es}(\HMF^{gr,-}_{R,\omega})$, denoted by $\mathcal{ES}$.
\begin{proposition}
The functor $\mathcal{ES}$ is a categorical equivalence.
\end{proposition}
\begin{proof}
We need to show that $(a)$ any object $\overline{N}$ of $\kom^b_{es}(\HMF^{gr,-}_{R,\omega})$ is isomorphic to $\mathcal{ES}(\overline{M})$ for some objects $\overline{M}$ of $\kom^b(\HMF^{gr,-}_{R,\omega})$ and $(b)$ the functor $\mathcal{ES}$ is fully faithful.
The fact $(a)$ is satisfied by definition of $\mathcal{ES}$.
$(b)$ The functor $\mathcal{ES}$ is a full functor since $\kom^b_{es}(\HMF^{gr,-}_{R,\omega})$ is the full subcategory of $\kom^b(\HMF^{gr,-}_{R,\omega})$.
$\mathcal{ES}$ is a faithful since any morphism from a contractible factorization or to a contractible factorization is null homotopic.
\end{proof}
\indent
The definition of the tensor product naturally adjusts to the category $\kom^b(\HMF^{gr,all})$ and $\kom^b(\HMF^{gr})$.
We also denote the tensor product in $\kom^b(\HMF^{gr,-})$ by $\boxtimes$;
\begin{eqnarray}
\nonumber
\boxtimes&:&\kom^b(\HMF_{R,\omega}^{gr,all})\times\kom^b(\HMF_{R',\omega'}^{gr,all})\longrightarrow\kom^b(\HMF_{R\otimes R',\omega+\omega'}^{gr,all}).
\end{eqnarray}
By Proposition \ref{preserve-fdc}, we also obtain bifunctors
\begin{eqnarray}
\nonumber
\boxtimes&:&\kom^b(\HMF^{gr,all}_{R\otimes R',\omega-\omega'})\times\kom^b(\HMF^{gr,all}_{R'\otimes R'',\omega'-\omega''})\to\kom^b(\HMF^{gr,all}_{R\otimes R'',\omega-\omega''}),\\
\nonumber
\boxtimes&:&\kom^b(\HMF^{gr}_{R\otimes R',\omega-\omega'})\times\kom^b(\HMF^{gr}_{R'\otimes R'',\omega'-\omega''})\to\kom^b(\HMF^{gr}_{R\otimes R'',\omega-\omega''}).
\end{eqnarray}
\indent
Finally, we show a proposition for Koszul factorization.
\begin{proposition}
Let $a_i$, $a_i'$ and $b_i$ $(i=1,\ldots,k)$ be sequences of homogeneous $\Z$-graded polynomials in $R$ such that 
\begin{itemize}
\item[(1)]$\displaystyle c a_1 b_1+\sum_{i=2}^{k}a_i b_i=c a_1' b_1+\sum_{i=2}^{k} a_i' b_i$, where $c$ is a homogeneous $\Z$-graded polynomial of $R$,
\item[(2)]$\mathbf{b}$ is a regular sequence.
\end{itemize}
Put $\overline{S}=K((a_2,\ldots,a_k);(b_2,\ldots,b_k))_R$ and $\overline{S}'=K((a_2',\ldots,a_k');(b_2,\ldots,b_k))_R$.
By Corollary \ref{induce-sq1}, we have isomorphisms
\begin{eqnarray}
\nonumber
&&\xymatrix{K(c a_1;b_1)_R\boxtimes\overline{S}\ar[r]^{\overline{\varphi}}&K(c a_1';b_1)_R\boxtimes\overline{S}'},\\
\nonumber
&&\xymatrix{K(a_1;c b_1)_R\boxtimes\overline{S}\ar[r]^{\overline{\psi}}&K(a_1';c b_1)_R\boxtimes\overline{S}'}.
\end{eqnarray}
$(1)$ We have the $\Z$-grading preserving morphisms between matrix factorizations
\begin{eqnarray}
\nonumber
&&\xymatrix{K(c a_1;b_1)_R\boxtimes\overline{S}\ar[rr]^(.45){(c,1)\boxtimes\id_{\overline{S}}}&&K(a_1;c b_1)_R\boxtimes\overline{S}\{-\deg c\}},\\
\nonumber
&&\xymatrix{K(c a_1;b_1)_R\boxtimes\overline{S}'\ar[rr]^(.45){(c,1)\boxtimes\id_{\overline{S}'}}&&K(a_1;c b_1)_R\boxtimes\overline{S}'\{-\deg c\}}.
\end{eqnarray}
Then, this morphism satisfies the condition
\begin{equation}
\nonumber
\left((c,1)\boxtimes\id_{\overline{S}}\right)\,\cdot\overline{\psi}=\overline{\varphi }\cdot \left((c,1)\boxtimes\id_{\overline{S}'}\right).
\end{equation}
That is, these natural morphisms between matrix factorization give the commute diagram
\begin{equation}
\nonumber
\xymatrix{
K(c a_1;b_1)_R\boxtimes\overline{S}\ar[d]^{\overline{\varphi}}\ar[rr]^(.45){(c,1)\boxtimes\id_{\overline{S}}}&&K(a_1;c b_1)_R\boxtimes\overline{S}\{-\deg c\}\ar[d]^{\overline{\psi}}\\
K(c a_1';b_1)_R\boxtimes\overline{S}'\ar[rr]^(.45){(c,1)\boxtimes\id_{\overline{S}'}}&&K(a_1';c b_1)_R\boxtimes\overline{S}'\{-\deg c\}.
}
\end{equation}
$(2)$ We have the $\Z$-grading preserving morphisms between matrix factorizations
\begin{eqnarray}
\nonumber
&&\xymatrix{K(a_1;c b_1)_R\boxtimes\overline{S}\ar[rr]^{(1,c)\boxtimes\id_{\overline{S}}}&&K(c a_1;b_1)_R\boxtimes\overline{S}},\\
\nonumber
&&\xymatrix{K(a_1;c b_1)_R\boxtimes\overline{S}'\ar[rr]^{(1,c)\boxtimes\id_{\overline{S}'}}&&K(c a_1;b_1)_R\boxtimes\overline{S}'}.
\end{eqnarray}
Then, this morphism satisfies the condition
\begin{equation}
\nonumber
 \left((1,c)\boxtimes\id_{\overline{S}}\right)\,\cdot\overline{\varphi}=\overline{\psi}\cdot \left((1,c)\boxtimes\id_{\overline{S}'}\right).
\end{equation}
That is, these natural morphisms between matrix factorization give the commute diagram
\begin{equation}
\nonumber
\xymatrix{
K(a_1;c b_1)_R\boxtimes\overline{S}\ar[d]^{\overline{\psi}}\ar[rr]^(.45){(1,c)\boxtimes\id_{\overline{S}}}&&K(c a_1;b_1)_R\boxtimes\overline{S}\ar[d]^{\overline{\varphi}}\\
K(a_1';c b_1)_R\boxtimes\overline{S}'\ar[rr]^(.45){(1,c)\boxtimes\id_{\overline{S}'}}&&K(c a_1';b_1)_R\boxtimes\overline{S}'.
}
\end{equation}
\end{proposition}
\begin{proof}
It suffices to show influence of isomorphisms of Proposition \ref{equiv} on the following morphism
\begin{equation}
\nonumber
\xymatrix{
{K\left(\left(\begin{array}{c}a_1\\a_2\end{array}\right);\left(\begin{array}{c}c b_1\\b_2\end{array}\right)\right)_R}
\ar[rr]^{(c,1)\boxtimes\id}
&&{K\left(\left(\begin{array}{c}c a_1\\a_2\end{array}\right);\left(\begin{array}{c}b_1\\b_2\end{array}\right)\right)_R}}
\end{equation}
We obtain the following morphism by direct calculation of morphism composition
\begin{equation}
\nonumber
\xymatrix{
{K\left(\left(\begin{array}{c}a_1'\\a_2'\end{array}\right);\left(\begin{array}{c}c b_1\\b_2\end{array}\right)\right)_R}
\ar[rr]^{(c,1)\boxtimes\id}
&&{K\left(\left(\begin{array}{c}c a_1'\\a_2'\end{array}\right);\left(\begin{array}{c}b_1\\b_2\end{array}\right)\right)_R}}
\end{equation}
\end{proof}
%
%
%
%
\subsection{Cohomology of complex of $\Z$-graded matrix factorizations}\label{sec2.12}
We define a cohomology of a complex of factorizations.\\
\begin{proposition}
A complex of factorizations $\overline{M}^{\bullet}$ of $\ob(\kom^b(\HMF^{gr,-}_{R,\omega}))$ induces the following complex of $\Z\oplus\Z_2$-graded $\Q$-vector spaces, denoted by $\H_{mf}(\overline{M}^{\bullet})$:
\begin{equation}
\nonumber
\xymatrix{\H_{mf}(\overline{M}^{\bullet})&=&
\ldots\ar[r]&
\H(\overline{M}^{i-1})\ar[rr]^{\H({d_c}_{\overline{M}^{i-1}})}&&
\H(\overline{M}^{i})\ar[rr]^{\H({d_c}_{\overline{M}^{i}})}&&
\H(\overline{M}^{i+1})\ar[r]&
\ldots 
}.
\end{equation}
\end{proposition}
\begin{proof}
By Proposition \ref{null-0}, a null-homotopic morphism between factorizations induces 0 between cohomologies of these factorizations.
Therefore, the condition ${d_c}_{\overline{M}^{i+1}}{d_c}_{\overline{M}^{i}}=0$ in $\mor(\HMF^{gr,-}_{R,\omega})$ induces $\H({d_c}_{\overline{M}^{i+1}})\H({d_c}_{\overline{M}^{i}})=0$ for every $i$.
\end{proof}
\begin{proposition}
A $\Z$-grading preserving morphism $\overline{f}^{\bullet}$ from $\overline{M}^{\bullet}$ to $\overline{N}^{\bullet}$ of $\mor(\kom^b(\HMF^{gr,-}_{R,\omega}))$ induces a $\Z$-grading preserving morphism $\H_{mf}(\overline{f}^{\bullet})$ from $\H_{mf}(\overline{M}^{\bullet})$ to $\H_{mf}(\overline{N}^{\bullet})$;
\begin{equation}
\nonumber
\xymatrix{
\H_{mf}(\overline{M}^{\bullet})\ar[d]^{\H_{mf}(\overline{f}^{\bullet})}&\ar@{}[d]|-{=}&
\ldots\ar[r]&
\H(\overline{M}^{i-1})\ar[rr]^{\H({d_c}_{\overline{M}^{i-1}})}\ar[d]^{\H(\overline{f}^{i-1})}&&
\H(\overline{M}^{i})\ar[rr]^{\H({d_c}_{\overline{M}^{i}})}\ar[d]^{\H(\overline{f}^{i})}&&
\H(\overline{M}^{i+1})\ar[r]\ar[d]^{\H(\overline{f}^{i+1}})&
\ldots \\
\H_{mf}(\overline{N}^{\bullet})&&
\ldots\ar[r]&
\H(\overline{N}^{i-1})\ar[rr]^{\H({d_c}_{\overline{N}^{i-1}})}&&
\H(\overline{N}^{i})\ar[rr]^{\H({d_c}_{\overline{N}^{i}})}&&
\H(\overline{N}^{i+1})\ar[r]&
\ldots 
}.
\end{equation}
\end{proposition}
\begin{proof}
The condition $\overline{f}^{i+1}{d_c}_{\overline{M}^{i}}={d_c}_{\overline{N}^{i}}\overline{f}^{i}$ in $\mor(\HMF^{gr,-}_{R,\omega})$ induces the condition $\H(\overline{f}^{i+1})\H({d_c}_{\overline{M}^{i}})=\H({d_c}_{\overline{N}^{i}})\H(\overline{f}^{i})$ for every $i$ by Proposition \ref{null-0}.
\end{proof}
\begin{corollary}
$\H_{mf}$ is a functor from $\kom^b(\HMF^{gr,-}_{R,\omega})$ to the category of complexes of $\Z\oplus\Z_2$-graded $\Q$-vector spaces and $\Z$-grading preserving morphisms, denoted by $\kom^b(\Q$-$\mathcal{V}ect)$.
\end{corollary}
\begin{proposition}
A complex null-homotopic morphism $\overline{f}^{\bullet}$ from $\overline{M}^{\bullet}$ to $\overline{N}^{\bullet}$ of $\mor(\kom^b(\HMF^{gr,-}_{R,\omega}))$ induces a complex null-hotopic morphism $\H_{mf}(\overline{f}^{\bullet})$ from $\H_{mf}(\overline{M}^{\bullet})$ to $\H_{mf}(\overline{N}^{\bullet})$ of $\mor(\kom^b(\Q$-$\mathcal{V}ect))$.
\end{proposition}
\begin{proof}
There exists a collection $\overline{h_c}^{\bullet}=(h_c^i:\overline{M}^i\to\overline{N}^{i-1})$ such that $\overline{f}^{i}= h_c^{i+1}{d_c}_{\overline{M}^{i}}+{d_c}_{\overline{N}^{i-1}}h_c^i$ in $\mor(\HMF^{gr,-}_{R,\omega})$ for any $i$.
The collection $\overline{h_c}^{\bullet}=(\overline{h}_c^i:\overline{M}^i\to\overline{N}^{i-1})$ induces the collection $\H_{mf}(\overline{h_c}^{\bullet})=(\H(\overline{h}_c^i):\H(\overline{M}^i)\to\H(\overline{N}^{i-1}))$ satisfying that the condition $\H(\overline{f}^{i})=\H(\overline{h}_c^{i+1})\H({d_c}_{\overline{M}^{i}})+\H({d_c}_{\overline{N}^{i-1}})\H(\overline{h}_c^i)$ for every $i$ by Proposition \ref{null-0}.
\end{proof}
\begin{corollary}
$\H_{mf}$ is a functor from $\k^b(\HMF^{gr,-}_{R,\omega})$ to the homotopy category $\k^b(\Q$-$\mathcal{V}ect)$.
\end{corollary}
We denote the cohomology of $\H_{mf}(\overline{M}^{\bullet})$ by $\H(\overline{M}^{\bullet})$ and call it a {\bf cohomology of a complex of factorizations}.
This is a $\Z\oplus\Z\oplus\Z_2$-graded $\Q$-vector space $\H(\overline{M}^{\bullet})=\bigoplus_{k\in\Z}\H^{k,0}(\overline{M}^{\bullet})\oplus\H^{k,1}(\overline{M}^{\bullet})$, where $k$ is corresponding to the complex grading.
\begin{definition}
Poincar\'e polynomial $P(\overline{M}^{\bullet})$ of a complex of matrix factorizations $\overline{M}^{\bullet}$ is defined to be
\begin{equation}
\nonumber
P(\overline{M}^{\bullet}):=\sum_{k\in\Z}t^k\left\{\gdim(\H^{k,0}(\overline{M}^{\bullet}))+\, s\,\gdim(\H^{k,1}(\overline{M}^{\bullet}))\right\}.
\end{equation}
\end{definition}

\begin{remark}\label{rem-shift}
We find the following equations
\begin{eqnarray}
\nonumber
&&P(\overline{M}^{\bullet}\{m\})=q^mP(\overline{M}^{\bullet}),\\
\nonumber
&&P(\overline{M}^{\bullet}\{f(q)\}_q)=f(q)P(\overline{M}^{\bullet}),\\
\nonumber
&&P(\overline{M}^{\bullet}[m])=t^mP(\overline{M}^{\bullet}).
\end{eqnarray}
\end{remark}
%
%
%
%
\section{Symmetric function and its generating function}\label{sec3}
In this section, we give a few special symmetric functions and their generating functions.
Using these functions, we define matrix factorizations for colored planar diagrams and show isomorphisms corresponding to some relations of the MOY bracket in next Section \ref{sec4}.\\
%
%
%
%
\subsection{Homogeneous $\Z$-graded polynomial}\label{sec3.1}
Let $x_{k,i}$ be a variable with $\Z$-grading $2k$($k\in \N$, $i$: a formal index) and we define $x_{0,i}=1$ for any $i$.
Let $\mathbb{X}^{(m)}_{(i)}$ be a sequence of $m$ variables $x_{l,i}$ $(1 \leq l \leq m)$;
\begin{equation}
\nonumber
\mathbb{X}^{(m)}_{(i)}=(x_{1,i},x_{2,i},\ldots,x_{m,i}).
\end{equation}
Let $(i_1,i_2,\ldots,i_k)$ be a sequence of indexes.
For a sequence of positive integers $(m_1,m_2,\ldots,m_k)$, we define $R_{(i_1,i_2,\ldots,i_k)}^{(m_1,m_2,\ldots,m_k)}$ to be a $\Z$-graded polynomial ring over $\Q$ generated by variables of sequences $\mathbb{X}^{(m_1)}_{(i_1)}$, $\mathbb{X}^{(m_2)}_{(i_2)}$, $\ldots$, $\mathbb{X}^{(m_k)}_{(i_k)}$;
\begin{equation}
\nonumber
R_{(i_1,i_2,\ldots,i_k)}^{(m_1,m_2,\ldots,m_k)}=\Q[x_{1,i_1},x_{2,i_1},\ldots,x_{m_1,i_1},x_{1,i_2},x_{2,i_2},\ldots,x_{m_2,i_2},\ldots,x_{1,i_k},x_{2,i_k},\ldots,x_{m_k,i_k}]
\end{equation}
whose $\Z$-grading is induced by the $\Z$-gradings $\deg (x_{l,i_j})=2l$ $(1\leq l \leq m_j,1\leq j \leq k)$.
Let $s(m)$ be a function which is $1$ if $m\geq 0$ and $-1$ if $m<0$.
For a sequence of integers $(m_1,m_2,\ldots,m_k)$, we define $X_{(i_1,i_2,\ldots ,i_l)}^{(m_1,m_2,\ldots ,m_l)}$ to be a generating function composed of symmetric polynomials $\mathbb{X}_{(i_k)}^{(m_k)}$ ($k=1,...,l$)
$$
\prod^{l}_{k=1}(1+x_{1,i_k}+\ldots +x_{\left|m_k\right|,i_k})^{s(m_k)}
$$
and define $X_{m,(i_1,i_2,\ldots ,i_l)}^{(m_1,m_2,\ldots ,m_l)}$ to be the homogeneous term with $\Z$-grading $2m$ of $X_{(i_1,i_2,\ldots ,i_l)}^{(m_1,m_2,\ldots ,m_l)}$.
For example,
\begin{eqnarray}
\nonumber
X_{(i_1,i_2,i_3)}^{(1,-2,3)}&=&\frac{(1+x_{1,i_1})(1+x_{1,i_3}+x_{2,i_3}+x_{3,i_3})}{(1+x_{1,i_2}+x_{2,i_2})},\\
\nonumber
X_{3,(i_1,i_2,i_3)}^{(1,-2,3)}&=&\mathrm{polynomial}\hspace{0.1cm}\mathrm{with}\hspace{0.1cm}\Z\mathrm{ -}\mathrm{grading}\hspace{0.1cm}6\hspace{0.1cm}\mathrm{of}\hspace{0.1cm}
\frac{(1+x_{1,i_1})(1+x_{1,i_3}+x_{2,i_3}+x_{3,i_3})}{(1+x_{1,i_2}+x_{2,i_2})}\\
\nonumber
&=&2x_{1,i_2}x_{2,i_2}-x_{1,i_2}^3 +(-x_{2,i_2}+x_{1,i_2}^2)(x_{1,i_1}+x_{1,i_3})-x_{1,i_2}(x_{2,i_3}+x_{1,i_1}x_{1,i_3})+x_{3,i_3}+x_{1,i_1}x_{2,i_3}.
\end{eqnarray}
In general, we denote the sequence of homogeneous $\Z$-graded $X_{m,(i_1,i_2,\ldots ,i_l)}^{(m_1,m_2,\ldots ,m_l)}$ $(m\in\N_{\geq 1})$ by $\mathbb{X}^{(m_1,m_2,\ldots ,m_l)}_{(i_1,i_2,\ldots ,i_l)}$;
\begin{equation}
\nonumber
\mathbb{X}^{(m_1,m_2,\ldots ,m_l)}_{(i_1,i_2,\ldots ,i_l)}=(X_{m,(i_1,i_2,\ldots ,i_l)}^{(m_1,m_2,\ldots ,m_l)})_{m\in\N_{\geq 1}}.
\end{equation}
\indent
These polynomials have the following properties.
\begin{proposition}\label{sym-prop}
\begin{itemize}
\item[(1)]For any $\sigma\in S_k$, where $S_k$ is symmetric group, 
$$
X_{m,(i_1,i_2,\ldots,i_k)}^{(m_1,m_2,\ldots,m_k)}=
X_{m,(i_{\sigma(1)},i_{\sigma(2)},\ldots,i_{\sigma(k)})}^{(m_{\sigma(1)},m_{\sigma(2)},\ldots,m_{\sigma(k)})}.
$$
\item[(2)]For any $l\in\{1,2,\ldots,k-1\}$,
$$
X_{m,(i_1,i_2,\ldots,i_k)}^{(m_1,m_2,\ldots,m_k)}=
\sum_{j=0}^{m}X_{m-j,(i_1,\ldots,i_l)}^{(m_1,\ldots ,m_l)}X_{j,(i_{l+1},\ldots,i_k)}^{(m_{l+1},\ldots,m_k)}.
$$
\item[(3)]For any positive integer $m_1$,
$$
X_{m,(i_1,i_2,\ldots,i_k)}^{(-m_1,m_2,\ldots,m_k)}=
X_{m,(i_2,\ldots,i_k)}^{(m_2,\ldots ,m_k)}
-x_{1,i_1}X_{m-1,(i_1,i_2,\ldots,i_k)}^{(-m_1,m_2,\ldots,m_k)}
-\ldots
-x_{m_1,i_1}X_{m-m_1,(i_1,i_2,\ldots,i_k)}^{(-m_1,m_2,\ldots,m_k)}.
$$
\item[(4)]For any positive integer $m$, we have
\begin{equation}
\nonumber
\sum_{l=0}^{m}X_{m-l,(i_1,\ldots,i_k)}^{(m_1,\ldots,m_k)}X_{l,(i_1,\ldots,i_k)}^{(-m_1,\ldots,-m_k)}=0.
\end{equation}
\item[(5)]For any number $l\in\{0,1,\ldots,k-1\}$, we find that $(X_{1,(i_1,\ldots,i_k)}^{(m_1,\ldots,m_k)},\ldots,X_{m,(i_1,\ldots,i_k)}^{(m_1,\ldots,m_k)})$ and $(X_{1,(i_1,\ldots,i_l)}^{(m_1,\ldots ,m_l)}-X_{1,(i_{l+1},\ldots,i_k)}^{(-m_{l+1},\ldots ,-m_k)},\ldots,X_{m,(i_1,\ldots,i_l)}^{(m_1,\ldots,m_l)}-X_{m,(i_{l+1},\ldots,i_k)}^{(-m_{l+1},\ldots,-m_k)})$ transform to each other by linear translations over $R_{(i_1,i_2,\ldots,i_k)}^{(m_1,m_2,\ldots,m_k)}$.
\end{itemize}
\end{proposition}
\begin{proof}
$(1)$: $X_{m,(i_1,i_2,\ldots,i_k)}^{(m_1,m_2,\ldots,m_k)}$ and $X_{m,(i_{\sigma(1)},i_{\sigma(2)},\ldots ,i_{\sigma(k)})}^{(m_{\sigma(1)},m_{\sigma(2)},\ldots,m_{\sigma(k)})}$ have the same generating function.
Then, we obtain $(1)$\\
$(2)$: The generating function of $X_{m,(i_1,i_2,\ldots,i_k)}^{(m_1,m_2,\ldots,m_k)}$ equals 
$$
(\prod_{s=0}^{\infty}X_{s,(i_1,\ldots,i_l)}^{(m_1,\ldots ,m_l)})
(\prod_{t=0}^{\infty}X_{t,(i_{l+1},\ldots,i_k)}^{(m_{l+1},\ldots,m_k)}).
$$
$(3)$ is obtained by the equality
$$
\frac{1}{1+x_{1,i_1}+\ldots+x_{m_1,i_1}}=1-x_{1,i_1}\frac{1}{1+x_{1,i_1}+\ldots+x_{m_1,i_1}}-\ldots-x_{m_1,i_1}\frac{1}{1+x_{1,i_1}+\ldots+x_{m_1,i_1}}.
$$
$(4)$: By Proposition \ref{sym-prop} $(2)$, the left-hand side is an $m$-graded polynomial of $X_{(i_1,\ldots,i_k)}^{(m_1,\ldots,m_k)}X_{(i_1,\ldots,i_k)}^{(-m_1,\ldots,-m_k)}$.
However, we have $X_{(i_1,\ldots,i_k)}^{(m_1,\ldots,m_k)}X_{(i_1,\ldots,i_k)}^{(-m_1,\ldots,-m_k)}=1$.
Then, we obtain the equation of $(4)$.\\
$(5)$: We use the induction to $m$.
We have
\begin{equation}
\nonumber
X_{1,(i_1,\ldots ,i_k)}^{(m_1,\ldots ,m_k)}=X_{1,(i_1,\ldots ,i_l)}^{(m_1,\ldots ,m_l)}-X_{1,(i_{l+1},\ldots ,i_k)}^{(-m_{l+1},\ldots ,-m_k)}.
\end{equation}
By assumption of the induction and Proposition \ref{sym-prop} $(2)$, we have
\begin{equation}
\nonumber
X_{s,(i_1,\ldots ,i_k)}^{(m_1,\ldots ,m_k)}=\sum_{t=0}^{s}X_{s-t,(i_1,\ldots,i_l)}^{(m_1,\ldots,m_l)}X_{t,(i_{l+1},\ldots,i_k)}^{(m_{l+1},\ldots,m_k)}.
\end{equation}
Moreover, by the assumption and Proposition \ref{sym-prop} $(4)$, we have
\begin{eqnarray}
\nonumber
\sum_{t=0}^{s}X_{s-t,(i_1,\ldots,i_l)}^{(m_1,\ldots,m_l)}X_{t,(i_{l+1},\ldots,i_k)}^{(m_{l+1},\ldots,m_k)}
&=&X_{s,(i_1,\ldots,i_l)}^{(m_1,\ldots,m_l)}+\sum_{t=1}^{s}X_{s-t,(i_{l+1},\ldots,i_k)}^{(-m_{l+1},\ldots,-m_k)}X_{t,(i_{l+1},\ldots,i_k)}^{(m_{l+1},\ldots,m_k)}\\
\nonumber
&=&X_{s,(i_1,\ldots,i_l)}^{(m_1,\ldots,m_l)}-X_{s,(i_{l+1},\ldots,i_k)}^{(-m_{l+1},\ldots,-m_k)}.
\end{eqnarray}
Therefore, $(X_{1,(i_{l+1},\ldots,i_k)}^{(m_{l+1},\ldots,m_k)},\ldots,X_{s,(i_{l+1},\ldots,i_k)}^{(m_{l+1},\ldots,m_k)})$ and $(X_{1,(i_1,\ldots,i_l)}^{(m_1,\ldots,m_l)}-X_{1,(i_{l+1},\ldots,i_k)}^{(-m_{l+1},\ldots,-m_k)},\ldots,X_{s,(i_1,\ldots,i_l)}^{(m_1,\ldots,m_l)}-X_{s,(i_{l+1},\ldots,i_k)}^{(-m_{l+1},\ldots,-m_k)})$ transform to each.
Thus, we obtain the claim $(5)$.
\end{proof}
\begin{remark}
Proposition \ref{sym-prop} $(5)$ is obtained by the following equivalent in the words of a generating function:
\begin{eqnarray}
\nonumber
&&X_{m,(i_1,\ldots ,i_k)}^{(m_1,\ldots ,m_k)}=0\hspace{.3cm}\mathrm{for}\hspace{.1cm}\mathrm{any}\hspace{.1cm}m\in\N_{\geq 1}
\\
\nonumber
&\Leftrightarrow&
\prod_{j=1}^{k}(1+x_{1,i_j}+\ldots +x_{\left|m_j\right|,i_j})^{s(m_j)}=1
\\
\nonumber
&\Leftrightarrow&
\prod_{j=1}^{l}(1+x_{1,i_j}+\ldots +x_{\left|m_j\right|,i_j})^{s(m_j)}
-\prod_{j=l+1}^{k}(1+x_{1,i_j}+\ldots +x_{\left|m_j\right|,i_j})^{-s(m_j)}=0
\\
\nonumber
&\Leftrightarrow&
X_{m,(i_1,\ldots,i_l)}^{(m_1,\ldots,m_l)}-X_{m,(i_{l+1},\ldots,i_k)}^{(-m_{l+1},\ldots,-m_k)}=0\hspace{.3cm}\mathrm{for}\hspace{.1cm}\mathrm{any}\hspace{.1cm}m\in\N_{\geq 1}.
\end{eqnarray}
\end{remark}
%
%
%
%
\subsection{Power sum, elementary and complete symmetric function}\label{sec3.2}
\indent
Hereinafter, we fix an integer $n$.
The integer $n$ means that we consider a homology theory corresponding to the quantum $\mathfrak{sl}_n$ link invariant.
We suppose that variables $t_{1,i}$, $t_{2,i}$, $\ldots$, $t_{m,i}$, where $i$ is a formal index, have $\Z$-grading 2.
We consider the power sum $t_{1,i}^{n+1}+t_{2,i}^{n+1}+\ldots +t_{m,i}^{n+1}$ in the polynomial ring $\Q[t_{1,i},\ldots,t_{m,i}]$.
The elementary symmetric functions $x_{j,i}=\sum_{1\leq k_1<\ldots<k_j\leq m}t_{k_1,i}\ldots t_{k_j,i}$ $(1\leq j\leq m)$ form a basis of symmetric functions (Its $\Z$-grading is naturally $2j$).
Then, the power sum is represented as a polynomial of the subring $\Q[x_{1,i},\ldots,x_{m,i}]$ generated by the elementary symmetric functions, denoted by $F_{m}(x_{1,i},x_{2,i},\ldots,x_{m,i})$ or $F_{m}(\mathbb{X}_{(i)}^{(m)})$ for short;
\begin{equation}
\nonumber
F_{m}(\mathbb{X}_{(i)}^{(m)})=F_{m}(x_{1,i},x_{2,i},\ldots,x_{m,i})=t_{1,i}^{n+1}+t_{2,i}^{n+1}+\ldots +t_{m,i}^{n+1}.
\end{equation}
We find that the elementary symmetric function $x_{k,i}$ naturally has $\Z$-grading $2k$.
\begin{proposition}
Put $x_{j,i}=\sum_{1\leq k_1<\ldots<k_j\leq m}t_{k_1,i}\ldots t_{k_j,i}$ $(1\leq j\leq m)$, which is the elementary symmetric functions of variables $t_{1,i}$, $t_{2,i}$, $\ldots$, $t_{m,i}$, and $y_{j,i}=\sum_{1\leq k_1\leq\ldots\leq k_j\leq m}t_{k_1,i}\ldots t_{k_j,i}$ $(1\leq j\leq m)$, which is the complete symmetric functions of variables $t_{1,i}$, $t_{2,i}$, $\ldots$, $t_{m,i}$.
\begin{itemize}
\item[(1)]
$X^{(m)}_{(i)}$ is a generating function of elementary symmetric functions $x_{j,i}$.
\item[(2)]
$X^{(-m)}_{(i)}$ is a generating function of complete symmetric functions up to $\pm 1$.
\item[(3)]
For $m\leq n$, we have
\begin{equation}
\nonumber
F_{m}(x_{1,i},x_{2,i},\ldots,x_{m,i})=\sum_{k=1}^{m}(-1)^{n+1-k}k\,x_{k,i}X_{n+1-k,(i)}^{(-m)}.
\end{equation}
\end{itemize}
\end{proposition}
\begin{proof}
(1): It is obvious by definition.\\
(2): We find that
\begin{equation}
\nonumber
(-1)^k\frac{1}{k!}\left(\frac{d}{d T}\right)^k\left.\left(\frac{1}{1+x_{1,i}T+x_{2,i}T^2+\ldots+x_{m,i}T^m}\right)\right|_{T=0}
=
\left|
\begin{array}{ccccc}
x_{1,i}&1&0&\cdots&0\\
x_{2,i}&x_{1,i}&1&\ddots&\vdots\\
x_{3,i}&x_{2,i}&x_{1,i}&\ddots&0\\
\vdots&\vdots&&\ddots&1\\
x_{k,i}&x_{k-1,i}&\cdots&x_{2,i}&x_{1,i}
\end{array}
\right|.
\end{equation}
We pick out a homogeneous polynomial with $\Z$-grading $2k$ from the rational function $X_{(i)}^{(-m)}$ on the left hand side of the equation, that is, $X_{k,(i)}^{(-m)}$.
On the other hand side, the determinant is the complete symmetric function with $\Z$-grading $2k$ described by elementary symmetric functions.\\
(3):We have
\begin{equation}
\nonumber
F_{m}(x_{1,i},x_{2,i},\ldots,x_{m,i})=\left|
\begin{array}{ccccc}
x_{1,i}&1&0&\cdots&0\\
2x_{2,i}&x_{1,i}&1&\ddots&\vdots\\
3x_{3,i}&x_{2,i}&x_{1,i}&\ddots&0\\
\vdots&\vdots&&\ddots&1\\
(n+1)x_{n+1,i}&x_{n,i}&\cdots&x_{2,i}&x_{1,i}
\end{array}
\right|,
\end{equation}
where some variables in the right-hand determinant satisfy $x_{m+1,i}=x_{m+2,i}=\ldots=x_{n+1,i}=0$.
We apply Laplace expansion to the first column of the determinant. 
Then, we obtain this proposition (3) by (2).
\end{proof}
\begin{proposition}
$(1)$The sum of the polynomials $F_{m_1}(\mathbb{X}^{(m_1)}_{(i_1)})$ and $F_{m_2}(\mathbb{X}^{(m_2)}_{(i_2)})$ equals to $F_{m_1+m_2}(\mathbb{X}^{(m_1,m_2)}_{(i_1,i_2)})$;
$$
F_{m_1}(\mathbb{X}^{(m_1)}_{(i_1)})+F_{m_2}(\mathbb{X}^{(m_2)}_{(i_2)})=F_{m_1+m_2}(\mathbb{X}^{(m_1,m_2)}_{(i_1,i_2)}).
$$
$(2)$The polynomial $F_m(\mathbb{X}^{(m)}_{(i)})$ is a potential of $R_{(i)}^{(m)}$.\\
\end{proposition}
\begin{proof}
$(1)$It is obvious by redescribing $x_{j,i}$ as $\sum_{1\leq k_1<\ldots<k_j\leq m}t_{k_1,i}\ldots t_{k_j,i}$.\\
$(2)$When $m\geq n+1$, $\frac{\partial F_m}{\partial x_{n+1,i}}=1$.
Then, the Jacobi ring $J_{F_m}\simeq\Q$.
Therefore, for $m\leq n$, we show that the Jacobi ring $J_{F_m}=R_{(i)}^{(m)}\left/\left<\frac{\partial F_m}{\partial x_{1,i}},\ldots,\frac{\partial F_m}{\partial x_{m,i}}\right>\right.$ is finite dimension over $Q$.
In other words, we show that the sequence $(\frac{\partial F_m}{\partial x_{1,i}}$, $\ldots$, $\frac{\partial F_m}{\partial x_{m,i}})$ forms regular in $R_{(i)}^{(m)}$.
We find 
\begin{eqnarray}
\nonumber
\frac{\partial F_m(\mathbb{X}_{(i)}^{(m)})}{\partial x_{j,i}}&=&
(-1)^{j-1}j\,X_{n+1-j,(i)}^{(-m)}+(-1)^{j-1}\sum_{k=1}^{n+1-j}F_k(\mathbb{X}_{(i)}^{(m)})X_{n-k,(i)}^{(-m)} \hspace{1cm}(j=1,\ldots,m)\\
\nonumber
&=&(-1)^{j-1}j\,X_{n+1-j,(i)}^{(-m)}+(-1)^{j-1}(n+1-j)X_{n+1-j,(i)}^{(-m)}=(-1)^{j-1}(n+1)X_{n+1-j,(i)}^{(-m)}.
\end{eqnarray}
The radical ideal of $\langle X_{n,(i)}^{(-m)},\ldots,X_{n+1-m,(i)}^{(-m)}\rangle$ is equal to the maximal ideal $\langle x_{1,i},\ldots,x_{m,i}\rangle$.
Thus, the sequence $(\frac{\partial F_m}{\partial x_{1,i}}$, $\ldots$, $\frac{\partial F_m}{\partial x_{m,i}})$ is regular.
\end{proof}
\begin{corollary}\label{potential-f}
$(1)$The sum of the polynomials $F_{m_k}(\mathbb{X}^{(m_k)}_{(i_k)})$ $(k=1,\ldots,j)$ equals to $F_{\sum_{k=1}^{j}m_k}(\mathbb{X}^{(m_1,m_2,\ldots,m_j)}_{(i_1,i_2,\ldots,i_j)})$;
$$
\sum_{k=1}^{j}F_{m_k}(\mathbb{X}^{(m_k)}_{(i_k)})=F_{\sum_{k=1}^{j}m_k}(\mathbb{X}^{(m_1,m_2,\ldots,m_j)}_{(i_1,i_2,\ldots,i_j)}).
$$
$(2)$The polynomial $\displaystyle \sum_{k=1}^{j}F_{m_k}(\mathbb{X}^{(m_k)}_{(i_k)})$ is a potential of $R_{(i_1,i_2,\ldots,i_j)}^{(m_1,m_2,\ldots,m_j)}$.
\end{corollary}
%
%
%
%
\section{Colored planar diagrams and matrix factorizations}\label{sec4}
\begin{figure}[hbt]
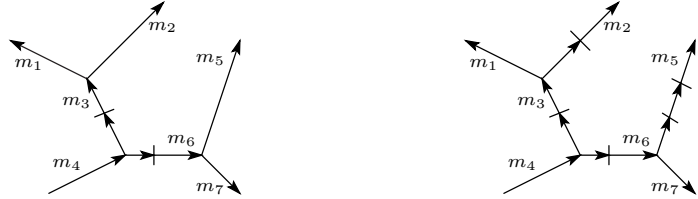

\input{figure/fig-decomp-eff}\hspace{3cm}\input{figure/fig-decomp-waste}
\caption{Effective decomposition and non-effective decomposition}\label{decomposition}
\end{figure}

\indent
In paper \cite{KR1}, Khovanov and Rozansky gave a potential for the vector representation $V_n$ of $U_q(\mathfrak{sl}_n)$ and defined matrix factorizations for intertwiners between tensor products of $V_n$.
Then, they showed that there exist isomorphisms of matrix factorizations corresponding to relations of intertwiners, see the MOY relations between planar diagrams with coloring $1$ and $2$ in Appendix \ref{NMOY}, and defined a complex for an oriented link diagram using these matrix factorizations.
Moreover, they discussed a potential for anti-symmetric tensor product of $V_n$, called the fundamental representation, in Section $11$ of \cite{KR1}.
H. Wu and the author independently defined matrix factorizations for intertwiners of the fundamental representations $\land^iV_n$ ($i=1,\ldots,n-1$) \cite{Wu1}\cite{Yone1}.
They independently showed that there exist isomorphisms of factorizations corresponding to most relations of the MOY bracket.\\
\indent
In this section, we give definition of the factorization for colored planar diagrams and isomorphisms between factorizations corresponding to most MOY relations.
\\
\indent
Before defining factorizations for colored planar diagrams, we show the structure of the colored planar diagrams derived from a colored oriented link diagram by using the MOY bracket.
The MOY bracket expands a single $[i,j]$-crossing into a linear combination of colored planar diagrams in Figure \ref{colored-planar}.
The colored planar diagram is locally composed of three types of oriented diagrams called {\bf essential}, see diagrams in Figure \ref{essential-planar}.
Therefore, colored planar diagrams obtained by applying the MOY bracket to a colored oriented link diagram also locally consist of the essential planar diagrams.
\\
\indent
For a colored planar diagram $\Gamma$, we consider a decomposition of $\Gamma$ into some essential diagrams using markings, see Figure \ref{decomposition}.
\begin{definition}
A decomposition into essential diagrams is {\bf effective} if there exists no marking such that the decomposition cleared the marking off still consists of essential diagrams.
A decomposition into essential diagrams is {\bf non-effective} if there exists such a marking.
\end{definition}
For a given colored planar diagram, its effective decomposition is uniquely determined up to isotopy.
\begin{definition}
A colored planar diagram is a {\bf cycle} if the diagram has a region encircled by edges of the diagram and is a {\bf tree} otherwise.
\end{definition}
The colored planar diagram produced from relations of the MOY bracket can be roughly divided into two types of cycle and tree. 
\begin{figure}[hbt]
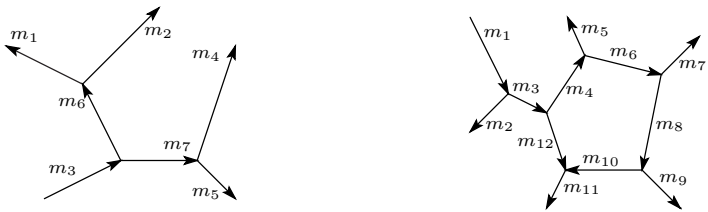

\input{figure/fig-3-valent1}\hspace{3cm} \input{figure/fig-3-valent-cycle} 
\caption{Tree diagram and cycle diagram}\label{tree-cycle}
\end{figure}\\
\indent
Matrix factorizations for the colored planar diagrams are defined using the expression of the power sum in the elementary symmetric functions and homogeneous $\Z$-graded polynomials in Section \ref{sec3.2}.
%
%
%
%
\subsection{Potential of colored planar diagram}\label{sec4.1}
\begin{figure}[htb]
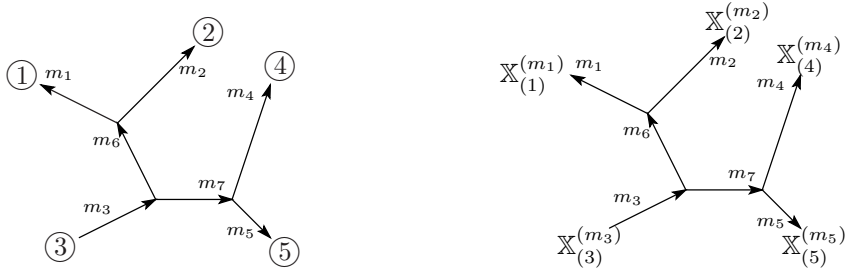

$$
\input{figure/fig-3-valent1-formal-index-en}
\hspace{2cm}
\hspace{2cm}
\input{figure/fig-3-valent1-formal-index}
$$
\caption{Planar diagram assigned formal indexes and diagram assigned sequences}\label{planar-assign-index}
\end{figure}

\indent
We define a potential for a colored planar diagram.
It is a power sum determined by coloring, orientation of the diagram and an additional data which is a formal index.\\
\indent
For a given colored planar diagram, we assign a distinct formal index $i$ to each end of the diagram and, then, assign a power sum to each end as follows.
When an edge including an $i$-assigned end has a coloring $m$ and an orientation from inside diagram to outside end, we assign the polynomial $+F_{m}(x_{1,i},x_{2,i},\ldots,x_{m,i})$ to the end, and when an edge has an opposite orientation from outside to inside, assign the polynomial $-F_{m}(x_{1,i},x_{2,i},\ldots,x_{m,i})$.
A {\bf potential} of a colored planar diagram is defined to be the sum of these assigned polynomials over every ends of the diagram.\\
\indent
To each end of the edge with coloring $m$ we simply assign only a formal index $i$ or a sequence of variables $\mathbb{X}^{(m)}_{(i)}$ for convenience, see Figure \ref{planar-assign-index}.
These datum are enough to seek a potential of a diagram.
\\
\indent
For instance, the potential of the diagram in Figure \ref{planar-assign-index} is
\begin{equation}
\nonumber
F_{m_1}(\mathbb{X}_{(1)}^{(m_1)})+F_{m_2}(\mathbb{X}_{(2)}^{(m_2)})-F_{m_3}(\mathbb{X}_{(3)}^{(m_3)})+F_{m_4}(\mathbb{X}_{(4)}^{(m_4)})+F_{m_5}(\mathbb{X}_{(5)}^{(m_5)}).
\end{equation}
%
%
%
%
\subsection{Essential planar diagrams and matrix factorizations}\label{sec4.2}
For an essential planar diagram, we define a matrix factorizations with the potential of the diagram. 
\begin{definition}
\indent
A matrix factorization for a colored planar line,
\begin{center}
\input{figure/figmoy3-mf}\hspace{2cm}${\rm (1 \leq m \leq n)},$
\end{center}
\hspace{1mm}\\\\
is defined to be
\begin{equation}
\label{line-mf}
\c \Biggl( \input{figure/figmoy3-mf} \Biggl)_n:=
\mathop{\boxtimes}_{j=1}^{m} K\Big( L^{[m]}_{j,(1;2)} ;X_{j,(1)}^{(m)}-X_{j,(2)}^{(m)} \Big)_{R_{(1,2)}^{(m,m)}} ,
\end{equation}
where 
$$
L^{[m]}_{j,(1;2)} =\frac{F_{m}(X_{1,(2)}^{(m)},...,X_{j-1,(2)}^{(m)} ,X_{j,(1)}^{(m)},...,X_{m,(1)}^{(m)})-F_{m}(X_{1,(2)}^{(m)},...,X_{j,(2)}^{(m)},X_{j+1,(1)}^{(m)},...,X_{m,(1)}^{(m)})}{X_{j,(1)}^{(m)}-X_{j,(2)}^{(m)}}.
$$
\end{definition}
It is obvious that this matrix factorization is a finite factorization of $\MF^{gr,fin}_{R_{(1,2)}^{(m,m)},F_m (\mathbb{X}^{(m)}_{(1)}) -F_m (\mathbb{X}^{(m)}_{(2)}) }$.
We denote this matrix factorization $\overline{L}^{[m]}_{(1;2)}$ for short.
\begin{remark}
For $m \geq n+1$, we can consider the matrix factorization for a line colored $m$, \input{figure/figmoy5} , as the above definition.
However, we find that such matrix factorizations are contractible, that is, isomorphic to the zero matrix factorization in $\HMF^{gr}$.
Because, in the case that $m\geq n+1$, the matrix factorization $\overline{L}^{[m]}_{(1;2)}$ includes the contractible matrix factorization
\begin{equation}
\label{trivial}
K(L^{[m]}_{n+1,(1;2)};x_{n+1,(1)}^{(m)}-x_{n+1,(2)}^{(m)})_{R_{(1,2)}^{(m,m)}}.
\end{equation}
The polynomial with $m$ variables $F_{m}$ is the expression of the power sum $t_1^{n+1}+t_2^{n+1}+\ldots +t_m^{n+1}$ with the elementary symmetric functions $x_j=\sum_{1\leq i_1<\ldots<i_j\leq m}t_{i_1}\ldots t_{i_j}$ $(1\leq j\leq m)$.
However, in the case of $m\geq n+1$, the power sum $t_1^{n+1}+t_2^{n+1}+\ldots +t_m^{n+1}$ is described as a polynomial of $n+1$ variables $x_1,x_2,\ldots,x_{n+1}$.
Thus, we find $L^{[m]}_{n+1,(1;2)}=(-1)^n (n+1)$.
Then, the matrix factorization (\ref{trivial}) is contractible.
\end{remark}
\begin{definition}
\indent
We define matrix factorizations for the following trivalent diagrams:
\begin{center}
\input{figure/figgluing-in-3valent-mf} \hspace{1cm} \input{figure/figgluing-out-3valent-mf}\hspace{2cm}$(2\leq m_3=m_1+m_2\leq n)$.\vspace{1cm}
\end{center}
The first one is defined to be
\begin{equation}
\label{n-mf}
\c \left(\input{figure/figgluing-in-3valent-mf} \right)_n := \mathop{\boxtimes}_{j=1}^{m_{3}} 
K\Big( \Lambda_{j,(3;1,2)}^{[m_1,m_2]} ;X_{j,(3)}^{(m_3)}-X^{(m_1,m_2)}_{j,(1,2)} \Big)_{R_{(1,2,3)}^{(m_1,m_2,m_3)}} ,
\end{equation}\\[-0.1em]
where 
\begin{eqnarray}
\nonumber
\Lambda_{j,(3;1,2)}^{[m_1,m_2]}=
\frac{
F_{m_{3}}(X^{(m_1,m_2)}_{1,(1,2)},...,X^{(m_1,m_2)}_{j-1,(1,2)} ,X_{j,(3)}^{(m_3)},...,X_{m_{3},(3)}^{(m_3)})
-F_{m_{3}}(X^{(m_1,m_2)}_{1,(1,2)},...,X^{(m_1,m_2)}_{j,(1,2)},X_{j+1,(3)}^{(m_3)},...,X_{m_{3},(3)}^{(m_3)})
}
{
X_{j,(3)}^{(m_3)}-X^{(m_1,m_2)}_{j,(1,2)}
},
\end{eqnarray}
denoted this matrix factorization by $\overline{\Lambda}_{(3;1,2)}^{[m_1,m_2]}$ for short.
The second one is defined to be
\begin{equation}
\label{v-mf}
\c \left(\input{figure/figgluing-out-3valent-mf} \right)_n := \mathop{\boxtimes}_{j=1}^{m_{3}} 
K\Big( V_{j,(1,2;3)}^{[m_1,m_2]} ;X^{(m_1,m_2)}_{j,(1,2)}-x_{j,(3)}^{(m_3)} \Big)_{R_{(1,2,3)}^{(m_1,m_2,m_3)}} \{ - m_1 m_2 \},
\end{equation}
where 
\begin{eqnarray}
\nonumber
V_{j,(1,2;3)}^{[m_1,m_2]}=
\frac{
F_{m_{3}}(X_{1,(3)}^{(m_3)},...,X_{j-1,(3)}^{(m_3)} ,X^{(m_1,m_2)}_{j,(1,2)},...,X^{(m_1,m_2)}_{m_3,(1,2)})
-F_{m_{3}}(X_{1,(3)}^{(m_3)},...,X_{j,(3)}^{(m_3)},X^{(m_1,m_2)}_{j+1,(1,2)},...,X^{(m_1,m_2)}_{m_3,(1,2)})
}{
X^{(m_1,m_2)}_{j,(1,2)}-X_{j,(3)}^{(m_3)}
},
\end{eqnarray}
denoted this matrix factorization by $\overline{V}_{(1,2;3)}^{[m_1,m_2]}$ for short.
\end{definition}
Put $\omega=F_{m_1} (\mathbb{X}^{(m_1)}_{(1)}) +F_{m_2} (\mathbb{X}^{(m_2)}_{(2)}) - F_{m_3} (\mathbb{X}^{(m_3)}_{(3)})$.
Two matrix factorizations (\ref{n-mf}) and (\ref{v-mf}) are finite factorizations of $\MF^{gr,fin}_{R_{(1,2,3)}^{(m_1,m_2,m_3)},-\omega}$ and $\MF^{gr,fin}_{R_{(1,2,3)}^{(m_1,m_2,m_3)}, \omega}$ respectively.\\
\begin{remark}
\begin{itemize}
\item[(1)]For $m_3 \geq n+1$, we can consider the matrix factorization for colored planar diagrams \,\input{figure/figmoy-out}\,,\, \input{figure/figmoy-in} as the above definition.
However, we find that such matrix factorizations are contractible, that is, isomorphic to the zero matrix factorization.
\item[(2)]By definition, we can describe matrix factorizations for essential trivalent diagrams as a matrix factorization for a colored line;
\begin{eqnarray}
\nonumber
\c\left(\hspace{0.1cm}\input{figure/figgluing-in-3valent-mf3}\hspace{0.1cm}\right)_n
&=&
\c\left(\input{figure/figgluing-in-3valent-mf2}\right)_n,\\
\nonumber
\c\left(\hspace{0.1cm}\input{figure/figgluing-out-3valent-mf3}\hspace{0.1cm}\right)_n
&=&
\c\left(\input{figure/figgluing-out-3valent-mf2}\right)_n.
\end{eqnarray}
\end{itemize}
\end{remark}
%
%
%
%
\subsection{Glued diagram and matrix factorization}\label{sec4.3}
\begin{figure}[htb]
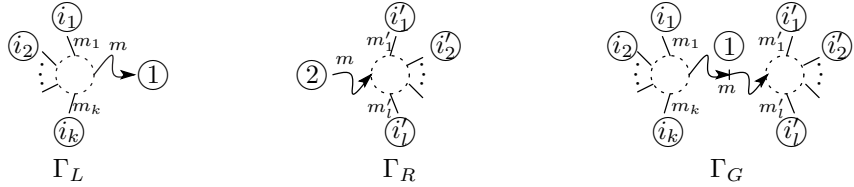

\hspace{1cm}\input{figure/figgluing-gen-planar1} \hspace{1cm} \input{figure/figgluing-gen-planar2} \hspace{2cm} \input{figure/figgluing-gen1}
\caption{Gluing planar diagrams}\label{fig-tree}
\end{figure}

\begin{definition}
For a colored planar diagram $\Gamma$ composed of the disjoint union of diagrams $\Gamma_1$ and $\Gamma_2$, we define a matrix factorization for $\Gamma$ to be tensor product of the matrix factorizations for $\Gamma_1$ and $\Gamma_2$;
\begin{equation}
\nonumber
\c(\ \Gamma\ )_{n} := \c(\ \Gamma_1\ )_{n}\boxtimes \c(\ \Gamma_2\ )_{n}.
\end{equation}
\end{definition}
We consider only a colored planar diagram locally composed of essential planar diagrams.
We inductively define a matrix factorization for the colored planar diagram obtained by gluing essential diagrams.\\
\indent
We consider two tree diagrams which have an $m$-colored edge and can be match with keeping the orientation on the edge, see the left and the middle diagrams in Figure \ref{fig-tree}.
These diagrams $\Gamma_L$ and $\Gamma_R$ can be glued at the markings $\hspace{.2cm}\en{1}\hspace{.2cm}$ and $\hspace{.2cm}\en{2}\hspace{.2cm}$ and, then, we obtain a tree diagram.
See the right diagram in Figure \ref{fig-tree}.
\begin{definition}
Let $\omega+F_{m}(\mathbb{X}^{(m)}_{(1)})$ be a potential of $\Gamma_L$ and $\omega'-F_{m}(\mathbb{X}^{(m)}_{(2)})$ be a potential of $\Gamma_R$.
We denote the factorization for $\Gamma_L$ in $\ob (\MF^{gr}_{R_{(i_1,\ldots,i_k,1)}^{(m_1,\ldots,m_k,m)},\omega + F_{m}(\mathbb{X}^{(m)}_{(1)})})$ by $\c(\Gamma_L)_{n}$ and $\Gamma_R$ in $\ob (\MF^{gr}_{R_{(i_1',\ldots,i_l',2)}^{(m_1',\ldots,m_l',m)},\omega' - F_{m}(\mathbb{X}^{(m)}_{(2)})})$ by $\c(\Gamma_R)_{n}$.
A matrix factorization for the glued diagram $\Gamma_G$ is defined to be
\begin{equation}
\nonumber
\c(\Gamma_G)_{n}:=\c(\Gamma_L)_{n} \boxtimes \c(\Gamma_R)_{n}\Big|_{\mathbb{X}^{(m)}_{(2)}=\mathbb{X}^{(m)}_{(1)}}.
\end{equation}
\end{definition}
The definition means that we identify the sequence $\mathbb{X}^{(m)}_{(1)}$ and the sequence $\mathbb{X}^{(m)}_{(2)}$ after taking the tensor product of these matrix factorizations.
Remark that the definition is essentially the same with the definition of gluing factorizations using a quotient factorization by Khovanov and Rozansky.
The glued factorization is an infinite-rank factorization but has finite-dimensional cohomology.
Therefore, the factorization is an object of $\MF^{gr}_{R_{(i_1,...,i_k,i_1',...,i_l')}^{(m_1,...,m_k,m_1',...,m_l')},\omega+\omega'}$.
\begin{proposition}\label{prop1}
The glued matrix factorization $\c(\Gamma_G)_{n}$ has finite-dimensional cohomology.
\end{proposition}
\begin{proof}
We can prove this proposition by Proposition \ref{preserve-fdc} since an essential factorization is finite and a glued diagram is decomposed into essential diagrams.
\end{proof}
\begin{figure}[hbt]
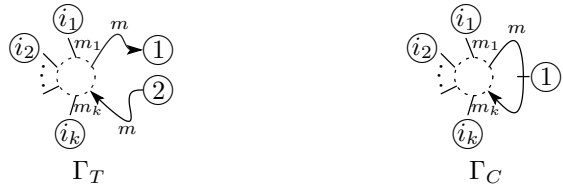

\input{figure/fig-tree-cycle1}\hspace{3cm}\input{figure/fig-tree-cycle2}
\caption{Diagram $\Gamma_T$ and cycle diagram $\Gamma_C$}\label{tree-cycle2}
\end{figure}
We consider a colored tree diagram $\Gamma_T$ and a cycle diagram $\Gamma_C$ obtained by joining ends of edges with the same coloring, see Figure \ref{tree-cycle2}.
\begin{definition}
Let $\omega+ F_{m}(\mathbb{X}^{(m)}_{(1)})- F_{m}(\mathbb{X}^{(m)}_{(2)})$ be a potential of the tree diagram $\Gamma_T$.
For factorization $\c(\Gamma_T)_{n}$ in $\ob (\MF^{gr}_{R_{(i_1,\ldots,i_k,1,2)}^{(m_1,\ldots,m_k,m,m)},\omega + F_{m}(\mathbb{X}^{(m)}_{(1)})- F_{m}(\mathbb{X}^{(m)}_{(2)})})$, a matrix factorization for the cycle diagram $\Gamma_C$ is defined to be
\begin{equation}
\nonumber
\c(\Gamma_C)_{n}:=\c(\Gamma_T)_{n}|_{\mathbb{X}_{(2)}^{(m)}=\mathbb{X}_{(1)}^{(m)}}.
\end{equation}
The factorization is an object of $\MF^{gr}_{R_{(i_1,\ldots,i_k)}^{(m_1,\ldots,m_k)},\omega}$.
\end{definition}
\begin{proposition}\label{prop2}
The glued matrix factorization $\c(\Gamma_C)_{n}$ has finite-dimensional cohomology.
\end{proposition}
\begin{proof}
The complex $\c(\Gamma_C)_{n}/\mathfrak{m}\c(\Gamma_C)_{n}$ contains the tensor product of a finite factorization with the potential $F_m(\mathbb{X}_{(1)}^{(m)})$ and a finite factorization with the potential $-F_m(\mathbb{X}_{(1)}^{(m)})$.
Then, the cohomology $\H(\c(\Gamma_C)_{n})$ is finitely dimensional by Proposition \ref{preserve-fdc}.
\end{proof}
\indent
We find that a glued matrix factorization loses potentials at glued ends.
Therefore, the potential as a colored planar diagram is compatible with the potential as a matrix factorization for a colored planar diagram.\\
\indent
For a given colored planar diagram, the matrix factorization for the diagram does not depend on a decomposition of the diagram in $\HMF^{gr}$.
\begin{proposition}\label{independent-decomp}
A matrix factorization for a colored planar diagram is independent of a decomposition of the diagram in the homotopy category $\HMF^{gr}$.
\end{proposition}
\begin{proof}
For a colored planar diagram, an effective decomposition of the diagram is uniquely determined.
Therefore, we show a factorization for any non-effective decomposition of the diagram is isomorphic to the factorization for the effective decomposition.
It suffices to show the following lemma.
\begin{lemma}\label{preserve-finiteness}
{\rm\bf (1)}We consider the following planar diagrams
\begin{center}
\input{figure/fig-gen-gluing-mf5-claim}\hspace{1cm}
\input{figure/fig-gen-gluing-mf3-claim}\hspace{1cm}
\input{figure/fig-gen-gluing-mf1-claim}.
\end{center}
There is the following canonical isomorphism in $\HMF^{gr}_{R_{(1)}^{(m)}\otimes R,\omega +F_m(\mathbb{X}^{(m)}_{(1)})}$, where the polynomial ring $R$ and the potential $\omega$ are determined by sequences of ends of the diagram except the sequence $\mathbb{X}^{(m)}_{(1)}$:
\begin{eqnarray}
\nonumber
\c\Big( \ \input{figure/fig-gen-gluing-mf3} \ \Big)_{n} \simeq \c\Big( \ \input{figure/fig-gen-gluing-mf1} \ \Big)_{n}
\end{eqnarray}
{\rm\bf (2)}We consider the following planar diagrams
\begin{center}
\input{figure/fig-gen-gluing-mf6-claim}\hspace{1cm}
\input{figure/fig-gen-gluing-mf4-claim}\hspace{1cm}
\input{figure/fig-gen-gluing-mf2-claim}.
\end{center}
There is the following canonical isomorphism in $\HMF^{gr}_{R_{(1)}^{(m)}\otimes R,\omega -F_m(\mathbb{X}^{(m)}_{(1)})}$, where the polynomial ring $R$ 
and the potential $\omega$ are determined by sequences of ends of the diagram except the sequence $\mathbb{X}^{(m)}_{(1)}$:
\begin{eqnarray}
\nonumber
\c\Big( \ \input{figure/fig-gen-gluing-mf4} \ \Big)_{n} \simeq \c\Big( \ \input{figure/fig-gen-gluing-mf2} \ \Big)_{n}
\end{eqnarray}
\end{lemma}
\noindent
{\it Proof of Lemma \ref{preserve-finiteness}}\\
We prove Lemma \ref{preserve-finiteness} (1).
By construction, a matrix factorization of the planar diagram $\Gamma_a$ forms as follows.
\begin{equation}
\nonumber
\c\big(\input{figure/fig-gen-gluing-mf5}\big)_n=\overline{M}_a\boxtimes\mathop{\boxtimes}_{k=1}^{m}K(q_k;x_{k,2}-p_k)_{R_{(2)}^{(m)}\otimes R}\in\ob(\MF^{gr}_{R_{(2)}^{(m)}\otimes R,\omega+F_m(\mathbb{X}_{(2)}^{(m)})}),
\end{equation}
where the factorization $\overline{M}_a$ and the polynomial $p_k$ are independent of variables $\mathbb{X}_{(2)}^{(m)}$.
Then, we have the following matrix factorization of the planar diagram $\Gamma_b$
\begin{eqnarray}
\label{glued-factorization}
\overline{M}_a\boxtimes\mathop{\boxtimes}_{k=1}^{m}K(q_k;x_{k,2}-p_k)_{R_{(2)}^{(m)}\otimes R}\boxtimes\mathop{\boxtimes}_{k=1}^{m}K(L_{k,(1;2)}^{[m]};x_{k,1}-x_{k,2})_{R_{(1,2)}^{(m,m)}}.
\end{eqnarray}
The potential of this factorization is $\omega+F_m(\mathbb{X}_{(1)}^{(m)})$.
We choose $x_{k,1}-x_{k,2}$ ($k=1,...,m$) as $b_j(\underline{x},\underline{y})$ of Corollary \ref{cor2-11}.
Then, the factorization (\ref{glued-factorization}) is isomorphic in $\HMF^{gr}_{R_{(1)}^{(m)}\otimes R,\omega+F_m(\mathbb{X}_{(1)}^{(m)})}$ to
\begin{equation}
\nonumber
\overline{M}_a\boxtimes\mathop{\boxtimes}_{k=1}^{m}K(r_k;x_{k,1}-p_k)_{R_{(1)}^{(m)}\otimes R},
\end{equation}
where $r_k$ is the polynomial $q_k\left|_{\mathbb{X}_{(2)}^{(k)}=\mathbb{X}_{(1)}^{(k)}}\right.$.
By Proposition \ref{reg-eq}, this is isomorphic to a factorization of the planar diagram $\Gamma_c$.
\\
\indent
We similarly prove Lemma \ref{preserve-finiteness} (2).
\end{proof}
\indent
By Proposition \ref{independent-decomp}, it suffices to obtain a factorization for a planar diagram that we consider the effective decomposition of the planar diagram.
The factorization is obtained by gluing factorizations for essential planar diagrams of the decomposition.
A matrix factorization obtained by gluing essential factorizations generally becomes an infinite factorization.
However, the glued factorization is isomorphic to a finite factorization in the homotopy category $\HMF^{gr}$ since it has finite-dimensional cohomology by Proposition \ref{prop1} and Proposition \ref{prop2}.
%
%
%
%
\subsection{MOY relations and isomorphisms between matrix factorizations}\label{sec4.4}
\indent
We show isomorphisms between factorizations for colored planar diagrams corresponding to the MOY relations in Appendix \ref{NMOY}. 
\begin{proposition}\label{mat-equiv1}
Let $\omega_1$ be a polynomial $F_{m_4}(\mathbb{X}^{(m_4)}_{(4)})-F_{m_1}(\mathbb{X}^{(m_1)}_{(1)})-F_{m_2}(\mathbb{X}^{(m_2)}_{(2)})-F_{m_3}(\mathbb{X}^{(m_3)}_{(3)})$.
\begin{enumerate}
\item[{\rm\bf (1)}] There is an isomorphism in $\HMF^{gr}_{R_{(1,2,3,4)}^{(m_1,m_2,m_3,m_4)},\omega_1}$
$$
\c\left( \input{figure/fig-ass-dia1-mf}\right)_n \simeq \c\left( \input{figure/fig-ass-dia3-mf}\right)_n \simeq \c\left( \input{figure/fig-ass-dia2-mf}\right)_n.
$$ 
\item[{\rm\bf (2)}] There is an isomorphism in $\HMF^{gr}_{R_{(1,2,3,4)}^{(m_1,m_2,m_3,m_4)},-\omega_1}$
$$
\c\left( \input{figure/fig-coass-dia1-mf}\right)_n \simeq \c\left( \input{figure/fig-coass-dia3-mf}\right)_n \simeq \c\left( \input{figure/fig-coass-dia2-mf}\right)_n,
$$
\end{enumerate}
where $1 \leq m_1,m_2,m_3 \leq n-2$, $m_5 = m_1 +m_2\leq n-1$, $m_6 = m_2 + m_3\leq n-1$ and $m_4 = m_1 + m_2 + m_3 \leq n$.
\end{proposition}
%
\begin{proposition}\label{mat-equiv2}
{\rm\bf (1)}There is an isomorphism in $\HMF^{gr}_{\Q ,0}$
\begin{eqnarray}
\nonumber
\c\Big( \ \input{figure/figcircle-weight-r-mf} \ \Big)_{n} &=& \mathop{\boxtimes}_{j=1}^{m} K\Big( L^{[m]}_{j,(1;2)} ;x_{j,1}-x_{j,2} \Big)_{R_{(1,2)}^{(m,m)}}\Big|_{\mathbb{X}^{(m)}_{(2)} = \mathbb{X}^{(m)}_{(1)}}\\[-0.1em]
\nonumber
&\simeq& (\, J_{F_{m}(\mathbb{X}^{(m)}_{(1)})} \to 0\to J_{F_{m}(\mathbb{X}^{(m)}_{(1)})}  \,)\left\{ -mn+m^2 \right\} \left< m \right>,
\end{eqnarray}
where $J_{F_{m}(\mathbb{X}^{(m)}_{(1)})}$ is Jacobi algebra for the polynomial $F_{m}(\mathbb{X}^{(m)}_{(1)})$,
$$
J_{F_{m}(\mathbb{X}^{(m)}_{(1)})}=R_{(1)}^{(m)} \left/ 
\left< \frac{\partial F_{m}}{\partial x_{1,1}}, \ldots , \frac{\partial F_{m}}{\partial x_{m,1}} \right>\right. .
$$
{\rm\bf (2)}There is an isomorphism in 
$\HMF^{gr}_{R_{(1,2)}^{(m_3,m_3)} ,F_{m_3}(\mathbb{X}^{(m_3)}_{(1)})-F_{m_3}(\mathbb{X}^{(m_3)}_{(2)})}$
$$
\c\left( \input{figure/fig-bubble-color-mf1}\right)_n \simeq 
\c\left( \input{figure/fig-line-color-mf}\right)_n\left\{ \left[m_3\atop m_1\right]_q \right\}_{q},
$$
{\rm\bf (3)}There is an isomorphism in 
$\HMF^{gr}_{R_{(1,2)}^{(m_1,m_1)},F_{m_1}(\mathbb{X}^{(m_1)}_{(1)})-F_{m_1}(\mathbb{X}^{(m_1)}_{(2)})}$
$$
\c\left( \input{figure/fig-bubble-color1-mf1}\right)_n 
\simeq \c\left( \input{figure/fig-line-color1-mf}\right)_n \left\{ \left[n-m_1\atop m_2\right]_q \right\}_{q}\left< m_2 \right> ,
$$ 
where $1 \leq m_1,m_2 \leq n-1$ and $m_3 = m_1 + m_2 \leq n$.
\end{proposition}
%
\begin{proposition}\label{mat-equiv3}
There are isomorphisms in $\HMF^{gr}_{R_{(1,2,3,4)}^{(1,j,1,j)},F_1(\mathbb{X}^{(1)}_{(1)})+F_j(\mathbb{X}^{(j)}_{(2)})-F_1(\mathbb{X}^{(1)}_{(3)})-F_j(\mathbb{X}^{(j)}_{(4)})}$
\begin{eqnarray}
\nonumber
{\bf(1)}\c\left( \input{figure/figsquare1j--1--2j-1--1--1j-mf1}\right)_n &\simeq & \c\left( \input{figure/figsquare1j--j+1--1j-mf}\right)_n\bigoplus \c\left( \input{figure/figsquare1j-mf}\right)_n\{[m-1]_q\}_{q}\\[-0.1em]
\nonumber
{\bf(2)}\c\left( \input{figure/figsquare1j--j+1--j1--j+1--1j-rev-mf1}\right)_n &\simeq & \c\left( \input{figure/figsquare1j-rev-mf}\right)_n\bigoplus\c\left( \input{figure/figsquare1j--j-1--1j-rev-mf}\right)_n\{[n-m-1]_q\}_{q} \left<1\right>.
\end{eqnarray}
\end{proposition}
\begin{proof}[{\bf Proof of Proposition \ref{mat-equiv1}}]
We prove this proposition $(1)$.
The left-hand side factorization forms
\begin{eqnarray}
\nonumber
&&\overline{\Lambda}^{[m_1+m_2,m_3]}_{(4;5,3)}
\boxtimes
\overline{\Lambda}^{[m_1,m_2]}_{(5;1,2)}\\
\nonumber
&=& \mathop{\boxtimes}_{j=1}^{m_1+m_2+m_3} 
K\Big( \Lambda_{j,(4;5,3)}^{[m_1+m_2,m_3]} ;x_{j,4}-X^{(m_1+m_2,m_3)}_{j,(5,3)} \Big)_{R_{(3,4,5)}^{(m_3,m_1+m_2+m_3,m_1+m_2)}}\\
\nonumber
&&\hspace{1cm}
\boxtimes\mathop{\boxtimes}_{j=1}^{m_1+m_2}
K\Big( \Lambda_{j,(5;1,2)}^{[m_1,m_2]} ;x_{j,5}-X^{(m_1,m_2)}_{j,(1,2)} \Big)_{R_{(1,2,5)}^{(m_1,m_2,m_1+m_2)}}.
\end{eqnarray}
Since the potential of this factorization does not include the variables of $\mathbb{X}^{(m_1+m_2)}_{(5)}$ and 
\begin{equation}
\nonumber
(x_{1,5}-X^{(m_1,m_2)}_{1,(1,2)},\cdots ,x_{m_1+m_2,5}-X^{(m_1,m_2)}_{m_1+m_2,(1,2)})|_{(\mathbb{X}^{(m_1)}_{(1)},\mathbb{X}^{(m_2)}_{(2)},\mathbb{X}^{(m_3)}_{(3)},\mathbb{X}^{(m_1+m_2+m_3)}_{(4)})=(\underline{0})}=(\mathbb{X}^{(m_1+m_2)}_{(5)})
\end{equation}
is a regular sequence, we can apply Corollary \ref{cor2-11} to the variables of $\mathbb{X}^{(m_1+m_2)}_{(5)}$. Then, the matrix factorization is isomorphic to
\begin{equation}
\nonumber
\mathop{\boxtimes}_{j=1}^{m_1+m_2+m_3} K\Big( \Lambda_{j,(4;5,3)}^{[m_1+m_2,m_3]} ;x_{j,4}-X^{(m_1+m_2,m_3)}_{j,(5,3)} \Big)_{R_{(1,2,3,4,5)}^{(m_1,m_2,m_3,m_1+m_2+m_3,m_1+m_2)}/
\left< x_{1,5}-X^{(m_1,m_2)}_{1,(1,2)},\cdots ,x_{m_1+m_2,5}-X^{(m_1,m_2)}_{m_1+m_2,(1,2)}\right> }. 
\end{equation}
In the quotient $R_{(1,2,3,4,5)}^{(m_1,m_2,m_3,m_1+m_2+m_3,m_1+m_2)}/\left< x_{1,5}-X^{(m_1,m_2)}_{1,(1,2)},\cdots ,x_{m_1+m_2,5}-X^{(m_1,m_2)}_{m_1+m_2,(1,2)} \right>$, the polynomial $X^{(m_1+m_2,m_3)}_{j,(5,3)}$ equals $X^{(m_1,m_2,m_3)}_{j,(1,2,3)}$.
Then, the polynomial $\Lambda_{j,(4;5,3)}^{[m_1+m_2,m_3]}$ equals to
$$
\frac{F_{m_1+m_2+m_3}(\cdots ,X^{(m_1,m_2,m_3)}_{j-1,(1,2,3)},x_{j,4},x_{j+1,4},\cdots )-F_{m_1+m_2+m_3}(\cdots ,X^{(m_1,m_2,m_3)}_{j-1,(1,2,3)},X^{(m_1,m_2,m_3)}_{j,(1,2,3)},x_{j+1,4},\cdots )}{x_{j,4}-X^{(m_1,m_2,m_3)}_{j,(1,2,3)}}.
$$
We denote this polynomial by $\Lambda_{j,(4;1,2,3)}^{[m_1,m_2,m_3]}$.  
Since we find that
\begin{equation}
\nonumber
R_{(1,2,3,4,5)}^{(m_1,m_2,m_3,m_1+m_2+m_3,m_1+m_2)}/\left< x_{1,5}-X^{(m_1,m_2)}_{1,(1,2)},\cdots ,x_{m_1+m_2,5}-X^{(m_1,m_2)}_{m_1+m_2,(1,2)} \right> \simeq R_{(1,2,3,4)}^{(m_1,m_2,m_3,m_1+m_2+m_3)}
\end{equation}
as a $\Z$-graded $R_{(1,2,3,4)}^{(m_1,m_2,m_3,m_1+m_2+m_3)}$-module, the matrix factorization is isomorphic to
$$
\mathop{\boxtimes}_{j=1}^{m_1+m_2+m_3} K\Big( \Lambda_{j,(4;1,2,3)}^{[m_1,m_2,m_3]} ;x_{j,4}-X^{(m_1,m_2,m_3)}_{j,(1,2,3)} \Big)_{R_{(1,2,3,4)}^{(m_1,m_2,m_3,m_1+m_2+m_3)}}.
$$
\indent
The right-hand side factorization forms
\begin{eqnarray}
\nonumber
&&\overline{\Lambda}^{[m_1,m_2+m_3]}_{(4;1,6)}
\boxtimes
\overline{\Lambda}^{[m_2,m_3]}_{(6;2,3)}\\
\nonumber
&=& \mathop{\boxtimes}_{j=1}^{m_1+m_2+m_3}
K\Big( \Lambda_{j,(4;1,6)}^{[m_1,m_2+m_3]} ;x_{j,4}-X^{(m_1,m_2+m_3)}_{j,(1,6)} \Big)_{R_{(1,4,6)}^{(m_1,m_1+m_2+m_3,m_2+m_3)}}\\
\nonumber
&&\hspace{1cm}\boxtimes\mathop{\boxtimes}_{j=1}^{m_2+m_3}
K\Big( \Lambda_{j,(6;2,3)}^{[m_2,m_3]} ;x_{j,6}-X^{(m_2,m_3)}_{j,(2,3)} \Big)_{R_{(2,3,6)}^{(m_2,m_3,m_2+m_3)}} 
\end{eqnarray}
Since the potential of this matrix factorization does not include the variables of $\mathbb{X}_{m_2+m_3,6}$ and
\begin{equation}
\nonumber
(x_{1,6}-X^{(m_2,m_3)}_{1,(2,3)},\cdots ,x_{m_2+m_3,6}-X^{(m_2,m_3)}_{m_2+m_3,(2,3)})|_{(\mathbb{X}_{m_1,1},\mathbb{X}_{i_2,2},\mathbb{X}_{i_3,3},\mathbb{X}_{m_1+m_2+m_3,4})=(\underline{0})}=(\mathbb{X}_{m_2+m_3,6})
\end{equation}
is a regular sequence in $R_{(1,2,3,4,6)}^{(m_1,m_2,m_3,m_1+m_2+m_3,m_2+m_3)}$, we can apply Corollary \ref{cor2-11} to these variables. We similarly obtain the result that the matrix factorization is isomorphic to
$$
\mathop{\boxtimes}_{j=1}^{m_1+m_2+m_3} K\Big( \Lambda_{j,(4;1,2,3)}^{[m_1,m_2,m_3]} ;x_{j,4}-X^{(m_1,m_2,m_3)}_{j,(1,2,3)}\Big)_{R_{(1,2,3,4)}^{(m_1,m_2,m_3,m_1+m_2+m_3)}}. 
$$
\indent
We similarly prove this proposition $(2)$.
\end{proof}
\begin{proof}[{\bf Proof of Proposition \ref{mat-equiv2}}]
$(1)$
We have
\begin{equation}
\nonumber
\mathop{\boxtimes}_{j=1}^{m} 
K\Big( L^{[m]}_{j,(1;2)} ;x_{j,1}-x_{j,2} \Big)_{R_{(1,2)}^{(m,m)}}\Big|_{\mathbb{X}_{m,2} = \mathbb{X}_{m,1}} =
\mathop{\boxtimes}_{j=1}^{m} 
\left( R_{(1)}^{(m)},R_{(1)}^{(m)}\{ 2j-1-n \},L^{[m]}_{j,(1;2)}|_{\mathbb{X}_{m,2} = \mathbb{X}_{m,1}},0 \right).
\end{equation}
The polynomial $L^{[m]}_{j,(1;2)}|_{\mathbb{X}_{m,2} = \mathbb{X}_{m,1}}$ is
$$
\frac{F_{m}(\cdots ,x_{j-1,2},x_{j,1},x_{j+1,1},\cdots )-F_{m}(\cdots ,x_{j-1,2},x_{j,2},x_{j+1,1},\cdots )}{x_{j,1}-x_{j,2}}
\Big|_{\mathbb{X}_{m,2} = \mathbb{X}_{m,1}} =\frac{\partial F_{m}(\mathbb{X}_{m,1})}{\partial x_{j,1}}.
$$
Hence, we apply Theorem \ref{exclude} to these polynomials of the matrix factorization after using Proposition \ref{functor1} and \ref{functor2};
\begin{eqnarray}
\nonumber
&&\mathop{\boxtimes}_{j=1}^{m} 
K\Big( L^{[m]}_{j,(1;2)} ;x_{j,1}-x_{j,2} \Big)_{R_{(1,2)}^{(m,m)}}\Big|_{\mathbb{X}_{m,2}=\mathbb{X}_{m,1}}\\
\nonumber
&\simeq&
\mathop{\boxtimes}_{j=1}^{m} 
\left( R_{(1)}^{(m)},R_{(1)}^{(m)}\{ n+1-2j \},0,\frac{\partial F_{m}(\mathbb{X}_{m,1})}{\partial x_{j,1}}\right) \{ -m n + m^2 \} \left< m \right> \\[-0.1em]
\nonumber
&\simeq&\left(J_{F_{m}(\mathbb{X}_{m,1})},0,0,0\right)
\left\{ -mn+m^2 \right\} \left< m \right>.
\end{eqnarray}
$(2)$ 
We have
\begin{eqnarray}
\nonumber
&&\c\left( \input{figure/fig-bubble-color-mf}\right)_n \\[-0.1em]
\nonumber
&=&
\mathop{\boxtimes}^{m_3}_{j=1} K\Big( \Lambda_{j,(1;3,4)}^{[m_1,m_2]} ;x_{j,1}-X^{(m_1,m_2)}_{j,(3,4)} \Big)_{R_{(1,3,4)}^{(m_3,m_1,m_2)}} 
\boxtimes\mathop{\boxtimes}^{m_3}_{j=1} K\Big( V_{j,(3,4;2)}^{[m_1,m_2]} ;X^{(m_1,m_2)}_{j,(3,4)}-x_{j,2} \Big)_{R_{(2,3,4)}^{(m_3,m_1,m_2)}}
\left\{ -m_1 m_2 \right\} .
\end{eqnarray}
The potential of this matrix factorization does not include the variables of $\mathbb{X}^{(m_1)}_{(3)}$, $\mathbb{X}^{(m_2)}_{(4)}$ and we find that the following sequence is regular:
\begin{equation}
\nonumber
(X^{(m_1,m_2)}_{1,(3,4)}-x_{1,2},... ,X^{(m_1,m_2)}_{m_3,(3,4)}-x_{m_3,2})
\left|_{(\mathbb{X}_{m_{3},1},\mathbb{X}_{m_{3},2})=(\underline{0})}\right.=(X^{(m_1,m_2)}_{1,(3,4)},... ,X^{(m_1,m_2)}_{m_3,(3,4)}).
\end{equation}
Therefore, we can apply Corollary \ref{cor2-11} to the matrix factorization. 
Then, we have
\begin{equation}
\nonumber
\mathop{\boxtimes}^{m_3}_{j=1} K\Big( \Lambda_{j,(1;3,4)}^{[m_1,m_2]} ;x_{j,1}-x_{j,2} \Big)_{R_{(1,2,3,4)}^{(m_3,m_3,m_1,m_2)}/\left< X^{(m_1,m_2)}_{1,(3,4)}-x_{1,2},... ,X^{(m_1,m_2)}_{m_3,(3,4)}-x_{m_3,2} \right> } \left\{ -m_1 m_2 \right\}.
\end{equation}
In the quotient $R_{(1,2,3,4)}^{(m_3,m_3,m_1,m_2)}/\left< X^{(m_1,m_2)}_{1,(3,4)}-x_{1,2},... ,X^{(m_1,m_2)}_{m_3,(3,4)}-x_{m_3,2}\right>$, the polynomial $\Lambda_{j,(1;3,4)}^{[m_1,m_2]}$ is equal to $L_{j,(1;2)}^{[m_3]}$.
We find that the quotient $R_{(1,2,3,4)}^{(m_3,m_3,m_1,m_2)}/\left< X^{(m_1,m_2)}_{1,(3,4)}-x_{1,2},... ,X^{(m_1,m_2)}_{m_3,(3,4)}-x_{m_3,2}\right> \left\{ -m_1 m_2 \right\}$ is isomorphic to $R_{(1,2)}^{(m_3,m_3)}\left\{\left[m_3\atop m_1\right]_q\right\}_{q}$ as a $\Z$-graded $R_{(1,2)}^{(m_3,m_3)}$-module.
Thus, we obtain the isomorphism of $(2)$.\\
\noindent
$(3)$ We have
\begin{eqnarray}
\nonumber
&&\c\left( \input{figure/fig-bubble-color1-mf}\right)_n \\[-0.1em]
\nonumber
&=&\mathop{\boxtimes}_{j=1}^{m_3} K(\Lambda_{j,(3;2,4)}^{[m_1,m_2]};x_{j,3}-X_{j,(2,4)}^{(m_1,m_2)})_{R_{(2,3,4)}^{(m_1,m_3,m_2)}}
\boxtimes\mathop{\boxtimes}_{j=1}^{m_3} K(V_{j,(1,4;3)}^{[m_1,m_2]};X_{j,(1,4)}^{(m_1,m_2)}-x_{j,3})_{R_{(1,3,4)}^{(m_1,m_3,m_2)}}\left\{ -m_1 m_2 \right\}.
\end{eqnarray}
The potential of this matrix factorization does not include the variables of $\mathbb{X}^{(m_3)}_{(3)}$, $\mathbb{X}^{(m_2)}_{(4)}$ and we find that the following sequence is regular:
\begin{equation}
\nonumber
(X_{1,(1,4)}^{(m_1,m_2)}-x_{1,3},... ,X_{m_3,(1,4)}^{(m_1,m_2)}-x_{m_3,3})\left|_{(\mathbb{X}^{(m_1)}_{(1)},\mathbb{X}^{(m_1)}_{(2)})=
(\underline{0})}\right. 
=(x_{1,4}-x_{1,3},...,x_{m_2,4}-x_{m_2,3},-x_{m_2+1,3},...,-x_{m_3,3}).
\end{equation}
is regular.
Therefore, we can apply Corollary \ref{cor2-11} 
to the matrix factorization.
Then we have
\begin{eqnarray}
\nonumber
&&\mathop{\boxtimes}_{j=1}^{m_3}K(\Lambda_{j,(3;2,4)}^{[m_1,m_2]};X_{j,(1,4)}^{(m_1,m_2)}-X_{j,(2,4)}^{(m_1,m_2)})_{R_{(1,2,3,4)}^{(m_1,m_1,m_3,m_2)}/\left< X_{1,(1,4)}^{(m_1,m_2)}-x_{1,3} ,... ,X_{m_3,(1,4)}^{(m_1,m_2)}-x_{m_3,3} \right> }\left\{ -m_1 m_2 \right\}\\[-0.1em]
\nonumber
\label{mf1}&\simeq&\mathop{\boxtimes}_{j=1}^{m_3}K(\widetilde{\Lambda_{j,(3;2,4)}^{[m_1,m_2]}}\hspace{1mm};X_{j,(1,4)}^{(m_1,m_2)}-X_{j,(2,4)}^{(m_1,m_2)})_{R_{(1,2,4)}^{(m_1,m_1,m_2)}}\left\{ -m_1 m_2 \right\} ,
\end{eqnarray}
where 
\begin{eqnarray}
\nonumber
&&\hspace{-0.1cm}\widetilde{\Lambda_{j,(3;2,4)}^{[m_1,m_2]}}\\
\nonumber
&&\hspace{-0.1cm}=\frac{F_{m_3}(X_{1,(2,4)}^{(m_1,m_2)}\hspace{-0.1cm},... ,X_{j-1,(2,4)}^{(m_1,m_2)},X_{j,(1,4)}^{(m_1,m_2)}\hspace{-0.1cm},... ,X_{m_3,(1,4)}^{(m_1,m_2)})-
F_{m_3}(X_{1,(2,4)}^{(m_1,m_2)}\hspace{-0.1cm},... ,X_{j,(2,4)}^{(m_1,m_2)},X_{j+1,(1,4)}^{(m_1,m_2)},... ,X_{m_3,(1,4)}^{(m_1,m_2)}
)}{X_{j,(1,4)}^{(m_1,m_2)}-X_{j,(2,4)}^{(m_1,m_2)}}.
\end{eqnarray}
$X_{j,(1,4)}^{(m_1,m_2)}-X_{j,(2,4)}^{(m_1,m_2)}$ is a polynomial with $\Z$-grading $2j$ of
\begin{equation}
\nonumber
((x_{1,1}-x_{1,2})+(x_{2,1}-x_{2,2})+... +(x_{m_1,1}-x_{m_1,2}))(1+x_{1,4}+x_{2,4}+... +x_{m_2,4}).
\end{equation}
Then, the polynomials $(X_{m_1+1,(1,4)}^{(m_1,m_2)}-X_{m_1+1,(2,4)}^{(m_1,m_2)},...,X_{m_3,(1,4)}^{(m_1,m_2)}-X_{m_3,(2,4)}^{(m_1,m_2)})$ can be described as the linear sum of the polynomials $(X_{1,(1,4)}^{(m_1,m_2)}-X_{1,(2,4)}^{(m_1,m_2)},...,X_{m_1,(1,4)}^{(m_1,m_2)}-X_{m_1,(2,4)}^{(m_1,m_2)})$.
Applying Proposition \ref{equiv} to the matrix factorization $(\ref{mf1})$, it is isomorphic to
\begin{eqnarray}
\label{factorization1}
&&\mathop{\boxtimes}_{j=1}^{m_1}K(\ast ;x_{j,1}-x_{j,2})_{R_{(1,2,4)}^{(m_1,m_1,m_2)}}\\[-0.1em]
\nonumber
&&\boxtimes\mathop{\boxtimes}_{k=m_1+1}^{m_3}
(R_{(1,2,4)}^{(m_1,m_1,m_2)},R_{(1,2,4)}^{(m_1,m_1,m_2)}\{2k-n-1\},
\widetilde{\Lambda_{k,(3;2,4)}^{[m_1,m_2]}},0)\left\{ -m_1 m_2 \right\}.
\end{eqnarray}
We find that the following sequence is regular:
\begin{equation}
\nonumber
(\widetilde{\Lambda_{m_1+1,(3;2,4)}^{[m_1,m_2]}},...,\widetilde{\Lambda_{m_3,(3;2,4)}^{[m_1,m_2]}})
\left|_{(\mathbb{X}^{(m_1)}_{(1)},\mathbb{X}^{(m_1)}_{(2)})=(\underline{0})}\right.
=((-1)^{m_1}(n+1)X_{n-m_1,(4)}^{(-m_3)},...,(-1)^{m_3-1}(n+1)X_{n+1-m_3,(4)}^{(-m_3)}).
\end{equation}
The potential of the factorization (\ref{factorization1}) does not include the variables of $\mathbb{X}^{(m_2)}_{(4)}$.
Then, the partial factorization of (\ref{factorization1}), $\mathop{\boxtimes}_{k=m_1+1}^{m_3}
(R_{(1,2,4)}^{(m_1,m_1,m_2)},R_{(1,2,4)}^{(m_1,m_1,m_2)}\{2k-n-1\},
\widetilde{\Lambda_{k,(3;2,4)}^{[m_1,m_2]}},0)\left\{ -m_1 m_2 \right\}$, is isomorphic to
\begin{equation}
\nonumber
\left( R_{(1,2,4)}^{(m_1,m_1,m_2)}/\langle X_{n-m_1,(4)}^{(-m_3)},...,X_{n+1-m_3,(4)}^{(-m_3)}\rangle,0,0,0\right)
\left\{ \sum_{m_1+1}^{m_3}2k-n-1\right\}\left\{ -m_1 m_2 \right\}\langle m_2\rangle .
\end{equation}
As a factorization of $\Z$-graded $R_{(1,2)}^{(m_1,m_1)}$-modules, this is isomorphic to
\begin{equation}
\nonumber
\left( R_{(1,2)}^{(m_1,m_1)}\left\{\left[n-m_1 \atop m_2\right]_q\right\}_{q},0,0,0\right)
\langle m_2\rangle .
\end{equation}
By Theorem \ref{reg-eq}, the other partial matrix factorization of (\ref{factorization1}), $\displaystyle \mathop{\boxtimes}_{j=1}^{m_1}K(\ast ;x_{j,1}-x_{j,2})_{R_{(1,2,4)}^{(m_1,m_1,m_2)}}$, is isomorphic to
\begin{equation}
\nonumber
\c\left( \input{figure/fig-line-color1-mf}\right)_n\boxtimes (R_{(4)}^{(m_2)},0,0,0).
\end{equation} 
Hence, we obtain the isomorphism of $(3)$.
\end{proof}
\begin{proof}[{\bf Proof of Proposition \ref{mat-equiv3}}]
$(1)$We have
\begin{eqnarray*}
&&\c\left( \input{figure/figsquare1j--1--2j-1--1--1j-mf}\right)_n\\[-0.1em]
&=&
K\left(\hspace{-0.1cm}
\left( V^{[1,1]}_{1,(1,6;5)} \atop V^{[1,1]}_{2,(1,6;5)}\right);
\left(X_{1,(1,6)}^{(1,1)}-x_{1,5} \atop X_{2,(1,6)}^{(1,1)}-x_{2,5}\right)\hspace{-0.1cm}
\right)_{R_{(1,5,6)}^{(1,2,1)}}\hspace{-1cm}\{-1\}
\boxtimes 
K\left(\hspace{-0.1cm}
\left(
\Lambda^{[1,1]}_{1,(5;3,8)} \atop \Lambda^{[1,1]}_{2,(5;3,8)}
\right)\hspace{-0.1cm};\hspace{-0.1cm}
\left(
x_{1,5}-X_{1,(3,8)}^{(1,1)} \atop x_{2,5}-X_{2,(3,8)}^{(1,1)}
\right)\hspace{-0.1cm}
\right)_{R_{(3,5,8)}^{(1,2,1)}}\\[-0.1em]
&&\boxtimes
K\left(\hspace{-0.1cm}
\left( 
\begin{array}{c}
\Lambda^{[1,m-1]}_{1,(2;6,7)}\\[.5em]
\vdots\\[.5em]
\Lambda^{[1,m-1]}_{m,(2;6,7)}
\end{array}
\right)\hspace{-0.1cm};\hspace{-0.1cm}
\left(
\begin{array}{c}
x_{1,2}-X_{1,(6,7)}^{(1,m-1)}\\[.6em]
\vdots\\[.6em]
x_{m,2}-X_{m,(6,7)}^{(1,m-1)}
\end{array}
\right)\hspace{-0.1cm}
\right)_{R_{(2,6,7)}^{(m,1,m-1)}}
\hspace{-1cm}
\boxtimes
\hspace{0.5cm}
K\left(\hspace{-0.1cm}
\left( 
\begin{array}{c}
V_{1,(8,7;4)}^{[1,m-1]}\\[.5em]
\vdots\\[.5em]
V_{m,(8,7;4)}^{[1,m-1]}
\end{array}
\right)\hspace{-0.1cm};\hspace{-0.1cm}
\left(
\begin{array}{c}
X_{1,(7,8)}^{(m-1,1)}-x_{1,4}\\[.6em]
\vdots\\[.6em]
X_{m,(7,8)}^{(m-1,1)}-x_{m,4}
\end{array}
\right)\hspace{-0.1cm}
\right)_{R_{(2,7,8)}^{(m,m-1,1)}}\hspace{-1.5cm}\{-m+1\}.
\end{eqnarray*}
We apply Corollary \ref{cor2-11} to the matrix factorization. Then we obtain
\begin{equation}
\label{factorization2}
K\left(
\left( 
\begin{array}{c}
\Lambda^{[1,1]}_{2,(5;3,8)}\\[.5em]
\Lambda^{[1,m-1]}_{m,(2;6,7)}\\[.5em]
V_{1,(8,7;4)}^{[1,m-1]}\\[.5em]
\vdots\\[.5em]
V_{m,(8,7;4)}^{[1,m-1]}
\end{array}
\right);
\left(
\begin{array}{c}
x_{2,5}-X_{2,(3,8)}^{(1,1)}\\[.6em]
x_{m,2}-X_{m,(6,7)}^{(1,m-1)}\\[.6em]
X_{1,(7,8)}^{(m-1,1)}-x_{1,4}\\[.6em]
\vdots\\[.6em]
X_{m,(7,8)}^{(m-1,1)}-x_{m,4}
\end{array}
\right)
\right)_{Q_1}\{-m\},
\end{equation}
where $Q_1=R_{(1,2,3,4,5,6,7,8)}^{(1,m,1,m,2,1,m-1,1)}\left/
\left< 
\begin{array}{c}
X_{1,(1,6)}^{(1,1)}-x_{1,5}, X_{2,(1,6)}^{(1,1)}-x_{2,5}, x_{1,5}-X_{1,(3,8)}^{(1,1)}, \\
x_{1,2}-X_{1,(6,7)}^{(1,m-1)},  \cdots , x_{m,2}-X_{m,(6,7)}^{(1,m-1)} 
\end{array}
\right>\right.$.\\
In the quotient, there are the following equations

\begin{eqnarray*}
&&
x_{1,5}=X_{1,(1,6)}^{(1,1)},\,
x_{2,5}=X_{2,(1,6)}^{(1,1)},\,
x_{1,8}=X_{1,(1,3,6)}^{(1,-1,1)},\\[-0.1em]
&&
x_{k,7}=X_{k,(2,6)}^{(m,-1)}\hspace{1cm}(k=1,2,\cdots ,m-1).
\end{eqnarray*}
Therefore, $Q_1$ is isomorphic to $R_{(1,2,3,4,6)}^{(1,m,1,m,1)}$ as a $\Z$-graded $R_{(1,2,3,4)}^{(1,m,1,m)}$-module.
That is, the variables $x_{1,5}$, $x_{2,5}$, $x_{1,8}$ and $x_{k,7}$ can be removed from the quotient $Q_1$ using the above equations.
Then, the matrix factorization (\ref{factorization2}) is isomorphic to
\begin{equation*}
K\left(
\left( 
\begin{array}{c}
\widetilde{\Lambda^{[1,1]}_{2,(5;3,8)}}\\[.5em]
\widetilde{\Lambda^{[1,m-1]}_{m,(2;6,7)}}\\[.5em]
\widetilde{V^{[1,m-1]}_{1,(8,7;4)}}\\[.5em]
\widetilde{V^{[1,m-1]}_{2,(8,7;4)}}\\[.5em]
\vdots\\[.5em]
\widetilde{V^{[1,m-1]}_{m-1,(8,7;4)}}\\[.5em]
\widetilde{V^{[1,m-1]}_{m,(8,7;4)}}
\end{array}
\right);
\left(
\begin{array}{c}
X_{1,(1,3)}^{(1,-1)}(x_{1,6}-x_{1,3})\\[1em]
X_{m,(2,6)}^{(m,-1)}\\[1em]
x_{1,1}+x_{1,2}-x_{1,3}-x_{1,4}\\[1em]
X_{1,(1,3)}^{(1,-1)}X_{1,(2,6)}^{(m,-1)}+x_{2,2}-x_{2,4}\\[1em]
\vdots\\[1em]
X_{1,(1,3)}^{(1,-1)}X_{m-2,(2,6)}^{(m,-1)}+x_{m-1,2}-x_{m-1,4}\\[1em]
X_{1,(1,3)}^{(1,-1)}X_{m-1,(2,6)}^{(m,-1)}-X_{m,(2,6)}^{(m,-1)}+x_{m,2}-x_{m,4}
\end{array}
\right)
\right)_{R_{(1,2,3,4,6)}^{(1,m,1,m,1)}.}\{-m\}
\end{equation*}
Moreover, by Theorem \ref{exclude}, we obtain
\begin{equation*}
K\left(
\left( 
\begin{array}{c}
\widehat{\Lambda^{[1,1]}_{2,(5;3,8)}}\\[.5em]
\widehat{V^{[1,m-1]}_{1,(8,7;4)}}\\[.5em]
\widehat{V^{[1,m-1]}_{2,(8,7;4)}}\\[.5em]
\vdots\\[.5em]
\widehat{V^{[1,m-1]}_{m,(8,7;4)}}
\end{array}
\right);
\left(
\begin{array}{c}
X_{1,(1,3)}^{(1,-1)}(x_{1,6}-x_{1,3})\\[.9em]
x_{1,1}+x_{1,2}-x_{1,3}-x_{1,4}\\[.9em]
X_{1,(1,3)}^{(1,-1)}X_{1,(2,6)}^{(m,-1)}+x_{2,2}-x_{2,4}\\[.9em]
\vdots\\[.9em]
X_{1,(1,3)}^{(1,-1)}X_{m-1,(2,6)}^{(m,-1)}+x_{m,2}-x_{m,4}
\end{array}
\right)
\right)_{R_{(1,2,3,4,6)}^{(1,m,1,m,1)}\left/\left<X_{m,(2,6)}^{(m,-1)}\right>\right.}\hspace{-2cm}\{-m\},
\end{equation*}
where $\widehat{\Lambda^{[1,1]}_{2,(5;3,8)}}=\widetilde{\Lambda^{[1,1]}_{2,(5;3,8)}}$ and $\widehat{V^{[1,m-1]}_{i,(8,7;4)}}=\widetilde{V^{[1,m-1]}_{i,(8,7;4)}}$ in the quotient $R_{(1,2,3,4,6)}^{(1,m,1,m,1)}\left/\left<X_{m,(2,6)}^{(m,-1)}\right>\right.$.
Since the polynomials $X_{k,(2,6)}^{(m,-1)}$ are described as
\begin{eqnarray*}
X_{k,(2,6)}^{(m,-1)}&=&-X_{k-1,(2,3,6)}^{(m,-1,-1)}(x_{1,6}-x_{1,3})+X_{k,(2,3)}^{(m,-1)}\,\, (k=1,...,m-1),
\end{eqnarray*}
the above matrix factorization is isomorphic to
\begin{equation*}
K\left(
\left( 
\begin{array}{c}
\widehat{\Lambda^{[1,1]}_{2,(5;3,8)}}-\sum_{k=2}^{m}X_{k-1,(2,3,6)}^{(m,-1,-1)}\widehat{V^{[1,m-1]}_{k,(8,7;4)}}\\[.5em]
\widehat{V^{[1,m-1]}_{1,(8,7;4)}}\\[.5em]
\widehat{V^{[1,m-1]}_{2,(8,7;4)}}\\[.5em]
\vdots\\[.5em]
\widehat{V^{[1,m-1]}_{m,(8,7;4)}}
\end{array}
\right);
\left(
\begin{array}{c}
X_{1,(1,3)}^{(1,-1)}(x_{1,6}-x_{1,3})\\[.9em]
x_{1,1}+x_{1,2}-x_{1,3}-x_{1,4}\\[.9em]
X_{1,(1,3)}^{(1,-1)}X_{1,(2,3)}^{(m,-1)}+x_{2,2}-x_{2,4}\\[.9em]
\vdots\\[.9em]
X_{1,(1,3)}^{(1,-1)}X_{m-1,(2,3)}^{(m,-1)}+x_{m,2}-x_{m,4}
\end{array}
\right)
\right)_{R_{(1,2,3,4,6)}^{(1,m,1,m,1)}\left/\left<X_{m,(2,6)}^{(m,-1)}\right>\right.}\hspace{-2cm}\{-m\}.
\end{equation*}
Applying Corollary \ref{induce-sq1} (2) to this matrix factorization, we find that there are polynomials $a_k\in R_{(1,2,3,4)}^{(1,m,1,m)}$ $(k=1,\cdots m)$ and $a_0$ $\in$ $R_{(1,2,3,4,6)}^{(1,m,1,m,1)}\left/\left<X_{m,(2,6)}^{(m,-1)}\right>\right.$ such that $a_0(x_{1,6}-x_{1,3})\equiv a \in R_{(1,2,3,4)}^{(1,m,1,m)}$ $(\mod R_{(1,2,3,4,6)}^{(1,m,1,m,1)}\left/\left<X_{m,(2,6)}^{(m,-1)}\right>\right.)$ and we have an isomorphism between the above matrix factorization and the following matrix factorization
\begin{equation}
\label{total-mat1}
K\left(
\left( 
\begin{array}{c}
a_0\\[.3em]
a_1\\[.3em]
a_2\\[.3em]
\vdots\\[.3em]
a_m
\end{array}
\right);
\left(
\begin{array}{c}
X_{1,(1,3)}^{(1,-1)}(x_{1,6}-x_{1,3})\\[.1em]
x_{1,1}+x_{1,2}-x_{1,3}-x_{1,4}\\[.1em]
X_{1,(1,3)}^{(1,-1)}X_{1,(2,3)}^{(m,-1)}+x_{2,2}-x_{2,4}\\[.1em]
\vdots\\[.1em]
X_{1,(1,3)}^{(1,-1)}X_{m-1,(2,3)}^{(m,-1)}+x_{m,2}-x_{m,4}
\end{array}
\right)
\right)_{R_{(1,2,3,4,6)}^{(1,m,1,m,1)}\left/\left<X_{m,(2,6)}^{(m,-1)}\right>\right.}\hspace{-2cm}\{-m\}.
\end{equation}
We choose a decompositions of $R_{(1,2,3,4,6)}^{(1,m,1,m,1)}\left/\left<X_{m,(2,6)}^{(m,-1)}\right>\right.$ to be
\begin{eqnarray}
\nonumber
R_1&\simeq& 
R_{(1,2,3,4)}^{(1,m,1,m)}\oplus
(x_{1,6}-x_{1,3})R_{(1,2,3,4)}^{(1,m,1,m)}\oplus
x_{1,6}(x_{1,6}-x_{1,3})R_{(1,2,3,4)}^{(1,m,1,m)}\oplus
\cdots\oplus
x_{1,6}^{m-2}(x_{1,6}-x_{1,3})R_{(1,2,3,4)}^{(1,m,1,m)},\\
\nonumber
R_2&\simeq& 
R_{(1,2,3,4)}^{(1,m,1,m)}\oplus
x_{1,6}R_{(1,2,3,4)}^{(1,m,1,m)}\oplus
x_{1,6}^2R_{(1,2,3,4)}^{(1,m,1,m)}\oplus
\cdots\oplus
x_{1,6}^{m-1}R_{(1,2,3,4)}^{(1,m,1,m)}.
\end{eqnarray}
Then, the partial factorization $K(a_0;X_{1,(1,3)}^{(1,-1)}(x_{1,6}-x_{1,3}))_{R_{(1,2,3,4,6)}^{(1,m,1,m,1)}\left/\left<X_{m,(2,6)}^{(m,-1)}\right>\right.}$ of (\ref{total-mat1}) is isomorphic to
\begin{eqnarray}
\nonumber
\left(R_1,R_2\{3-n\},f_0,f_1\right).
\end{eqnarray}
where 
\begin{eqnarray*}
f_0&=&\left(
\begin{array}{cc}
{}^t\mathfrak{0}_{m-1}&E_{m-1}(a)\\[0.1em]
\displaystyle\frac{a_0}{X_{m-1,(2,3,6)}^{(m-1,-1,-1)}}&\mathfrak{0}_{m-1}
\end{array}
\right),\\
f_1&=&\left(
\begin{array}{cc}
\mathfrak{0}_{m-1}&X_{m,(2,3)}^{(m-1,-1)}(x_{1,1}-x_{1,3})\\[0.1em]
E_{m-1}(x_{1,1}-x_{1,3})&{}^{t}\mathfrak{0}_{m-1}
\end{array}
\right).
\end{eqnarray*}
$E_{k}(f)$ is the $k\times k$ diagonal matrix of $f$ and $\mathfrak{0}_k$ is the zero low vector with length $m$.
Remark that $\displaystyle\frac{a_0}{X_{m-1,(2,3,6)}^{(m-1,-1,-1)}}$ naturally become a polynomial of $R_{(1,2,3,4)}^{(1,m,1,m)}$ in the quotient $R_{(1,2,3,4,6)}^{(1,m,1,m,1)}\left/\left<X_{m,(2,6)}^{(m,-1)}\right>\right.$ by the structure of matrix factorization.
This polynomial $\displaystyle\frac{a_0}{X_{m-1,(2,3,6)}^{(m-1,-1,-1)}}$ in $R_{(1,2,3,4)}^{(1,m,1,m)}$ is denoted by $b$.\\
\indent
Hence, the partial matrix factorization splits into the following direct sum
\begin{eqnarray*}
&&\bigoplus_{i=1}^{m-1}K(a;X_{1,(1,3)}^{(1,1)})_{R_{(1,2,3,4)}^{(1,m,1,m)}}\{2i\}
\oplus K(b;X_{1,(1,3)}^{(1,1)}X_{m,(2,3)}^{(m-1,-1)})_{R_{(1,2,3,4)}^{(1,m,1,m)}}.
\end{eqnarray*}
The other partial factorization of (\ref{total-mat1}) consists of polynomials which do not include variable $x_{1,6}$.
Thus, using Theorem \ref{reg-eq}, the factorization (\ref{total-mat1}) is isomorphic to
\begin{eqnarray*}
&&\hspace{-1cm}\bigoplus_{i=1}^{m-1}K\left(
\left( 
\begin{array}{c}
a\\[.1em]
a_1\\[.1em]
a_2\\[.1em]
\vdots\\[.1em]
a_m
\end{array}
\right);
\left(
\begin{array}{c}
X_{1,(1,3)}^{(1,1)}\\[.1em]
x_{1,1}+x_{1,2}-x_{1,3}-x_{1,4}\\[.1em]
X_{1,(1,3)}^{(1,1)}X_{1,(2,3)}^{(m-1,-1)}+x_{2,2}-x_{2,4}\\[.1em]
\vdots\\[.1em]
X_{1,(1,3)}^{(1,1)}X_{m-1,(2,3)}^{(m-1,-1)}+x_{m,2}-x_{m,4}
\end{array}
\right)
\right)_{R_{(1,2,3,4)}^{(1,m,1,m)}}\{2i-m\}\\[-0.1em]
&&\bigoplus
K\left(
\left( 
\begin{array}{c}
a_1\\[.1em]
a_2\\[.1em]
\vdots\\[.1em]
a_m\\[.1em]
b
\end{array}
\right);
\left(
\begin{array}{c}
x_{1,1}+x_{1,2}-x_{1,3}-x_{1,4}\\[.1em]
X_{1,(1,3)}^{(1,1)}X_{1,(2,3)}^{(m-1,-1)}+x_{2,2}-x_{2,4}\\[.1em]
\vdots\\[.1em]
X_{1,(1,3)}^{(1,1)}X_{m-1,(2,3)}^{(m-1,-1)}+x_{m,2}-x_{m,4}\\[.1em]
X_{1,(1,3)}^{(1,1)}X_{m,(2,3)}^{(m-1,-1)}
\end{array}
\right)
\right)_{R_{(1,2,3,4)}^{(1,m,1,m)}}\{-m\}\\[-0.1em]
&&\hspace{-1cm}\simeq
\bigoplus_{i=1}^{m-1}K\left(
\left( 
\begin{array}{c}
a+\sum_{l=1}^{m}X_{l-1,(2,3)}^{(m-1,-1)}a_l\\[.1em]
a_1\\[.1em]
\vdots\\[.1em]
a_m
\end{array}
\right);
\left(
\begin{array}{c}
x_{1,1}-x_{1,3}\\[.1em]
x_{1,2}-x_{1,4}\\[.1em]
\vdots\\[.1em]
x_{m,2}-x_{m,4}
\end{array}
\right)
\right)_{R_{(1,2,3,4)}^{(1,m,1,m)}}\{2i-m\}\\[-0.1em]
&&\bigoplus
K\left(
\left( 
\begin{array}{c}
\sum_{k=1}^{m+1}a_k (-x_{1,3})^{k-1}\\[.1em]
\sum_{k=2}^{m+1}a_k (-x_{1,3})^{k-2}\\[.1em]
\vdots\\[.1em]
a_m-x_{1,3}a_{m+1}\\[.1em]
b
\end{array}
\right);
\left(
\begin{array}{c}
X_{1,(1,2)}^{(1,m)}-X_{1,(3,4)}^{(1,m)}\\[.1em]
X_{2,(1,2)}^{(1,m)}-X_{2,(3,4)}^{(1,m)}\\[.1em]
\vdots\\[.1em]
X_{m,(1,2)}^{(1,m)}-X_{m,(3,4)}^{(1,m)}\\[.1em]
X_{m+1,(1,2)}^{(1,m)}-X_{m+1,(3,4)}^{(1,m)}
\end{array}
\right)
\right)_{R_{(1,2,3,4)}^{(1,m,1,m)}}\{-m\}.
\end{eqnarray*}
Using Theorem \ref{reg-eq}, we find the above matrix factorization is isomorphic to
\begin{equation*}
\c\left( \input{figure/figsquare1j--j+1--1j-mf}\right)\oplus \c\left( \input{figure/figsquare1j-mf}\right)_n\{[m-1]_q\}_{q}
\end{equation*}
$(2)$
We have
\begin{eqnarray*}
&&\hspace{-1cm}\c\left( \input{figure/figsquare1j--j+1--j1--j+1--1j-rev-mf}\right)_n=\\
&&\hspace{-1cm}K\left(\hspace{-0.1cm}\left(\hspace{-0.1cm}
\begin{array}{c}
\Lambda_{1,(6;1,5)}^{[1,m]}\\[.5em]
\vdots\\[.5em]
\Lambda_{m+1,(6;1,5)}^{[1,m]}
\end{array}\hspace{-0.1cm}
\right)\hspace{-0.1cm}
;\hspace{-0.1cm}
\left(\hspace{-0.1cm}
\begin{array}{c}
x_{1,6}-X_{1,(1,5)}^{(1,m)}\\[.6em]
\vdots\\[.6em]
x_{m+1,6}-X_{m+1,(1,5)}^{(1,m)}
\end{array}\hspace{-0.1cm}
\right)\hspace{-0.1cm}
\right)_{R_{(1,5,6)}^{(1,m,m+1)}}
\boxtimes
K\left(\hspace{-0.1cm}\left(\hspace{-0.1cm}
\begin{array}{c}
V_{1,(3,5;8)}^{[1,m]}\\[.5em]
\vdots\\[.5em]
V_{m+1,(3,5;8)}^{[1,m]}
\end{array}\hspace{-0.1cm}
\right)\hspace{-0.1cm}
;\hspace{-0.1cm}
\left(\hspace{-0.1cm}
\begin{array}{c}
X_{1,(3,5)}^{(1,m)}-x_{1,8}\\[.6em]
\vdots\\[.6em]
X_{m+1,(3,5)}^{(1,m)}-x_{m+1,8}
\end{array}\hspace{-0.1cm}
\right)\hspace{-0.1cm}
\right)_{R_{(3,5,8)}^{(1,m,m+1)}}\hspace{-1cm}\{-m\}\\
&&\hspace{-0.5cm}\boxtimes
K\left(\hspace{-0.1cm}\left(\hspace{-0.1cm}
\begin{array}{c}
V_{1,(7,2;6)}^{[1;m]}\\[.5em]
\vdots\\[.5em]
V_{m+1,(7,2;6)}^{[1;m]}
\end{array}\hspace{-0.1cm}
\right)\hspace{-0.1cm}
;\hspace{-0.1cm}
\left(\hspace{-0.1cm}
\begin{array}{c}
X_{1,(2,7)}^{(m,1)}-x_{1,6}\\[.6em]
\vdots\\[.6em]
X_{m+1,(2,7)}^{(m,1)}-x_{m+1,6}
\end{array}\hspace{-0.1cm}
\right)\hspace{-0.1cm}
\right)_{R_{(2,6,7)}^{(m,m+1,1)}}\hspace{-1cm}\{-m\}
\boxtimes
K\left(\hspace{-0.1cm}\left(\hspace{-0.1cm}
\begin{array}{c}
\Lambda_{1,(8;7,4)}^{[1,m]}\\[.5em]
\vdots\\[.5em]
\Lambda_{m+1,(8;7,4)}^{[1,m]}
\end{array}\hspace{-0.1cm}
\right)\hspace{-0.1cm}
;\hspace{-0.1cm}
\left(\hspace{-0.1cm}
\begin{array}{c}
x_{1,8}-X_{1,(4,7)}^{(m,1)}\\[.6em]
\vdots\\[.6em]
x_{m+1,8}-X_{m+1,(4,7)}^{(m,1)}
\end{array}\hspace{-0.1cm}
\right)\hspace{-0.1cm}
\right)_{R_{(4,7,8)}^{(m,1,m+1)}.}
\end{eqnarray*}
We apply Corollary \ref{cor2-11} to this matrix factorization. Then we obtain
\begin{equation}
\label{factorization3}
K\left(\left(
\begin{array}{c}
\Lambda_{1,(6;1,5)}^{[1,m]}\\[.5em]
\vdots\\[.5em]
\Lambda_{m+1,(6;1,5)}^{[1,m]}\\[.5em]
V_{m+1,(3,5;8)}^{[1,m]}
\end{array}
\right)
;
\left(
\begin{array}{c}
x_{1,6}-X_{1,(1,5)}^{(1,m)}\\[.6em]
\vdots\\[1em]
x_{m+1,6}-X_{m+1,(1,5)}^{(1,m)}\\[.6em]
X_{m+1,(3,5)}^{(1,m)}-x_{m+1,8}
\end{array}
\right)
\right)_{Q_2}\{-2m\},
\end{equation}
where 
$
Q_2=R_{(1,2,3,4,5,6,7,8)}^{(1,m,1,m,m,m+1,1,m+1)}
\left/\left<
\begin{array}{c}
X_{1,(3,5)}^{(1,m)}-x_{1,8},
\cdots,
X_{m,(1,5)}^{(1,m)}-x_{m,8},\\
X_{1,(2,7)}^{(m,1)}-x_{1,6},
\cdots,
X_{m+1,(2,7)}^{(m,1)}-x_{m+1,6},\\
x_{1,8}-X_{1,(4,7)}^{(m,1)},
\cdots,
x_{m+1,8}-X_{m+1,(4,7)}^{(m,1)}
\end{array}
\right>\right. .
$\\
In the quotient, we have equations
\begin{eqnarray}
\label{eq1}
\left\{
\begin{array}{l}
x_{k,5}= X_{k,(3,4,7)}^{(-1,m,1)}\hspace{1cm}(k=1,\cdots,m),\\[-0.1em]
x_{k,6}= X_{k,(2,7)}^{(m,1)}\hspace{1cm}(k=1,\cdots,m+1),\\[-0.1em]
x_{k,8}= X_{k,(4,7)}^{(m,1)}\hspace{1cm}(k=1,\cdots,m+1).
\end{array}
\right.
\end{eqnarray}
Using the equations, we find the factorization $(\ref{factorization3})$ is isomorphic to
\begin{eqnarray*}
&&
K\left(\left(
\begin{array}{c}
\Lambda_{1,(6;1,5)}^{[1,m]}\\[.5em]
\vdots\\[.5em]
\Lambda_{m,(6;1,5)}^{[1,m]}\\[.5em]
\Lambda_{m+1,(6;1,5)}^{[1,m]}\\[.5em]
V_{m+1,(3,5;8)}^{[1,m]}
\end{array}
\right)
;
\left(
\begin{array}{c}
X_{1,(2,7)}^{(m,1)}-X_{1,(1,3,4,7)}^{(1,-1,m,1)}\\[.5em]
\vdots\\[.6em]
X_{m,(2,7)}^{(m,1)}-X_{m,(1,3,4,7)}^{(-1,m,1)}\\[.6em]
X_{m+1,(2,7)}^{(m,1)}-x_{1,1}X_{m,(3,4,7)}^{(-1,m,1)}\\[.6em]
x_{1,3}X_{m,(3,4,7)}^{(-1,m,1)}-X_{m+1,(4,7)}^{(m,1)}
\end{array}
\right)
\right)_{Q_2}\{-2m\},\\
&\simeq&
K\left(\left(
\begin{array}{c}
\Lambda_{1,(6;1,5)}^{[1,m]}+x_{1,7}\Lambda_{2,(6;1,5)}^{[1,m]}\\[.3em]
\vdots\\[.3em]
\Lambda_{m,(6;1,5)}^{[1,m]}+x_{1,7}\Lambda_{m+1,(6;1,5)}^{[1,m]}\\[.3em]
\Lambda_{m+1,(6;1,5)}^{[1,m]}\\[.3em]
V_{m+1,(3,5;8)}^{[1,m]}
\end{array}
\right)
;
\left(
\begin{array}{c}
x_{1,2}-X_{1,(1,3,4)}^{(1,-1,m)}\\[.3em]
\vdots\\[.4em]
x_{m,2}-X_{m,(1,3,4)}^{(1,-1,m)}\\[.4em]
(x_{1,7}-x_{1,1})X_{m,(3,4)}^{(-1,m)}\\[.4em]
(x_{1,3}-x_{1,7})X_{m,(3,4)}^{(-1,m)}
\end{array}
\right)
\right)_{Q_2}\{-2m\}.
\end{eqnarray*}
By Eq. (\ref{eq1}) in the quotient $Q_2$, the polynomial $\Lambda_{m+1,(6;1,5)}^{[1,m]}$ can be described by
\begin{eqnarray}
\label{eq2}
\Lambda_{m+1,(6;1,5)}^{[1,m]}&=&\frac{F_{m+1}(X_{1,(1,5)}^{(1,m)},\cdots,X_{m,(1,5)}^{(1,m)},x_{m+1,6})-F_{m+1}(X_{1,(1,5)}^{(1,m)},\cdots,X_{m,(1,5)}^{(1,m)},X_{m+1,(1,5)}^{(1,m)})}{x_{m+1,6}-X_{m+1,(1,5)}^{(1,m)}}\\[-0.1em]
\nonumber
&=&c_0\left(X_{1,(1,5)}^{(1,m)}\right)^{n-m}+\cdots\\[-0.1em]
\nonumber
&\equiv&c_0(x_{1,7}+x_{1,1}-x_{1,3}+x_{1,4})^{n-m}+\cdots \hspace{1cm}(\mod Q_2)\\[-0.1em]
\nonumber
&=&c_0 x_{1,7}^{n-m}+r_1 x_{1,7}^{n-m-1} + \cdots +r_{n-m},
\end{eqnarray}
where $c_0 \in \Q$ and $r_k \in R_{(1,2,3,4)}^{(1,m,1,m)}$ ($k=1,\cdots,n-m$). Using Theorem \ref{exclude}, we have
\begin{eqnarray*}
&&
K\left(\left(
\begin{array}{c}
\Lambda_{1,(6;1,5)}^{[1,m]}+x_{1,7}\Lambda_{2,(6;1,5)}^{[1,m]}\\[.3em]
\vdots\\[.5em]
\Lambda_{m,(6;1,5)}^{[1,m]}+x_{1,7}\Lambda_{m+1,(6;1,5)}^{[1,m]}\\[.5em]
(x_{1,7}-x_{1,1})X_{m,(3,4)}^{(-1,m)}\\[.5em]
V_{m+1,(3,5;8)}^{[1,m]}
\end{array}
\right)
;
\left(
\begin{array}{c}
x_{1,2}-X_{1,(1,3,4)}^{(1,-1,m)}\\[.5em]
\vdots\\[.6em]
x_{m,2}-X_{m,(1,3,4)}^{(1,-1,m)}\\[.6em]
\Lambda_{m+1,(6;1,5)}^{[1,m]}\\[.6em]
(x_{1,3}-x_{1,7})X_{m,(3,4)}^{(-1,m)}
\end{array}
\right)
\right)_{Q_2}\{-2m\}\{2m+1-n\}\left<1\right>,\\[-0.1em]
&\simeq&
K\left(\left(
\begin{array}{c}
\Lambda_{1,(6;1,5)}^{[1,m]}+x_{1,7}\Lambda_{2,(6;1,5)}^{[1,m]}\\[.3em]
\vdots\\[.5em]
\Lambda_{m,(6;1,5)}^{[1,m]}+x_{1,7}\Lambda_{m+1,(6;1,5)}^{[1,m]}\\[.5em]
V_{m+1,(3,5;8)}^{[1,m]}
\end{array}
\right)
;
\left(
\begin{array}{c}
x_{1,2}-X_{1,(1,3,4)}^{(1,-1,m)}\\[.6em]
\vdots\\[.6em]
x_{m,2}-X_{m,(1,3,4)}^{(1,-1,m)}\\[.6em]
(x_{1,3}-x_{1,7})X_{m,(3,4)}^{(-1,m)}
\end{array}
\right)
\right)_{Q_2/\Lambda_{m+1,(6;1,5)}^{[1,m]}}\{1-n\}\left<1\right>.
\end{eqnarray*}
Applying Corollary \ref{induce-sq1} to this factorization, there are polynomials $b_k \in R_{(1,2,3,4)}^{(1,m,1,m)}$ ($k=1,\cdots,m$) and $b_0\in Q_2/\Lambda_{m+1,(6;1,5)}^{[1,m]}$ such that $b_0(x_{1,3}-x_{1,7})\equiv B_1 \in R_{(1,2,3,4)}^{(1,m,1,m)} (\mod Q_2/\Lambda_{m+1,(6;1,5)}^{[1,m]})$ and we have an isomorphism between the above factorization and the following factorization
\begin{equation}
\label{factorization4}
K\left(\left(
\begin{array}{c}
b_1\\[.5em]
\vdots\\[.5em]
b_m\\[.6em]
b_0
\end{array}
\right)
;
\left(
\begin{array}{c}
x_{1,2}-X_{1,(1,3,4)}^{(1,-1,m)}\\[.5em]
\vdots\\[.5em]
x_{m,2}-X_{m,(1,3,4)}^{(1,-1,m)}\\[.5em]
(x_{1,3}-x_{1,7})X_{m,(3,4)}^{(-1,m)}
\end{array}
\right)
\right)_{Q_2/\Lambda_{m+1,(6;1,5)}^{[1,m]}}\{1-n\}\left<1\right>.
\end{equation}
We choose decompositions of $Q_2/\Lambda_{m+1,(6;1,5)}^{[1,m]}$ to be
\begin{eqnarray}
\nonumber
R_3&:=&R_{(1,2,3,4)}^{(1,m,1,m)}\oplus
(x_{1,3}-x_{1,7})R_{(1,2,3,4)}^{(1,m,1,m)}\oplus
\cdots\oplus
x_{1,7}^{n-m-2}(x_{1,3}-x_{1,7})R_{(1,2,3,4)}^{(1,m,1,m)},\\
\nonumber
R_4&:=&R_{(1,2,3,4)}^{(1,m,1,m)}\oplus
x_{1,7} R_{(1,2,3,4)}^{(1,m,1,m)}\oplus
\cdots\oplus
x_{1,7}^{n-m-2} R_{(1,2,3,4)}^{(1,m,1,m)}\oplus
\beta R_{(1,2,3,4)}^{(1,m,1,m)},
\end{eqnarray}
where $\beta=\sum_{k=0}^{n-m-1}x_{1,7}^{n-m-1-k}\left(c_0x_{1,3}^k+\sum_{l=1}^{k}x_{1,3}^{k-l}r_l\right)$.
Then, we consider the partial matrix factorization of (\ref{factorization4}) $K(b_0;(x_{1,3}-x_{1,7})X_{m,(3,4)}^{(-1,m)})_{Q_2/\Lambda_{m+1,(6;1,5)}^{[1,m]}}$. This is isomorphic to
\begin{eqnarray*}
\left(R_3,R_4\{3-n\},g_0,g_1\right),
\end{eqnarray*}
where 
\begin{eqnarray*}
g_0&=&\left(
\begin{array}{cc}
{}^t\mathfrak{0}_{n-m-1}&E_{n-m-1}(B_1)\\
\frac{b_0}{\beta}&\mathfrak{0}_{n-m-1}
\end{array}
\right),\\[-0.1em]
g_1&=&\left(
\begin{array}{cc}
\mathfrak{0}_{n-m-1}&\beta (x_{1,3}-x_{1,7})X_{m,(3,4)}^{(-1,m)}\\
E_{n-m-1}(X_{m,(3,4)}^{(-1,m)})&{}^t\mathfrak{0}_{n-m-1}
\end{array}
\right).
\end{eqnarray*}
Remark that $\displaystyle\frac{b_0}{\beta}$ is a polynomial, denoted by $B_2$, in $R_{(1,2,3,4)}^{(1,m,1,m)}$ and we have
\begin{eqnarray*}
\beta(x_{1,3}-x_{1,7})&\equiv&c_0 x_{1,3}^{n-m}+r_1 x_{1,3}^{n-m-1}+\cdots +r_{n-m}\hspace{1cm}(\mod Q_2/\Lambda_{m+1,(6;1,5)}^{[1,m]})\\
&=:&B_3
.
\end{eqnarray*}
Therefore, the partial factorization has a direct decomposition
\begin{eqnarray*}
&&K\left(B_1;X_{m,(3,4)}^{(-1,m)}\right)_{R_{(1,2,3,4)}^{(1,m,1,m)}}\{[n-m-1]_q\}_{q}\{n-m\}\\[-0.1em]
&&\hspace{1cm}\bigoplus K\left(B_2;B_3 X_{m,(3,4)}^{(-1,m)}\right)_{R_{(1,2,3,4)}^{(1,m,1,m)}}.
\end{eqnarray*}
Then, the factorization (\ref{factorization4}) is isomorphic to
\begin{eqnarray}\label{total-mat2}
&&
K\left(\left(
\begin{array}{c}
b_1\\[.5em]
\vdots\\[.5em]
b_m\\[.8em]
B_1
\end{array}
\right)
;
\left(
\begin{array}{c}
x_{1,2}-X_{1,(1,3,4)}^{(1,-1,m)}\\[.5em]
\vdots\\[.5em]
x_{m,2}-X_{m,(1,3,4)}^{(1,-1,m)}\\[.5em]
X_{m,(3,4)}^{(-1,m)}
\end{array}
\right)
\right)_{R_{(1,2,3,4)}^{(1,m,1,m)}}\{[n-m-1]_q\}_{q}\{1-m\}\left<1\right>\\[-0.1em]
\label{total-mat3}&&\hspace{1cm}\bigoplus
K\left(\left(
\begin{array}{c}
b_1\\[.5em]
\vdots\\[.5em]
b_m\\[.8em]
B_2
\end{array}
\right)
;
\left(
\begin{array}{c}
x_{1,2}-X_{1,(1,3,4)}^{(1,-1,m)}\\[.5em]
\vdots\\[.5em]
x_{m,2}-X_{m,(1,3,4)}^{(1,-1,m)}\\[.5em]
B_3 X_{m,(3,4)}^{(-1,m)}
\end{array}
\right)
\right)_{R_{(1,2,3,4)}^{(1,m,1,m)}}\{1-n\}\left<1\right>.
\end{eqnarray}
Applying Theorem \ref{reg-eq} to the partial factorization (\ref{total-mat2}), we have a factorization
\begin{equation*}
\c\left(\input{figure/figsquare1j--j-1--1j-rev-mf}\right)_n
\{[m-n-1]_q\}_{q}\left<1\right>.
\end{equation*}
\indent
To show an isomorphism between the partial factorization (\ref{total-mat3}) and
\begin{equation}
\nonumber
\c\left(\input{figure/figsquare1j-rev-mf}\right)_n,
\end{equation} 
we calculate a specific form of $B_3$ and $B_3X_{m,(3,4)}^{(-1,m)}$.
First, we have
\begin{eqnarray}
\nonumber
B_3&=&c_0 x_{1,3}^{n-m}+r_1 x_{1,3}^{n-m-1}+\cdots +r_{n-m}\\
\nonumber
&=&(c_0 x_{1,7}^{n-m}+r_1 x_{1,7}^{n-m-1}+\cdots +r_{n-m})|_{x_{1,7}= x_{1,3}}\\
\nonumber
&\stackrel{Eqs.(\ref{eq2})}{\equiv}&\Lambda_{m+1,(6;1,5)}^{[1,m]}\left|_{x_{1,7}= x_{1,3}}\right.(\mod Q_2)
.
\end{eqnarray}
In the quotient $Q_2$, we have $\left.X_{k,(1,5)}^{(1,m)}\right|_{x_{1,7}= x_{1,3}}\equiv X_{k,(1,4)}^{(1,m)}$ ($k=1,...,m+1$) and $x_{m+1,6}\equiv X_{m+1,(2,3)}^{(m,1)}$.
Then, we have
\begin{eqnarray*}
B_3&\equiv&\Lambda_{m+1,(6;1,5)}^{[1,m]}\left|_{x_{1,7}= x_{1,3}}\right.\\
&\equiv&
\frac{F_{m+1}(X_{1,(1,4)}^{(1,m)},\cdots,X_{m,(1,4)}^{(1,m)},X_{m+1,(2,3)}^{(m,1)})
-F_{m+1}(X_{1,(1,4)}^{(1,m)},\cdots,X_{m,(1,4)}^{(1,m)},X_{m+1,(1,4)}^{(1,m)})}
{X_{m+1,(2,3)}^{(m,1)}-X_{m+1,(1,4)}^{(1,m)}}\\
&\equiv&\frac{F_{m+1}(X_{1,(1,4)}^{(1,m)},\cdots,X_{m,(1,4)}^{(1,m)},x_{1,3}X_{m,(1,3,4)}^{(1,-1,m)})-F_{m+1}(X_{1,(1,4)}^{(1,m)},\cdots,X_{m,(1,4)}^{(1,m)},X_{m+1,(1,4)}^{(1,m)})}{x_{1,3}X_{m,(1,3,4)}^{(1,-1,m)}-X_{m+1,(1,4)}^{(1,m)}}\\
&&\hspace{8cm}(\mod \left<x_{m,2}-X_{m,(1,3,4)}^{1,-1,m}\right>_{R_{(1,2,3,4)}^{(1,m,1,m)}})\\
&=&\frac{F_{m+1}(X_{1,(1,4)}^{(1,m)},\cdots,X_{m,(1,4)}^{(1,m)},x_{1,3}X_{m,(1,3,4)}^{(1,-1,m)})-F_{m+1}(X_{1,(1,4)}^{(1,m)},\cdots,X_{m,(1,4)}^{(1,m)},X_{m+1,(1,4)}^{(1,m)})}{(x_{1,3}-x_{1,1})X_{m,(3,4)}^{(-1,m)}}\\
&=&\frac{F_1(x_{1,3})-F_1(x_{1,1})}{(x_{1,3}-x_{1,1})X_{m,(3,4)}^{(-1,m)}}+\frac{F_{m}(X_{1,(1,3,4)}^{(1,-1,m)},\cdots,X_{m,(1,3,4)}^{(1,-1,m)})-F_{m}(x_{1,4},\cdots,x_{m,4})}{(x_{1,3}-x_{1,1})X_{m,(3,4)}^{(-1,m)}}\\
&=&\frac{F_1(x_{1,3})-F_1(x_{1,1})}{(x_{1,3}-x_{1,1})X_{m,(3,4)}^{(-1,m)}}+
\frac{F_{m}(X_{1,(1,3,4)}^{(1,-1,m)},X_{2,(1,3,4)}^{(1,-1,m)},\cdots,X_{m,(1,3,4)}^{(1,-1,m)})-F_{m}(x_{1,4},X_{2,(1,3,4)}^{(1,-1,m)},\cdots,X_{m,(1,3,4)}^{(1,-1,m)})}{(x_{1,3}-x_{1,1})X_{m,(3,4)}^{(-1,m)}}\\
&&+\frac{F_{m}(x_{1,4},X_{2,(1,3,4)}^{(1,-1,m)},X_{3,(1,3,4)}^{(1,-1,m)},\cdots,X_{m,(1,3,4)}^{(1,-1,m)})-F_{m}(x_{1,4},x_{2,4},X_{3,(1,3,4)}^{(1,-1,m)},\cdots,X_{m,(1,3,4)}^{(1,-1,m)})}{(x_{1,3}-x_{1,1})X_{m,(3,4)}^{(-1,m)}}\\
&&+\cdots +\frac{F_{m}(x_{1,4},\cdots,x_{m-1,4},X_{m,(1,3,4)}^{(1,-1,m)})-F_{m}(x_{1,4},\cdots,x_{m-1,4},x_{m,4})}{(x_{1,3}-x_{1,1})X_{m,(3,4)}^{(-1,m)}}.
\end{eqnarray*}
Using the equation $X_{k,(1,3,4)}^{(1,-1,m)}-x_{k,4}=(x_{1,1}-x_{1,3})X_{k-1,(3,4)}^{(-1,m)}$ $(k=1,...,m)$, we find
\begin{eqnarray*}
B_3 X_{m,(3,4)}^{(-1,m)}&\equiv&\frac{F_1(x_{1,1})-F_1(x_{1,3})}{x_{1,1}-x_{1,3}}\\
&&-\frac{F_{m}(X_{1,(1,3,4)}^{(1,-1,m)},X_{2,(1,3,4)}^{(1,-1,m)},\cdots,X_{m,(1,3,4)}^{(1,-1,m)})-F_{m}(x_{1,4},X_{2,(1,3,4)}^{(1,-1,m)},\cdots,X_{m,(1,3,4)}^{(1,-1,m)})}{X_{1,(1,3,4)}^{(1,-1,m)}-x_{1,4}}\\
&&-\frac{F_{m}(x_{1,4},X_{2,(1,3,4)}^{(1,-1,m)},X_{3,(1,3,4)}^{(1,-1,m)},\cdots,X_{m,(1,3,4)}^{(1,-1,m)})-F_{m}(x_{1,4},x_{2,4},X_{3,(1,3,4)}^{(1,-1,m)},\cdots,X_{m,(1,3,4)}^{(1,-1,m)})}{X_{2,(1,3,4)}^{(1,-1,m)}-x_{2,4}}
X_{1,(3,4)}^{(-1,m)}\\
&&-\cdots -\frac{F_{m}(x_{1,4},\cdots,x_{m-1,4},X_{m,(1,3,4)}^{(1,-1,m)})-F_{m}(x_{1,4},\cdots,x_{m-1,4},x_{m,4})}{X_{m,(1,3,4)}^{(1,-1,m)}-x_{m,4}}
X_{m-1,(3,4)}^{(-1,m)}\\
&\equiv&\frac{F_1(x_{1,1})-F_1(x_{1,3})}{x_{1,1}-x_{1,3}}\\
&&-\frac{F_{m}(x_{1,2},x_{2,2},\cdots,x_{m,2})-F_{m}(x_{1,4},x_{2,2},\cdots,x_{m,2})}{x_{1,2}-x_{1,4}}\\
&&-\frac{F_{m}(x_{1,4},x_{2,2},x_{3,2},\cdots,x_{m,2})-F_{m}(x_{1,4},x_{2,4},x_{3,2},\cdots,x_{m,2})}{x_{2,2}-x_{2,4}}
X_{1,(3,4)}^{(-1,m)}\\
&&-\cdots -\frac{F_{m}(x_{1,4},\cdots,x_{m-1,4},x_{m,2})-F_{m}(x_{1,4},\cdots,x_{m-1,4},x_{m,4})}{x_{m,2}-x_{m,4}}
X_{m-1,(3,4)}^{(-1,m)}\\
&&\hspace{5.5cm}(\mod \left<x_{1,2}-X_{1,(1,3,4)}^{1,-1,m},...,x_{m,2}-X_{m,(1,3,4)}^{1,-1,m}\right>_{R_{(1,2,3,4)}^{(1,m,1,m)}})\\
&=&L_{1,(1;3)}^{[1]}-L_{1,(2;4)}^{[m]}-L_{2,(2;4)}^{[m]}X_{1,(3,4)}^{(-1,m)}-\cdots -L_{m,(2;4)}^{[m]}X_{m-1,(3,4)}^{(-1,m)}.
\end{eqnarray*}

Hence using Theorem \ref{reg-eq} and Proposition \ref{equiv}, the matrix factorization (\ref{total-mat3}) is isomorphic to
\begin{eqnarray*}
&&
K\left(\left(
\begin{array}{c}
\ast\\[.5em]
\vdots\\[.5em]
\ast\\[.8em]
B_2
\end{array}
\right)
;
\left(
\begin{array}{c}
x_{1,2}-X_{1,(1,3,4)}^{(1,-1,m)}\\[.5em]
\vdots\\[.5em]
x_{m,2}-X_{m,(1,3,4)}^{(1,-1,m)}\\[.5em]
L_{1,(1;3)}^{[1]}-\sum_{k=1}^{m}L_{k,(2;4)}^{[m]}X_{k-1,(3,4)}^{(-1,m)}
\end{array}
\right)
\right)_{R_{(1,2,3,4)}^{(1,m,1,m)}}\{1-n\}\left<1\right>\\
&\simeq&
K\left(\left(
\begin{array}{c}
L_{1,(2;4)}^{[m]}\\[.5em]
L_{2,(2;4)}^{[m]}\\[.8em]
\vdots\\[.5em]
L_{m,(2;4)}^{[m]}\\[.8em]
x_{1,3}-x_{1,1}
\end{array}
\right)
;
\left(
\begin{array}{c}
x_{1,2}-x_{1,4}+x_{1,3}-x_{1,1}\\[.5em]
x_{2,2}-x_{2,4}+(x_{1,3}-x_{1,1})X_{1,(3,4)}^{(-1,m)}\\[.5em]
\vdots\\[.5em]
x_{m,2}-x_{m,4}+(x_{1,3}-x_{1,1})X_{m-1,(3,4)}^{(-1,m)}\\[.5em]
L_{1,(1;3)}^{[1]}-\sum_{k=1}^{m}L_{k,(2;4)}^{[m]}X_{k-1,(3,4)}^{(-1,m)}
\end{array}
\right)
\right)_{R_{(1,2,3,4)}^{(1,m,1,m)}}\{1-n\}\left<1\right>\\
&\simeq&
K\left(\left(
\begin{array}{c}
L_{1,(2;4)}^{[m]}\\[.5em]
\vdots\\[.5em]
L_{m,(2;4)}^{[m]}\\[.5em]
L_{1,(3;1)}^{[1]}
\end{array}
\right)
;
\left(
\begin{array}{c}
x_{1,2}-x_{1,4}\\[.8em]
\vdots\\[.5em]
x_{m,2}-x_{m,4}\\[.8em]
x_{1,3}-x_{1,1}
\end{array}
\right)
\right)_{R_{(1,2,3,4)}^{(1,m,1,m)}}\simeq \c\left( \input{figure/figsquare1j-rev-mf}\right)_n.
\end{eqnarray*}
\end{proof}
\indent
We can obtain isomorphisms corresponding to the other MOY relations in \cite{MOY} by properties in the category $\HMF^{gr}_{R,\omega}$ as a Krull-Schmidt category.
%
\begin{corollary}\label{cor-square}
\begin{eqnarray}
\nonumber
\c\left( \input{figure/figsquare1j--k--k+1j-k--k--1j-mf3}\right)_n
\simeq 
\c\left( \input{figure/figsquare1j1--j1+1--1j1-mf}\right)_n\left\{\left[m_1-1\atop m_2-1\right]_q\right\}_{q}
\bigoplus\c\left( \input{figure/figsquare1j1-mf}\right)_n\left\{\left[m_1-1\atop m_2\right]_q\right\}_{q}.
\end{eqnarray}
\end{corollary}
\begin{proof}[{\bf Proof of Corollary \ref{cor-square}}]
We consider the following matrix factorization
\begin{equation}\label{gene-cor}
\c\left(\input{figure/fig-gene-cor}\right)_n.
\end{equation}
Using Proposition \ref{mat-equiv3} (2) and Proposition \ref{mat-equiv2} (1), the matrix factorization is isomorphic to
\begin{eqnarray}
\nonumber&&\c\left( \input{figure/figsquare1j--k--k+1j-k--k--1j-mf3}\right)_n\oplus\c\left(\input{figure/figsquare1j-bubble-mf}\right)_n\{[m_2-1]_q\}_{q}\\[-0.5em]
\label{part-mat1}
&\simeq&\c\left( \input{figure/figsquare1j--k--k+1j-k--k--1j-mf3}\right)_n
\oplus\c\left(\input{figure/figsquare1j1-mf}\right)_n\left\{[m_2]_q\left[m_1\atop m_2\right]_q\right\}_{q}.
\end{eqnarray}
On the other hand, using Proposition \ref{mat-equiv1}, Proposition \ref{mat-equiv2} (1) and Proposition \ref{mat-equiv3} (2), the matrix factorization (\ref{gene-cor}) is isomorphic to
\begin{eqnarray}
\nonumber&&\c\left(\,\input{figure/figsquare1j1--1--2j1-1-bubble-mf}\,\,\right)_n\simeq\c\left(\input{figure/figsquare1j1--1--2j1-1--1--1j1-mf}\right)_n\left\{\left[m_1-1\atop m_2-1\right]_q\right\}_{q}\\
\label{part-mat2}
&\simeq&
\c\left( \input{figure/figsquare1j1--j1+1--1j1-mf}\right)_n\left\{\left[m_1-1\atop m_2-1\right]_q\right\}_{q}\oplus\c\left( \input{figure/figsquare1j1-mf}\right)_n\left\{\left[m_1-1\atop m_2-1\right]_q[m_1-1]_q\right\}_{q}.
\end{eqnarray}
We have the equation
\begin{equation}
\nonumber
[m_1-1]_{q}\left[m_1-1\atop m_2-1\right]_{q}-[m_2-1]_{q}\left[m_1\atop m_2\right]_{q}=\left[m_1-1\atop m_2\right]_{q}.
\end{equation}
Therefore, Krull-Schmidt property turns the isomorphism between the factorization (\ref{part-mat1}) and the factorization (\ref{part-mat2}) into the isomorphism of the corollary.
\end{proof}
%
%
%
%
%
%
%
%
%
\section{Complexes of matrix factorizations for $[1,k]$-crossing and $[k,1]$-crossing}\label{sec5}
This section includes the author's new result that we define complexes of matrix factorizations for $[1,k]$-crossing and $[k,1]$-crossing $(k=1,\ldots,n-1)$ and, for tangle diagrams with $[1,k]$-crossing and $[k,1]$-crossing only, we show that there exists an isomorphism in $\k^b(\HMF^{gr}_{R,\omega})$ between complexes of matrix factorizations for these tangle diagrams in Section \ref{sec5.1}, \ref{sec5.2}, \ref{sec5.3} and \ref{sec5.4}.
Remark that this construction is a generalization of a complex of matrix factorizations for $[1,2]$-crossing given by Rozansky \cite{Roz}.
Although we can calculate the homological invariant for a given link diagram, it is not easy the calculation of a link diagram with many crossings.
We calculate the homological invariant for Hopf link with $[1,k]$-crossing and $[k,1]$-crossing in Section \ref{sec5.5}. 
\\[-2em]
\begin{figure}[htb]
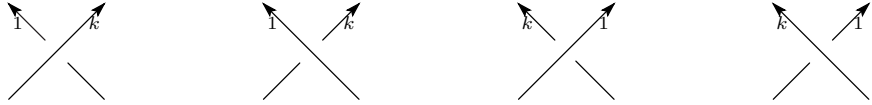
\label{1k-crossing}
\input{figure/figplus-1-k}\hspace{2cm} \input{figure/figminus-1-k}\hspace{2cm} \input{figure/figplus-k-1}\hspace{2cm} \input{figure/figminus-k-1}\\[1em]
\caption{$[1,k]$-plus crossing, $[1,k]$-minus crossing, $[k,1]$-plus crossing and $[k,1]$-minus crossing }
\end{figure}\\
The $[1,k]$-crossing and $[k,1]$-crossing are locally expanded into a linear sum by the normalized MOY calculus of the quantum $(\mathfrak{sl}_n , \land  V_n)$ link invariant as follows,
\begin{eqnarray*}
\left<\hspace{.1cm}\input{figure/figplus-1-k-intro}\hspace{.1cm}\right>_n&=&
(-1)^{1-k}q^{kn-1}\left<\hspace{.3cm}\input{figure/figsquarek1--k-1--1k}\hspace{.3cm}\right>_n+(-1)^{-k}q^{kn}\left<\hspace{.3cm}\input{figure/figsquarek1--k+1--1k}\hspace{.3cm}\right>_n,\\
\left<\hspace{.1cm}\input{figure/figplus-k-1-intro}\hspace{.1cm}\right>_n&=&
(-1)^{1-k}q^{kn-1}\left<\hspace{.3cm}\input{figure/figsquare1k--k-1--k1}\hspace{.3cm}\right>_n+(-1)^{-k}q^{kn}\left<\hspace{.3cm}\input{figure/figsquare1k--k+1--k1}\hspace{.3cm}\right>_n,\\
\left<\hspace{.1cm}\input{figure/figminus-1-k-intro}\hspace{.1cm}\right>_n&=&
(-1)^{k-1}q^{-kn+1}\left<\hspace{.3cm}\input{figure/figsquarek1--k-1--1k}\hspace{.3cm}\right>_n+(-1)^{k}q^{-kn}\left<\hspace{.3cm}\input{figure/figsquarek1--k+1--1k}\hspace{.3cm}\right>_n,\\
\left<\hspace{.1cm}\input{figure/figminus-k-1-intro}\hspace{.1cm}\right>_n&=&
(-1)^{k-1}q^{-kn+1}\left<\hspace{.3cm}\input{figure/figsquare1k--k-1--k1}\hspace{.3cm}\right>_n+(-1)^{k}q^{-kn}\left<\hspace{.3cm}\input{figure/figsquare1k--k+1--k1}\hspace{.3cm}\right>_n.
\end{eqnarray*}
%
%
%
%
\subsection{Definition of complexes for $[1,k]$-crossing and $[k,1]$-crossing}\label{sec5.1}
\begin{figure}[htb]
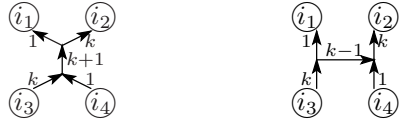

\input{figure/figsquarek1--k+1--1k-mf}\hspace{3cm}\input{figure/figsquarek1--k-1--1k-mf}\\[1em]
\caption{Colored planar diagrams assigned indexes}\label{1k-planar}
\end{figure}

\indent
First, we consider matrix factorizations for colored planar diagrams appearing in the MOY bracket for $[k,1]$-crossing and $[k,1]$-crossing, see Figure \ref{1k-planar}.\\[-2em]
\\
\begin{proposition}
The matrix factorization $\c\left(\hspace{0.5cm}\input{figure/figsquarek1--k+1--1k-mf}\hspace{0.5cm}\right)_n$ $(1\leq k\leq n-1)$ is isomorphic to the following finite matrix factorization
\begin{equation}
\nonumber
\overline{M}^{[1,k]}_{(i_1,i_2,i_3,i_4)}:=\overline{S}_{(i_1,i_2,i_3,i_4)}^{[1,k]}\boxtimes K(u^{[1,k]}_{k+1,(i_1,i_2,i_3,i_4)};(x_{1,i_1}-x_{1,i_4})X_{k,(i_2,i_4)}^{(k,-1)})_{R_{(i_1,i_2,i_3,i_4)}^{(1,k,k,1)}}\{ -k\}
\end{equation}
and the matrix factorization $\c\left(\hspace{0.5cm}\input{figure/figsquarek1--k-1--1k-mf}\hspace{0.5cm}\right)_n$ $(1\leq k\leq n-1)$ is isomorphic to the following finite matrix factorization
\begin{equation}
\nonumber
\overline{N}^{[1,k]}_{(i_1,i_2,i_3,i_4)}:=\overline{S}_{(i_1,i_2,i_3,i_4)}^{[1,k]}\boxtimes K(u^{[1,k]}_{k+1,(i_1,i_2,i_3,i_4)}(x_{1,i_1}-x_{1,i_4});X_{k,(i_2,i_4)}^{(k,-1)})_{R_{(i_1,i_2,i_3,i_4)}^{(1,k,k,1)}}\{ -k+1\},
\end{equation}
where $\overline{S}^{[1,k]}_{(i_1,i_2,i_3,i_4)}$ is the matrix factorization
\begin{eqnarray}
&&K\left(
\left(
\begin{array}{c}
A_{1,(i_1,i_2,i_3,i_4)}^{[1,k]}\vspace{0.1cm}\\
\vdots\vspace{0.1cm}\\
A_{k,(i_1,i_2,i_3,i_4)}^{[1,k]}
\end{array}
\right);
\left(
\begin{array}{c}
\vspace{0.1cm}X_{1,(i_1,i_2)}^{(1,k)}-X_{1,(i_3,i_4)}^{(k,1)}\\
\vspace{0.1cm}\vdots\\
X_{k,(i_1,i_2)}^{(1,k)}-X_{k,(i_3,i_4)}^{(k,1)}
\end{array}
\right)
\right)_{R_{(i_1,i_2,i_3,i_4)}^{(1,k,k,1)}},\\
\nonumber
&&A_{j,(i_1,i_2,i_3,i_4)}^{[1,k]}=u^{[1,k]}_{j,(i_1,i_2,i_3,i_4)}-(-x_{1,i_4})^{k+1-j}u^{[1,k]}_{k+1,(i_1,i_2,i_3,i_4)}\hspace{1cm}(1\leq j\leq k),
\end{eqnarray}
$u^{[1,k]}_{j,(i_1,i_2,i_3,i_4)}$ $(1\leq j \leq k+1)$ is the polynomial
\begin{equation}
\nonumber
\frac{
F_{k+1}(X_{1,(i_3,i_4)}^{(k,1)},\ldots,X_{j-1,(i_3,i_4)}^{(k,1)},X_{j,(i_1,i_2)}^{(1,k)},\ldots,X_{k+1,(i_1,i_2)}^{(1,k)})
-F_{k+1}(X_{1,(i_3,i_4)}^{(k,1)},\ldots,X_{j,(i_3,i_4)}^{(k,1)},X_{j+1,(i_1,i_2)}^{(1,k)},\ldots,X_{k+1,(i_1,i_2)}^{(1,k)})
}
{
X_{j,(i_1,i_2)}^{(1,k)}-X_{j,(i_3,i_4)}^{(k,1)}
}.
\end{equation}
\end{proposition}
\begin{proof}
\indent
By definition and Corollary \ref{cor2-11}, we have
\begin{eqnarray}
\nonumber
\c\left(\hspace{0.5cm}\input{figure/figsquarek1--k+1--1k-mf}\hspace{0.5cm}\right)_n
&=&\overline{V}_{(i_1,i_2;i_5)}^{[1,k]}\boxtimes\overline{\Lambda}_{(i_5;i_3,i_4)}^{[k,1]}
\\
\label{mf1}
&\simeq&K\left(
\left(
\begin{array}{c}
u^{[1,k]}_{1,(i_1,i_2,i_3,i_4)}\\
\vdots\\
u^{[1,k]}_{k+1,(i_1,i_2,i_3,i_4)}
\end{array}
\right);
\left(
\begin{array}{c}
X_{1,(i_1,i_2)}^{(1,k)}-X_{1,(i_3,i_4)}^{(k,1)}\\
\vdots\\
X_{k+1,(i_1,i_2)}^{(1,k)}-X_{k+1,(i_3,i_4)}^{(k,1)}
\end{array}
\right)
\right)_{R_{(i_1,i_2,i_3,i_4)}^{(1,k,k,1)}}\{-k\}
\end{eqnarray}
and
\begin{eqnarray}
\nonumber
\c\left(\hspace{0.5cm}\input{figure/figsquarek1--k-1--1k-mf}\hspace{0.5cm}\right)_n
&=&\overline{V}_{(i_1,i_5;i_3)}^{[1,k-1]}\boxtimes\overline{\Lambda}_{(i_2;i_5,i_4)}^{[k-1,1]}
\\
\label{mf2}
&\simeq&K\left(
\left(
\begin{array}{c}
v^{[1,k]}_{1,(i_1,i_2,i_3,i_4)}\\
\vdots\\
v^{[1,k]}_{k,(i_1,i_2,i_3,i_4)}\\
v^{[1,k]}_{k+1,(i_1,i_2,i_3,i_4)}
\end{array}
\right);
\left(
\begin{array}{c}
x_{1,i_2}-a_1\vspace{0.2cm}\\
\vdots\\
x_{k,i_2}-a_k\vspace{0.2cm}\\
b_k-x_{k,i_3}
\end{array}
\right)
\right)_{R_{(i_1,i_2,i_3,i_4)}^{(1,k,k,1)}}\{-k+1\},
\end{eqnarray}
where $a_j=X_{j,(i_1,i_3)}^{(-1,k)}+x_{1,i_4}X_{j-1,(i_1,i_3)}^{(-1,k)}$ $(1 \leq j\leq k-1)$, $a_k=x_{1,i_4}X_{k-1,(i_1,i_3)}^{(-1,k)}$, $b_k=x_{1,i_1}X_{k-1,(i_1,i_3)}^{(-1,k)}$,
\begin{align}
\nonumber
&v^{[1,k]}_{j,(i_1,i_2,i_3,i_4)}=\frac{F_{k}(a_1,\ldots,a_{j-1},x_{j,i_2},\ldots,x_{k,i_2})-F_{k}(a_1,\ldots,a_{j},x_{j+1,i_2},\ldots,x_{k,i_2})}
{x_{j,i_2}-a_j}&& (1\leq j\leq k),\\
\nonumber
&v^{[1,k]}_{k+1,(i_1,i_2,i_3,i_4)}=\frac{F_{k}(x_{1,i_3},\ldots,x_{k-1,i_3},b_k)-F_{k}(x_{1,i_3},\ldots,x_{k,i_3})}{b_k-x_{k,i_3}}.&&
\end{align}
\indent

By Proposition \ref{sym-prop} $(5)$, $x_{j,i_3}-X_{j,(i_1,i_2,i_4)}^{(1,k,-1)}$ $(j=1,...,k)$ is in the ideal $I=\langle X_{j,(i_1,i_2)}^{(1,k)}-X_{j,(i_3,i_4)}^{(k,1)} \rangle_{j=1,...,k}$.
By the polynomial $x_{k,i_3}-X_{k,(i_1,i_2,i_4)}^{(1,k,-1)}$ in the ideal $I$ and the equation $X_{k,(i_1,i_2,i_4)}^{(1,k,-1)}=x_{1,i_1}X_{k,(i_2,i_4)}^{(k,-1)}+X_{k,(i_1,i_2,i_4)}^{(k,-1)}$ of Proposition \ref{sym-prop} (2), we have
\begin{eqnarray}
\nonumber
X_{k+1,(i_1,i_2)}^{(1,k)}-X_{k+1,(i_3,i_4)}^{(k,1)} &= & x_{1,i_1}x_{k,i_2}-x_{k,i_3}x_{1,i_4}\\
\nonumber
&\equiv &x_{1,i_1}x_{k,i_2}-X_{k,(i_1,i_2,i_4)}^{(1,k,-1)}x_{1,i_4}\hspace{0.5cm}(\mod I)\\
\nonumber
&=&(x_{1,i_1}-x_{1,i_4})X_{k,(i_2,i_4)}^{(k,-1)}.
\end{eqnarray}
By Theorem \ref{equiv}, there exist an isomorphism $\overline{\varphi_1}_{(i_1,i_2,i_3,i_4)}$ from the matrix factorization (\ref{mf1}) to
\begin{eqnarray}
&&\overline{S}_{(i_1,i_2,i_3,i_4)}^{[1,k]}\boxtimes K(u^{[1,k]}_{k+1,(i_1,i_2,i_3,i_4)};(x_{1,i_1}-x_{1,i_4})X_{k,(i_2,i_4)}^{(k,-1)})_{R_{(i_1,i_2,i_3,i_4)}^{(1,k,k,1)}}\{-k\}.
\end{eqnarray} 

\indent
By Proposition \ref{sym-prop} $(2)$, we have
\begin{eqnarray}
\nonumber
x_{j,i_2}-a_j&=&x_{j,i_2}-X_{j,(i_1,i_3,i_4)}^{(-1,k,1)}\hspace{0.5cm}(1 \leq j\leq k-1),\\
\nonumber
x_{k,i_2}-a_k&=&x_{k,i_2}-X_{k,(i_1,i_3,i_4)}^{(-1,k,1)}.
\end{eqnarray}
By Proposition \ref{sym-prop} $(5)$, we find $\langle x_{j,i_2}-X_{j,(i_1,i_3,i_4)}^{(-1,k,1)}\rangle_{j=1,...,k}=I$.
By these polynomials $x_{j,i_3}-X_{j,(i_1,i_2,i_4)}^{(1,k,-1)}$ $(j=1,...,k)$ in the ideal $I$ and the equation $X_{k,(i_1,i_2,i_4)}^{(1,k,-1)}=x_{1,i_1}X_{k-1,(i_2,i_4)}^{(k,-1)}+X_{k,(i_2,i_4)}^{(k,-1)}$ of Proposition \ref{sym-prop} (2), we have
\begin{eqnarray}
\nonumber
b_k-x_{k,i_3}&=&x_{1,i_1}X_{k-1,(i_1,i_3)}^{(-1,k)}-x_{k,i_3}\\
\nonumber
&\equiv &x_{1,i_1}X_{k-1,(i_2,i_4)}^{(k,-1)}-X_{k,(i_1,i_2,i_4)}^{(1,k,-1)}\hspace{0.5cm}(\mod I)\\
\nonumber
&=&X_{k,(i_2,i_4)}^{(k,-1)}.
\end{eqnarray}
By Proposition \ref{equiv} and Theorem \ref{reg-eq}, there exists an isomorphism $\overline{\varphi_2}_{(i_1,i_2,i_3,i_4)}$ from the matrix factorization $(\ref{mf2})$ to
\begin{eqnarray}
&&\overline{S}_{(i_1,i_2,i_3,i_4)}^{[1,k]}\boxtimes K(u^{[1,k]}_{k+1,(i_1,i_2,i_3,i_4)}(x_{1,i_1}-x_{1,i_4});X_{k,(i_2,i_4)}^{(k,-1)})_{R_{(i_1,i_2,i_3,i_4)}^{(1,k,k,1)}}\{-k+1\}.
\end{eqnarray}
\end{proof}
\indent
We find that there are two $\Z$-grading preserving morphisms between the matrix factorizations $\overline{M}^{[1,k]}_{(i_1,i_2,i_3,i_4)}$ and $\overline{N}^{[1,k]}_{(i_1,i_2,i_3,i_4)}$.
\begin{corollary}
There exist the following natural $\Z$-grading preserving morphisms
\begin{equation}
\id_{\overline{S}^{[1,k]}_{(i_1,i_2,i_3,i_4)}}\boxtimes (1,x_{1,i_1}-x_{1,i_4}):\overline{M}^{[1,k]}_{(i_1,i_2,i_3,i_4)}\longrightarrow\overline{N}^{[1,k]}_{(i_1,i_2,i_3,i_4)}\{-1\}
\end{equation}
and
\begin{equation}
\id_{\overline{S}^{[1,k]}_{(i_1,i_2,i_3,i_4)}}\boxtimes (x_{1,i_1}-x_{1,i_4},1):\overline{N}^{[1,k]}_{(i_1,i_2,i_3,i_4)}\longrightarrow\overline{M}^{[1,k]}_{(i_1,i_2,i_3,i_4)}\{-1\}.
\end{equation}
\end{corollary}
\indent
We define complexes of matrix factorizations corresponding to $\pm$-crossings with coloring $[1,k]$ and $[k,1]$.
\begin{definition}\label{def-crossing-1-k}
We define complexes of matrix factorization for $\pm$-crossings with coloring $[1,k]$ and $[k,1]$ as follows,
\begin{eqnarray}
\nonumber
%
%
&&\xymatrix{
&-k\ar@{.}[d]
&&1-k\ar@{.}[d]&\\
\c\left(\hspace{.5cm}\input{figure/figplus-1-k-mf}\hspace{.5cm}\right)_n\hspace{-.2cm}:=
0\ar[r]&
\overline{M}^{[1,k]}_{(i_1,i_2,i_3,i_4)}\{kn\}\left<k\right>
\ar[rr]^(.5){\chi_{+,(i_1,i_2,i_3,i_4)}^{[1,k]}}
&&\overline{N}^{[1,k]}_{(i_1,i_2,i_3,i_4)}\{kn-1\}\left<k\right>\ar[r]&0,
}
%
%
%
%
\\
\nonumber
%
%
&&\xymatrix{
&-k\ar@{.}[d]
&&1-k\ar@{.}[d]&\\
\c\left(\hspace{.5cm}\input{figure/figplus-k-1-mf}\hspace{.5cm}\right)_n\hspace{-.2cm}:=
0\ar[r]&
\overline{M}^{[1,k]}_{(i_2,i_1,i_4,i_3)}\{kn\}\left<k\right>
\ar[rr]^(.5){\chi_{+,(i_2,i_1,i_4,i_3)}^{[1,k]}}
&&\overline{N}^{[1,k]}_{(i_2,i_1,i_4,i_3)}\{kn-1\}\left<k\right>\ar[r]&0,
}
%
%
%
%
\\
\nonumber
%
%
&&\xymatrix{
&k-1\ar@{.}[d]
&&k\ar@{.}[d]\\
\c\left(\hspace{.5cm}\input{figure/figminus-1-k-mf}\hspace{.5cm}\right)_n\hspace{-.2cm}:=
0\ar[r]&
\overline{N}^{[1,k]}_{(i_1,i_2,i_3,i_4)}\{1-kn\}\left<k\right>
\ar[rr]^(.5){\chi_{-,(i_1,i_2,i_3,i_4)}^{[1,k]}}
&&\overline{M}^{[1,k]}_{(i_1,i_2,i_3,i_4)}\{-kn\}\left<k\right>
\ar[r]&0,
}
%
%
%
%
\\
\nonumber
%
%
&&\xymatrix{
&k-1\ar@{.}[d]
&&k\ar@{.}[d]\\
\c\left(\hspace{.5cm}\input{figure/figminus-k-1-mf}\hspace{.5cm}\right)_n\hspace{-.2cm}:=
0\ar[r]&
\overline{N}^{[1,k]}_{(i_2,i_1,i_4,i_3)}\{1-kn\}\left<k\right>
\ar[rr]^(.5){\chi_{-,(i_2,i_1,i_4,i_3)}^{[1,k]}}
&&\overline{M}^{[1,k]}_{(i_2,i_1,i_4,i_3)}\{ -kn\}\left<k\right>
\ar[r]&0,
}
%
%
%
%
\end{eqnarray}
where 
\begin{eqnarray}
\nonumber
\chi_{+,(i_1,i_2,i_3,i_4)}^{[1,k]}&:=&\id_{\overline{S}^{[1,k]}_{(i_1,i_2,i_3,i_4)}}\boxtimes(1,x_{1,i_1}-x_{1,i_4}),\\
\nonumber
\chi_{-,(i_1,i_2,i_3,i_4)}^{[1,k]}&:=&\id_{\overline{S}^{[1,k]}_{(i_1,i_2,i_3,i_4)}}\boxtimes(x_{1,i_1}-x_{1,i_4},1).
\end{eqnarray}
\end{definition}
%
%
%
%
\subsection{Decomposition of a tangle diagram and a complex for a tangle diagram}\label{sec5.15}
\indent
We consider a decomposition of a colored tangle diagram $T$ into colored crossings and colored planar lines using markings.
\begin{definition}
A decomposition of a colored tangle diagram $T$ is {\bf effective} if the decomposition consists of colored single crossings only.
A decomposition of a colored tangle diagram $T$ is {\bf non-effective} if the decomposition consists of colored crossings and not less than one colored line.\\
\end{definition}
\begin{figure}[hbt]
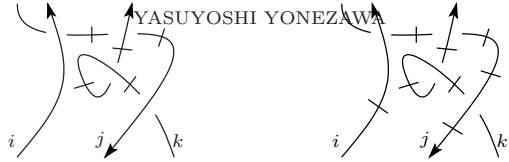

\input{figure/fig-tangle-diag-effective}
\hspace{2cm}
\input{figure/fig-tangle-diag-wasted}
\caption{Effective decomposition and non-effective decomposition of a tangle diagram}
\end{figure}
\indent
For a tangle diagram with $[1,k]$-crossings and $[1,k]$-crossings, we define a complex of matrix factorizations as follows:
We decompose the tangle diagram into $[1,k]$-crossings, $[1,k]$-crossings and colored lines using markings and assign different indexes to the markings and end points.
Then, we take tensor product of these complexes of matrix factorizations for $[1,k]$-crossings, $[1,k]$-crossings and colored lines in the decomposition.
In the categories $\kom^b(\MF^{gr})$ and $\k^b(\MF^{gr})$ (resp. $\kom^b(\HMF^{gr})$ and $\k^b(\HMF^{gr})$), an object $\overline{M}$ in $\MF^{gr}$ (resp. $\HMF^{gr}$) is defined by the following complex
\begin{equation}
\nonumber
\xymatrix{
&-1\ar@{.}[d]&0\ar@{.}[d]&1\ar@{.}[d]&\\
\cdots\ar[r]&0\ar[r]&\overline{M}\ar[r]&0\ar[r]&\cdots.
}
\end{equation}
\indent
By Lemma \ref{preserve-finiteness}, it suffices to obtain a complex for a tangle diagram with $[1,k]$-crossings and $[1,k]$-crossings that we consider the effective decomposition of the tangle diagram.
\begin{definition}
For a colored tangle diagram with $[1,k]$-crossings and $[1,k]$-crossings $T$, we define a complex of matrix factorizations to be tensor product of complexes for $[1,k]$-crossings and $[1,k]$-crossings of the effective decomposition of $T$.
\end{definition}
%
%
%
%
\subsection{Invariance under Reidemeister moves}\label{sec5.2}
\begin{theorem}[In the case $k=1$, Khovanov-Rozansky\cite{KR1}]\label{main1}
\indent
We consider tangle diagrams with $[1,k]$-crossings and $[k,1]$-crossings transforming to each other under colored Reidemeister moves composed of $[1,k]$-crossings and $[k,1]$-crossings.
Complexes of factorizations for these tangle diagrams are isomorphic in $\k^{b}(\HMF^{gr}_{R,\omega})$:
\begin{eqnarray*}
(I_1)&&\c\left(\input{figure/r1p-1}\right)_n\simeq\c\left(\input{figure/r1c-1}\right)_n\simeq\c\left(\input{figure/r1m-1}\right)_n,\\
(IIa_{1k})&&\c\left(\input{figure/r2-1kl}\right)_n\simeq\c\left(\input{figure/r2-1kc}\right)_n\simeq\c\left(\input{figure/r2-1kr}\right)_n,\,
\c\left(\input{figure/r2-k1l}\right)_n\simeq\c\left(\,\input{figure/r2-k1c}\right)_n\simeq\c\left(\,\input{figure/r2-k1r}\right)_n,\\
(IIb_{1k})&&\c\left(\input{figure/r2-1klrev}\right)_n\simeq\c\left(\input{figure/r2-1kcrev}\right)_n\simeq\c\left(\input{figure/r2-1krrev}\right)_n,\,
\c\left(\input{figure/r2-1klrevori}\right)_n\simeq\c\left(\input{figure/r2-1kcrevori}\right)_n\simeq\c\left(\input{figure/r2-1krrevori}\right)_n,\\
(III_{11k})&&\c\left(\input{figure/r3-k11}\right)_n\simeq\c\left(\input{figure/r3-k11rev}\right)_n,
\c\left(\input{figure/r3-1k1}\right)_n\simeq\c\left(\input{figure/r3-1k1rev}\right)_n,
\c\left(\input{figure/r3-11k}\right)_n\simeq\c\left(\input{figure/r3-11krev}\right)_n.
\end{eqnarray*}
\end{theorem}
\begin{proof}
The invariance of $(I_1)$ is proved by M. Khovanov and L. Rozansky in \cite{KR1}. The other invariance is proved in the following section.
\end{proof}
%
%
%
%
\subsection{Proof of invariance under Reidemeister moves IIa and IIb}\label{sec5.3}
By definition of a complex of tangle diagram in Section \ref{sec5.15}, the complex $\c\left(\input{figure/r2-1kl-mf}\right)_n$ is the tensor product of these complexes $\c\left(\input{figure/r2-1kl-decomp-mf1}\right)_n$ and $\c\left(\input{figure/r2-1kl-decomp-mf2}\right)_n$.
The complex is an object of $\k^b(\HMF^{gr}_{R_{(1,2,3,4)}^{(1,k,1,k)},\omega_1})$, where $\omega_1=F_1(\mathbb{X}^{(1)}_{(1)})+F_k(\mathbb{X}^{(k)}_{(2)})-F_1(\mathbb{X}^{(1)}_{(3)})-F_k(\mathbb{X}^{(k)}_{(4)})$.
Then, we have
\begin{equation}
\label{reid2}
\xymatrix{
&-1\ar@{.}[d]&&&0\ar@{.}[d]&&&1\ar@{.}[d]\\
\c\left(\input{figure/r2-1kl-mf}\right)_n=&
\overline{M}_{00}\{1\}
\ar[rrr]^{\left(\begin{array}{c}\overline{\mu}_1\\ \overline{\mu}_2\end{array}\right)}&&&
\text{$\begin{array}{c}\overline{M}_{10}\\ \oplus \\ \overline{M}_{01}\end{array}$}
\ar[rrr]^{\txt{$(\overline{\mu}_3,\overline{\mu}_4)$}}&&&
\overline{M}_{11}\{-1\},
}
\end{equation}
where 
\begin{align}
&\nonumber \overline{M}_{00}=\overline{N}^{[1,k]}_{(1,2,5,6)}\boxtimes\overline{M}^{[1,k]}_{(6,5,4,3)},
&&\overline{M}_{10}=\overline{M}^{[1,k]}_{(1,2,5,6)}\boxtimes\overline{M}^{[1,k]}_{(6,5,4,3)},\\
&\nonumber \overline{M}_{01}=\overline{N}^{[1,k]}_{(1,2,5,6)}\boxtimes\overline{N}^{[1,k]}_{(6,5,4,3)},
&&\overline{M}_{11}=\overline{M}^{[1,k]}_{(1,2,5,6)}\boxtimes\overline{N}^{[1,k]}_{(6,5,4,3)},\\
&\nonumber \overline{\mu}_1=
\left(\id_{\overline{S}^{[1,k]}_{(1,2,5,6)}}\boxtimes(x_{1,1}-x_{1,6},1)\right)\boxtimes\id_{\overline{M}^{[1,k]}_{(6,5,4,3)}},
&&\overline{\mu}_2=
\id_{\overline{N}^{[1,k]}_{(1,2,5,6)}}\boxtimes\left(\id_{\overline{S}^{[1,k]}_{(6,5,4,3)}}\boxtimes(1,x_{1,6}-x_{1,3})\right),\\
&\overline{\mu}_3=
\id_{\overline{M}^{[1,k]}_{(1,2,5,6)}}\boxtimes\left(\id_{\overline{S}^{[1,k]}_{(6,5,4,3)}}\boxtimes(1,x_{1,6}-x_{1,3})\right), 
&&\nonumber \overline{\mu}_4=
-\left(\id_{\overline{S}^{[1,k]}_{(1,2,5,6)}}\boxtimes(x_{1,1}-x_{1,6},1)\right)\boxtimes\id_{\overline{N}^{[1,k]}_{(6,5,4,3)}}.
\end{align}

\begin{lemma}\label{lem-IIa}
There exist isomorphisms in $\HMF^{gr}_{R_{(1,2,3,4)}^{(1,k,1,k)},\omega_1}$
\begin{eqnarray}
\nonumber
\overline{M}_{00}\{ 1\} & \simeq & \overline{M}^{[1,k]}_{(1,2,4,3)}\{[k]_q\}_q \{ 1\} , \\
\nonumber
\overline{M}_{10} & \simeq & \overline{M}^{[1,k]}_{(1,2,4,3)} \{[k+1]_q\}_q , \\
\nonumber
\overline{M}_{01} & \simeq & \overline{M}^{[1,k]}_{(1,2,4,3)} \{[k-1]_q\}_q \oplus \overline{L}^{[1,k]}_{(1,2,4,3)} , \\
\nonumber
\overline{M}_{11}\{ -1\} & \simeq & \overline{M}^{[1,k]}_{(1,2,4,3)} \{[k]_q\}_q \{ -1\},
\end{eqnarray}
,where 
\begin{equation}
\nonumber
\overline{L}^{[1,k]}_{(1,2,4,3)}=\overline{S}^{[1,k]}_{(1,2,4,3)}\boxtimes K(u^{[1,k]}_{k+1,(1,2,4,3)}X_{k,(2,3)}^{(k,-1)};x_{1,1}-x_{1,3})_{R_{(1,2,3,4)}^{(1,k,1,k)}},
\end{equation}
such that the matrix forms of $\overline{\mu}_1$, $\overline{\mu}_2$, $\overline{\mu}_3$ and $\overline{\mu}_4$ in the complex $(\ref{reid2})$ are a $(k+1) \times k$ matrix, a $k \times k$ matrix, a $k \times (k+1)$ matrix and a $k \times k$ matrix as follows
\begin{eqnarray}
\nonumber
\overline{\mu}_1&=&\left(
\begin{array}{cc}
\mathfrak{0}_{k-1} &-X_{k,(2,3)}^{(k,-1)}\id_{\overline{M}_{(1,2,4,3)}^{[1,k]}}\\
E_{k-1}(\id_{\overline{M}_{(1,2,4,3)}^{[1,k]}})&{}^{t}\mathfrak{0}_{k-1}\\
\mathfrak{0}_{k-1}&\id_{\overline{M}_{(1,2,4,3)}^{[1,k]}}
\end{array}
\right),\\
\nonumber
\overline{\mu}_2&=&\left(
\begin{array}{cc}
E_{k-1}(\id_{\overline{M}_{(1,2,4,3)}^{[1,k]}})&{}^{t}\mathfrak{0}_{k-1}\\
\mathfrak{0}_{k-1}&\id_{\overline{S}^{[1,k]}_{(1,2,4,3)}}\boxtimes(1,X_{k,(2,3)}^{(k,-1)})
\end{array}
\right),\\
\nonumber
\overline{\mu}_3&=&\left(
\begin{array}{cc}
E_{k}(\id_{\overline{M}_{(1,2,4,3)}^{[1,k]}})&{}^{t}\mathfrak{0}_{k}
\end{array}
\right),\\
\nonumber
\overline{\mu}_4&=&-\left(
\begin{array}{cc}
\mathfrak{0}_{k-1}&\id_{\overline{S}^{[1,k]}_{(1,2,4,3)}}\boxtimes(-X_{k,(2,3)}^{(k,-1)},-1)\\
E_{k-1}(\id_{\overline{M}_{(1,2,4,3)}^{[1,k]}})&{}^{t}\mathfrak{0}_{k-1}
\end{array}
\right),
\end{eqnarray}
where $E_m(f)$ is the diagonal matrix of $f$ with the order $m$ and $\mathfrak{0}_m$ is the zero low vector with length $m$.
\end{lemma}
\begin{remark}\label{remark1}
We have
\begin{equation}
\nonumber
\overline{L}^{[1,k]}_{(1,2,4,3)}\simeq\c\left(\input{figure/r2-1kc-mf}\right)_n.
\end{equation}
Thus, by this Lemma \ref{lem-IIa}, we obtain the following isomorphism in $\k^b(\HMF^{gr}_{R_{(1,2,3,4)}^{(1,k,1,k)},\omega_1})$
\begin{equation}
\nonumber
\c\left(\input{figure/r2-1kl-mf}\right)_n\simeq\c\left(\input{figure/r2-1kc-mf}\right)_n.
\end{equation}
\end{remark}
\begin{proof}[{\bf Proof of Lemma \ref{lem-IIa}}]
By Corollary \ref{cor2-11}, we have the following direct sum decompositions of $\overline{M}_{00}\{1\}$, $\overline{M}_{10}$, $\overline{M}_{01}$ and $\overline{M}_{11}\{-1\}$. Firstly, we have
\begin{eqnarray}
\nonumber
\overline{M}_{00}\{1\}
\nonumber&=&
\overline{S}_{(1,2,5,6)}^{[1,k]}
\boxtimes
K(u^{[1,k]}_{k+1,(1,2,5,6)}(x_{1,1}-x_{1,6});X_{k,(2,6)}^{(k,-1)})_{R_{(1,2,5,6)}^{(1,k,k,1)}}
\{-k+1\}\\
\nonumber
&&\boxtimes
\overline{S}_{(6,5,4,3)}^{[1,k]}
\boxtimes
K(u^{[1,k]}_{k+1,(6,5,4,3)};(x_{1,6}-x_{1,3})X_{k,(5,3)}^{(k,-1)})_{R_{(3,4,5,6)}^{(1,k,k,1)}}
\{-k\}\{1\}
\\
\nonumber
&\simeq&
\overline{S}_{(1,2,4,3)}^{[1,k]}\boxtimes
K\left((x_{1,1}-x_{1,6})u^{[1,k]}_{k+1,(1,2,5,6)};X_{k,(2,6)}^{(k,-1)}\right)_{Q_1}\\
\nonumber
&&\boxtimes
K\left(u^{[1,k]}_{k+1,(6,5,4,3)};(x_{1,6}-x_{1,3})X_{k,(5,3)}^{(k,-1)}\right)_{Q_1}\{-2k+2\},
\end{eqnarray}
where 
\begin{equation}
\nonumber
Q_1=R_{(1,2,3,4,5,6)}^{(1,k,1,k,k,1)}\left/\left<X_{1,(1,2)}^{(1,k)}-X_{1,(5,6)}^{(k,1)},\ldots,X_{k,(1,2)}^{(1,k)}-X_{k,(5,6)}^{(k,1)}\right>\right. .
\end{equation} 
Moreover, using Corollary \ref{cor2-10}, we have
\begin{eqnarray}
\label{subcom00-1}
&\simeq&
\overline{S}_{(1,2,4,3)}^{[1,k]}\boxtimes
K\left(u^{[1,k]}_{k+1,(6,5,4,3)};(x_{1,6}-x_{1,3})X_{k,(5,3)}^{(k,-1)}\right)_{Q_1\left/\langle X_{k,(2,6)}^{(k,-1)}\rangle\right.}
\{-2k+2\}.
\end{eqnarray}
In the quotient $Q_1\left/\langle X_{k,(2,6)}^{(k,-1)}\rangle\right.$, we have the following equalities
\begin{eqnarray}
\nonumber
(x_{1,6}-x_{1,3})X_{k,(5,3)}^{(k,-1)}&=&(x_{1,1}-x_{1,3})X_{k,(2,3)}^{(k,-1)},\\
\nonumber
u^{[1,k]}_{k+1,(6,5,4,3)}&=&u^{[1,k]}_{k+1,(1,2,4,3)}.
\end{eqnarray}
Then, the matrix factorization (\ref{subcom00-1}) equals to
\begin{equation}
\label{subcom00-2}
\overline{S}^{[1,k]}_{(1,2,4,3)}\boxtimes K\left(u^{[1,k]}_{k+1,(1,2,4,3)};(x_{1,1}-x_{1,3})X_{k,(2,3)}^{(k,-1)}\right)_{Q_1\left/\langle X_{k,(2,6)}^{(k,-1)}\rangle\right.}
\{-2k+2\}.
\end{equation}
The quotient $Q_1\left/\langle X_{k,(2,6)}^{(k,-1)}\rangle\right.$ is isomorphic as a $\Z$-graded $R_{(1,2,3,4)}^{(1,k,1,k)}$-module to
\begin{equation}
\nonumber
R_1:= R_{(1,2,3,4)}^{(1,k,1,k)}\oplus x_{1,6}R_{(1,2,3,4)}^{(1,k,1,k)}\oplus\ldots\oplus x_{1,6}^{k-2}R_{(1,2,3,4)}^{(1,k,1,k)}\oplus X_{k-1,(2,3,6)}^{(k,-1,-1)} R_{(1,2,3,4)}^{(1,k,1,k)}.
\end{equation}
Since the polynomials $u^{[1,k]}_{k+1,(1,2,4,3)}$ and $(x_{1,1}-x_{1,3})X_{k,(2,3)}^{(k,-1)}$ do not include the variables of $\mathbb{X}^{(k)}_{(5)}$ and $\mathbb{X}^{(1)}_{(6)}$, then the partial matrix factorization $K\left(u^{[1,k]}_{k+1,(1,2,4,3)};(x_{1,1}-x_{1,3})X_{k,(2,3)}^{(k,-1)}\right)_{Q_1\left/\langle X_{k,(2,6)}^{(k,-1)}\rangle\right.}\{-2k+2\}$ is isomorphic to
\begin{eqnarray}
\nonumber
&&\xymatrix{
R_1\{-2k+2\}\ar[rrr]_{E_{k}(u^{[1,k]}_{k+1,(1,2,4,3)})}&&&
R_1\{3-n\}\ar[rrr]_{E_{k}((x_{1,1}-x_{1,3})X_{k,(2,3)}^{(k,-1)})}&&&
R_1\{-2k+2\}
}\\%
\label{subcom00-3}
&&\bigoplus_{i=0}^{k-1}\hspace{1cm}\simeq K\left(u^{[1,k]}_{k+1,(1,2,4,3)};(x_{1,1}-x_{1,3})X_{k,(2,3)}^{(k,-1)}\right)_{R_{(1,2,3,4)}^{(1,k,1,k)}}\{2i-2k+2\}.
\end{eqnarray}
Thus, the matrix factorization (\ref{subcom00-2}) is isomorphic to
\begin{equation}
\label{subcom00-4}
\overline{M}^{[1,k]}_{(1,2,4,3)}\{[k]_q\}_q\{1\}.
\end{equation}
\indent
Secondly, we have
\begin{eqnarray}
\nonumber\overline{M}_{10}
&=&
\overline{S}_{(1,2,5,6)}^{[1,k]}
\boxtimes
K(u^{[1,k]}_{k+1,(1,2,5,6)};(x_{1,1}-x_{1,6})X_{k,(2,6)}^{(k,-1)})_{R_{(1,2,5,6)}^{(1,k,k,1)}}
\{-k\}\\
\nonumber
&&\boxtimes
\overline{S}_{(5,6,4,3)}^{[1,k]}
\boxtimes
K(u^{[1,k]}_{k+1,(6,5,4,3)};(x_{1,6}-x_{1,3})X_{k,(5,3)}^{(k,-1)})_{R_{(6,5,4,3)}^{(k,1,1,k)}}
\{-k\}\\
\label{subcom10-1}
&\simeq&
\overline{S}_{(1,2,4,3)}^{[1,k]}\boxtimes
K\left(u^{[1,k]}_{k+1,(6,5,4,3)};(x_{1,6}-x_{1,3})X_{k,(5,3)}^{(k,-1)}\right)_{Q_1\left/\langle (x_{1,1}-x_{1,6})X_{k,(2,6)}^{(k,-1)}\rangle\right.}
\{-2k\}.
\end{eqnarray}
In the quotient $Q_1\left/\langle (x_{1,1}-x_{1,6})X_{k,(2,6)}^{(k,-1)}\rangle\right.$, we have the following equalities
\begin{eqnarray}
\nonumber
(x_{1,6}-x_{1,3})X_{k,(5,3)}^{(k,-1)}&=&(x_{1,1}-x_{1,3})X_{k,(2,3)}^{(k,-1)},\\
\nonumber
u^{[1,k]}_{k+1,(6,5,4,3)}&=&u^{[1,k]}_{k+1,(1,2,4,3)}.
\end{eqnarray}
Then, the matrix factorization (\ref{subcom10-1}) equals to
\begin{equation}
\label{subcom10-2}
\overline{S}^{[1,k]}_{(1,2,4,3)}\boxtimes K\left(u^{[1,k]}_{k+1,(1,2,4,3)};(x_{1,1}-x_{1,3})X_{k,(2,3)}^{(k,-1)}\right)_{Q_1\left/\langle (x_{1,1}-x_{1,6})X_{k,(2,6)}^{(k,-1)}\rangle\right.}
\{-2k\}.
\end{equation}
The quotient $Q_1\left/\langle (x_{1,1}-x_{1,6})X_{k,(2,6)}^{(k,-1)}\rangle\right.$ is isomorphic as a $\Z$-graded $R_{(1,2,3,4)}^{(1,k,1,k)}$-module to
\begin{equation}
\nonumber
R_2:= R_{(1,2,3,4)}^{(1,k,1,k)}\oplus (x_{1,1}-x_{1,6})R_{(1,2,3,4)}^{(1,k,1,k)}\oplus\ldots\oplus x_{1,6}^{k-2}(x_{1,1}-x_{1,6})R_{(1,2,3,4)}^{(1,k,1,k)}\oplus X_{k,(3,5)}^{(-1,k)} R_{(1,2,3,4)}^{(1,k,1,k)}.
\end{equation}
Then, the factorization $K\left(u^{[1,k]}_{k+1,(1,2,4,3)};(x_{1,1}-x_{1,3})X_{k,(2,3)}^{(k,-1)}\right)_{Q_1\left/\langle (x_{1,1}-x_{1,6})X_{k,(2,6)}^{(k,-1)}\rangle\right.}\{-2k\}$ equals to
\begin{eqnarray}
\nonumber
&&\xymatrix{
R_2\{-2k\}\ar[rrr]_{E_{k+1}(u^{[1,k]}_{k+1,(1,2,4,3)})}&&&
R_2\{1-n\}\ar[rrrr]_{E_{k+1}((x_{1,1}-x_{1,3})X_{k,(2,3)}^{(k,-1)})}&&&&
R_2\{-2k\}
}\\
\label{subcom10-3}
&\simeq&\bigoplus_{i=0}^{k}K\left(u^{[1,k]}_{k+1,(1,2,4,3)};(x_{1,1}-x_{1,3})X_{k,(2,3)}^{(k,-1)}\right)_{R_{(1,2,3,4)}^{(1,k,1,k)}}\{2i-2k\}.
\end{eqnarray}
Thus, the matrix factorization (\ref{subcom10-2}) is isomorphic to
\begin{equation}
\label{subcom10-4}
\overline{M}^{[1,k]}_{(1,2,4,3)}\{[k+1]_q\}_q.
\end{equation}
\indent
Thirdly, we have
\begin{eqnarray}
\nonumber\overline{M}_{01}
&=&
\overline{S}_{(1,2,5,6)}^{[1,k]}
\boxtimes
K(u^{[1,k]}_{k+1,(1,2,5,6)}(x_{1,1}-x_{1,6});X_{k,(2,6)}^{(k,-1)})_{R_{(1,2,5,6)}^{(1,k,k,1)}}
\{-k+1\}\\
\nonumber
&&
\boxtimes
\overline{S}_{(5,6,4,3)}^{[1,k]}
\boxtimes
K(u^{[1,k]}_{k+1,(6,5,4,3)}(x_{1,6}-x_{1,3});X_{k,(5,3)}^{(k,-1)})_{R_{(6,5,4,3)}^{(k,1,1,k)}}
\{-k+1\}\\
\label{subcom01-1}
&\simeq&
\overline{S}_{(1,2,4,3)}^{[1,k]}
\boxtimes
K\left(u^{[1,k]}_{k+1,(6,5,4,3)}(x_{1,6}-x_{1,3});X_{k,(5,3)}^{(k,-1)}\right)_{Q_1\left/\langle X_{k,(2,6)}^{(k,-1)}\rangle\right.}
\{-2k+2\}.
\end{eqnarray}
In the quotient $Q_1\left/\langle X_{k,(2,6)}^{(k,-1)}\rangle\right.$, we have equalities
\begin{eqnarray}
\nonumber
X_{k,(5,3)}^{(k,-1)}=(x_{1,1}-x_{1,3})X_{k-1,(2,3,6)}^{(k,-1,-1)},\\
\nonumber
u^{[1,k]}_{k+1,(6,5,4,3)}(x_{1,6}-x_{1,3})=u^{[1,k]}_{k+1,(1,2,4,3)}(x_{1,6}-x_{1,3}).
\end{eqnarray}
Then, the matrix factorization (\ref{subcom01-1}) equals to
\begin{equation}
\label{subcom01-2}
\overline{S}^{[1,k]}_{(1,2,4,3)}\boxtimes K\left(u^{[1,k]}_{k+1,(1,2,4,3)}(x_{1,6}-x_{1,3});(x_{1,1}-x_{1,3})X_{k-1,(2,3,6)}^{(k,-1,-1)}\right)_{Q_1\left/\langle X_{k,(2,6)}^{(k,-1)}\rangle\right.}
\{-2k+2\}.
\end{equation}
The quotient $Q_1\left/\langle X_{k,(2,6)}^{(k,-1)}\rangle\right.$ is isomorphic as an $R_{(1,2,3,4)}^{(1,k,1,k)}$-module to $R_1$ and
\begin{equation}
\nonumber
R_3:= R_{(1,2,3,4)}^{(1,k,1,k)}\oplus (x_{1,6}-x_{1,3})R_{(1,2,3,4)}^{(1,k,1,k)}\oplus\ldots\oplus x_{1,6}^{k-2}(x_{1,6}-x_{1,3})R_{(1,2,3,4)}^{(1,k,1,k)}.
\end{equation}
Then, the partial factorization $K\left(u^{[1,k]}_{k+1,(1,2,4,3)}(x_{1,6}-x_{1,3});(x_{1,1}-x_{1,3})X_{k-1,(2,3,6)}^{(k,-1,-1)}\right)_{Q_1\left/\langle X_{k,(2,6)}^{(k,-1)}\rangle\right.}\{-2k+2\}$ is isomorphic to
\begin{eqnarray}
\nonumber&&\xymatrix{
R_1\{-2k+2\}\ar[rrrr]_{
\left(
\begin{array}{cc}
\hspace{-0.3cm}\mathfrak{0}_{k-1}&\hspace{-0.3cm}X_{k,(2,3)}^{k,-1}u^{[1,k]}_{k+1,(1,2,4,3)}\\
\hspace{-0.3cm}E_{k-1}(u^{[1,k]}_{k+1,(1,2,4,3)})&\hspace{-0.3cm}{}^{t}\mathfrak{0}_{k-1}
\end{array}\hspace{-0.3cm}
\right)
}&&&&
R_3\{1-n\}\ar[rrrrr]_{
\left(
\begin{array}{cc}
\hspace{-0.3cm}{}^{t}\mathfrak{0}_{k-1}&\hspace{-0.3cm}E_{k-1}((x_{1,1}-x_{1,3})X_{k,(2,3)}^{(k,-1)})\\
\hspace{-0.3cm}x_{1,1}-x_{1,3}&\hspace{-0.3cm}\mathfrak{0}_{k-1}
\end{array}\hspace{-0.3cm}
\right)
}&&&&&
R_1\{-2k+2\}
}\\
\label{subcom01-3}
&&\hspace{1cm}\simeq\bigoplus_{i=0}^{k-2}K\left(u^{[1,k]}_{k+1,(1,2,4,3)};(x_{1,1}-x_{1,3})X_{k,(2,3)}^{(k,-1)}\right)_{R_{(1,2,3,4)}^{(1,k,1,k)}}\{2i-2k+2\}\\
\nonumber&&\hspace{2cm}\oplus K\left(u^{[1,k]}_{k+1,(1,2,4,3)}X_{k,(2,3)}^{(k,-1)};x_{1,1}-x_{1,3}\right)_{R_{(1,2,3,4)}^{(1,k,1,k)}}.
\end{eqnarray}
Thus, the matrix factorization (\ref{subcom01-2}) is isomorphic to
\begin{equation}
\label{subcom01-4}
\overline{M}^{[1,k]}_{(1,2,4,3)}\{[k-1]_q\}_q\oplus\overline{L}^{[1,k]}_{(1,2,4,3)}.
\end{equation}
\indent
Finally, we have
\begin{eqnarray}
\nonumber
\overline{M}_{11}\{-1\}&=&
\overline{S}_{(1,2,5,6)}^{[1,k]}
\boxtimes
K(u^{[1,k]}_{k+1,(1,2,5,6)};(x_{1,1}-x_{1,6})X_{k,(2,6)}^{(k,-1)})_{R_{(1,2,5,6)}^{(1,k,k,1)}}
\{-k\}\\
\nonumber
&&
\boxtimes
\overline{S}_{(5,6,4,3)}^{[1,k]}
\boxtimes
K(u^{[1,k]}_{k+1,(6,5,4,3)}(x_{1,6}-x_{1,3});X_{k,(5,3)}^{(k,-1)})_{R_{(6,5,4,3)}^{(k,1,1,k)}}
\{-k+1\}\{-1\}\\
\label{subcom11-1}
&\simeq&
\overline{S}_{(1,2,4,3)}^{[1,k]}
\boxtimes
K\left(u^{[1,k]}_{k+1,(6,5,4,3)}(x_{1,6}-x_{1,3});X_{k,(5,3)}^{(k,-1)}\right)_{Q_1\left/\langle (x_{1,1}-x_{1,6})X_{k,(2,6)}^{(k,-1)}\rangle\right.}
\{-2k\}.
\end{eqnarray}
Then, by Corollary \ref{induce-sq1} (1), the matrix factorization (\ref{subcom11-1}) is isomorphic to
\begin{equation}
\label{subcom11-2}
\overline{S}^{[1,k]}_{(1,2,4,3)}\boxtimes K\left(u^{[1,k]}_{k+1,(1,2,4,3)}(x_{1,6}-x_{1,3});X_{k,(5,3)}^{(k,-1)}\right)_{Q_1\left/\langle (x_{1,1}-x_{1,6})X_{k,(2,6)}^{(k,-1)}\rangle\right.}
\{-2k\}.
\end{equation}
The quotient $Q_1\left/\langle (x_{1,1}-x_{1,6})X_{k,(2,6)}^{(k,-1)}\rangle\right.$ is isomorphic as an $R_{(1,2,3,4)}^{(1,k,1,k)}$-module to $R_2$ and
\begin{equation}
\nonumber
R_4:=R_{(1,2,3,4)}^{(1,k,1,k)}\oplus (x_{1,6}-x_{1,3})R_{(1,2,3,4)}^{(1,k,1,k)}\oplus (x_{1,1}-x_{1,6})(x_{1,6}-x_{1,3})R_{(1,2,3,4)}^{(1,k,1,k)}\oplus\ldots\oplus x_{1,6}^{k-2}(x_{1,1}-x_{1,6})(x_{1,6}-x_{1,3})R_{(1,2,3,4)}^{(1,k,1,k)}.
\end{equation}
Then, the partial matrix factorization $K\left(u^{[1,k]}_{k+1,(1,2,4,3)}(x_{1,6}-x_{1,3});X_{k,(5,3)}^{(k,-1)}\right)_{Q_1\left/\langle (x_{1,1}-x_{1,6})X_{k,(2,6)}^{(k,-1)}\rangle\right.}\{-2k\}$ is isomorphic to
\begin{eqnarray}
\nonumber&&\xymatrix{
R_2\{-2k\}\ar[rrrrrr]_{
\left(
\begin{array}{cc}
\hspace{-0.3cm}\mathfrak{0}_{k-1}&\hspace{-0.3cm}X_{k,(2,3)}^{k,-1}u^{[1,k]}_{k+1,(1,2,4,3)}(x_{1,1}-x_{1,3})\\
\hspace{-0.3cm}E_{k-1}(u^{[1,k]}_{k+1,(1,2,4,3)})&\hspace{-0.3cm}{}^{t}\mathfrak{0}_{k-1}
\end{array}\hspace{-0.3cm}
\right)
}&&&&&&
R_4\{-1-n\}\ar[rrrrr]_{
\left(
\begin{array}{cc}
\hspace{-0.2cm}{}^{t}\mathfrak{0}_{k-1}&\hspace{-0.2cm}E_{k-1}((x_{1,1}-x_{1,3})X_{k,(2,3)}^{(k,-1)})\\
\hspace{-0.2cm}1&\hspace{-0.2cm}\mathfrak{0}_{k-1}
\end{array}\hspace{-0.2cm}
\right)
}&&&&&
R_2\{-2k\}
}\\
\nonumber&&\hspace{1cm}\simeq\bigoplus_{i=0}^{k-1}K\left(u^{[1,k]}_{k+1,(1,2,4,3)};(x_{1,1}-x_{1,3})X_{k,(2,3)}^{(k,-1)}\right)_{R_{(1,2,3,4)}^{(1,k,1,k)}}\{2i-2k\}\\
\nonumber&&\hspace{2cm}\oplus K\left(u^{[1,k]}_{k+1,(1,2,4,3)}X_{k,(2,3)}^{(k,-1)}(x_{1,1}-x_{1,3});1\right)_{R_{(1,2,3,4)}^{(1,k,1,k)}}\\
\label{subcom11-3}
&&\hspace{1cm}\simeq\bigoplus_{i=0}^{k-1}K\left(u^{[1,k]}_{k+1,(1,2,4,3)};(x_{1,1}-x_{1,3})X_{k,(2,3)}^{(k,-1)}\right)_{R_{(1,2,3,4)}^{(1,k,1,k)}}\{2i-2k\}.
\end{eqnarray}
Thus, the matrix factorization (\ref{subcom11-2}) is isomorphic to
\begin{equation}
\label{subcom11-4}
\overline{M}^{[1,k]}_{(1,2,4,3)}\{[k]_q\}_q.
\end{equation}
\indent
We show how the morphisms $\overline{\mu}_1$, $\overline{\mu}_2$, $\overline{\mu}_3$ and $\overline{\mu}_4$ of the complex (\ref{reid2}) transform by the above isomorphisms.
Since a matrix representation of $\xymatrix{R_1\{2\}\ar[rr]^{x_{1,1}-x_{1,6}}&&R_2}$ as an $R_{(1,2,3,4)}^{(1,k,1,k)}$-module morphism is a $(k+1) \times k$ matrix
\begin{equation}
\nonumber
\left(
\begin{array}{cc}
\mathfrak{0}_{k-1}&-X_{k,(2,3)}^{(k,-1)}\\
E_{k-1}(1)&{}^{t}\mathfrak{0}_{k-1}\\
\mathfrak{0}_{k-1}&1	
\end{array}
\right),
\end{equation}
then the morphism $(x_{1,1}-x_{1,6},x_{1,1}-x_{1,6})$ from $K\left(u^{[1,k]}_{k+1,(1,2,4,3)};(x_{1,1}-x_{1,3})X_{k,(2,3)}^{(k,-1)}\right)_{Q_1\left/\langle X_{k,(2,6)}^{(k,-1)}\rangle\right.} \{-2k+2\}$ to $K\left(u^{[1,k]}_{k+1,(1,2,4,3)};(x_{1,1}-x_{1,3})X_{k,(2,3)}^{(k,-1)}\right)_{Q_1\left/\langle (x_{1,1}-x_{1,6})X_{k,(2,6)}^{(k,-1)}\rangle\right.} \{-2k\}$ transforms, for the decomposition (\ref{subcom00-3}) and (\ref{subcom10-3}), into a $(k+1) \times k$ matrix
\begin{equation}
\nonumber
\left(
\begin{array}{cc}
\mathfrak{0}_{k-1}&(-X_{k,(2,3)}^{(k,-1)},-X_{k,(2,3)}^{(k,-1)})\\
E_{k-1}((1,1))&{}^{t}\mathfrak{0}_{k-1}\\
\mathfrak{0}_{k-1}&(1,1)
\end{array}
\right).
\end{equation}
Since the tensor product of $\id_{\overline{S}_{(1,2,4,3)}^{[1,k]}}$ and the above morphism is $\overline{\eta_1}$ 
, thus the morphism $\overline{\mu}_1$ transforms into $\overline{\eta}_1$.
\\
\indent
We show how the morphism $\overline{\mu}_2$ transforms for decompositions (\ref{subcom00-4}) and (\ref{subcom01-4}).
Since a matrix representation of $\xymatrix{R_1\ar[rr]^1&&R_1}$ is a matrix $E_k(1)$
and $\xymatrix{R_1\ar[rr]^{x_{1,6}-x_{1,3}}&&R_3}$ is a $k \times k$ matrix
\begin{equation}
\left(
\begin{array}{cc}
\mathfrak{0}_{k-1}&X_{k,(2,3)}^{(k,-1)}\\
E_{k-1}(1)&{}^{t}\mathfrak{0}_{k-1}
\end{array}
\right)
,
\end{equation}
then the morphism $(\, 1 \, ,\, x_{1,6}\, -\, x_{1,3}\, )$ from $K\left(\, u^{[1,k]}_{k+1,(1,2,4,3)}\, ;\, (x_{1,1}-x_{1,3})X_{k,(2,3)}^{(k,-1)}\, \right)_{Q_1\left/\langle X_{k,(2,6)}^{(k,-1)}\rangle\right.}\{\, -2k+2\, \}$ to $K\left(u^{[1,k]}_{k+1,(1,2,4,3)}(x_{1,6}-x_{1,3});X_{k,(2,3)}^{(k,-1)}\right)_{Q_1\left/\langle (x_{1,1}-x_{1,6})X_{k,(2,6)}^{(k,-1)}\rangle\right.}\{-2k+2\}$ transforms, for the decomposition (\ref{subcom00-3}) and (\ref{subcom01-3}), into a $k \times k$ matrix
\begin{equation}
\nonumber
\left(
\begin{array}{cc}
E_{k-1}((1,1))&{}^{t}\mathfrak{0}_{k-1}\\
\mathfrak{0}_{k-1}&(1,X_{k,(2,3)}^{(k,-1)})
\end{array}
\right).
\end{equation}
Thus, the morphism $\overline{\mu}_2$ transforms into $\overline{\eta}_2$.\\
\indent
We show how the morphism $\overline{\mu}_3$ transforms for decompositions (\ref{subcom10-4}) and (\ref{subcom11-4}).
Since a matrix representation of $\xymatrix{R_2\ar[rr]^1&&R_2}$ is a matrix $E_{k+1}(1)$
and $\xymatrix{R_2\ar[rr]^{x_{1,6}-x_{1,3}}&&R_4}$ is a $(k+1) \times (k+1)$ matrix
\begin{equation}
\left(
\begin{array}{cc}
\mathfrak{0}_{k}&(x_{1,1}-x_{1,3})X_{k,(2,3)}^{(k,-1)}\\
E_{k}(1)&{}^{t}\mathfrak{0}_{k}
\end{array}
\right)
,
\end{equation}
then the morphism $(1,x_{1,6}-x_{1,3})$ from $K\left(u^{[1,k]}_{k+1,(1,2,4,3)};(x_{1,1}-x_{1,3})X_{k,(2,3)}^{(k,-1)}\right)_{Q_1\left/\langle (x_{1,1}-x_{1,6})X_{k,(2,6)}^{(k,-1)}\rangle\right.}\{-2k\}$ to  $K\left(u^{[1,k]}_{k+1,(1,2,4,3)}(x_{1,6}-x_{1,3});X_{k,(2,3)}^{(k,-1)}\right)_{Q_1\left/\langle (x_{1,1}-x_{1,6})X_{k,(2,6)}^{(k,-1)}\rangle\right.}\{-2k\}$ transforms, for the decomposition (\ref{subcom10-3}) and (\ref{subcom11-3}), into a $k \times (k+1)$ matrix
\begin{equation}
\nonumber
\left(
\begin{array}{cc}
E_{k}((1,1))&{}^{t}\mathfrak{0}_{k}
\end{array}
\right).
\end{equation}
Thus, the morphism $\overline{\mu}_3$ transforms into $\overline{\eta}_3$.\\
\indent
We show how the morphism $\overline{\mu}_4$ transforms for decompositions (\ref{subcom01-4}) and (\ref{subcom11-4}).
Since a matrix representation of $\xymatrix{R_1\{2\}\ar[rr]^{-x_{1,1}+x_{1,6}}&&R_2}$ is a $(k+1) \times k$ matrix
\begin{equation}
\nonumber
\left(
\begin{array}{cc}
\mathfrak{0}_{k-1}&X_{k,(2,3)}^{(k,-1)}\\
-E_{k-1}(1)&{}^{t}\mathfrak{0}_{k-1}
\end{array}
\right)
\end{equation}
and $\xymatrix{R_3\{2\}\ar[rr]^{-x_{1,1}+x_{1,6}}&&R_4}$ is a $(k+1) \times k$ matrix
\begin{equation}	
\nonumber
\left(
\begin{array}{cc}
-x_{1,1}+x_{1,3}&\mathfrak{0}_{k-1}\\
1&\mathfrak{0}_{k-1}\\
{}^{t}\mathfrak{0}_{k-1}&-E_{k-1}(1)
\end{array}
\right)
,
\end{equation}
then the morphism $(1,x_{1,6}-x_{1,3})$ from $K\left(u^{[1,k]}_{k+1,(1,2,4,3)}(x_{1,6}-x_{1,3});(x_{1,1}-x_{1,3})X_{k-1,(2,3,6)}^{(k,-1,-1)}\right)_{Q_1\left/\langle X_{k,(2,6)}^{(k,-1)}\rangle\right.}\{-2k+2\}$ to $K\left(u^{[1,k]}_{k+1,(1,2,4,3)}(x_{1,6}-x_{1,3});X_{k,(2,3)}^{(k,-1)}\right)_{Q_1\left/\langle (x_{1,1}-x_{1,6})X_{k,(2,6)}^{(k,-1)}\rangle\right.}\{-2k\}$ transforms, for the decomposition (\ref{subcom10-3}) and (\ref{subcom11-3}), into a $k \times (k+1)$ matrix
\begin{equation}
\nonumber
\left(
\begin{array}{ccc}
(1,1)&\mathfrak{0}_{k-1}&(-X_{k,(2,3)}^{(k,-1)},-1)\\
{}^{t}\mathfrak{0}_{k-1}&E_{k-1}((1,1))&{}^{t}\mathfrak{0}_{k-1}
\end{array}
\right).
\end{equation}
Thus, we find that the morphism $\overline{\mu}_4$ transforms into $\overline{\eta}_4$.
\end{proof}
We obtain the following isomorphism in $\k^b(\HMF^{gr}_{R_{(1,2,3,4)}^{1,k,1,k},\omega_1})$ by this lemma as we have already indicated the isomorphism in Remark \ref{remark1}
\begin{equation}
\nonumber
\c\left(\input{figure/r2-1kl-mf}\right)_n\simeq\c\left(\input{figure/r2-1kc-mf}\right)_n.
\end{equation}
The unproved isomorphism of corresponding to invariance under the Reidemeister moves $(IIa_{1k})$ and the isomorphisms corresponding to the Reidemeister moves $(IIb_{1k})$ are proved in Section \ref{remain-inv-IIa} and Section \ref{remain-inv-IIb}.
%
%
%
%
\subsection{Proof of invariance under Reidemeister move III}\label{sec5.4}
\begin{proposition}\label{prop-r3}
The following isomorphisms exist in $\k^b(\HMF^{gr}_{R,\omega})$:
\begin{itemize}
\item[(1)]$\c\left(\input{figure/r3-lem-plus-1k-1}\right)_n\simeq\c\left(\input{figure/r3-lem-plus-1k-2}\right)_n,$
\hspace{0.5cm}\rm{(2)}\,\,$\c\left(\input{figure/r3-lem-minus-1k-1}\right)_n\simeq\c\left(\input{figure/r3-lem-minus-1k-2}\right)_n,$
\item[(3)]$\c\left(\input{figure/r3-lem-plus-k1-1}\right)_n\simeq\c\left(\input{figure/r3-lem-plus-k1-2}\right)_n,$
\hspace{0.5cm}\rm{(4)}\,\,$\c\left(\input{figure/r3-lem-minus-k1-1}\right)_n\simeq\c\left(\input{figure/r3-lem-minus-k1-2}\right)_n,$
\item[(5)]$\c\left(\input{figure/r3-lem-plus-k1-3}\right)_n\simeq\c\left(\input{figure/r3-lem-plus-k1-4}\right)_n,$
\hspace{0.5cm}\rm{(6)}\,\,$\c\left(\input{figure/r3-lem-minus-k1-3}\right)_n\simeq\c\left(\input{figure/r3-lem-minus-k1-4}\right)_n,$
\item[(7)]$\c\left(\input{figure/r3-lem-plus-1k-3}\right)_n\simeq\c\left(\input{figure/r3-lem-plus-1k-4}\right)_n,$
\hspace{0.5cm}\rm{(8)}\,\,$\c\left(\input{figure/r3-lem-minus-1k-3}\right)_n\simeq\c\left(\input{figure/r3-lem-minus-1k-4}\right)_n.$
\end{itemize}
\end{proposition}
\begin{proof}
We prove this proposition in Section \ref{prop-r3-proof}.
\end{proof}
We immediately find the following corollary by this proposition.
\begin{corollary}\label{cor-r3}
The following isomorphisms exist in $\k^b(\HMF^{gr}_{R,\omega})$
\begin{align}
\nonumber
&\c\left(\input{figure/r3-lem-plus-1k-5}\right)_n\simeq
\c\left(\input{figure/r3-lem-plus-1k-6}\right)_n,
&&
\c\left(\input{figure/r3-lem-plus-k1-5}\right)_n\simeq
\c\left(\input{figure/r3-lem-plus-k1-6}\right)_n
,\\
\nonumber
&\c\left(\input{figure/r3-lem-plus-1k-1-5}\right)_n\simeq
\c\left(\input{figure/r3-lem-plus-1k-1-6}\right)_n,
&&
\c\left(\input{figure/r3-lem-plus-kk-1-5}\right)_n\simeq
\c\left(\input{figure/r3-lem-plus-kk-1-6}\right)_n
.
\end{align}\\[.1em]
\end{corollary}
\begin{proof}[{\bf Proof of invariance under Reidemeister move $III_{11k}$ (Theorem \ref{main1} $(III_{11k})$)}]
The diagram \input{figure/r3-k11-text} is represented as an object of $\k^b(\HMF^{gr}_{R_{(1,2,3,4,5,6)}^{(1,1,k,k,1,1)},\omega_{III}})$ $(\omega_{III}=F_{1}(\mathbb{X}^{(1)}_{(1)})+F_{1}(\mathbb{X}^{(1)}_{(2)})+F_{k}(\mathbb{X}^{(k)}_{(3)})-F_{k}(\mathbb{X}^{(k)}_{(4)})-F_{1}(\mathbb{X}^{(1)}_{(5)})-F_{1}(\mathbb{X}^{(1)}_{(6)}))$. By Corollary \ref{cor-r3}, we have an isomorphism, in $\k^b(\HMF^{gr}_{R_{(1,2,3,4,5,6)}^{(1,1,k,k,1,1)},\omega_{III}})$,
\begin{eqnarray}
\nonumber
&&
\xymatrix{
\c\left(\input{figure/r3-lem-plus-k1-5}\right)_n
\ar[rrr]^{\overline{\chi}_{+}^{[k,1]}\boxtimes\id}&&&
\c\left(\input{figure/r3-lem-plus-kk-1-5}\right)_n
}\\
\nonumber
&\simeq&
\xymatrix{
\c\left(\input{figure/r3-lem-plus-k1-6}\right)_n
\ar[rrr]^{\widetilde{\overline{\chi}_{+}^{[k,1]}\boxtimes\id}}&&&
\c\left(\input{figure/r3-lem-plus-kk-1-6}\right)_n
.}
\end{eqnarray}
In general, the morphism $\widetilde{\overline{\chi}_{+}^{[k,1]}\boxtimes\id}$ is different from the morphism $\overline{\chi}_{+}^{[k,1]}\boxtimes\id$.
However, we find that the morphism $\widetilde{\overline{\chi}_{+}^{[k,1]}\boxtimes\id}$ is the same with the morphism $\overline{\chi}_{+}^{[k,1]}\boxtimes\id$ as follows.\\
\indent
We put
\begin{eqnarray}
\nonumber
&&\c\left(\input{figure/r3-k11rev-cut1}\right)_n\\
\label{complex-r3}&&=
\xymatrix{
0\ar[r]&C^{-k-1}
\ar[rrr]_{\left(\begin{array}{c}\id_{\overline{M}}\boxtimes\chi^{[1,1]}_{+,(2,3,9,8)}\\ \chi^{[1,k]}_{+,(9,1,7,4)}\boxtimes\id_{\overline{M}}\end{array}\right)}
&&&
{\begin{array}{c}C_1^{-k}\\ \oplus\\C_2^{-k}\end{array}}
\ar[rrrrr]_{\left(\chi^{[1,k]}_{+,(9,1,7,4)}\boxtimes\id_{\overline{N}} ,\,\,\,-\id_{\overline{N}}\boxtimes\chi^{[1,1]}_{+,(2,3,9,8)}\right)}
&&&&&C^{-k+1}\ar[r]&0
},
\end{eqnarray}
where
\begin{eqnarray}
\nonumber
C^{-k-1}&=&\overline{M}_{(9,1,7,4)}^{[1,k]}\boxtimes\overline{M}_{(2,3,9,8)}^{[1,1]}\{(k+1)(n-1)\} \langle k+1 \rangle,\\
\nonumber
C_1^{-k}&=&\overline{M}_{(9,1,7,4)}^{[1,k]}\boxtimes\overline{N}_{(2,3,9,8)}^{[1,1]}\{(k+1)(n-1)\} \langle k+1 \rangle,\\
\nonumber
C_2^{-k}&=&\overline{N}_{(9,1,7,4)}^{[1,k]}\boxtimes\overline{M}_{(2,3,9,8)}^{[1,1]}\{(k+1)(n-1)\} \langle k+1 \rangle,\\
\nonumber
C^{-k+1}&=&\overline{N}_{(9,1,7,4)}^{[1,k]}\boxtimes\overline{N}_{(2,3,9,8)}^{[1,1]}\{(k+1)(n-1)\} \langle k+1 \rangle.
\end{eqnarray}
The morphism $\widetilde{\overline{\chi}_{+}^{[k,1]}\boxtimes\id}$ consists of a tensor product of an endomorphism $\Phi $ of the complex (\ref{complex-r3}) and a morphism of $\Hom_{\HMF}\left(\c\left(\input{figure/figsquare1k--k+1--k1}\right)_n,\c\left(\input{figure/figsquare1k--k-1--k1}\right)_n\{-1\}\right)$, where the endomorphism $\Phi$ denotes
\begin{equation}
\nonumber
\xymatrix{
\ar@{}[d]^{\Phi:}
&0\ar[r]&C^{-k-1}
\ar[d]_{\overline{f}}
\ar[rr]
&&
{\begin{array}{c}C_1^{-k}\\ \oplus\\C_2^{-k}\end{array}}
\ar[d]_{\left(\begin{array}{cc}\overline{g_{00}}&\overline{g_{01}}\\ \overline{g_{10}}&\overline{g_{11}}\end{array}\right)}
\ar[rr]
&&C^{-k+1}
\ar[d]_{\overline{h}}
\ar[r]&0
\\
&0\ar[r]&C^{-k-1}
\ar[rr]
&&
{\begin{array}{c}C_1^{-k}\\ \oplus\\C_2^{-k}\end{array}}
\ar[rr]
&&C^{-k+1}\ar[r]&0.
}
\end{equation}
Since the isomorphisms of Corollary \ref{cor-r3} transform $\overline{\chi}_{+}^{[k,1]}\boxtimes\id$ into a morphism $\widetilde{\overline{\chi}_{+}^{[k,1]}\boxtimes\id}$, we have
\begin{eqnarray}
\nonumber
\overline{f}\not=0,\,\,
\overline{g_{00}}\not=0,\,\,
\overline{g_{11}}\not=0,\,\,
\overline{h}\not=0.
\end{eqnarray}
Moreover, by Corollary \ref{hom-dim}, we have
\begin{eqnarray}
\nonumber
&&
\dim_{\Q} \Hom_{\HMF}(C^{-k-1},C^{-k-1})=
\dim_{\Q} \Hom_{\HMF}(C_1^{-k},C_1^{-k})=1,\\
\nonumber
&&
\dim_{\Q} \Hom_{\HMF}(C_2^{-k},C_2^{-k})=
\dim_{\Q} \Hom_{\HMF}(C^{-k+1},C^{-k+1})=1,\\
\nonumber
&&
\dim_{\Q} \Hom_{\HMF}(C_1^{-k},C_2^{-k})=
\dim_{\Q} \Hom_{\HMF}(C_1^{-k},C_2^{-k})=0.
\end{eqnarray}
Therefore, the morphism $\Phi$ is the following morphism up to homotopy equivalence
\begin{equation}
\nonumber
\left(...,0,\id_{C^{-k-1}},\left(\begin{array}{cc}\id_{C_1^{-k}}&0\\0&\id_{C_2^{-k}}\end{array}\right),\id_{C^{-k+1}},0,...\right).
\end{equation}
We find $\dim_{\Q}\Hom_{\HMF}\left(\c\left(\input{figure/figsquare1k--k+1--k1}\right)_n,\c\left(\input{figure/figsquare1k--k-1--k1}\right)_n\{-1\}\right)=1$.
Thus, we obtain the isomorphism
\begin{equation}
\nonumber
\c\left(\input{figure/r3-k11}\right)_n\simeq
\c\left(\input{figure/r3-k11rev}\right)_n.
\end{equation}
We similarly obtain the isomorphisms,
\begin{equation}
\nonumber
\c\left(\input{figure/r3-1k1}\right)_n\simeq
\c\left(\input{figure/r3-1k1rev}\right)_n,
\c\left(\input{figure/r3-11k}\right)_n\simeq
\c\left(\input{figure/r3-11krev}\right)_n.
\end{equation}
\end{proof}
%
%
%
%
\subsection{Example: Homology of Hopf link with $[1,k]$-coloring}\label{sec5.5}
We show Poincar\'e polynomial of the link homology of $[1,k]$-colored Hopf link by Definition \ref{def-crossing-1-k}.
\begin{eqnarray}
\nonumber
P\left(\input{figure/hopf-k-1}\right)_n&=&t^{-2k}s^{k+1}q^{2kn+k}\left[ n\atop k\right]_q [n-k]_q+t^{-2k+2}s^{k+1}q^{2kn-n+k-2}\left[ n\atop k\right]_q [k]_q\\
\nonumber
&=&t^{-2k}s^{k+1}\left[ n\atop k\right]_qq^{2kn}([k]_q q^{-n+k-2} t^{2} + [n-k]_q q^{k}),\\
\nonumber
P\left(\input{figure/hopf-1-k}\right)_n&=&t^{2k-2}s^{k+1}q^{-2kn+n-k+2}\left[ n\atop k\right]_q [k]_q+t^{2k}s^{k+1}q^{-2kn-k}\left[ n\atop k\right]_q [n-k]_q\\
\nonumber
&=&t^{2k}s^{k+1}\left[ n\atop k\right]_qq^{-2kn}([k]_q q^{n-k+2} t^{-2} + [n-k]_q q^{-k}).
\end{eqnarray}
H. Awata and H. Kanno calculated a homological Hopf link invariant by refined topological vertex\cite{AK}.
The evaluation for a $[1,k]$-colored Hopf link is
\begin{equation}
\overline{\mathcal{P}}_{(k,1)}(q',t')=q'^{-2n+k^2-k}(-t')^{k}\left[ n\atop k\right]_{q'}([k]_{q'} q'^{n+k-2} t'^{-2} + [n-k]_{q'} q'^{2n+k}).
\end{equation}
Therefore, we find a relation between the evaluations as follows.
\begin{eqnarray}
P\left(\input{figure/hopf-1-k}\right)_n&=&\overline{\mathcal{P}}_{(k,1)}(q^{-1},-t)s^{k+1}t^k q^{-2kn+k^2-k},\\
P\left(\input{figure/hopf-k-1}\right)_n&=&\overline{\mathcal{P}}_{(k,1)}(q,-t^{-1})s^{k+1}t^{-k}q^{2kn-k^2+k}.
\end{eqnarray}
%
%
%
%
\section{Complexes of matrix factorizations for $[i,j]$-crossing}\label{sec6}
%
%
There are properties between the factorizations for colored planar diagrams and complexes for $[1,k]$-crossing and $[k,1]$-crossing in $\k^b(\HMF^{gr}_{R,\omega})$ as follows:
\begin{proposition}\label{prop-twist}
The following isomorphisms exist in $\k^b(\HMF^{gr}_{R,\omega})$
\begin{eqnarray}
&(1)&
\nonumber
\c\left(\input{figure/fig-plus-1k-tri1}\right)_n\simeq\c\left(\input{figure/fig-planar-1k-tri1}\right)_n
\{kn+k\}\left<k\right>\left[-k\right],
\hspace{.5cm}
\c\left(\input{figure/fig-plus-k1-tri1}\right)_n\simeq\c\left(\input{figure/fig-planar-k1-tri1}\right)_n
\{kn+k\}\left<k\right>\left[-k\right],\\
&(2)&
\nonumber
\c\left(\input{figure/fig-minus-1k-tri1}\right)_n\simeq\c\left(\input{figure/fig-planar-1k-tri1}\right)_n
\{-kn-k\}\left<k\right>\left[k\right],
\hspace{.5cm}
\c\left(\input{figure/fig-minus-k1-tri1}\right)_n\simeq\c\left(\input{figure/fig-planar-k1-tri1}\right)_n
\{-kn-k\}\left<k\right>\left[k\right],\\
&(3)&
\nonumber
\c\left(\input{figure/fig-plus-1k-tri2}\right)_n\simeq\c\left(\input{figure/fig-planar-1k-tri2}\right)_n
\{kn+k\}\left<k\right>\left[-k\right],
\hspace{.5cm}
\c\left(\input{figure/fig-plus-k1-tri2}\right)_n\simeq\c\left(\input{figure/fig-planar-k1-tri2}\right)_n
\{kn+k\}\left<k\right>\left[-k\right],\\
&(4)&
\nonumber
\c\left(\input{figure/fig-minus-1k-tri2}\right)_n\simeq\c\left(\input{figure/fig-planar-1k-tri2}\right)_n
\{-kn-k\}\left<k\right>\left[k\right],
\hspace{.5cm}
\c\left(\input{figure/fig-minus-k1-tri2}\right)_n\simeq\c\left(\input{figure/fig-planar-k1-tri2}\right)_n
\{-kn-k\}\left<k\right>\left[k\right].
\end{eqnarray}
\end{proposition}
\begin{proof}
We prove this proposition in Section \ref{prop-twist-proof}.
\end{proof}
Using this proposition, we construct a polynomial link invariant associated to a homological link invariant whose Euler characteristic is the quantum $(\mathfrak{sl}_n,\land  V_{n})$ link invariant.
%
%
%
%
\subsection{Wide edge and propositions}\label{sec6.1}
\begin{figure}[htb]
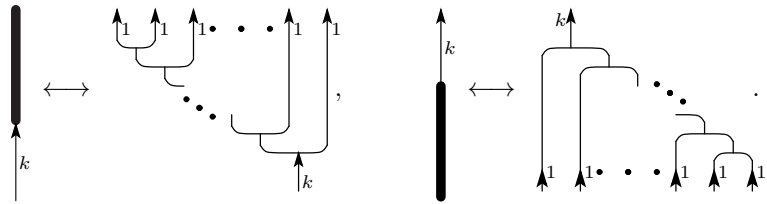

\begin{eqnarray}
\nonumber
\input{figure/fig-decomp2-text} \longleftrightarrow \input{figure/fig-decomp1-text} \hspace{.1cm} ,\hspace{1cm} 
\input{figure/fig-decomp4-text} \longleftrightarrow \input{figure/fig-decomp3-text}.
\end{eqnarray}\\[3em]
\caption{Wide edge and a bunch of $k$ $1$-colored lines}\label{wide-edge}
\end{figure}

We introduce a wide edge to define a complex of matrix factorizations for $[i,j]$-crossing.
The wide edge represents a bunch of $1$-colored lines with the same orientation.
We suppose that a $k$-colored edge branches into a bunch of $k$ $1$-colored lines and a bunch of $k$ $1$-colored lines joins into a $k$-colored edge, see Figure \ref{wide-edge}.
We naturally consider a crossing of a wide edge and colored edge and a crossing of wide edges. For example,
\begin{eqnarray}
\nonumber
\input{figure/fig-decomp-crossing2}\longleftrightarrow \input{figure/fig-decomp-crossing1}\hspace{.1cm} ,\hspace{1cm}
\input{figure/fig-decomp-crossing4}\longleftrightarrow \input{figure/fig-decomp-crossing3}.
\end{eqnarray}\\[6em]
%
%
\begin{proposition}\label{prop-wide}
We have isomorphisms in $\k^b(\HMF^{gr}_{R,\omega})$
\begin{eqnarray}
\nonumber
\c\left(
\input{figure/fig-wide-noori1}
\right)_n
\simeq
\c\left(
\input{figure/fig-wide-noori2}
\right)_n,
\hspace{1cm}
\c\left(
\input{figure/fig-wide-noori3}
\right)_n
\simeq
\c\left(
\input{figure/fig-wide-noori4}
\right)_n.
\end{eqnarray}
For diagrams with the other crossing, their complexes are isomorphic in $\k^b(\HMF^{gr}_{R,\omega})$.
\end{proposition}
\begin{proof}
We find this property by Proposition \ref{prop-r3}.
\end{proof}
%
%
\begin{corollary}\label{cor-6.4}
We have isomorphisms in $\k^b(\HMF^{gr}_{R,\omega})$
\begin{equation}
\nonumber
\c\left(
\input{figure/fig-wide5}
\right)_n
\simeq
\c\left(
\input{figure/fig-wide6}
\right)_n,
\hspace{1cm}
\c\left(
\input{figure/fig-wide7}
\right)_n
\simeq
\c\left(
\input{figure/fig-wide8}
\right)_n.
\end{equation}
For diagrams with the other crossing, their complexes are isomorphic in $\k^b(\HMF^{gr}_{R,\omega})$.
\end{corollary}
\begin{proof}
We immediately find this corollary by Proposition \ref{prop-wide}.
\end{proof}
%
%
Proposition \ref{prop-r3} and Proposition \ref{prop-twist} give the following properties of colored planar diagrams with a wide edge.
\begin{corollary}\label{cor-twist}
The following isomorphisms exist in $\k^b(\HMF^{gr}_{R,\omega})$:
\begin{eqnarray}
\nonumber
&(1)&
\c\left(\input{figure/fig-plus-1k-tri1-wide}\right)_n\simeq\c\left(\input{figure/fig-planar-1k-tri1-wide}\right)_n
\{kn+k\}\left<k\right>\left[-k\right],
\hspace{.5cm}
\c\left(\input{figure/fig-plus-k1-tri1-wide}\right)_n\simeq\c\left(\input{figure/fig-planar-k1-tri1-wide}\right)_n
\{kn+k\}\left<k\right>\left[-k\right],\\
\nonumber
&(2)&
\c\left(\input{figure/fig-minus-1k-tri1-wide}\right)_n\simeq\c\left(\input{figure/fig-planar-1k-tri1-wide}\right)_n
\{-kn-k\}\left<k\right>\left[k\right],
\hspace{.5cm}
\c\left(\input{figure/fig-minus-k1-tri1-wide}\right)_n\simeq\c\left(\input{figure/fig-planar-k1-tri1-wide}\right)_n
\{-kn-k\}\left<k\right>\left[k\right],\\
\nonumber
&(3)&
\c\left(\input{figure/fig-plus-1k-tri2-wide}\right)_n\simeq\c\left(\input{figure/fig-planar-1k-tri2-wide}\right)_n
\{kn+k\}\left<k\right>\left[-k\right],
\hspace{.5cm}
\c\left(\input{figure/fig-plus-k1-tri2-wide}\right)_n\simeq\c\left(\input{figure/fig-planar-k1-tri2-wide}\right)_n
\{kn+k\}\left<k\right>\left[-k\right],\\
\nonumber
&(4)&
\c\left(\input{figure/fig-minus-1k-tri2-wide}\right)_n\simeq\c\left(\input{figure/fig-planar-1k-tri2-wide}\right)_n
\{-kn-k\}\left<k\right>\left[k\right],
\hspace{.5cm}
\c\left(\input{figure/fig-minus-k1-tri2-wide}\right)_n\simeq\c\left(\input{figure/fig-planar-k1-tri2-wide}\right)_n
\{-kn-k\}\left<k\right>\left[-k\right].
\end{eqnarray}
\end{corollary}
%
%
\subsection{Approximate complex for $[i,j]$-crossing}\label{sec6.2}
We consider an approximate crossing of an $[i,j]$-crossing in Figure \ref{approximate-diag}.
This approximate crossing has only $[i,1]$-crossings.
Thus, we can define a complex for matrix factorization for the approximate crossing using definition of $[i,1]$-crossing in Section \ref{sec5.1}.
\begin{definition}
We define a complex of matrix factorization for an $[i,j]$-crossing as an object of $\k^b(\HMF^{gr}_{R,\omega})$ for its approximate crossing:
\begin{eqnarray}
\nonumber
&&
\overline{\c}\left(
\input{figure/figplus-color}
\right)_n:=
\c\left(
\input{figure/figplus-approximate-color}
\right)_n\{-i(i-1)(n+1)\}\left[i(i-1)\right]
\hspace{.5cm}
(i\geq j),
\\
\nonumber
&&
\overline{\c}\left(
\input{figure/figplus-color}
\right)_n:=
\c\left(
\input{figure/figplus-approximate-color}
\right)_n\{-j(j-1)(n+1)\}\left[j(j-1)\right]
\hspace{.5cm}
(i<j),
\\
\nonumber
&&
\overline{\c}\left(
\input{figure/figminus-color}
\right)_n:=
\c\left(
\input{figure/figminus-approximate-color}
\right)_n\{j(j-1)(n+1)\}\left[-j(j-1)\right]
\hspace{.5cm}
(i\leq j),
\\
\nonumber
&&
\overline{\c}\left(
\input{figure/figminus-color}
\right)_n:=
\c\left(
\input{figure/figminus-approximate-color}
\right)_n\{i(i-1)(n+1)\}\left[-i(i-1)\right]
\hspace{.5cm}
(i>j).
\end{eqnarray}
\end{definition}
%
%
\begin{theorem}\label{main2}
We have the following isomorphisms in $\k^b(\HMF^{gr}_{R,\omega})$
\begin{eqnarray}
\nonumber
(\overline{I})
&&
\overline{\c}\left(\input{figure/r1p-i}\right)_n
\hspace{.2cm}\simeq\hspace{.2cm}
\overline{\c}\left(\input{figure/r1c-i}\right)_n\{[i]_{q}!\}_{q}\hspace{.2cm}\simeq\hspace{.2cm}
\overline{\c}\left(\input{figure/r1m-i}\right)_n,\\
\nonumber
(\overline{IIa})
&&
\overline{\c}\left(\input{figure/r2-ijl}\right)_n
\hspace{.2cm}\simeq\hspace{.2cm}
\overline{\c}\left(\input{figure/r2-ijc}\right)_n\{[i]_{q}![j]_{q}!\}_{q}
\hspace{.2cm}\simeq\hspace{.2cm}
\overline{\c}\left(\input{figure/r2-ijr}\right)_n,\\
\nonumber
(\overline{IIb})
&&
\overline{\c}\left(\input{figure/r2-ijlrev}\right)_n
\hspace{.2cm}\simeq\hspace{.2cm}
\overline{\c}\left(\input{figure/r2-ijcrev}\right)_n\{[i]_{q}![j]_{q}!\}_{q},
\hspace{.5cm}
\overline{\c}\left(\input{figure/r2-ijlrevori}\right)_n
\hspace{.2cm}\simeq\hspace{.2cm}
\overline{\c}\left(\input{figure/r2-ijcrevori}\right)_n\{[i]_{q}![j]_{q}!\}_{q},\\
\nonumber
(\overline{III})&&
\overline{\c}\left(\input{figure/r3-ijk}\right)_n
\hspace{.2cm}\simeq\hspace{.2cm}
\overline{\c}\left(\input{figure/r3-ijkrev}\right)_n.
\end{eqnarray}
\end{theorem}
For a colored oriented link diagram $D$, we obtain the homology associated to the complex $\overline{\c}(D)$.
We consider the Poincar\'e polynomial $\overline{P}(D)$ in $\Q[t^{\pm 1},q^{\pm 1},s]/\langle s^2-1\rangle$ of the homology associated to the complex $\overline{\c}(D)$.
We have the following properties.
\begin{corollary}\label{main-cor}
For colored oriented link diagrams transforming to each other under colored Reidemeister moves, we have the following equations of the evaluation of Poincar\'e polynomial $\overline{P}$ for diagrams:
\begin{eqnarray}
\nonumber
(\overline{I})
&&
\overline{P}\left(\input{figure/r1p-i}\right)_n
\hspace{.2cm}=\hspace{.2cm}
\overline{P}\left(\input{figure/r1c-i}\right)_n[i]_{q}!
\hspace{.2cm}=\hspace{.2cm}
\overline{P}\left(\input{figure/r1m-i}\right)_n,\\
\nonumber
(\overline{IIa})
&&
\overline{P}\left(\input{figure/r2-ijl}\right)_n
\hspace{.2cm}=\hspace{.2cm}
\overline{P}\left(\input{figure/r2-ijc}\right)_n[i]_{q}![j]_{q}!
\hspace{.2cm}=\hspace{.2cm}
\overline{P}\left(\input{figure/r2-ijr}\right)_n,\\
\nonumber
(\overline{IIb})
&&
\overline{P}\left(\input{figure/r2-ijlrev}\right)_n
\hspace{.2cm}=\hspace{.2cm}
\overline{P}\left(\input{figure/r2-ijcrev}\right)_n[i]_{q}![j]_{q}!,
\hspace{.5cm}
\overline{P}\left(\input{figure/r2-ijlrevori}\right)_n
\hspace{.2cm}=\hspace{.2cm}
\overline{P}\left(\input{figure/r2-ijcrevori}\right)_n[i]_{q}![j]_{q}!,\\
\nonumber
(\overline{III})&&
\overline{P}\left(\input{figure/r3-ijk}\right)_n
\hspace{.2cm}=\hspace{.2cm}
\overline{P}\left(\input{figure/r3-ijkrev}\right)_n.
\end{eqnarray}
\end{corollary}
We normalize the Poincar\'e polynomial $\overline{P}$ associated to $\overline{\c}(D)$.
We define the function $\mathrm{Cr}_k$ $(k=1,...,n-1)$ on a colored oriented link diagram $D$ as follows,
\begin{center}
$\mathrm{Cr}_k(D):=$ the number of $[\ast,k]$-crossings in $D$.
\end{center}
We define a normalized Poincar\'e polynomial $P(D)$ to be
\begin{equation}
\nonumber
P(D):=\overline{P}(D)\prod_{k=1}^{n-1}\frac{1}{\left([k]_q!\right)^{\mathrm{Cr}_k(D)}}.
\end{equation}
By Corollary \ref{main-cor} the following corollary is obtained.
\begin{corollary}\label{main3}
For two colored oriented link diagrams $D$ and $D'$ transforming to each other under colored Reidemeister moves, these evaluations of $P$ are the same,
\begin{equation}
\nonumber
P(D)=P(D').
\end{equation}
That is, we have equations for evaluations of colored oriented link diagrams,
\begin{eqnarray}
\nonumber
&&
P\left(\input{figure/r1p-i}\right)
\hspace{.1cm}=\hspace{.1cm}
P\left(\input{figure/r1c-i}\right)
\hspace{.1cm}=\hspace{.1cm}
P\left(\input{figure/r1m-i}\right),
\hspace{.4cm}
P\left(\input{figure/r2-ijl}\right)
\hspace{.1cm}=\hspace{.1cm}
P\left(\input{figure/r2-ijc}\right)
\hspace{.1cm}=\hspace{.1cm}
P\left(\input{figure/r2-ijr}\right),\\
\nonumber
&&
P\left(\input{figure/r2-ijlrev}\right)
\hspace{.1cm}=\hspace{.1cm}
P\left(\input{figure/r2-ijcrev}\right),
\hspace{.4cm}
P\left(\input{figure/r2-ijlrevori}\right)
\hspace{.1cm}=\hspace{.1cm}
P\left(\input{figure/r2-ijcrevori}\right),
\hspace{.4cm}
P\left(\input{figure/r3-ijk}\right)
\hspace{.1cm}=\hspace{.1cm}
P\left(\input{figure/r3-ijkrev}\right).
\end{eqnarray}
where outsides of colored tangle diagrams in each equation have the same picture.
\end{corollary}
The polynomial $P(D)$ is a refined link invariant of the quantum $(\mathfrak{sl}_n,\land V_n)$ link invariant since $P(D)$ is the quantum $(\mathfrak{sl}_n,\land V_n)$ link invariant by specializing $t$ to $-1$ and $s$ to $1$.
\begin{proof}
\noindent
{\bf Proof of Theorem \ref{main2} $(\overline{I})$}
By Corollary \ref{cor-6.4}, We have
\begin{eqnarray}
\nonumber
\overline{\c}\left(\input{figure/r1p-i}\right)_n&=&\c\left(\input{figure/r1p-i-wide}\right)_n\{-i(i-1)(n+1)\}\left[i(i-1)\right]\\
\nonumber
&\simeq&\c\left(\input{figure/r1p-i-wide-proof1}\right)_n\{-i(i-1)(n+1)\}\left[i(i-1)\right]
\end{eqnarray}
We show the following lemma.
\begin{lemma}
We have the following isomorphism in $\k^b(\HMF^{gr})$
\begin{equation}
\c\left(\input{figure/r1p-i-wide-proof1}\right)_n\{-i(i-1)(n+1)\}\left[i(i-1)\right]
\simeq
\c\left(\input{figure/r1c-i-wide-proof1}\right)_n.
\end{equation}
\end{lemma}
\begin{proof}
We prove the lemma by induction over $i$.
If $i = 2$, then we have the following isomorphism by Theorem \ref{main1} and Proposition \ref{prop-twist}:
\begin{eqnarray}
\nonumber
\c\left(\input{figure/r1p-2-wide-proof1}\right)_n\{-2(n+1)\}\left[2\right]&=&\c\left(\input{figure/r1p-2-wide-proof2}\right)_n\{-2(n+1)\}\left[2\right]\\
\nonumber
&\simeq&\c\left(\input{figure/fig-full-twist-wide}\right)_n\{-2(n+1)\}\left[2\right]\\
\nonumber
&\simeq&\c\left(\input{figure/fig-full-twist-wide1}\right)_n=\c\left(\input{figure/r1c-2-wide-proof1}\right)_n.
\end{eqnarray}
We assume that the lemma holds for $i=k-1$ and consider the case $i=k$.
We have the following isomorphism by Theorem \ref{main1} and Proposition \ref{prop-twist}:
\begin{eqnarray}
\nonumber
\c\left(\input{figure/r1p-k-wide-proof1}\right)_n\{-k(k-1)(n+1)\}\left[k(k-1)\right]&=&\c\left(\input{figure/r1p-k-wide-proof2}\right)_n\{-k(k-1)(n+1)\}\left[k(k-1)\right]\\
\nonumber
&\simeq&\c\left(\input{figure/r1p-k-wide-proof3}\right)_n\{-k(k-1)(n+1)\}\left[k(k-1)\right]\\
\label{twist-complex}
&\simeq&\c\left(\input{figure/r1p-k-wide-proof4}\right)_n\{-(k-1)(k-2)(n+1)\}\left[(k-1)(k-2)\right].
\end{eqnarray}
By the assumption of induction, the complex (\ref{twist-complex}) is isomorphic to
\begin{equation}
\nonumber
\c\left(\input{figure/fig-full-twist-wide2}\right)_n=\c\left(\input{figure/r1c-k-wide-proof1}\right)_n.
\end{equation}
\indent
We can similarly prove the other isomorphism for a minus $[i]$-curl.\\\\
\end{proof}
\noindent
{\bf Proof of Theorem \ref{main2} $(\overline{II})$}
\begin{eqnarray}
\nonumber
\overline{\c}\left(\input{figure/r2-ijl}\right)_n&=&\c\left(\input{figure/r2-ijl-wide}\right)_n\simeq\c\left(\input{figure/r2-ijl-wide-proof1}\right)_n \\
\nonumber
&\simeq&\c\left(\input{figure/r2-ijl-wide-proof2}\right)_n\simeq\c\left(\input{figure/r2-ijl-wide-proof3}\right)_n\{[i]_q!\}_{q} \\
\nonumber
&\simeq&\c\left(\input{figure/r2-ijc-wide-proof1}\right)_n\{[i]_q!\}_{q}\simeq\c\left(\input{figure/r2-ijc}\right)_n\{[i]_q![j]_q!\}_{q}.
\end{eqnarray}
We can similarly prove the other isomorphism of $(\overline{IIa})$ and isomorphisms of $(\overline{IIb})$.\\\\
\noindent
{\bf Proof of Theorem \ref{main2} $(\overline{III})$}
It is sufficient that we consider the case $i<j<k$.
We can similarly prove invariance of the Reidemeister moves $(\overline{III})$ for the other coloring case.
\begin{eqnarray}
\nonumber
\overline{\c}\left(\input{figure/r3-ijk}\right)_n&=&
\c\left(\input{figure/r3-ijk-wide}\right)_n\{\alpha\}\left[\beta\right]\hspace{1cm}(\alpha=(-2k(k-1)-j(j-1))(n+1),\beta=2k(k-1)+j(j-1))\\
\nonumber
&\simeq&\c\left(\input{figure/r3-ijk-wide-proof1}\right)_n\{\alpha\}\left[\beta\right]
\simeq\c\left(\input{figure/r3-ijk-wide-proof2}\right)_n\{\alpha\}\left[\beta\right]
\simeq\c\left(\input{figure/r3-ijk-wide-proof3}\right)_n\{\alpha\}\left[\beta\right]\{[k]_q!\}_{q}\\
\nonumber
&\simeq&\c\left(\input{figure/r3-ijk-wide-proof4}\right)_n\{\alpha\}\left[\beta\right]\{[k]_q!\}_{q}
\simeq\c\left(\input{figure/r3-ijkrev-wide-proof2}\right)_n\{\alpha\}\left[\beta\right]\{[k]_q!\}_{q}.
\end{eqnarray}
On the other side, we have
\begin{eqnarray}
\nonumber
\overline{\c}\left(\input{figure/r3-ijkrev}\right)_n
&=&\c\left(\input{figure/r3-ijkrev-wide}\right)_n\{\alpha\}\left[\beta\right]
\simeq\c\left(\input{figure/r3-ijkrev-wide-proof1}\right)_n \{\alpha\}\left[\beta\right]\\
\nonumber
&\simeq&
\c\left(\input{figure/r3-ijkrev-wide-proof3}\right)_n \{ \alpha \}\left[\beta\right]
\simeq\c\left(\input{figure/r3-ijkrev-wide-proof4}\right)_n\{ \alpha \}\left[ \beta\right] \{ [k]_q! \}_{q}
\\
\nonumber
&\simeq&\c\left(\input{figure/r3-ijkrev-wide-proof2}\right)_n\{\alpha\}\left[\beta\right]\{[k]_q!\}_{q}.
\end{eqnarray}
Thus, we have
\begin{equation}
\nonumber
\overline{\c}\left(\input{figure/r3-ijk}\right)_n\simeq\overline{\c}\left(\input{figure/r3-ijkrev}\right)_n.
\end{equation}
\end{proof}
%
%
%
%
\section{Proof of Theorem and Proposition}\label{sec7}
%
%
%
%
\subsection{Invariance under Reidemeister move IIa}\label{remain-inv-IIa}
\indent
We can similarly prove the remains of invariance under the Reidemeister moves $(IIa_{1k})$. 
The complex of matrix factorization $\c\left(\input{figure/r2-1kr-mf}\right)_n$ is described as follows,
\begin{equation}
\label{com-2}
\xymatrix{
-1\ar@{.}[d]&&0\ar@{.}[d]&&1\ar@{.}[d]\\
\overline{M}_{11}\{1\}\ar[rr]^{
\left(
\begin{array}{c}\overline{\mu}_5\\\overline{\mu}_6 \end{array}
\right)
}&&
\text{$\begin{array}{c}\overline{M}_{10}\\ \oplus \\ \overline{M}_{01}\end{array}$}
\ar[rr]^{\txt{$(\overline{\mu}_7 ,\overline{\mu}_8)$}}
&&
\overline{M}_{00}\{-1\},
}
\end{equation}
where
\begin{align}
\nonumber 
&\overline{\mu}_5=\id_{\overline{M}_{(1,2,5,6)}^{[1,k]}}\boxtimes
\left(\id_{\overline{S}^{[1,k]}_{(6,5,4,3)}}\boxtimes(x_{1,6}-x_{1,3},1)\right),
&& \overline{\mu}_6=\left(\id_{\overline{S}^{[1,k]}_{(1,2,5,6)}}\boxtimes(1,x_{1,1}-x_{1,6})\right)
\boxtimes\id_{\overline{N}_{(6,5,4,3)}^{[1,k]}},\\
\nonumber 
&\overline{\mu}_7=\left(\id_{\overline{S}^{[1,k]}_{(1,2,5,6)}}\boxtimes(1,x_{1,1}-x_{1,6})\right)
\boxtimes\id_{\overline{M}_{(6,5,4,3)}^{[1,k]}},
&& \overline{\mu}_8=-\id_{\overline{N}_{(1,2,5,6)}^{[1,k]}}\boxtimes
\left(\id_{\overline{S}^{[1,k]}_{(6,5,4,3)}}\boxtimes(x_{1,6}-x_{1,3},1)\right).
\end{align}
As we discussed the homotopy equivalence of the complex $\c\left(\input{figure/r2-1kl-mf}\right)_n$, we similarly have the following isomorphisms (see isomorphisms (\ref{subcom00-2}), (\ref{subcom10-2}), (\ref{subcom01-2}) and (\ref{subcom11-2})):
\begin{eqnarray}
\nonumber \overline{M}_{00}\{-1\}&\simeq&\overline{S}^{[1,k]}_{(1,2,4,3)}\boxtimes K\left(u^{[1,k]}_{k+1,(1,2,4,3)};(x_{1,1}-x_{1,3})X_{k,(2,3)}^{(k,-1)}\right)_{Q_1\left/\langle X_{k,(2,6)}^{(k,-1)}\rangle\right.}
\{-2k\},
\\
\nonumber \overline{M}_{10}&\simeq&\overline{S}^{[1,k]}_{(1,2,4,3)}\boxtimes K\left(u^{[1,k]}_{k+1,(1,2,4,3)};(x_{1,1}-x_{1,3})X_{k,(2,3)}^{(k,-1)}\right)_{Q_1\left/\langle (x_{1,1}-x_{1,6})X_{k,(2,6)}^{(k,-1)}\rangle\right.}
\{-2k\},
\\
\nonumber \overline{M}_{01}&\simeq&\overline{S}^{[1,k]}_{(1,2,4,3)}\boxtimes K\left(u^{[1,k]}_{k+1,(1,2,4,3)}(x_{1,6}-x_{1,3});(x_{1,1}-x_{1,3})X_{k-1,(2,3,6)}^{(k,-1,-1)}\right)_{Q_1\left/\langle X_{k,(2,6)}^{(k,-1)}\rangle\right.}
\{-2k+2\},
\\
\nonumber \overline{M}_{11}\{1\}&\simeq&\overline{S}^{[1,k]}_{(1,2,4,3)}\boxtimes K\left(u^{[1,k]}_{k+1,(1,2,4,3)}(x_{1,6}-x_{1,3});X_{k,(5,3)}^{(k,-1)}\right)_{Q_1\left/\langle (x_{1,1}-x_{1,6})X_{k,(2,6)}^{(k,-1)}\rangle\right.}
\{-2k+2\},
\end{eqnarray}

where
\begin{equation}
\nonumber
Q_1=R_{(1,2,3,4,5,6)}^{(1,k,1,k,k,1)}\left/\left<X_{1,(1,2)}^{(1,k)}-X_{1,(5,6)}^{(k,1)},\ldots,X_{k,(1,2)}^{(1,k)}-X_{k,(5,6)}^{(k,1)}\right>\right. .
\end{equation} 
For these isomorphisms, the morphisms $\overline{\mu}_5$, $\overline{\mu}_6$, $\overline{\mu}_7$ and $\overline{\mu}_8$ transform into
\begin{align}
&\nonumber\overline{\mu}_5\simeq\id_{\overline{S}^{[1,k]}_{(1,2,4,3)}}\boxtimes(x_{1,6}-x_{1,3},1),
&&\overline{\mu}_7\simeq\id_{\overline{S}^{[1,k]}_{(1,2,4,3)}}\boxtimes(1,1),\\
&\nonumber\overline{\mu}_6\simeq\id_{\overline{S}^{[1,k]}_{(1,2,4,3)}}\boxtimes(1,1),
&&\overline{\mu}_8\simeq-\id_{\overline{S}^{[1,k]}_{(1,2,4,3)}}\boxtimes(x_{1,6}-x_{1,3},1).
\end{align} 
We also consider the isomorphisms $R_1$ and $R_3$ of $Q_1\left/\langle X_{k,(2,6)}^{(k,-1)}\rangle\right.$ and $R_2$ of $Q_1\left/\langle (x_{1,1}-x_{1,6})X_{k,(2,6)}^{(k,-1)}\rangle\right.$. Moreover, we consider an isomorphism $R_5$ of $Q_1\left/\langle (x_{1,1}-x_{1,6})X_{k,(2,6)}^{(k,-1)}\rangle\right.$,
\begin{equation}
\nonumber
R_5= R_{(1,2,3,4)}^{(1,k,1,k)}\oplus (x_{1,6}-x_{1,3})R_{(1,2,3,4)}^{(1,k,1,k)}\oplus\ldots\oplus x_{1,6}^{k-1}(x_{1,6}-x_{1,3})R_{(1,2,3,4)}^{(1,k,1,k)}.
\end{equation}
By these isomorphism, we obtain
\begin{eqnarray}
\nonumber
\overline{M}_{00}\{-1\}&\simeq&\overline{S}^{[1,k]}_{(1,2,4,3)}\boxtimes
K\left(u^{[1,k]}_{k+1,(1,2,4,3)};(x_{1,1}-x_{1,3})X_{k,(2,3)}^{(k,-1)}\right)_{Q_1\left/\langle X_{k,(2,6)}^{(k,-1)}\rangle\right.}
\{-2k\}
\\
\nonumber
&\simeq&\overline{S}^{[1,k]}_{(1,2,4,3)}\boxtimes
\left(
\xymatrix{
R_3\{-2k\}
\ar[rrr]_{
E_{k}\left(
{u^{[1,k]}_{k+1,(1,2,4,3)}}
\right)
}&&&
R_3\{1-n\}
\ar[rrr]_{
E_{k}\left(
{(x_{1,1}-x_{1,3})X_{k,(2,3)}^{(k,-1)}}
\right)
}&&&
R_3\{-2k\}}
\right)\{-2k\}\\
\nonumber
&\simeq&\bigoplus_{i=0}^{k-1}\overline{M}^{[1,k]}_{(1,2,4,3)}\{2i-k\},\\
\nonumber
\overline{M}_{10}&\simeq&\overline{S}^{[1,k]}_{(1,2,4,3)}\boxtimes
K\left(u^{[1,k]}_{k+1,(1,2,4,3)};(x_{1,1}-x_{1,3})X_{k,(2,3)}^{(k,-1)}\right)_{Q_1\left/\langle (x_{1,1}-x_{1,6})X_{k,(2,6)}^{(k,-1)}\rangle\right.}
\{-2k\}
\\
\nonumber
&\simeq&\overline{S}^{[1,k]}_{(1,2,4,3)}\boxtimes
\left(
\xymatrix{
R_5\{-2k\}
\ar[rrr]_{
E_{k+1}\left(
{u^{[1,k]}_{k+1,(1,2,4,3)}}
\right)
}&&&
R_5\{1-n\}
\ar[rrrr]_{
E_{k+1}\left(
{(x_{1,1}-x_{1,3})X_{k,(2,3)}^{(k,-1)}}
\right)
}&&&&
R_5\{-2k\}}
\right)\{-2k\}\\
\nonumber
&\simeq&\bigoplus_{i=0}^{k}\overline{M}^{[1,k]}_{(1,2,4,3)}\{2i-k\},
\end{eqnarray}%
\begin{eqnarray}%
\nonumber
\overline{M}_{01}&\simeq&\overline{S}^{[1,k]}_{(1,2,4,3)}\boxtimes
K\left(u^{[1,k]}_{k+1,(1,2,4,3)}(x_{1,6}-x_{1,3});(x_{1,1}-x_{1,3})X_{k-1,(2,3,6)}^{(k,-1,-1)}\right)_{Q_1\left/\langle X_{k,(2,6)}^{(k,-1)}\rangle\right.}
\{-2k+2\}
\\
\nonumber
&\simeq&\overline{S}^{[1,k]}_{(1,2,4,3)}\boxtimes
\left(
\xymatrix{
R_1\{-2k+2\}
\ar[rr]_{f_1}&&
R_3\{1-n\}
\ar[rr]_{f_2}&&
R_1\{-2k+2\}}
\right)\{-2k+2\}\\
\nonumber
&\simeq&\bigoplus_{i=0}^{k-2}\overline{M}^{[1,k]}_{(1,2,4,3)}\{2i-k+2\}\oplus \overline{L}^{[1,k]}_{(1,2,4,3)},\\
\nonumber
\overline{M}_{11}\{1\}&\simeq&\overline{S}^{[1,k]}_{(1,2,4,3)}\boxtimes
K\left(u^{[1,k]}_{k+1,(1,2,4,3)}(x_{1,6}-x_{1,3});X_{k,(5,3)}^{(k,-1)}\right)_{Q_1\left/\langle (x_{1,1}-x_{1,6})X_{k,(2,6)}^{(k,-1)}\rangle\right.}
\{-2k+2\}\\
\nonumber
&\simeq&\overline{S}^{[1,k]}_{(1,2,4,3)}\boxtimes
\left(\xymatrix{
R_2\{-2k+2\}
\ar[rr]_{f_3}&&
R_5\{1-n\}
\ar[rr]_{f_4}&&
R_2\{-2k+2\}}
\right)\{-2k+2\}\\
\nonumber
&\simeq&\bigoplus_{i=0}^{k-1}\overline{M}^{[1,k]}_{(1,2,4,3)}\{2i-k+2\},\\
\nonumber
f_1&=&\left(
\begin{array}{cc}
\mathfrak{0}_{k-1}&u^{[1,k]}_{k+1,(1,2,4,3)}X_{k,(2,3)}^{(k,-1)}\\
E_{k-1}\left(u^{[1,k]}_{k+1,(1,2,4,3)}\right)&{}^{t}\mathfrak{0}_{k-1}
\end{array}
\right),\\
\nonumber
f_2&=&\left(
\begin{array}{cc}
{}^{t}\mathfrak{0}_{k-1}&E_{k-1}\left((x_{1,1}-x_{1,3})X_{k,(2,3)}^{(k,-1)}\right)\\
x_{1,1}-x_{1,3}&\mathfrak{0}_{k-1}
\end{array}
\right),\\
\nonumber
f_3&=&\left(
\begin{array}{cc}
\mathfrak{0}_{k-1}&u^{[1,k]}_{k+1,(1,2,4,3)}(x_{1,1}-x_{1,3})X_{k,(2,3)}^{(k,-1)}\\
E_{k-1}\left(u^{[1,k]}_{k+1,(1,2,4,3)}\right)&{}^{t}\mathfrak{0}_{k-1}
\end{array}
\right),\\
\nonumber
f_4&=&\left(
\begin{array}{cc}
{}^{t}\mathfrak{0}_{k-1}&E_{k-1}\left((x_{1,1}-x_{1,3})X_{k,(2,3)}^{(k,-1)}\right)\\
1&\mathfrak{0}_{k-1}
\end{array}
\right).
\end{eqnarray}
Then, for these decompositions of matrix factorizations, the morphisms $\overline{\mu}_5$, $\overline{\mu}_6$, $\overline{\mu}_7$ and $\overline{\mu}_8$ transform into
\begin{eqnarray}
\nonumber
\overline{\mu}_5&\simeq&
\left(
\begin{array}{c}
\mathfrak{0}_{k}\\
E_{k}\left(\id_{\overline{M}_{(1,2,4,3)}^{[1,k]}}\right)
\end{array}
\right),\\
\nonumber
\overline{\mu}_6&\simeq&
\left(
\begin{array}{cccc}
&&&\hspace{-0.5cm}-(-1)^{k-1}X_{k-1,(2,3)}^{(k,-1)}\id_{\overline{M}_{(1,2,4,3)}^{[1,k]}}\\
&\hspace{-0.5cm}E_{k}\left(\id_{\overline{M}_{(1,2,4,3)}^{[1,k]}}\right)\hspace{-0.5cm}&&\vdots\\
&&&\hspace{-0.5cm}-(-1)^{1}X_{1,(2,3)}^{(k,-1)}\id_{\overline{M}_{(1,2,4,3)}^{[1,k]}}\\
&\mathfrak{0}_{k}&&\hspace{-0.5cm}\id_{\overline{S}_{(1,2,4,3)}^{[1,k]}}\boxtimes\left((-1)^{k},-(-1)^{k}X_{k,(2,3)}^{(k,-1)}\right)
\end{array}\hspace{-0.2cm}
\right),
\end{eqnarray}
\begin{eqnarray}
\nonumber
\overline{\mu}_7&\simeq&
-\left(
\begin{array}{cccc}
&&&\hspace{-0.5cm}-(-1)^{k}X_{k,(2,3)}^{(k,-1)}\id_{\overline{M}_{(1,2,4,3)}^{[1,k]}}\\
&\hspace{-0.5cm}E_{k}\left(\id_{\overline{M}_{(1,2,4,3)}^{[1,k]}}\right)\hspace{-0.5cm}&&\vdots\\
&&&\hspace{-0.5cm}-(-1)^{1}X_{1,(2,3)}^{(k,-1)}\id_{\overline{M}_{(1,2,4,3)}^{[1,k]}}
\end{array}
\right),\\
\nonumber
\overline{\mu}_8&\simeq&
\left(
\begin{array}{cc}
\mathfrak{0}_{k-1}&\id_{\overline{S}_{(1,2,4,3)}^{[1,k]}}\boxtimes\left(X_{k,(2,3)}^{(k,-1)},1\right)\\
E_{k-1}\left(\id_{\overline{M}_{(1,2,4,3)}^{[1,k]}}\right)&{}^{t}\mathfrak{0}_{k-1}
\end{array}
\right).
\end{eqnarray}
Thus, the complex (\ref{com-2}) is isomorphic, in $\k^b(\HMF^{gr}_{R_{(1,2,3,4)}^{(1,k,1,k)},\omega_1})$, to $\overline{L}^{[1,k]}_{(1,2,4,3)}$:
\begin{equation*}
\c\left(\input{figure/r2-1kr-mf}\right)_n\simeq\c\left(\input{figure/r2-1kc-mf}\right)_n.
\end{equation*}
We can similarly prove the following isomorphisms for the Reidemeister moves $(IIa_{1k})$ with another coloring:
\begin{equation*}
\c\left(\input{figure/r2-k1l-mf}\right)_n\simeq\c\left(\input{figure/r2-k1c-mf}\right)_n\simeq\c\left(\input{figure/r2-k1r-mf}\right)_n.
\end{equation*}
%
%
%
%
\subsection{Invariance under Reidemeister move IIb}\label{remain-inv-IIb}
We show the following isomorphisms
\begin{equation}
\nonumber
\c\left(\input{figure/r2-1klrev-mf}\right)_n\simeq\c\left(\input{figure/r2-1kcrev-mf}\right)_n\simeq\c\left(\input{figure/r2-1krrev-mf}\right)_n.
\end{equation}
The complex $\c\left(\input{figure/r2-1klrev-mf}\right)_n$ is an object of $\k^b(\HMF^{gr}_{R_{(1,2,3,4)}^{(1,k,1,k)},\omega_2})$, $\omega_2=F_1(\mathbb{X}^{(1)}_{(1)})-F_k(\mathbb{X}^{(k)}_{(2)})-F_1(\mathbb{X}^{(1)}_{(3)})+F_k(\mathbb{X}^{(k)}_{(4)})$, 
and this object is isomorphic to
\begin{equation}
\label{com-rev-1}
\xymatrix{
-1\ar@{.}[d]&&0\ar@{.}[d]&&1\ar@{.}[d]\\
\overline{N}_{00}\{1\}\ar[rr]^{
\left(
\begin{array}{c}\overline{\nu}_1\\\overline{\nu}_2 \end{array}
\right)
}&&
\text{$\begin{array}{c}\overline{N}_{10}\\ \oplus \\ \overline{N}_{01}\end{array}$}
\ar[rr]^{\txt{$(\overline{\nu}_3,\overline{\nu}_4)$}
}
&&
\overline{N}_{11}\{-1\},
}
\end{equation}
where
\begin{align}
\nonumber&\overline{N}_{00}=\overline{M}^{[1,k]}_{(1,5,2,6)}\boxtimes\overline{N}^{[1,k]}_{(6,4,5,3)},&&
\overline{N}_{10}=\overline{M}^{[1,k]}_{(1,5,2,6)}\boxtimes\overline{M}^{[1,k]}_{(6,4,5,3)},\\
\nonumber&\overline{N}_{01}=\overline{N}^{[1,k]}_{(1,5,2,6)}\boxtimes\overline{N}^{[1,k]}_{(6,4,5,3)},&&
\overline{N}_{11}=\overline{N}^{[1,k]}_{(1,5,2,6)}\boxtimes\overline{M}^{[1,k]}_{(6,4,5,3)},\\
\nonumber&\overline{\nu}_1=\id_{\overline{M}^{[1,k]}_{(1,5,2,6)}}\boxtimes \id_{\overline{S}^{[1,k]}_{(6,4,5,3)}}\boxtimes (x_{1,6}-x_{1,3},1),&&
\overline{\nu}_2=\id_{\overline{S}^{[1,k]}_{(1,5,2,6)}}\boxtimes (1,x_{1,1}-x_{1,6})\boxtimes\id_{\overline{N}^{[1,k]}_{(6,4,5,3)}},\\
\nonumber&\overline{\nu}_3=\id_{\overline{S}^{[1,k]}_{(1,5,2,6)}}\boxtimes (1,x_{1,1}-x_{1,6})\boxtimes\id_{\overline{M}^{[1,k]}_{(6,4,5,3)}},&&
\overline{\nu}_4=-\id_{\overline{N}^{[1,k]}_{(1,5,2,6)}}\boxtimes\id_{\overline{S}^{[1,k]}_{(6,4,5,3)}}\boxtimes  (x_{1,6}-x_{1,3},1).
\end{align}
By Corollary \ref{cor2-11}, we have 
\begin{eqnarray}
\nonumber
\overline{N}_{00}\{1\}&=&
K\left(
\left(
\begin{array}{c}
A_{1,(1,5,2,6)}^{[1,k]}\\
\vspace{0.2cm}\vdots\\
\vspace{0.2cm}A_{k,(1,5,2,6)}^{[1,k]}\\
u_{k+1,(1,5,2,6)}^{[1,k]}
\end{array}
\right);
\left(
\begin{array}{c}
X_{1,(1,5)}^{(1,k)}-X_{1,(2,6)}^{(k,1)}\\
\vdots\\
X_{k,(1,5)}^{(1,k)}-X_{k,(2,6)}^{(k,1)}\\
(x_{1,1}-x_{1,6})X_{k,(5,6)}^{(k,-1)}
\end{array}
\right)
\right)_{R_{(1,5,2,6)}^{(1,k,k,1)}}\\
&&
\nonumber
\hspace{0.5cm}\boxtimes
K\left(
\left(
\begin{array}{c}
A_{1,(6,4,5,3)}^{[1,k]}\\
\vspace{0.2cm}\vdots\\
\vspace{0.2cm}A_{k,(6,4,5,3)}^{[1,k]}\\
u_{k+1,(6,4,5,3)}^{[1,k]}(x_{1,6}-x_{1,3})
\end{array}
\right);
\left(
\begin{array}{c}
X_{1,(6,4)}^{(1,k)}-X_{1,(5,3)}^{(k,1)}\\
\vdots\\
X_{k,(6,4)}^{(1,k)}-X_{k,(5,3)}^{(k,1)}\\
X_{k,(3,4)}^{(-1,k)}
\end{array}
\right)
\right)_{R_{(6,4,5,3)}^{(1,k,k,1)}}\hspace{-1cm}\{-2k+1\}\{1\}\\
\nonumber
&\simeq&K\left(
\left(
\begin{array}{c}
A_{1,(6,4,5,3)}^{[1,k]}\\
\vspace{0.2cm}\vdots\\
\vspace{0.2cm}A_{k,(6,4,5,3)}^{[1,k]}\\
u_{k+1,(6,4,5,3)}^{[1,k]}(x_{1,6}-x_{1,3})
\end{array}
\right);
\left(
\begin{array}{c}
X_{1,(6,4)}^{(1,k)}-X_{1,(1,2,3,6)}^{(-1,k,1,1)}\\
\vdots\\
X_{k,(6,4)}^{(1,k)}-X_{k,(1,2,3,6)}^{(-1,k,1,1)}\\
X_{k,(3,4)}^{(-1,k)}
\end{array}
\right)
\right)_{Q_2\left/\langle u_{k+1,(1,5,2,6)}^{[1,k]}\rangle\right.}\hspace{-2cm}\{3-n\}\left<1\right>,
\end{eqnarray}
where $Q_2=R_{(1,2,3,4,5,6)}^{(1,k,1,k,k,1)}\left/\left<X_{1,(1,5)}^{(1,k)}-X_{1,(2,6)}^{(k,1)},
\ldots,X_{k,(1,5)}^{(1,k)}-X_{k,(2,6)}^{(k,1)}\right>\right.$.
Moreover, by Theorem \ref{equiv}, the matrix factorization is isomorphic to
\begin{eqnarray}
\nonumber
\overline{N}_{00}\{1\}&\simeq&
K\left(
\left(
\begin{array}{c}
A_{1,(6,4,5,3)}^{[1,k]}+(x_{1,6}-x_{1,1})A_{2,(6,4,5,3)}^{[1,k]}\\
\vspace{0.1cm}\vdots\\
\vspace{0.1cm}A_{k-1,(6,4,5,3)}^{[1,k]}+(x_{1,6}-x_{1,1})A_{k,(6,4,5,3)}^{[1,k]}\\
\vspace{0.1cm}A_{k,(6,4,5,3)}^{[1,k]}\\
u_{k+1,(6,4,5,3)}^{[1,k]}(x_{1,6}-x_{1,3})
\end{array}
\right);
\left(
\begin{array}{c}
X_{1,(1,4)}^{(1,k)}-X_{1,(2,3)}^{(k,1)}\\
\vspace{0.1cm}\vdots\\
\vspace{0.1cm}X_{k-1,(1,4)}^{(1,k)}-X_{k-1,(2,3)}^{(k,1)}\\
\vspace{0.1cm}X_{k,(1,4)}^{(1,k)}-X_{k,(2,3)}^{(k,1)}\\
X_{k,(3,4)}^{(-1,k)}
\end{array}
\right)
\right)_{Q_2\left/\langle u_{k+1,(1,5,2,6)}^{[1,k]}\rangle\right.}\hspace{-2cm}\{3-n\}\left<1\right>\\
\label{subcom-rev-1}
&\simeq&
\overline{S}^{[1,k]}_{(1,4,2,3)}\boxtimes
K\left(u_{k+1,(1,4,2,3)}^{[1,k]}(x_{1,1}-x_{1,3});X_{k,(3,4)}^{(-1,k)}\right)_{Q_2\left/\langle u_{k+1,(1,5,2,6)}^{[1,k]}\rangle\right.}\{3-n\}\left<1\right>.
\end{eqnarray}
By Theorem \ref{equiv} and Corollary \ref{induce-sq1} (2), we also have
\begin{eqnarray}
\label{subcom-rev-2}
\overline{N}_{10}&\simeq&
\overline{S}^{[1,k]}_{(1,4,2,3)}\boxtimes
K\left(u_{k+1,(6,4,5,3)}^{[1,k]}+\alpha;(x_{1,6}-x_{1,3})X_{k,(3,4)}^{(-1,k)}\right)_{Q_2\left/\langle u_{k+1,(1,5,2,6)}^{[1,k]}\rangle\right.}\{1-n\}\left<1\right>,
\end{eqnarray}
where $\alpha$ is a linear combination of $X_{i,(1,4)}^{(1,k)}-X_{i,(2,3)}^{(k,1)}$ $(1 \leq i \leq k)$ with $\Z$-grading $2n-2k$ satisfying that 
\begin{equation}
\label{condition}
(u_{k+1,(6,4,5,3)}^{[1,k]}+\alpha)(x_{1,6}-x_{1,3})\equiv u_{k+1,(1,4,2,3)}^{[1,k]}(x_{1,1}-x_{1,3}) \hspace{1cm}\mod Q_2\left/\langle u_{k+1,(1,5,2,6)}^{[1,k]}\rangle\right. .
\end{equation}
Since, in the quotient $Q_2\left/\langle u_{k+1,(1,5,2,6)}^{[1,k]}\rangle\right.$, we have equations
\begin{equation}
\nonumber
X_{i,(5)}^{(k)}=X_{i(1,2,6)}^{(-1,k,1)}\hspace{1cm} (1 \leq i \leq k),
\end{equation}
then the polynomial $u_{k+1,(1,5,2,6)}^{[1,k]}$ is described as
\begin{eqnarray}
\nonumber
&&\frac{F_{k+1}(X_{1,(2,6)}^{(k,1)},\ldots,X_{k,(2,6)}^{(k,1)},X_{k+1,(1,5)}^{(1,k)})-F_{k+1}(X_{1,(2,6)}^{(k,1)},\ldots,X_{k,(2,6)}^{(k,1)},X_{k+1,(2,6)}^{(k,1)})}{X_{k+1,(1,5)}^{(1,k)}-X_{k+1,(2,6)}^{(k,1)}}\\
\nonumber
&=&c_1(X_{1,(2,6)}^{(k,1)})^{n-k}+c_2(X_{1,(2,6)}^{(k,1)})^{n-k-2}X_{2,(2,6)}^{(k,1)}+\ldots\\
\nonumber
&=&c_1 x_{1,6}^{n-k}+c_3 x_{1,2}x_{1,6}^{n-k-1}+\ldots,
\end{eqnarray}
where $c_1$ and $c_2$ are the coefficients of $F_{k+1}(x_1,x_2,...,x_{k+1})=c_1x_1^{n-k} x_{k+1}+c_2x_1^{m-k-2} e_2 x_{k+1}+...$ and $c_3=c_1(n-k)+c_2$.
Then, in the quotient $Q_2\left/\langle u_{k+1,(1,5,2,6)}^{[1,k]}\rangle\right.$, the polynomial $u_{k+1,(6,4,5,3)}^{[1,k]}$ is described as
\begin{eqnarray}
\nonumber
&&\frac{F_{k+1}(X_{1,(1,2,3,6)}^{(-1,k,1,1)},\ldots,X_{k,(1,2,3,6)}^{(-1,k,1,1)},X_{k+1,(4,6)}^{(k,1)})-F_{k+1}(X_{1,(1,2,3,6)}^{(-1,k,1,1)},\ldots,X_{k,(1,2,3,6)}^{(-1,k,1,1)},X_{k+1,(1,2,3,6)}^{(-1,k,1,1)})}{X_{k+1,(4,6)}^{(k,1)}-X_{k+1,(1,2,3,6)}^{(-1,k,1,1)}}\\
\nonumber
&=&c_1(X_{1,(1,2,3,6)}^{(-1,k,1,1)})^{n-k}+c_2(X_{1,(1,2,3,6)}^{(-1,k,1,1)})^{n-k-2}X_{2,(1,2,3,6)}^{(-1,k,1,1)}+\ldots\\
\nonumber
&=&c_1 x_{1,6}^{n-k}+c_3(-x_{1,1}+x_{1,2}+x_{1,3})x_{1,6}^{n-k-1}+\ldots\\
\nonumber
&\equiv&-c_3(x_{1,1}-x_{1,3})x_{1,6}^{n-k-1}+\ldots\hspace{1cm}\mod Q_2\left/\langle u_{k+1,(1,5,2,6)}^{[1,k]}\rangle\right. .
\end{eqnarray}
By the condition (\ref{condition}), we find
\begin{equation}
u_{k+1,(6,4,5,3)}^{[1,k]}+\alpha\equiv -c_3(x_{1,1}-x_{1,3})(x_{1,6}^{n-k-1}+\beta)\hspace{1cm}\mod Q_2\left/\langle u_{k+1,(1,5,2,6)}^{[1,k]}\rangle\right. ,
\end{equation}
where $\beta$ is a polynomial with $\Z$-grading $2n-2k-2$ such that the degree as a polynomial of variable $x_{1,6}$ is less than $n-k-1$ and $-(x_{1,6}-x_{1,3})c_3(x_{1,6}^{n-k-1}+\beta) \equiv u_{k+1,(1,4,2,3)}^{[1,k]} \hspace{0.1cm}(\mod Q_2\left/\langle u_{k+1,(1,5,2,6)}^{[1,k]}\rangle\right.)$.
Thus, the matrix factorization (\ref{subcom-rev-2}) forms into
\begin{eqnarray}
\label{subcom-rev-22}
\overline{N}_{10}&\simeq&
\overline{S}^{[1,k]}_{(1,4,2,3)}\boxtimes
K\left(-c_3(x_{1,1}-x_{1,3})(x_{1,6}^{n-k-1}+\beta);(x_{1,6}-x_{1,3})X_{k,(3,4)}^{(-1,k)}\right)_{Q_2\left/\langle u_{k+1,(1,5,2,6)}^{[1,k]}\rangle\right.}\{1-n\}\left<1\right>.
\end{eqnarray}
By Theorem \ref{equiv}, the matrix factorizations $\overline{N}_{01}$ and $\overline{N}_{11}\{-1\}$ are isomorphic to the following matrix factorizations
\begin{eqnarray}
\label{subcom-rev-3}
\overline{N}_{01}
&\simeq&
\overline{S}^{[1,k]}_{(1,4,2,3)}\boxtimes
K\left(u_{k+1,(1,4,2,3)}^{[1,k]}(x_{1,1}-x_{1,3});X_{k,(3,4)}^{(-1,k)}\right)_{Q_2\left/\langle u_{k+1,(1,5,2,6)}^{[1,k]}(x_{1,1}-x_{1,6})\rangle\right.}\{1-n\}\left<1\right>,\\
\label{subcom-rev-4}
\overline{N}_{11}\{-1\}&\simeq&
\overline{S}^{[1,k]}_{(1,4,2,3)}\boxtimes
K\left(u_{k+1,(6,4,5,3)}^{[1,k]}+\alpha;(x_{1,6}-x_{1,3})X_{k,(3,4)}^{(-1,k)}\right)_{Q_2\left/\langle u_{k+1,(1,5,2,6)}^{[1,k]}(x_{1,1}-x_{1,6})\rangle\right.}\{-1-n\}\left<1\right>
\end{eqnarray}
\indent
We consider the following isomorphisms of $Q_2\left/\langle u_{k+1,(1,5,2,6)}^{[1,k]}\rangle\right.$:
\begin{eqnarray}
\nonumber
R_6&:=&R_{(1,2,3,4)}^{(1,k,1,k)}\oplus x_{1,6}R_{(1,2,3,4)}^{(1,k,1,k)}\oplus \ldots\oplus x_{1,6}^{n-k-2}R_{(1,2,3,4)}^{(1,k,1,k)}\oplus -c_3(x_{1,6}^{n-k-1}+\beta)R_{(1,2,3,4)}^{(1,k,1,k)},\\
\nonumber
R_7&:=&R_{(1,2,3,4)}^{(1,k,1,k)}\oplus (x_{1,6}-x_{1,3})R_{(1,2,3,4)}^{(1,k,1,k)}\oplus x_{1,6}(x_{1,6}-x_{1,3})R_{(1,2,3,4)}^{(1,k,1,k)}\oplus \ldots\oplus x_{1,6}^{n-k-2}(x_{1,6}-x_{1,3})R_{(1,2,3,4)}^{(1,k,1,k)},
\end{eqnarray}
and the following isomorphisms of $Q_2\left/\langle u_{k+1,(1,5,2,6)}^{[1,k]}(x_{1,1}-x_{1,6})\rangle\right.$:
\begin{eqnarray}
\nonumber
R_8&:=&R_{(1,2,3,4)}^{(1,k,1,k)}\oplus (x_{1,1}-x_{1,6})R_{(1,2,3,4)}^{(1,k,1,k)}\oplus \ldots \oplus x_{1,6}^{n-k-2}(x_{1,1}-x_{1,6})R_{(1,2,3,4)}^{(1,k,1,k)}\oplus (u_{k+1,(1,4,2,3)}^{[1,k]}+\alpha)R_{(1,2,3,4)}^{(1,k,1,k)},\\
\nonumber
R_9&:=&R_{(1,2,3,4)}^{(1,k,1,k)}\oplus (x_{1,6}-x_{1,3})R_{(1,2,3,4)}^{(1,k,1,k)}\oplus (x_{1,1}-x_{1,6})(x_{1,6}-x_{1,3})R_{(1,2,3,4)}^{(1,k,1,k)}\oplus \\
\nonumber
&&\hspace{1cm}x_{1,6}(x_{1,1}-x_{1,6})(x_{1,6}-x_{1,3})R_{(1,2,3,4)}^{(1,k,1,k)}\oplus\ldots \oplus x_{1,6}^{n-k-2}(x_{1,1}-x_{1,6})(x_{1,6}-x_{1,3})R_{(1,2,3,4)}^{(1,k,1,k)}.
\end{eqnarray}
Then, the partial matrix factorization $K\left(u_{k+1,(1,4,2,3)}^{[1,k]}(x_{1,1}-x_{1,3});X_{k,(3,4)}^{(-1,k)}\right)_{Q_2\left/\langle u_{k+1,(1,5,2,6)}^{[1,k]}\rangle\right.}\{3-n\}\left<1\right>$ of (\ref{subcom-rev-1}) forms into
\begin{equation}
\nonumber
\left(\xymatrix{R_6\ar[rrrrr]_(.4){E_{n-k}\left(u_{k+1,(1,4,2,3)}^{[1,k]}(x_{1,1}-x_{1,3})\right)}&&&&&R_6\{2k-n-1\}\ar[rrr]_(.6){E_{n-k}\left(X_{k,(3,4)}^{(-1,k)}\right)}&&&R_6}\right)\{3-n\}\left<1\right>,
\end{equation}
the partial matrix factorization $K\left(-c_3(x_{1,1}-x_{1,3})(x_{1,6}^{n-k-1}+\beta);(x_{1,6}-x_{1,3})X_{k,(3,4)}^{(-1,k)}\right)_{Q_2\left/\langle u_{k+1,(1,5,2,6)}^{[1,k]}\rangle\right.}\{1-n\}\left<1\right>$ of (\ref{subcom-rev-22}) forms into
\begin{eqnarray}
\nonumber
&&\left(\xymatrix{R_7\ar[rr]_(.4){g_1}&&R_6\{2k-n+1\}\ar[rr]_(.7){g_2}&&R_6}\right)\{1-n\}\left<1\right>,\\
\nonumber
&&g_1=\left(
\begin{array}{cc}
{}^{t}\mathfrak{0}_{n-k-2}&E_{n-k-2}\left(u_{k+1,(1,4,2,3)}^{[1,k]}(x_{1,1}-x_{1,3})\right)\\
x_{1,1}-x_{1,3}&\mathfrak{0}_{n-k-2}
\end{array}
\right)
,\\
\nonumber
&&
g_2=\left(
\begin{array}{cc}
\mathfrak{0}_{n-k-2}&X_{k,(3,4)}^{(-1,k)}u_{k+1,(1,4,2,3)}^{[1,k]}\\
E_{n-k-2}\left(X_{k,(3,4)}^{(-1,k)}\right)&{}^{t}\mathfrak{0}_{n-k-2}
\end{array}
\right)
,
\end{eqnarray}
the partial matrix factorization $K\left(u_{k+1,(1,4,2,3)}^{[1,k]}(x_{1,1}-x_{1,3});X_{k,(3,4)}^{(-1,k)}\right)_{Q_2\left/\langle u_{k+1,(1,5,2,6)}^{[1,k]}(x_{1,1}-x_{1,6})\rangle\right.}\{1-n\}\left<1\right>$ of (\ref{subcom-rev-3}) forms into
\begin{equation}
\nonumber
\left(\xymatrix{R_8\ar[rrrrr]_(.4){E_{n-k+1}\left(u_{k+1,(1,4,2,3)}^{[1,k]}(x_{1,1}-x_{1,3})\right)}&&&&&R_8\{2k-n-1\}\ar[rrr]_(.6){E_{n-k+1}\left(X_{k,(3,4)}^{(-1,k)}\right)}&&&R_8}\right)\{1-n\}\left<1\right>
\end{equation}
and the partial matrix factorization $K\left(u_{k+1,(6,4,5,3)}^{[1,k]};(x_{1,6}-x_{1,3})X_{k,(3,4)}^{(-1,k)}\right)_{Q_2\left/\langle u_{k+1,(1,5,2,6)}^{[1,k]}(x_{1,1}-x_{1,6})\rangle\right.}\{3-n\}\left<1\right>$ of (\ref{subcom-rev-4}) forms into
\begin{eqnarray}
\nonumber
&&\left(\xymatrix{R_9\ar[rr]_(.4){g_3}&&R_8\{2k-n+1\}\ar[rr]_(.6){g_4}&&R_9}\right)\{-1-n\}\left<1\right>,\\
\nonumber
&&
g_3=\left(
\begin{array}{cc}
{}^{t}\mathfrak{0}_{n-k}&E_{n-k}\left(u_{k+1,(1,4,2,3)}^{[1,k]}(x_{1,1}-x_{1,3})\right)\\
1&\mathfrak{0}_{n-k}
\end{array}
\right)
,\\
\nonumber
&&
g_4=\left(
\begin{array}{cc}
\mathfrak{0}_{n-k}&X_{k,(3,4)}^{(-1,k)}u_{k+1,(1,4,2,3)}^{[1,k]}(x_{1,1}-x_{1,3})\\
E_{n-k}\left(X_{k,(3,4)}^{(-1,k)}\right)&{}^{t}\mathfrak{0}_{n-k}
\end{array}
\right)
.
\end{eqnarray}
By these decompositions, we obtain the following isomorphisms
\begin{eqnarray}
\nonumber
\overline{N}_{00}&\simeq&\bigoplus_{i=0}^{n-k-1}\overline{N}^{[1,k]}_{(1,4,2,3)}\{2i-n+k+2\}\left<1\right>,\\
\nonumber
\overline{N}_{10}&\simeq&\bigoplus_{i=0}^{n-k-2}\overline{N}^{[1,k]}_{(1,4,2,3)}\{2i-n+k+2\}\left<1\right>\oplus \overline{L}_{(1,4,2,3)}^{[1,k]},\\
\nonumber
\overline{N}_{01}&\simeq&\bigoplus_{i=0}^{n-k}\overline{N}^{[1,k]}_{(1,4,2,3)}\{2i-n+k\}\left<1\right>,\\
\nonumber
\overline{N}_{11}&\simeq&\bigoplus_{i=0}^{n-k-1}\overline{N}^{[1,k]}_{(1,4,2,3)}\{2i-n+k\}\left<1\right>.
\end{eqnarray}
For these decompositions, the morphisms $\overline{\nu}_1$, $\overline{\nu}_2$, $\overline{\nu}_3$ and $\overline{\nu}_4$ transform as follows,
\begin{eqnarray}
\nonumber
\overline{\nu}_1&\simeq&
\left(
\begin{array}{cc}
E_{n-k-1}\left(\id_{\overline{N}^{[1,k]}_{(1,4,2,3)}}\right)&{}^{t}\mathfrak{0}_{n-k-1}\\
\mathfrak{0}_{n-k-1}&(1,u_{k+1,(1,4,2,3)}^{[1,k]})
\end{array}
\right),\\
\nonumber
\overline{\nu}_2&\simeq&
\left(
\begin{array}{cc}
\mathfrak{0}_{n-k-1}&-u_{k+1,(1,4,2,3)}^{[1,k]}\id_{\overline{N}^{[1,k]}_{(1,4,2,3)}}\\
E_{n-k-1}\left(\id_{\overline{N}^{[1,k]}_{(1,4,2,3)}}\right)&{}^{t}\mathfrak{0}_{n-k-1}\\
\mathfrak{0}_{n-k-1}&\id_{\overline{N}^{[1,k]}_{(1,4,2,3)}}
\end{array}
\right),\\
\nonumber
\overline{\nu}_3&\simeq&
-\left(
\begin{array}{cc}
\mathfrak{0}_{n-k-1}&(-u_{k+1,(1,4,2,3)}^{[1,k]},-1)\\
E_{n-k-1}\left(\id_{\overline{N}^{[1,k]}_{(1,4,2,3)}}\right)&{}^{t}\mathfrak{0}_{n-k-1}
\end{array}
\right),\\
\nonumber
\overline{\nu}_4&\simeq&
\left(
\begin{array}{cc}
E_{n-k-1}\left(\id_{\overline{N}^{[1,k]}_{(1,4,2,3)}}\right)&{}^{t}\mathfrak{0}_{n-k-1}
\end{array}
\right).
\end{eqnarray}
Thus, the complex (\ref{com-rev-1}) is isomorphic, in $\k^b(\HMF^{gr}_{R_{(1,2,3,4)}^{(1,k,1,k)},\omega_2})$, to $\overline{L}^{[1,k]}_{(1,4,2,3)}$:
\begin{equation}
\nonumber
\c\left(\input{figure/r2-1klrev-mf}\right)_n\simeq\c\left(\input{figure/r2-1kcrev-mf}\right)_n.
\end{equation}
\begin{remark}
The above isomorphism
\begin{equation}
\nonumber
\overline{N}_{10}\simeq\overline{L}_{(1,4,2,3)}^{[1,k]}\oplus\bigoplus_{i=0}^{n-k-2}\overline{N}^{[1,k]}_{(1,4,2,3)}\{2i-n+k+2\}\left<1\right>
\end{equation}
is corresponding to the MOY relation
\begin{equation}
\nonumber
\left<\input{figure/figsquare1k--k+1--k1--k+1--1k-rev}\right>_n=
\left<\input{figure/figsquare1k-rev}\right>_n+
\left[n-k-1\right]_q\left<\input{figure/figsquare1k--k-1--1k-rev}\right>_n.
\end{equation}
\end{remark}
\indent
We show the remains of invariance under the Reidemeister moves $(IIb_{1k})$. 
The complex of matrix factorization $\c\left(\input{figure/r2-1krrev-mf}\right)_n$ is described as follows,
\begin{equation}
\label{com-rev-2}
\xymatrix{
-1\ar@{.}[d]&&0\ar@{.}[d]&&1\ar@{.}[d]\\
\overline{N}_{11}\{1\}
\ar[rr]^{
\left(
\begin{array}{c}\overline{\nu}_5\\\overline{\nu}_6 \end{array}
\right)
}&&
\text{$\begin{array}{c}\overline{N}_{10}\\ \oplus \\ \overline{N}_{01}\end{array}$}
\ar[rr]^{\txt{$(\overline{\nu}_7,\overline{\nu}_8)$}
}
&&
\overline{N}_{00}\{-1\},
}
\end{equation}
where
\begin{align}
\nonumber&\overline{\nu}_5=\id_{\overline{S}^{[1,k]}_{(1,5,2,6)}}\boxtimes (x_{1,1}-x_{1,6},1)\boxtimes \id_{\overline{M}^{[1,k]}_{(6,4,5,3)}},&&
\overline{\nu}_6=\id_{\overline{N}^{[1,k]}_{(1,5,2,6)}}\boxtimes\id_{\overline{S}^{[1,k]}_{(6,4,5,3)}}\boxtimes (1,x_{1,6}-x_{1,3}),\\
\nonumber&\overline{\nu}_7=\id_{\overline{M}^{[1,k]}_{(1,5,2,6)}}\boxtimes \id_{\overline{S}^{[1,k]}_{(6,4,5,3)}}\boxtimes(1,x_{1,6}-x_{1,3}),&&
\overline{\nu}_8=-\id_{\overline{S}^{[1,k]}_{(1,5,2,6)}}\boxtimes (x_{1,1}-x_{1,6},1)\boxtimes\id_{\overline{N}^{[1,k]}_{(6,4,5,3)}}.
\end{align}
We have isomorphisms (\ref{subcom-rev-1}), (\ref{subcom-rev-22}), (\ref{subcom-rev-3}) and (\ref{subcom-rev-4}),
\begin{eqnarray}
\nonumber
\overline{N}_{00}\{-1\}&\simeq&
\overline{S}^{[1,k]}_{(1,4,2,3)}\boxtimes
K\left(u_{k+1,(1,4,2,3)}^{[1,k]}(x_{1,1}-x_{1,3});X_{k,(3,4)}^{(-1,k)}\right)_{Q_2\left/\langle u_{k+1,(1,5,2,6)}^{[1,k]}\rangle\right.}\{1-n\}\left<1\right>,\\
\nonumber
\overline{N}_{10}&\simeq&
\overline{S}^{[1,k]}_{(1,4,2,3)}\boxtimes
K\left(-c_3(x_{1,1}-x_{1,3})(x_{1,6}^{n-k-1}+\beta);(x_{1,6}-x_{1,3})X_{k,(3,4)}^{(-1,k)}\right)_{Q_2\left/\langle u_{k+1,(1,5,2,6)}^{[1,k]}\rangle\right.}\{1-n\}\left<1\right>,\\
\nonumber
\overline{N}_{01}
&\simeq&
\overline{S}^{[1,k]}_{(1,4,2,3)}\boxtimes
K\left(u_{k+1,(1,4,2,3)}^{[1,k]}(x_{1,1}-x_{1,3});X_{k,(3,4)}^{(-1,k)}\right)_{Q_2\left/\langle u_{k+1,(1,5,2,6)}^{[1,k]}(x_{1,1}-x_{1,6})\rangle\right.}\{1-n\}\left<1\right>,\\
\nonumber
\overline{N}_{11}\{1\}&\simeq&
\overline{S}^{[1,k]}_{(1,4,2,3)}\boxtimes
K\left(u_{k+1,(6,4,5,3)}^{[1,k]}+\alpha;(x_{1,6}-x_{1,3})X_{k,(3,4)}^{(-1,k)}\right)_{Q_2\left/\langle u_{k+1,(1,5,2,6)}^{[1,k]}(x_{1,1}-x_{1,6})\rangle\right.}\{1-n\}\left<1\right>,
\end{eqnarray}
where $Q_2=R_{(1,2,3,4,5,6)}^{(1,k,1,k,k,1)}\left/\left<X_{1,(1,5)}^{(1,k)}-X_{1,(2,6)}^{(k,1)},\ldots,X_{k,(1,5)}^{(1,k)}-X_{k,(2,6)}^{(k,1)}\right>\right.$.\\
\indent
We consider $R_6$ and $R_7$, which are isomorphisms of $Q_2\left/\langle u_{k+1,(1,5,2,6)}^{[1,k]}\rangle\right.$, and $R_8$, which is an isomorphism of $Q_2\left/\langle u_{k+1,(1,5,2,6)}^{[1,k]}(x_{1,1}-x_{1,6})\rangle\right.$.
Moreover, we consider the following isomorphism of $Q_2\left/\langle u_{k+1,(1,5,2,6)}^{[1,k]}(x_{1,1}-x_{1,6})\rangle\right.$:
\begin{equation}
\nonumber
R_{10}=R_{(1,2,3,4)}^{(1,k,1,k)}\oplus (x_{1,6}-x_{1,3})R_{(1,2,3,4)}^{(1,k,1,k)}\oplus x_{1,6}(x_{1,6}-x_{1,3})R_{(1,2,3,4)}^{(1,k,1,k)}\oplus \ldots \oplus x_{1,6}^{n-k-1}(x_{1,6}-x_{1,3})R_{(1,2,3,4)}^{(1,k,1,k)}.
\end{equation}
Then, we have isomorphisms
\begin{eqnarray}
\nonumber
\overline{N}_{00}\{-1\}&\simeq&\overline{S}^{[1,k]}_{(1,4,2,3)}\boxtimes
K\left(u_{k+1,(1,4,2,3)}^{[1,k]}(x_{1,1}-x_{1,3});X_{k,(3,4)}^{(-1,k)}\right)_{Q_2\left/\langle u_{k+1,(1,5,2,6)}^{[1,k]}\rangle\right.}\{1-n\}\left<1\right>\\
\nonumber
&\simeq&\overline{S}^{[1,k]}_{(1,4,2,3)}\boxtimes
\left(\xymatrix{R_7\ar[rrrrr]_{E_{n-k}\left(u_{k+1,(1,4,2,3)}^{[1,k]}(x_{1,1}-x_{1,3})\right)}&&&&&R_7\ar[rrr]_{E_{n-k}\left(X_{k,(3,4)}^{(-1,k)}\right)}&&&R_7}\right)\{1-n\}\left<1\right>\\
\nonumber
&\simeq&\bigoplus_{i=0}^{n-k-1}\overline{N}_{(1,4,2,3)}^{[1,k]}\{2i-n+k\},\\
\nonumber
\overline{N}_{10}&\simeq&\overline{S}^{[1,k]}_{(1,4,2,3)}\boxtimes
K\left(-c_3(x_{1,1}-x_{1,3})(x_{1,6}^{n-k-1}+\beta);(x_{1,6}-x_{1,3})X_{k,(3,4)}^{(-1,k)}\right)_{Q_2\left/\langle u_{k+1,(1,5,2,6)}^{[1,k]}\rangle\right.}\{1-n\}\left<1\right>\\
\nonumber
&\simeq&\overline{S}^{[1,k]}_{(1,4,2,3)}\boxtimes
\left(\xymatrix{R_7\ar[rr]_{h_1}&&R_6\ar[rr]_{h_2}&&R_7}\right)\{1-n\}\left<1\right>,\\
\nonumber
&&h_1=\left(
\begin{array}{cc}
{}^{t}\mathfrak{0}_{n-k-1}&E_{n-k-1}\left(u_{k+1,(1,4,2,3)}^{[1,k]}(x_{1,1}-x_{1,3})\right)\\
x_{1,1}-x_{1,3}&\mathfrak{0}_{n-k-1}
\end{array}
\right),\\
\nonumber
&&h_2=\left(
\begin{array}{cc}
\mathfrak{0}_{n-k-1}&u_{k+1,(1,4,2,3)}^{[1,k]}X_{k,(3,4)}^{(-1,k)}\\
E_{n-k-1}\left(X_{k,(3,4)}^{(-1,k)}\right)&{}^{t}\mathfrak{0}_{n-k-1}
\end{array}
\right),\\
\nonumber
&\simeq&\bigoplus_{i=0}^{n-k-2}\overline{N}_{(1,4,2,3)}^{[1,k]}\{2i-n+k+2\}\oplus\overline{L}_{(1,4,2,3)}^{[1,k]},
\end{eqnarray}
\begin{eqnarray}
\nonumber
\overline{N}_{01}&\simeq&\overline{S}^{[1,k]}_{(1,4,2,3)}\boxtimes
K\left(u_{k+1,(1,4,2,3)}^{[1,k]}(x_{1,1}-x_{1,3});X_{k,(3,4)}^{(-1,k)}\right)_{Q_2\left/\langle u_{k+1,(1,5,2,6)}^{[1,k]}(x_{1,1}-x_{1,6})\rangle\right.}\{1-n\}\left<1\right>\\
\nonumber
&\simeq&\overline{S}^{[1,k]}_{(1,4,2,3)}\boxtimes
\left(\xymatrix{R_{10}\ar[rrrrr]_{E_{n-k+1}\left(u_{k+1,(1,4,2,3)}^{[1,k]}(x_{1,1}-x_{1,3})\right)}&&&&&R_{10}\ar[rrr]_{E_{n-k+1}\left(X_{k,(3,4)}^{(-1,k)}\right)}&&&R_{10}}\right)\{1-n\}\left<1\right>\\
\nonumber
&\simeq&\bigoplus_{i=0}^{n-k}\overline{N}_{(1,4,2,3)}^{[1,k]}\{2i-n+k\},\\ 
\nonumber
\overline{N}_{11}\{1\}&\simeq&\overline{S}^{[1,k]}_{(1,4,2,3)}\boxtimes
K\left(u_{k+1,(6,4,5,3)}^{[1,k]};(x_{1,6}-x_{1,3})X_{k,(3,4)}^{(-1,k)}\right)_{Q_2\left/\langle u_{k+1,(1,5,2,6)}^{[1,k]}(x_{1,1}-x_{1,6})\rangle\right.}\{1-n\}\left<1\right>\\
\nonumber
&\simeq&\overline{S}^{[1,k]}_{(1,4,2,3)}\boxtimes
\left(\xymatrix{R_{10}\ar[rr]_{h_3}&&R_8\ar[rr]_{h_4}&&R_{10}}\right)\{1-n\}\left<1\right>,\\
\nonumber
&&h_3=\left(
\begin{array}{cc}
{}^{t}\mathfrak{0}_{n-k}&E_{n-k}\left(u_{k+1,(1,4,2,3)}^{[1,k]}(x_{1,1}-x_{1,3})\right)\\
1&\mathfrak{0}_{n-k}
\end{array}
\right),\\
\nonumber
&&h_4=\left(
\begin{array}{cc}
\mathfrak{0}_{n-k}&u_{k+1,(1,4,2,3)}^{[1,k]}(x_{1,1}-x_{1,3})X_{k,(3,4)}^{(-1,k)}\\
E_{n-k}\left(X_{k,(3,4)}^{(-1,k)}\right)&{}^{t}\mathfrak{0}_{n-k}
\end{array}
\right),\\
\nonumber
&\simeq&\bigoplus_{i=0}^{n-k-1}\overline{N}_{(1,4,2,3)}^{[1,k]}\{2i-n+k+2\}.
\end{eqnarray}
For these decompositions, the morphisms $\overline{\nu}_5$, $\overline{\nu}_6$, $\overline{\nu}_7$ and $\overline{\nu}_8$ form into
\begin{eqnarray}
\nonumber
\overline{\nu}_5&\simeq&
\left(
\begin{array}{cc}
&a_1\id_{\overline{N}_{(1,4,2,3)}^{[1,k]}}\\
E_{n-k-1}\left(\id_{\overline{N}_{(1,4,2,3)}^{[1,k]}}\right)&\vdots\\
&a_{n-k-1}\id_{\overline{N}_{(1,4,2,3)}^{[1,k]}}\\
\mathfrak{0}_{n-k-1}&(a_0 u_{k+1,(1,4,2,3)}^{[1,k]},a_0)
\end{array}
\right)
,\\
\nonumber
\overline{\nu}_6&\simeq&
\left(
\begin{array}{c}
\mathfrak{0}_{n-k}\\
E_{n-k}\left(\id_{\overline{N}_{(1,4,2,3)}^{[1,k]}}\right)
\end{array}
\right),
\end{eqnarray}
\begin{eqnarray}
\nonumber
\overline{\nu}_7&\simeq&
\left(
\begin{array}{cc}
\mathfrak{0}_{n-k-1}&(1,u_{k+1,(1,4,2,3)}^{[1,k]})\\
E_{n-k-1}\left(\id_{\overline{N}_{(1,4,2,3)}^{[1,k]}}\right)&{}^{t}\mathfrak{0}_{n-k-1}
\end{array}
\right),
\\
\nonumber
\overline{\nu}_8&\simeq&
-\left(
\begin{array}{cc}
&a_0 u_{k+1,(1,4,2,3)}^{[1,k]} \id_{\overline{N}_{(1,4,2,3)}^{[1,k]}}\\
E_{n-k}\left(\id_{\overline{N}_{(1,4,2,3)}^{[1,k]}}\right)&\vdots\\
&a_{n-k-1}\id_{\overline{N}_{(1,4,2,3)}^{[1,k]}}
\end{array}
\right),
\end{eqnarray}
where $a_0$, $a_1$, $\ldots$, $a_{n-k-1}$ are polynomials derived from $x_{1,6}^{n-k-1}$ expanded by the basis for the isomorphism $R_3$ of $Q_2\left/\langle u_{k+1,(1,5,2,6)}^{[1,k]}\rangle\right.$, $\left<1,x_{1,6},\ldots,x_{1,6}^{n-k-2},-c_3(x_{1,6}^{n-k-1}+\beta)\right>$,
\begin{equation}
\nonumber
x_{1,6}^{n-k-1}= a_0 \left(-c_3 (x_{1,6}^{n-k-1}+\beta)\right) +a_1x_{1,6}^{n-k-2}+\ldots+a_{n-k-1}. 
\end{equation}
Thus, the complex (\ref{com-rev-2}) is isomorphic, in $\k^b(\HMF^{gr}_{R_{(1,2,3,4)}^{(1,k,1,k)},\omega_2})$, to $\overline{L}^{[1,k]}_{(1,4,2,3)}$:
\begin{equation}
\nonumber
\c\left(\input{figure/r2-1krrev-mf}\right)_n\simeq\c\left(\input{figure/r2-1kcrev-mf}\right)_n.
\end{equation}
It is obvious that we can similarly prove the following isomorphisms for the Reidemeister moves $(IIb_{1k})$:
\begin{equation}
\nonumber
\c\left(\input{figure/r2-1klrevori-mf}\right)_n\simeq\c\left(\input{figure/r2-1kcrevori-mf}\right)_n\simeq\c\left(\input{figure/r2-1krrevori-mf}\right)_n.
\end{equation}
%
%
%
%
\subsection{Proof of Proposition \ref{prop-r3}}\label{prop-r3-proof}
\begin{proof}[{\bf Proof of Proposition \ref{prop-r3} (1)}]
The complex for the diagram $\input{figure/r3-lem-plus-1k-1-text}$ is described as a complex of factorizations of $\k^b(\HMF^{gr}_{R_{(1,2,3,4,5)}^{(k+1,1,1,1,k)},\omega_3})$ $(\omega_3=F_{k+1}(\mathbb{X}^{(k+1)}_{(1)})+F_{1}(\mathbb{X}^{(1)}_{(2)})-F_{1}(\mathbb{X}^{(1)}_{(3)})-F_{1}(\mathbb{X}^{(1)}_{(4)})-F_{k}(\mathbb{X}^{(k)}_{(5)}))$,
\begin{eqnarray}
\label{r3-lem-planar-comp-1}
&&\hspace{1cm}
\c\left(\input{figure/r3-lem-plus-1k-1-mf}\right)_n=\\
\nonumber
&&
\xymatrix{
-k-1\ar@{.}[d]&-k\ar@{.}[d]&-k+1\ar@{.}[d]\\
{\c\left(\input{figure/r3-lem-planar-1k-1-mf}\right)_n
\hspace{-0.5cm}
\begin{array}{c}
{}_{\{(k+1)n\}}\\
\left<k+1\right>
\end{array}}
\ar[r]^{\left(\hspace{-0.2cm}\begin{array}{c}
\overline{\zeta}_{+,1}
\\
\overline{\zeta}_{+,2}
\end{array}\hspace{-0.2cm}\right)}&
{
\begin{array}{c}
\c\left(\input{figure/r3-lem-planar-1k-2-mf}\right)_n
\hspace{-0.5cm}
\begin{array}{c}
{}_{\{(k+1)n-1\}}\\
\left<k+1\right>
\end{array}
\\
\bigoplus\\
\c\left(\input{figure/r3-lem-planar-1k-3-mf}\right)_n
\hspace{-0.5cm}
\begin{array}{c}
{}_{\{(k+1)n-1\}}\\
\left<k+1\right>
\end{array}
\end{array}
}
\ar[r]^(.45){\txt{$(\overline{\zeta}_{+,3},\overline{\zeta}_{+,4})$}}&
{\c\left(\input{figure/r3-lem-planar-1k-4-mf}\right)_n
\hspace{-0.5cm}
\begin{array}{c}
{}_{\{(k+1)n-2\}}\\
\left<k+1\right>
\end{array}
},
}
\end{eqnarray}
where
\begin{eqnarray}
\nonumber
\overline{\zeta}_{+,1}&=&\id_{\overline{\Lambda}_{(1;6,7)}^{[1,k]}}\boxtimes\id_{\overline{S}_{(8,6,4,3)}^{[1,1]}}\boxtimes(1,x_{1,8}-x_{1,3})\boxtimes\id_{\overline{M}_{(2,7,5,8)}^{[1,k]}},
\\
\nonumber
\overline{\zeta}_{+,2}&=&\id_{\overline{\Lambda}_{(1;6,7)}^{[1,k]}}\boxtimes\id_{\overline{M}_{(8,6,4,3)}^{[1,1]}}\boxtimes\id_{\overline{S}_{(2,7,5,8)}^{[1,k]}}\boxtimes(1,x_{1,2}-x_{1,8}),\\
\nonumber
\overline{\zeta}_{+,3}&=&\id_{\overline{\Lambda}_{(1;6,7)}^{[1,k]}}\boxtimes\id_{\overline{N}_{(8,6,4,3)}^{[1,1]}}\boxtimes\id_{\overline{S}_{(2,7,5,8)}^{[1,k]}}\boxtimes(1,x_{1,2}-x_{1,8}),
\\
\nonumber
\overline{\zeta}_{+,4}&=&-\id_{\overline{\Lambda}_{(1;6,7)}^{[1,k]}}\boxtimes\id_{\overline{S}_{(8,6,4,3)}^{[1,1]}}\boxtimes(1,x_{1,8}-x_{1,3})\boxtimes\id_{\overline{N}_{(2,7,5,8)}^{[1,k]}}.
\end{eqnarray}
First, we have
\begin{eqnarray}
\nonumber
\c\left(\input{figure/r3-lem-planar-1k-1-mf}\right)_n
&=&\overline{\Lambda}_{(1;6,7)}^{[1,k]}\boxtimes \overline{M}_{(8,6,4,3)}^{[1,1]} \boxtimes \overline{M}_{(2,7,5,8)}^{[1,k]}\\
\nonumber
&\simeq&
K\left(
\left(
\begin{array}{c}
\Lambda_{1,(1;6,7)}^{[1,k]}\\
\vspace{0.1cm}\vdots\\
\Lambda_{k+1,(1;6,7)}^{[1,k]}
\end{array}
\right);
\left(
\begin{array}{c}
x_{1,1}-X_{1,(6,7)}^{(1,k)}\\
\vspace{0.1cm}\vdots\\
x_{k+1,1}-X_{k+1,(6,7)}^{(1,k)}
\end{array}
\right)
\right)_{R_{(1,6,7)}^{(1,1,k)}}\\
&&
\nonumber
\hspace{0.5cm}\boxtimes
K\left(
\left(
\begin{array}{c}
A_{1,(8,6,4,3)}^{[1,1]}\\
u_{2,(8,6,4,3)}^{[1,1]}
\end{array}
\right);
\left(
\begin{array}{c}
X_{1,(8,6)}^{(1,1)}-X_{1,(4,3)}^{(1,1)}\\
(x_{1,8}-x_{1,3})X_{1,(6,3)}^{(1,-1)}
\end{array}
\right)
\right)_{R_{(8,6,4,3)}^{(1,1,1,1)}}\{-1\}\\
&&
\nonumber
\hspace{0.5cm}\boxtimes
K\left(
\left(
\begin{array}{c}
\vspace{0.2cm}A_{1,(2,7,5,8)}^{[1,k]}\\
\vspace{0.2cm}\vdots\\
\vspace{0.2cm}A_{k,(2,7,5,8)}^{[1,k]}\\
u_{k+1,(2,7,5,8)}^{[1,k]}
\end{array}
\right);
\left(
\begin{array}{c}
\vspace{0.2cm}X_{1,(2,7)}^{(1,k)}-X_{1,(5,8)}^{(k,1)}\\
\vspace{0.2cm}\vdots\\
\vspace{0.2cm}X_{k,(2,7)}^{(1,k)}-X_{k,(5,8)}^{(k,1)}\\
(x_{1,2}-x_{1,8})X_{k,(7,8)}^{(k,-1)}
\end{array}
\right)
\right)_{R_{(2,7,5,8)}^{(1,k,k,1)}}\{-k\}\\
&\simeq&
\label{mf-r3-1-00}
K\left(
\left(
\begin{array}{c}
\vspace{0.1cm}\Lambda_{1,(1;6,7)}^{[1,k]}\\
\vspace{0.1cm}\vdots\\
\vspace{0.1cm}\Lambda_{k+1,(1;6,7)}^{[1,k]}\\
u_{k+1,(2,7,5,8)}^{[1,k]}
\end{array}
\right);
\left(
\begin{array}{c}
\vspace{0.1cm}x_{1,1}-X_{1,(6,7)}^{(1,k)}\\
\vspace{0.1cm}\vdots\\
\vspace{0.1cm}x_{k+1,1}-X_{k+1,(6,7)}^{(1,k)}\\
(x_{1,2}-x_{1,8})X_{k,(7,8)}^{(k,-1)}
\end{array}
\right)
\right)_{Q_3\left/\langle (x_{1,8}-x_{1,3})X_{1,(6,3)}^{(1,-1)}\rangle\right.}\hspace{-2cm}\{-k-1\},
\end{eqnarray}
where 

\begin{eqnarray}
\nonumber
Q_3:=R_{(1,2,3,4,5,6,7,8)}^{(k+1,1,1,1,k,1,k,1)}\left/\left<X_{1,(8,6)}^{(1,1)}-X_{1,(4,3)}^{(1,1)},X_{1,(2,7)}^{(1,k)}-X_{1,(5,8)}^{(k,1)},\ldots,X_{k,(2,7)}^{(1,k)}-X_{k,(5,8)}^{(k,1)}\right>\right. .
\end{eqnarray}
The quotient $Q_3\left/\langle (x_{1,8}-x_{1,3})X_{1,(6,3)}^{(1,-1)}\rangle\right.$ has equations
\begin{eqnarray}
\nonumber
&&x_{1,6}=X_{1,(3,4,8)}^{(1,1,-1)},\\
\nonumber
&&(x_{1,8}-x_{1,3})(x_{1,8}-x_{1,4})=0,\\
\nonumber
&&x_{j,7}=X_{j,(2,5,8)}^{(-1,k,1)} \hspace{1cm}(1 \leq j \leq k).
\end{eqnarray}
In the quotient $Q_3\left/\langle (x_{1,8}-x_{1,3})X_{1,(6,3)}^{(1,-1)}\rangle\right.$, $X_{j,(6,7)}^{(1,k)}$ $(1 \leq j\leq k)$ equals to
\begin{eqnarray}
\nonumber x_{j,7}+x_{1,6}x_{j-1,7}&\equiv&X_{j,(2,5,8)}^{(-1,k,1)}+X_{1,(3,4,8)}^{(1,1,-1)}X_{j-1,(2,5,8)}^{(-1,k,1)}\\
\nonumber&\equiv&X_{j,(2,5)}^{(-1,k)}+X_{1,(3,4)}^{(1,1)}X_{j-1,(2,5)}^{(-1,k)}
+X_{2,(3,4)}^{(1,1)}X_{j-2,(2,5)}^{(-1,k)}=X_{j,(2,3,4,5)}^{(-1,1,1,k)},
\end{eqnarray}
and $X_{k+1,(6,7)}^{(1,k)}$ equals to
\begin{eqnarray}
\nonumber x_{1,6}x_{k,7}&\equiv&X_{1,(3,4,8)}^{(1,1,-1)}X_{k,(2,5,8)}^{(-1,k,1)}\\
\nonumber&\equiv&
X_{k+1,(2,3,4,5)}^{(-1,1,1,k)}+(x_{1,2}-x_{1,8})X_{k,(2,5)}^{(-1,k)}.
\end{eqnarray}
Then, the matrix factorization (\ref{mf-r3-1-00}) is isomorphic to 
\begin{eqnarray}
\label{mf-r3-2-00}
&&K\left(
\left(
\begin{array}{c}
\vspace{0.1cm}\Lambda_{1,(1;6,7)}^{[1,k]}\\
\vspace{0.1cm}\vdots\\
\vspace{0.1cm}\Lambda_{k+1,(1;6,7)}^{[1,k]}\\
u_{k+1,(2,7,5,8)}^{[1,k]}-\Lambda_{k+1,(1;6,7)}^{[1,k]}
\end{array}
\right);
\left(
\begin{array}{c}
\vspace{0.1cm}x_{1,1}-X_{1,(2,3,4,5)}^{(-1,1,1,k)}\\
\vspace{0.1cm}\vdots\\
\vspace{0.1cm}x_{k+1,1}-X_{k+1,(2,3,4,5)}^{(-1,1,1,k)}\\
(x_{1,2}-x_{1,8})X_{k,(2,5)}^{(-1,k)}
\end{array}
\right)
\right)_{Q_3\left/\langle (x_{1,8}-x_{1,3})X_{1,(6,3)}^{(1,-1)}\rangle\right.}\{-k-1\}.
\end{eqnarray}
By Corollary \ref{induce-sq1}, there exist polynomials $B_1$, $B_2$, $\ldots$, $B_{k+1}$ $\in R_{(1,2,3,4,5)}^{(k+1,1,1,1,k)}$ and $B_0\in Q_3\left/\langle (x_{1,8}-x_{1,3})X_{1,(6,3)}^{(1,-1)}\rangle\right.$ satisfying $(x_{1,2}-x_{1,8})B_0\equiv B\in R_{(1,2,3,4,5)}^{(k+1,1,1,1,k)}$ $(\mod Q_3\left/\langle (x_{1,8}-x_{1,3})X_{1,(6,3)}^{(1,-1)}\rangle\right.)$ and we have an isomorphism between the factorization (\ref{mf-r3-2-00}) and the following factorization
\begin{eqnarray}
\label{mf-r3-3-00}
\overline{S}_{(1,2;3,4,5)}^{[k+1,1;1,1,k]}\boxtimes
K\left(B_0;(x_{1,2}-x_{1,8})X_{k,(2,5)}^{(-1,k)}\right)_{Q_3\left/\langle (x_{1,8}-x_{1,3})X_{1,(6,3)}^{(1,-1)}\rangle\right.}\{-k-1\},
\end{eqnarray}
where
\begin{equation}
\label{common-mf}
\overline{S}_{(1,2;3,4,5)}^{[k+1,1;1,1,k]}
:=K\left(
\left(
\begin{array}{c}
\vspace{0.1cm}B_1\\
\vspace{0.1cm}\vdots\\
\vspace{0.1cm}B_{k+1}
\end{array}
\right);
\left(
\begin{array}{c}
x_{1,1}-X_{1,(2,3,4,5)}^{(-1,1,1,k)}\\
\vdots\\
x_{k+1,1}-X_{k+1,(2,3,4,5)}^{(-1,1,1,k)}
\end{array}
\right)
\right)_{R_{(1,2,3,4,5)}^{(k+1,1,1,1,k)}}.
\end{equation}
We consider isomorphisms of $Q_3\left/\langle (x_{1,8}-x_{1,3})X_{1,(6,3)}^{(1,-1)}\rangle\right.$ to be
\begin{eqnarray}
\nonumber
R_{11}&\simeq&R_{(1,2,3,4,5)}^{(k+1,1,1,1,k)}\oplus(x_{1,2}-x_{1,8})R_{(1,2,3,4,5)}^{(k+1,1,1,1,k)},\\
\nonumber
R_{12}&\simeq&R_{(1,2,3,4,5)}^{(k+1,1,1,1,k)}\oplus(x_{1,8}+x_{1,2}-x_{1,3}-x_{1,4})R_{(1,2,3,4,5)}^{(k+1,1,1,1,k)}.
\end{eqnarray}
Then, the partial matrix factorization $K(B_0;(x_{1,2}-x_{1,8})X_{k,(2,5)}^{(-1,k)})_{Q_3\left/\langle (x_{1,8}-x_{1,3})X_{1,(6,3)}^{(1,-1)}\rangle\right.}$ is isomorphic to
\begin{eqnarray}
\nonumber
&&\xymatrix{R_{11}
\ar[rrrrr]_(.42){
\left(
\begin{array}{cc}
0&B\\
\frac{B}{(x_{1,2}-x_{1,3})(x_{1,2}-x_{1,4})}&0
\end{array}
\right)
}
&&&&&R_{12}\{2k-n+1\}
\ar[rrrrrrr]_(.55){
\left(
\begin{array}{cc}
0&(x_{1,2}-x_{1,3})(x_{1,2}-x_{1,4})X_{k,(2,5)}^{(-1,k)}\\
X_{k,(2,5)}^{(-1,k)}&0
\end{array}
\right)
}
&&&&&&&R_{11}}\\
\nonumber
&\simeq&K\left(
\frac{B}{(x_{1,2}-x_{1,3})(x_{1,2}-x_{1,4})};(x_{1,2}-x_{1,3})(x_{1,2}-x_{1,4})X_{k,(2,5)}^{(-1,k)}
\right)_{R_{(1,2,3,4,5)}^{(k+1,1,1,1,k)}}\oplus 
K\left(
B;X_{k,(2,5)}^{(-1,k)}
\right)_{R_{(1,2,3,4,5)}^{(k+1,1,1,1,k)}}\{2\}.
\end{eqnarray}
Thus, the matrix factorization (\ref{mf-r3-3-00}) is decomposed into
\begin{eqnarray}
\label{mf-r3-fin-00}
&&\overline{S}_{(1,2;3,4,5)}^{[k+1,1;1,1,k]}\boxtimes
K\left(
\frac{B}{(x_{1,2}-x_{1,3})(x_{1,2}-x_{1,4})};(x_{1,2}-x_{1,3})(x_{1,2}-x_{1,4})X_{k,(2,5)}^{(-1,k)}
\right)_{R_{(1,2,3,4,5)}^{(k+1,1,1,1,k)}}\{-k-1\}\\
\nonumber
&&\hspace{1cm}\oplus 
\overline{S}_{(1,2;3,4,5)}^{[k+1,1;1,1,k]}\boxtimes
K\left(
B;X_{k,(2,5)}^{(-1,k)}
\right)_{R_{(1,2,3,4,5)}^{(k+1,1,1,1,k)}}\{-k+1\}.
\end{eqnarray}
Secondly, we have
\begin{eqnarray}
\nonumber
\c\left(\input{figure/r3-lem-planar-1k-2-mf}\right)_n
&=&\overline{\Lambda}_{(1;6,7)}^{[1,k]}\boxtimes \overline{N}_{(8,6,4,3)}^{[1,1]} \boxtimes \overline{M}_{(2,7,5,8)}^{[1,k]}\\
\nonumber
&\simeq&
K\left(
\left(
\begin{array}{c}
\Lambda_{1,(1;6,7)}^{[1,k]}\\
\vspace{0.1cm}\vdots\\
\Lambda_{k+1,(1;6,7)}^{[1,k]}
\end{array}
\right);
\left(
\begin{array}{c}
x_{1,1}-X_{1,(6,7)}^{(1,k)}\\
\vspace{0.1cm}\vdots\\
x_{k+1,1}-X_{k+1,(6,7)}^{(1,k)}
\end{array}
\right)
\right)_{R_{(1,6,7)}^{(1,1,k)}}\\
&&
\nonumber
\hspace{0.5cm}\boxtimes
K\left(
\left(
\begin{array}{c}
A_{1,(8,6,4,3)}^{[1,1]}\\
u_{2,(8,6,4,3)}^{[1,1]}(x_{1,8}-x_{1,3})
\end{array}
\right);
\left(
\begin{array}{c}
X_{1,(8,6)}^{(1,1)}-X_{1,(4,3)}^{(1,1)}\\
X_{1,(6,3)}^{(1,-1)}
\end{array}
\right)
\right)_{R_{(8,6,4,3)}^{(1,1,1,1)}}\\
&&
\nonumber
\hspace{0.5cm}\boxtimes
K\left(
\left(
\begin{array}{c}
\vspace{0.2cm}A_{1,(2,7,5,8)}^{[1,k]}\\
\vspace{0.2cm}\vdots\\
\vspace{0.2cm}A_{k,(2,7,5,8)}^{[1,k]}\\
u_{k+1,(2,7,5,8)}^{[1,k]}
\end{array}
\right);
\left(
\begin{array}{c}
\vspace{0.2cm}X_{1,(2,7)}^{(1,k)}-X_{1,(5,8)}^{(k,1)}\\
\vspace{0.2cm}\vdots\\
\vspace{0.2cm}X_{k,(2,7)}^{(1,k)}-X_{k,(5,8)}^{(k,1)}\\
(x_{1,2}-x_{1,8})X_{k,(7,8)}^{(k,-1)}
\end{array}
\right)
\right)_{R_{(2,7,5,8)}^{(1,k,k,1)}}\{-k\}
\end{eqnarray}
\begin{eqnarray}
&\simeq&
\label{mf-r3-1-10}
K\left(
\left(
\begin{array}{c}
\vspace{0.1cm}\Lambda_{1,(1;6,7)}^{[1,k]}\\
\vspace{0.1cm}\vdots\\
\vspace{0.1cm}\Lambda_{k+1,(1;6,7)}^{[1,k]}\\
u_{k+1,(2,7,5,8)}^{[1,k]}
\end{array}
\right);
\left(
\begin{array}{c}
\vspace{0.1cm}x_{1,1}-X_{1,(6,7)}^{(1,k)}\\
\vspace{0.1cm}\vdots\\
\vspace{0.1cm}x_{k+1,1}-X_{k+1,(6,7)}^{(1,k)}\\
(x_{1,2}-x_{1,8})X_{k,(7,8)}^{(k,-1)}
\end{array}
\right)
\right)_{Q_3\left/\langle X_{1,(6,3)}^{(1,-1)}\rangle\right.}\{-k\}.
\end{eqnarray}
The quotient $Q_3\left/\langle X_{1,(6,3)}^{(1,-1)}\rangle\right.$ has equations
\begin{eqnarray}
\nonumber
&&x_{1,6}=x_{1,3},\\
\nonumber
&&x_{1,8}=x_{1,4},\\
\nonumber
&&x_{j,7}=X_{j,(2,4,5)}^{(-1,1,k)} \hspace{1cm}(1\leq j\leq k).
\end{eqnarray}
In the quotient $Q_3\left/\langle X_{1,(6,3)}^{(1,-1)}\rangle\right.$, $X_{j,(6,7)}^{(1,k)}$ $(1\leq j\leq k)$ equals to
\begin{eqnarray}
\nonumber
x_{j,7}+x_{1,6}x_{j-1,7}&\equiv& X_{j,(2,4,5)}^{(-1,1,k)}+x_{1,3}X_{j-1,(2,4,5)}^{(-1,1,k)}\\
&=&X_{j,(2,3,4,5)}^{(-1,1,1,k)}
\end{eqnarray}
and $X_{k+1,(6,7)}^{(1,k)}$ equals to
\begin{eqnarray}
\nonumber
x_{k,6}x_{1,7}&\equiv& x_{1,3}X_{k,(2,4,5)}^{(-1,1,k)}\\
\nonumber
&=&X_{k+1,(2,3,4,5)}^{(-1,1,1,k)}-(x_{1,2}-x_{1,4})X_{k,(2,5)}^{(-1,k)}.
\end{eqnarray}
Thus, we find $Q_3\left/\langle X_{1,(6,3)}^{(1,-1)}\rangle\right. \simeq R_{(1,2,3,4,5)}^{(k+1,1,1,1,k)}$. Then, the matrix factorization (\ref{mf-r3-1-10}) is isomorphic to 
\begin{eqnarray}
\nonumber
&&K\left(
\left(
\begin{array}{c}
\Lambda_{1,(1;6,7)}^{[1,k]}\\
\vspace{0.1cm}\vdots\\
\vspace{0.1cm}\Lambda_{k+1,(1;6,7)}^{[1,k]}\\
u_{k+1,(2,7,5,8)}^{[1,k]}-\Lambda_{k+1,(1;6,7)}^{[1,k]}
\end{array}
\right);
\left(
\begin{array}{c}
x_{1,1}-X_{1,(2,3,4,5)}^{(-1,1,1,k)}\\
\vspace{0.1cm}\vdots\\
\vspace{0.1cm}x_{k+1,1}-X_{k+1,(2,3,4,5)}^{(-1,1,1,k)}\\
(x_{1,2}-x_{1,4})X_{k,(2,5)}^{(-1,k)}
\end{array}
\right)
\right)_{Q_3\left/\langle X_{1,(6,3)}^{(1,-1)}\rangle\right.}\{-k\}\\
\label{mf-r3-fin-10}
&\simeq&
\overline{S}_{(1,2;3,4,5)}^{[k+1,1;1,1,k]}\boxtimes
K\left(\frac{B}{(x_{1,2}-x_{1,4})};(x_{1,2}-x_{1,4})X_{k,(2,5)}^{(-1,k)}\right)_{R_{(1,2,3,4,5)}^{(k+1,1,1,1,k)}}\{-k\}.
\end{eqnarray}
Thirdly, we have
\begin{eqnarray}
\nonumber
\c\left(\input{figure/r3-lem-planar-1k-3-mf}\right)_n
&=&\overline{\Lambda}_{(1;6,7)}^{[1,k]}\boxtimes \overline{M}_{(8,6,4,3)}^{[1,1]} \boxtimes \overline{N}_{(2,7,5,8)}^{[1,k]}\\
\nonumber
&\simeq&
K\left(
\left(
\begin{array}{c}
\Lambda_{1,(1;6,7)}^{[1,k]}\\
\vspace{0.1cm}\vdots\\
\Lambda_{k+1,(1;6,7)}^{[1,k]}
\end{array}
\right);
\left(
\begin{array}{c}
x_{1,1}-X_{1,(6,7)}^{(1,k)}\\
\vspace{0.1cm}\vdots\\
x_{k+1,1}-X_{k+1,(6,7)}^{(1,k)}
\end{array}
\right)
\right)_{R_{(1,6,7)}^{(1,1,k)}}\\
&&
\nonumber
\hspace{0.5cm}\boxtimes
K\left(
\left(
\begin{array}{c}
A_{1,(8,6,4,3)}^{[1,1]}\\
u_{2,(8,6,4,3)}^{[1,1]}
\end{array}
\right);
\left(
\begin{array}{c}
X_{1,(8,6)}^{(1,1)}-X_{1,(4,3)}^{(1,1)}\\
(x_{1,8}-x_{1,3})X_{1,(6,3)}^{(1,-1)}
\end{array}
\right)
\right)_{R_{(8,6,4,3)}^{(1,1,1,1)}}\{-1\}\\
&&
\nonumber
\hspace{0.5cm}\boxtimes
K\left(
\left(
\begin{array}{c}
A_{1,(2,7,5,8)}^{[1,k]}\\
\vspace{0.2cm}\vdots\\
\vspace{0.2cm}A_{k,(2,7,5,8)}^{[1,k]}\\
u_{k+1,(2,7,5,8)}^{[1,k]}(x_{1,2}-x_{1,8})
\end{array}
\right);
\left(
\begin{array}{c}
X_{1,(2,7)}^{(1,k)}-X_{1,(5,8)}^{(k,1)}\\
\vspace{0.2cm}\vdots\\
\vspace{0.2cm}X_{k,(2,7)}^{(1,k)}-X_{k,(5,8)}^{(k,1)}\\
X_{k,(7,8)}^{(k,-1)}
\end{array}
\right)
\right)_{R_{(2,7,5,8)}^{(1,k,k,1)}}\{-k+1\}
\end{eqnarray}
\begin{eqnarray}
&\simeq&
\nonumber
K\left(
\left(
\begin{array}{c}
\vspace{0.1cm}\Lambda_{1,(1;6,7)}^{[1,k]}\\
\vspace{0.1cm}\vdots\\
\vspace{0.1cm}\Lambda_{k+1,(1;6,7)}^{[1,k]}\\
u_{k+1,(2,7,5,8)}^{[1,k]}(x_{1,2}-x_{1,8})
\end{array}
\right);
\left(
\begin{array}{c}
\vspace{0.1cm}x_{1,1}-X_{1,(6,7)}^{(1,k)}\\
\vspace{0.1cm}\vdots\\
\vspace{0.1cm}x_{k+1,1}-X_{k+1,(6,7)}^{(1,k)}\\
X_{k,(7,8)}^{(k,-1)}
\end{array}
\right)
\right)_{Q_3\left/\langle (x_{1,8}-x_{1,3})X_{1,(6,3)}^{(1,-1)}\rangle\right.}\{-k\}\\
\label{mf-r3-1-01}
&\simeq&
\overline{S}_{(1,2;3,4,5)}^{[k+1,1;1,1,k]}\boxtimes
K\left(B;X_{k,(2,5)}^{(-1,k)}\right)_{Q_3\left/\langle (x_{1,8}-x_{1,3})X_{1,(6,3)}^{(1,-1)}\rangle\right.}\{-k\}.
\end{eqnarray}
Since the partial matrix factorization $K\left(B;X_{k,(2,5)}^{(-1,k)}\right)_{Q_3\left/\langle (x_{1,8}-x_{1,3})X_{1,(6,3)}^{(1,-1)}\rangle\right.}$ is decomposed into
\begin{eqnarray}
\nonumber
&&\xymatrix{R_{11}\ar[r]_(.4){
\left(\hspace{-0.2cm}
\begin{array}{cc}
B&0\\
0&B
\end{array}\hspace{-0.2cm}
\right)
}
&R_{11}\{2k-n-1\}\ar[rrrr]_(.6){
\left(
\begin{array}{cc}
X_{k,(2,5)}^{(-1,k)}&0\\
0&X_{k,(2,5)}^{(-1,k)}
\end{array}
\right)
}
&&&&R_{11}}\\
\nonumber
&\simeq&
K\left(B;X_{k,(2,5)}^{(-1,k)}\right)_{R_{(1,2,3,4,5)}^{(k+1,1,1,1,k)}}
\oplus K\left(B;X_{k,(2,5)}^{(-1,k)}\right)_{R_{(1,2,3,4,5)}^{(k+1,1,1,1,k)}}\{2\},
\end{eqnarray} 
then the matrix factorization (\ref{mf-r3-1-01}) is isomorphic to
\begin{eqnarray}
\label{mf-r3-fin-01}
&&
\overline{S}_{(1,2;3,4,5)}^{[k+1,1;1,1,k]}\boxtimes
K\left(B;X_{k,(2,5)}^{(-1,k)}\right)_{R_{(1,2,3,4,5)}^{(k+1,1,1,1,k)}}\{-k\}\\
\nonumber
&&\oplus \overline{S}_{(1,2;3,4,5)}^{[k+1,1;1,1,k]}\boxtimes
K\left(B;X_{k,(2,5)}^{(-1,k)}\right)_{R_{(1,2,3,4,5)}^{(k+1,1,1,1,k)}}\{-k+2\}.
\end{eqnarray}
Finally, we have
\begin{eqnarray}
\nonumber
\c\left(\input{figure/r3-lem-planar-1k-4-mf}\right)_n
&=&\overline{\Lambda}_{(1;6,7)}^{[1,k]}\boxtimes \overline{N}_{(8,6,4,3)}^{[1,1]} \boxtimes \overline{N}_{(2,7,5,8)}^{[1,k]}\\
\nonumber
&\simeq&
K\left(
\left(
\begin{array}{c}
\Lambda_{1,(1;6,7)}^{[1,k]}\\
\vspace{0.1cm}\vdots\\
\Lambda_{k+1,(1;6,7)}^{[1,k]}
\end{array}
\right);
\left(
\begin{array}{c}
x_{1,1}-X_{1,(6,7)}^{(1,k)}\\
\vspace{0.1cm}\vdots\\
x_{k+1,1}-X_{k+1,(6,7)}^{(1,k)}
\end{array}
\right)
\right)_{R_{(1,6,7)}^{(1,1,k)}}\\
&&
\nonumber
\hspace{0.5cm}\boxtimes
K\left(
\left(
\begin{array}{c}
A_{1,(8,6,4,3)}^{[1,1]}\\
u_{2,(8,6,4,3)}^{[1,1]}(x_{1,8}-x_{1,3})
\end{array}
\right);
\left(
\begin{array}{c}
X_{1,(8,6)}^{(1,1)}-X_{1,(4,3)}^{(1,1)}\\
X_{1,(6,3)}^{(1,-1)}
\end{array}
\right)
\right)_{R_{(8,6,4,3)}^{(1,1,1,1)}}\\
&&
\nonumber
\hspace{0.5cm}\boxtimes
K\left(
\left(
\begin{array}{c}
A_{1,(2,7,5,8)}^{[1,k]}\\
\vspace{0.2cm}\vdots\\
\vspace{0.2cm}A_{k,(2,7,5,8)}^{[1,k]}\\
u_{k+1,(2,7,5,8)}^{[1,k]}(x_{1,2}-x_{1,8})
\end{array}
\right);
\left(
\begin{array}{c}
X_{1,(2,7)}^{(1,k)}-X_{1,(5,8)}^{(k,1)}\\
\vspace{0.2cm}\vdots\\
\vspace{0.2cm}X_{k,(2,7)}^{(1,k)}-X_{k,(5,8)}^{(k,1)}\\
X_{k,(7,8)}^{(k,-1)}
\end{array}
\right)
\right)_{R_{(2,7,5,8)}^{(1,k,k,1)}}\{-k+1\}\\
&\simeq&
\nonumber
K\left(
\left(
\begin{array}{c}
\vspace{0.1cm}\Lambda_{1,(1;6,7)}^{[1,k]}\\
\vspace{0.1cm}\vdots\\
\vspace{0.1cm}\Lambda_{k+1,(1;6,7)}^{[1,k]}\\
u_{k+1,(2,7,5,8)}^{[1,k]}(x_{1,2}-x_{1,8})
\end{array}
\right);
\left(
\begin{array}{c}
\vspace{0.1cm}x_{1,1}-X_{1,(6,7)}^{(1,k)}\\
\vspace{0.1cm}\vdots\\
\vspace{0.1cm}x_{k+1,1}-X_{k+1,(6,7)}^{(1,k)}\\
X_{k,(7,8)}^{(k,-1)}
\end{array}
\right)
\right)_{Q_3\left/\langle X_{1,(6,3)}^{(1,-1)}\rangle\right.}\{-k+1\}\\
\label{mf-r3-fin-11}
&\simeq&
\overline{S}_{(1,2;3,4,5)}^{[k+1,1;1,1,k]}\boxtimes
K\left(B;X_{k,(2,5)}^{(-1,k)}\right)_{R_{(1,2,3,4,5)}^{(k+1,1,1,1,k)}}\{-k+1\}.
\end{eqnarray}
For these decompositions (\ref{mf-r3-fin-00}), (\ref{mf-r3-fin-10}), (\ref{mf-r3-fin-01}) and (\ref{mf-r3-fin-11}), the morphisms $\overline{\zeta}_{+,1}$, $\overline{\zeta}_{+,2}$, $\overline{\zeta}_{+,3}$ and $\overline{\zeta}_{+,4}$ of the complex (\ref{r3-lem-planar-comp-1}) transform into
\begin{eqnarray}
\nonumber
&&
\overline{\zeta}_{+,1}\simeq
\left(\id_{\overline{S}_{(1,2;3,4,5)}^{[k+1,1;1,1,k]}}\boxtimes(1,x_{1,2}-x_{1,3}),
\id_{\overline{S}_{(1,2;3,4,5)}^{[k+1,1;1,1,k]}}\boxtimes(x_{1,2}-x_{1,4},1)\right),\\
\nonumber
&&
\overline{\zeta}_{+,2}\simeq
\left(
\begin{array}{cc}
\id_{\overline{S}_{(1,2;3,4,5)}^{[k+1,1;1,1,k]}}\boxtimes(1,(x_{1,2}-x_{1,3})(x_{1,2}-x_{1,4}))&0\\
0&\id_{\overline{S}_{(1,2;3,4,5)}^{[k+1,1;1,1,k]}}\boxtimes(1,1)
\end{array}
\right)
,\\
\nonumber
&&
\overline{\zeta}_{+,3}\simeq
\id_{\overline{S}_{(1,2;3,4,5)}^{[k+1,1;1,1,k]}}\boxtimes\left(1,x_{1,2}-x_{1,4}\right),\\
\nonumber
&&
\overline{\zeta}_{+,4}\simeq
-\left(\id_{\overline{S}_{(1,2;3,4,5)}^{[k+1,1;1,1,k]}}\boxtimes(1,1),
\id_{\overline{S}_{(1,2;3,4,5)}^{[k+1,1;1,1,k]}}\boxtimes(x_{1,2}-x_{1,4},x_{1,2}-x_{1,4})\right).
\end{eqnarray}
Then, the complex (\ref{r3-lem-planar-comp-1}) is isomorphic, in $\k^b(\HMF^{gr}_{R_{(1,2,3,4,5)}^{(k+1,1,1,1,k)},\omega_3})$, to
\begin{equation}
\nonumber
\xymatrix{
-k-1\ar@{.}[d]&&&-k\ar@{.}[d]\\
{\overline{M}_1
\{(k+1)n\}
\left<k+1\right>
}
\ar[rrr]^{\id_{\overline{S}}\boxtimes(1,x_{1,2}-x_{1,3})}&&&
{\overline{M}_2
\{(k+1)n-1\}
\left<k+1\right>
}
,}
\end{equation}
where
\begin{eqnarray}
\label{mf-planar-1}
&&\overline{M}_1=
K\left(
\left(
\begin{array}{c}
B_1\\
\vspace{0.2cm}\vdots\\
\vspace{0.2cm}B_{k+1}\\
\frac{B}{(x_{1,2}-x_{1,3})(x_{1,2}-x_{1,4})}
\end{array}
\right);
\left(
\begin{array}{c}
x_{1,1}-X_{1,(2,3,4,5)}^{(-1,1,1,k)}\\
\vdots\\
x_{k+1,1}-X_{k+1,(2,3,4,5)}^{(-1,1,1,k)}\\
(x_{1,2}-x_{1,3})(x_{1,2}-x_{1,4})X_{k,(2,5)}^{(-1,k)}
\end{array}
\right)
\right)_{R_{(1,2,3,4,5)}^{(k+1,1,1,1,k)}}\hspace{-1cm}\{-k-1\},\\
\label{mf-planar-2}
&&\overline{M}_2=
K\left(
\left(
\begin{array}{c}
B_1\\
\vspace{0.2cm}\vdots\\
\vspace{0.2cm}B_{k+1}\\
\frac{B}{(x_{1,2}-x_{1,4})}
\end{array}
\right);
\left(
\begin{array}{c}
x_{1,1}-X_{1,(2,3,4,5)}^{(-1,1,1,k)}\\
\vdots\\
x_{k+1,1}-X_{k+1,(2,3,4,5)}^{(-1,1,1,k)}\\
(x_{1,2}-x_{1,4})X_{k,(2,5)}^{(-1,k)}
\end{array}
\right)
\right)_{R_{(1,2,3,4,5)}^{(k+1,1,1,1,k)}}\hspace{-1cm}\{-k\}.
\end{eqnarray}
By the way, we have
\begin{eqnarray}
\nonumber
&&\c\left(\input{figure/r3-lem-plus-1k-2-mf}\right)_n=\\
&&\xymatrix{
-k-1\ar@{.}[d]&&&-k\ar@{.}[d]\\
{
\c\left(\input{figure/r3-lem-planar-1k-5-mf}\right)_n
\begin{array}{c}
\{(k+1)n\}\\
\left<k+1\right>
\end{array}
}
\ar[rrr]^{\chi_{+,(2,1,6,3)}^{[1,k+1]}\boxtimes\id_{\overline{\Lambda}_{(6;4,5)}^{[1,k]}}}&&&
{
\c\left(\input{figure/r3-lem-planar-1k-6-mf}\right)_n
\begin{array}{c}
\{(k+1)n-1\}\\
\left<k+1\right>
\end{array}
},
}
\end{eqnarray}
\begin{eqnarray}
\nonumber
&&\c\left(\input{figure/r3-lem-planar-1k-5-mf}\right)_n\\
\label{mf-planar-3}
&&\simeq
K\left(
\left(
\begin{array}{c}
\vspace{0.1cm}A_{1,(2,1,6,3)}^{[1,k+1]}\\
\vspace{0.1cm}\vdots\\
\vspace{0.1cm}A_{k+1,(2,1,6,3)}^{[1,k+1]}\\
u_{k+2,(2,1,6,3)}^{[1,k+1]}\\
\end{array}
\right);
\left(
\begin{array}{c}
\vspace{0.1cm}X_{1,(1,2)}^{(k+1,1)}-X_{1,(3,4,5)}^{(1,1,k)}\\
\vspace{0.1cm}\vdots\\
\vspace{0.1cm}X_{k+1,(1,2)}^{(k+1,1)}-X_{k+1,(3,4,5)}^{(1,1,k)}\\
(x_{1,2}-x_{1,3})X_{k+1,(1,3)}^{(k+1,-1)}
\end{array}
\right)
\right)_{Q_4}\{-k-1\},\\
\nonumber
&&\c\left(\input{figure/r3-lem-planar-1k-6-mf}\right)_n\\
\label{mf-planar-4}
&&\simeq
K\left(
\left(
\begin{array}{c}
A_{1,(2,1,6,3)}^{[1,k+1]}\\
\vspace{0.1cm}\vdots\\
A_{k+1,(2,1,6,3)}^{[1,k+1]}\\
u_{k+2,(2,1,6,3)}^{[1,k+1]}(x_{1,2}-x_{1,3})
\end{array}
\right);
\left(
\begin{array}{c}
X_{1,(1,2)}^{(k+1,1)}-X_{1,(3,4,5)}^{(1,1,k)}\\
\vspace{0.1cm}\vdots\\
X_{k+1,(1,2)}^{(k+1,1)}-X_{k+1,(3,4,5)}^{(1,1,k)}\\
X_{k+1,(1,3)}^{(k+1,-1)}
\end{array}
\right)
\right)_{Q_4}\{-k\},
\end{eqnarray}
where

\begin{eqnarray}
\nonumber	
Q_4&:=&R_{(1,2,3,4,5,6)}^{(k+1,1,1,1,k,k+1)}\left/\left<x_{1,6}-X_{1,(4,5)}^{(1,k)},\ldots,x_{k+1,6}-X_{k+1,(4,5)}^{(1,k)}\right>\right. \\
\nonumber
&\simeq&R_{(1,2,3,4,5)}^{(k+1,1,1,1,k)}.
\end{eqnarray}
The right-hand side sequences 
$\left(
x_{1,1}-X_{1,(2,3,4,5)}^{(-1,1,1,k)},\ldots,x_{1,k+1}-X_{k+1,(2,3,4,5)}^{(-1,1,1,k)},(x_{1,2}-x_{1,4})X_{k,(2,5)}^{(-1,k)}
\right)$ of the matrix factorization $\overline{M}_2$ and 
$\left(
X_{1,(1,2)}^{(k+1,1)}-X_{1,(3,4,5)}^{(1,1,k)},\ldots,X_{k+1,(1,2)}^{(k+1,1)}-X_{k+1,(3,4,5)}^{(1,1,k)},X_{k+1,(1,3)}^{(k+1,-1)}
\right)$
of the matrix factorization (\ref{mf-planar-4}) transform to each other by a linear transformation over $R_{(1,2,3,4,5)}^{(k+1,1,1,1,k)}$.
Then, by Proposition \ref{equiv} and Theorem \ref{reg-eq}, we have
\begin{equation}
\label{isom1}
\overline{M}_2\simeq\c\left(\input{figure/r3-lem-planar-1k-6-mf}\right)_n.
\end{equation}
The sequences 
$\left(
x_{1,1}-X_{1,(2,3,4,5)}^{(-1,1,1,k)},\ldots,x_{1,k+1}-X_{k+1,(2,3,4,5)}^{(-1,1,1,k)},(x_{1,2}-x_{1,3})(x_{1,2}-x_{1,4})X_{k,(2,5)}^{(-1,k)}
\right)$ of the matrix factorization $\overline{M}_1$ and 
$\left(
X_{1,(1,2)}^{(k+1,1)}-X_{1,(3,4,5)}^{(1,1,k)},\ldots,X_{k+1,(1,2)}^{(k+1,1)}-X_{k+1,(3,4,5)}^{(1,1,k)},(x_{1,2}-x_{1,3})X_{k+1,(1,3)}^{(k+1,-1)}
\right)$
of the matrix factorization (\ref{mf-planar-3}) also transform to each other.
We have
\begin{equation}
\label{isom2}
\overline{M}_1\simeq\c\left(\input{figure/r3-lem-planar-1k-5-mf}\right)_n.
\end{equation}
Thus, in $\k^b(\HMF^{gr}_{R_{(1,2,3,4,5)}^{(k+1,1,1,1,k)},\omega_3})$, we obtain
\begin{equation}
\nonumber
\c\left(\input{figure/r3-lem-plus-1k-1}\right)_n\simeq\c\left(\input{figure/r3-lem-plus-1k-2}\right)_n.
\end{equation}
\end{proof}
\begin{proof}[{\bf Proof of Proposition \ref{prop-r3} (2)}]
The complex for the diagram $\input{figure/r3-lem-minus-1k-1-text}$ is described as a complex of factorizations of $\k^b(\HMF^{gr}_{R_{(1,2,3,4,5)}^{(k+1,1,1,1,k)},\omega_3})$,
\begin{eqnarray}
\label{r3-lem-planar-comp-2}
&&\hspace{1.5cm}
\c\left(\input{figure/r3-lem-minus-1k-1-mf}\right)_n=\\
\nonumber
&&
\xymatrix{
k-1\ar@{.}[d]&k\ar@{.}[d]&k+1\ar@{.}[d]\\
{\c\left(\input{figure/r3-lem-planar-1k-4-mf}\right)_n
\hspace{-0.5cm}
\begin{array}{c}
{}_{\{-(k+1)n+2\}}\\
\left<k+1\right>
\end{array}}
\ar[r]^{\left(\hspace{-0.2cm}\begin{array}{c}\overline{\zeta}_{-,1}\\ \overline{\zeta}_{-,2} \end{array}\hspace{-0.2cm}\right)}&
{
\begin{array}{c}
\c\left(\input{figure/r3-lem-planar-1k-2-mf}\right)_n
\hspace{-0.5cm}
\begin{array}{c}
{}_{\{-(k+1)n+1\}}\\
\left<k+1\right>
\end{array}
\\
\bigoplus\\
\c\left(\input{figure/r3-lem-planar-1k-3-mf}\right)_n
\hspace{-0.5cm}
\begin{array}{c}
{}_{\{-(k+1)n+1\}}\\
\left<k+1\right>
\end{array}
\end{array}
}
\ar[r]^{\txt{$(\overline{\zeta}_{-,3},\overline{\zeta}_{-,4})$}}&
{\c\left(\input{figure/r3-lem-planar-1k-1-mf}\right)_n
\hspace{-0.5cm}
\begin{array}{c}
{}_{\{-(k+1)n\}}\\
\left<k+1\right>
\end{array}
}.
}
\end{eqnarray}
By the discussion of Proof of lemma \ref{prop-r3} (1), we also have
\begin{eqnarray}
\label{minus-mf1}
\c\left(\input{figure/r3-lem-planar-1k-4-mf}\right)_n
&\simeq&\overline{S}_{(1,2;3,4,5)}^{[k+1,1;1,1,k]}\boxtimes
K\left(B;X_{k,(2,5)}^{(-1,k)}\right)_{R_{(1,2,3,4,5)}^{(k+1,1,1,1,k)}}\{-k+1\},\\
\label{minus-mf2}
\c\left(\input{figure/r3-lem-planar-1k-2-mf}\right)_n
&\simeq&\overline{S}_{(1,2;3,4,5)}^{[k+1,1;1,1,k]}\boxtimes
K\left(\frac{B}{(x_{1,2}-x_{1,4})};(x_{1,2}-x_{1,4})X_{k,(2,5)}^{(-1,k)}
\right)_{R_{(1,2,3,4,5)}^{(k+1,1,1,1,k)}}\{-k\},\\
\label{minus-mf3}
\c\left(\input{figure/r3-lem-planar-1k-3-mf}\right)_n
&\simeq&\overline{S}_{(1,2;3,4,5)}^{[k+1,1;1,1,k]}\boxtimes
K\left(B;X_{k,(2,5)}^{(-1,k)}\right)_{Q_3\left/\langle (x_{1,8}-x_{1,3})X_{1,(6,3)}^{(1,-1)}\rangle\right.}\{-k\},\\
\label{minus-mf4}
\c\left(\input{figure/r3-lem-planar-1k-1-mf}\right)_n
&\simeq&\overline{S}_{(1,2;3,4,5)}^{[k+1,1;1,1,k]}\boxtimes
K\left(B_0;(x_{1,2}-x_{1,8})X_{k,(2,5)}^{(-1,k)}\right)_{Q_3\left/\langle (x_{1,8}-x_{1,3})X_{1,(6,3)}^{(1,-1)}\rangle\right.}\{-k-1\},
\end{eqnarray}
where 

\begin{eqnarray}
\nonumber
Q_3:=R_{(1,2,3,4,5,6,7,8)}^{(k+1,1,1,1,k,1,k,1)}\left/\left<X_{1,(8,6)}^{(1,1)}-X_{1,(4,3)}^{(1,1)},X_{1,(2,7)}^{(1,k)}-X_{1,(5,8)}^{(k,1)},\ldots,X_{k,(2,7)}^{(1,k)}-X_{k,(5,8)}^{(k,1)}\right>\right. .
\end{eqnarray}
We consider the decomposition of the partial matrix factorization $K\left(B;X_{k,(2,5)}^{(-1,k)}\right)_{Q_3\left/\langle (x_{1,8}-x_{1,3})X_{1,(6,3)}^{(1,-1)}\rangle\right.}$ of $(\ref{minus-mf3})$
\begin{eqnarray}
\nonumber
&&\xymatrix{R_{12}\ar[r]_(.4){
\left(\hspace{-0.2cm}
\begin{array}{cc}
B&0\\
0&B
\end{array}\hspace{-0.2cm}
\right)
}
&R_{12}\{2k-n-1\}\ar[rrrr]_(.6){
\left(
\begin{array}{cc}
X_{k,(2,5)}^{(-1,k)}&0\\
0&X_{k,(2,5)}^{(-1,k)}
\end{array}
\right)
}
&&&&R_{12}}\\
\nonumber
&\simeq&
K\left(B;X_{k,(2,5)}^{(-1,k)}\right)_{R_{(1,2,3,4,5)}^{(k+1,1,1,1,k)}}
\oplus K\left(B;X_{k,(2,5)}^{(-1,k)}\right)_{R_{(1,2,3,4,5)}^{(k+1,1,1,1,k)}}\{2\}.
\end{eqnarray}
and the decomposition of the partial matrix factorization $K\left(B_0;(x_{1,2}-x_{1,8})X_{k,(2,5)}^{(-1,k)}\right)_{Q_3\left/\langle (x_{1,8}-x_{1,3})X_{1,(6,3)}^{(1,-1)}\rangle\right.}$ of $(\ref{minus-mf4})$
\begin{eqnarray}
\nonumber
&&\xymatrix{R_{11}
\ar[rrrrr]_(.42){
\left(
\begin{array}{cc}
0&B\\
\frac{B}{(x_{1,2}-x_{1,3})(x_{1,2}-x_{1,4})}&0
\end{array}
\right)
}
&&&&&R_{12}\{2k-n+1\}
\ar[rrrrrrr]_(.55){
\left(
\begin{array}{cc}
0&(x_{1,2}-x_{1,3})(x_{1,2}-x_{1,4})X_{k,(2,5)}^{(-1,k)}\\
X_{k,(2,5)}^{(-1,k)}&0
\end{array}
\right)
}
&&&&&&&R_{11}}\\
\nonumber
&\simeq&K\left(
\frac{B}{(x_{1,2}-x_{1,3})(x_{1,2}-x_{1,4})};(x_{1,2}-x_{1,3})(x_{1,2}-x_{1,4})X_{k,(2,5)}^{(-1,k)}
\right)_{R_{(1,2,3,4,5)}^{(k+1,1,1,1,k)}}\oplus 
K\left(
B;X_{k,(2,5)}^{(-1,k)}
\right)_{R_{(1,2,3,4,5)}^{(k+1,1,1,1,k)}}\{2\}
\end{eqnarray}
Then, we obtain isomorphisms of the matrix factorizations (\ref{minus-mf3}) and (\ref{minus-mf4})
\begin{eqnarray}
\nonumber
&&\c\left(\input{figure/r3-lem-planar-1k-3-mf}\right)_n
\\
\label{minus-mf6}
&&\simeq
\overline{S}_{(1,2;3,4,5)}^{[k+1,1;1,1,k]}\boxtimes
K\left(B;X_{k,(2,5)}^{(-1,k)}\right)_{R_{(1,2,3,4,5)}^{(k+1,1,1,1,k)}}\{-k\}
\oplus\overline{S}_{(1,2;3,4,5)}^{[k+1,1;1,1,k]}\boxtimes
K\left(B;X_{k,(2,5)}^{(-1,k)}\right)_{R_{(1,2,3,4,5)}^{(k+1,1,1,1,k)}}\{-k+2\},\\
\nonumber
&&\c\left(\input{figure/r3-lem-planar-1k-4-mf}\right)_n
\\
\label{minus-mf5}
&&\simeq
\overline{S}_{(1,2;3,4,5)}^{[k+1,1;1,1,k]}\boxtimes
K\left(
\frac{B}{(x_{1,2}-x_{1,3})(x_{1,2}-x_{1,4})};(x_{1,2}-x_{1,3})(x_{1,2}-x_{1,4})X_{k,(2,5)}^{(-1,k)}
\right)_{R_{(1,2,3,4,5)}^{(k+1,1,1,1,k)}}\{-k-1\}\\
\nonumber
&&\hspace{1cm}\oplus\overline{S}_{(1,2;3,4,5)}^{[k+1,1;1,1,k]}\boxtimes
K\left(B;X_{k,(2,5)}^{(-1,k)}\right)_{R_{(1,2,3,4,5)}^{(k+1,1,1,1,k)}}\{-k+1\}.
\end{eqnarray}
For these decompositions (\ref{minus-mf1}), (\ref{minus-mf2}), (\ref{minus-mf5}) and (\ref{minus-mf6}), the morphisms $\overline{\zeta}_{-,1}$, $\overline{\zeta}_{-,2}$, $\overline{\zeta}_{-,3}$ and $\overline{\zeta}_{-,4}$ of the complex (\ref{r3-lem-planar-comp-2}) transform into
\begin{eqnarray}
\nonumber
\overline{\zeta}_{-,1}&\simeq&
\id_{\overline{S}_{(1,2;3,4,5)}^{[k+1,1;1,1,k]}}\boxtimes(x_{1,2}-x_{1,4},1),\\
\nonumber
\overline{\zeta}_{-,2}&\simeq&
\left(
\begin{array}{c}
\id_{\overline{S}_{(1,2;3,4,5)}^{[k+1,1;1,1,k]}}\boxtimes(-x_{1,2}+x_{1,4},-x_{1,2}+x_{1,4})\\
\id_{\overline{S}_{(1,2;3,4,5)}^{[k+1,1;1,1,k]}}\boxtimes(1,1)
\end{array}
\right),\\
\nonumber
\overline{\zeta}_{-,3}&\simeq&
\left(
\begin{array}{c}
\id_{\overline{S}_{(1,2;3,4,5)}^{[k+1,1;1,1,k]}}\boxtimes(x_{1,2}-x_{1,3},1)\\
\id_{\overline{S}_{(1,2;3,4,5)}^{[k+1,1;1,1,k]}}\boxtimes(-1,-x_{1,2}+x_{1,4})
\end{array}
\right),\\
\nonumber
\overline{\zeta}_{-,4}&\simeq&
-\left(
\begin{array}{cc}
0&\id_{\overline{S}_{(1,2;3,4,5)}^{[k+1,1;1,1,k]}}\boxtimes((x_{1,2}-x_{1,3})(x_{1,2}-x_{1,4}),1)\\
\id_{\overline{S}_{(1,2;3,4,5)}^{[k+1,1;1,1,k]}}\boxtimes(1,1)&0
\end{array}
\right).
\end{eqnarray}
Then, the complex (\ref{r3-lem-planar-comp-2}) is isomorphic, in $\k^b(\HMF^{gr}_{R_{(1,2,3,4,5)}^{[k+1,1,1,1,k]},\omega_3})$, to
\begin{equation}
\nonumber
\xymatrix{
k\ar@{.}[d]&&&k+1\ar@{.}[d]\\
{\overline{M}_2
\{1-(k+1)n\}
\left<k+1\right>
}
\ar[rrr]^{\id_{\overline{S}}\boxtimes(x_{1,2}-x_{1,3},1)}&&&
{\overline{M}_1
\{-(k+1)n\}
\left<k+1\right>
}
.}
\end{equation}
Since we have (the isomorphisms (\ref{isom1}) and (\ref{isom2}))
\begin{equation}
\nonumber
\overline{M}_1\simeq\c\left(\input{figure/r3-lem-planar-1k-5-mf}\right)_n,
\overline{M}_2\simeq\c\left(\input{figure/r3-lem-planar-1k-6-mf}\right)_n,
\end{equation}
thus we obtain
\begin{equation}
\nonumber
\c\left(\input{figure/r3-lem-minus-1k-1}\right)_n\simeq\c\left(\input{figure/r3-lem-minus-1k-2}\right)_n.
\end{equation}
\end{proof}
\begin{proof}[{\bf Proof of Proposition \ref{prop-r3} (3)}]
The complex for the diagram $\input{figure/r3-lem-plus-k1-1-text}$ is described as a complex of factorizations of  $\k^b(\HMF^{gr}_{R_{(1,2,3,4,5)}^{(k+1,1,1,k,1)},\omega_4})$ $(\omega_4=F_{k+1}(\mathbb{X}^{(k+1)}_{(1)})+F_{1}(\mathbb{X}^{(1)}_{(2)})-F_{1}(\mathbb{X}^{(1)}_{(3)})-F_{k}(\mathbb{X}^{(k)}_{(4)})-F_{1}(\mathbb{X}^{(1)}_{(5)}))$,
\begin{eqnarray}
\label{r3-lem-planar-comp-3}
&&\hspace{1cm}
\c\left(\input{figure/r3-lem-plus-k1-1-mf}\right)_n=\\
\nonumber
&&
\xymatrix{
-k-1\ar@{.}[d]&-k\ar@{.}[d]&-k+1\ar@{.}[d]\\
{\c\left(\input{figure/r3-lem-planar-k1-1-mf}\right)_n
\hspace{-0.5cm}
\begin{array}{c}
{}_{\{(k+1)n\}}\\
\left<k+1\right>
\end{array}}
\ar[r]^{\left(\hspace{-0.2cm}
\begin{array}{c}
\overline{\xi}_{+,1}\\
\overline{\xi}_{+,2} 
\end{array}\hspace{-0.2cm}\right)}&
{
\begin{array}{c}
\c\left(\input{figure/r3-lem-planar-k1-2-mf}\right)_n
\hspace{-0.5cm}
\begin{array}{c}
{}_{\{(k+1)n-1\}}\\
\left<k+1\right>
\end{array}
\\
\bigoplus\\
\c\left(\input{figure/r3-lem-planar-k1-3-mf}\right)_n
\hspace{-0.5cm}
\begin{array}{c}
{}_{\{(k+1)n-1\}}\\
\left<k+1\right>
\end{array}
\end{array}
}
\ar[r]^{\txt{$(\overline{\xi}_{+,3},\overline{\xi}_{+,4})$}}&
{\c\left(\input{figure/r3-lem-planar-k1-4-mf}\right)_n
\hspace{-0.5cm}
\begin{array}{c}
{}_{\{(k+1)n-2\}}\\
\left<k+1\right>
\end{array}
},
}
\end{eqnarray}
where
\begin{eqnarray}
\nonumber
\overline{\xi}_{+,1}&=&\id_{\overline{\Lambda}_{(1;6,7)}^{[k,1]}}\boxtimes\id_{\overline{S}_{(8,6,4,3)}^{[1,k]}}\boxtimes(1,x_{1,8}-x_{1,3})\boxtimes\id_{\overline{M}_{(2,7,5,8)}^{[1,1]}},
\\
\nonumber
\overline{\xi}_{+,2}&=&\id_{\overline{\Lambda}_{(1;6,7)}^{[k,1]}}\boxtimes\id_{\overline{M}_{(8,6,4,3)}^{[1,k]}}\boxtimes\id_{\overline{S}_{(2,7,5,8)}^{[1,1]}}\boxtimes(1,x_{1,2}-x_{1,8}),\\
\nonumber
\overline{\xi}_{+,3}&=&\id_{\overline{\Lambda}_{(1;6,7)}^{[k,1]}}\boxtimes\id_{\overline{N}_{(8,6,4,3)}^{[1,k]}}\boxtimes\id_{\overline{S}_{(2,7,5,8)}^{[1,1]}}\boxtimes(1,x_{1,2}-x_{1,8}),
\\
\nonumber
\overline{\xi}_{+,4}&=&-\id_{\overline{\Lambda}_{(1;6,7)}^{[k,1]}}\boxtimes\id_{\overline{S}_{(8,6,4,3)}^{[1,k]}}\boxtimes(1,x_{1,8}-x_{1,3})\boxtimes\id_{\overline{N}_{(2,7,5,8)}^{[1,1]}}.
\end{eqnarray}
First, we have
\begin{eqnarray}
\nonumber
\c\left(\input{figure/r3-lem-planar-k1-1-mf}\right)_n
&=&\overline{\Lambda}_{(1;6,7)}^{[k,1]}\boxtimes \overline{M}_{(8,6,4,3)}^{[1,k]} \boxtimes \overline{M}_{(2,7,5,8)}^{[1,1]}\\
\nonumber
&\simeq&
K\left(
\left(
\begin{array}{c}
\Lambda_{1,(1;6,7)}^{[k,1]}\\
\vspace{0.1cm}\vdots\\
\Lambda_{k+1,(1;6,7)}^{[k,1]}
\end{array}
\right);
\left(
\begin{array}{c}
x_{1,1}-X_{1,(6,7)}^{(k,1)}\\
\vspace{0.1cm}\vdots\\
x_{k+1,1}-X_{k+1,(6,7)}^{(k,1)}
\end{array}
\right)
\right)_{R_{(1,6,7)}^{(1,k,1)}}\\
&&
\nonumber
\hspace{0.5cm}\boxtimes
K\left(
\left(
\begin{array}{c}
A_{1,(8,6,4,3)}^{[1,k]}\\
\vdots\\
A_{k,(8,6,4,3)}^{[1,k]}\\
u_{k+1,(8,6,4,3)}^{[1,k]}
\end{array}
\right);
\left(
\begin{array}{c}
X_{1,(8,6)}^{(1,k)}-X_{1,(4,3)}^{(k,1)}\\
\vdots\\
X_{k,(8,6)}^{(1,k)}-X_{k,(4,3)}^{(k,1)}\\
(x_{1,8}-x_{1,3})X_{k,(6,3)}^{(k,-1)}
\end{array}
\right)
\right)_{R_{(8,6,4,3)}^{(1,k,k,1)}}\{-k\}\\
&&
\nonumber
\hspace{0.5cm}\boxtimes
K\left(
\left(
\begin{array}{c}
A_{1,(2,7,5,8)}^{[1,1]}\\
u_{2,(2,7,5,8)}^{[1,1]}
\end{array}
\right);
\left(
\begin{array}{c}
X_{1,(2,7)}^{(1,1)}-X_{1,(5,8)}^{(1,1)}\\
(x_{1,2}-x_{1,8})X_{1,(7,8)}^{(1,-1)}
\end{array}
\right)
\right)_{R_{(2,7,5,8)}^{(1,1,1,1)}}\{-1\}\\
&\simeq&
\label{mf-r3-1-00-3}
K\left(
\left(
\begin{array}{c}
\vspace{0.1cm}\Lambda_{1,(1;6,7)}^{[k,1]}\\
\vspace{0.1cm}\vdots\\
\vspace{0.1cm}\Lambda_{k+1,(1;6,7)}^{[k,1]}\\
u_{2,(2,7,5,8)}^{[1,1]}
\end{array}
\right);
\left(
\begin{array}{c}
\vspace{0.1cm}x_{1,1}-X_{1,(6,7)}^{(k,1)}\\
\vspace{0.1cm}\vdots\\
\vspace{0.1cm}x_{k+1,1}-X_{k+1,(6,7)}^{(k,1)}\\
(x_{1,2}-x_{1,8})X_{1,(7,8)}^{(1,-1)}
\end{array}
\right)
\right)_{Q_5\left/\langle (x_{1,8}-x_{1,3})X_{k,(6,3)}^{(k,-1)}\rangle\right.}\{-k-1\},
\end{eqnarray}
where 

\begin{equation}
\nonumber
Q_5:=R_{(1,2,3,4,5,6,7,8)}^{(k+1,1,1,k,1,k,1,1)}\left/\left<
X_{1,(8,6)}^{(1,k)}-X_{1,(4,3)}^{(k,1)},
\ldots,
X_{k,(8,6)}^{(1,k)}-X_{k,(4,3)}^{(k,1)},
X_{1,(2,7)}^{(1,1)}-X_{1,(5,8)}^{(1,1)}\right>\right. .
\end{equation}
The quotient $Q_5\left/\langle (x_{1,8}-x_{1,3})X_{k,(6,3)}^{(k,-1)}\rangle\right.$ has equations
\begin{eqnarray}
\nonumber
&&x_{j,6}=X_{j,(3,4,8)}^{(1,k,-1)} \hspace{1cm}(1 \leq j \leq k),\\
\nonumber
&&(x_{1,8}-x_{1,3})X_{k,(4,8)}^{(k,-1)}=0,\\
\nonumber
&&x_{1,7}=X_{1,(2,5,8)}^{(-1,1,1)}.
\end{eqnarray}
In the quotient $Q_5\left/\langle (x_{1,8}-x_{1,3})X_{k,(6,3)}^{(k,-1)}\rangle\right.$, $X_{j,(6,7)}^{(k,1)}$ $(1\leq j \leq k)$ equals to
\begin{eqnarray}
\nonumber
x_{j,6}+x_{j-1,6}x_{1,7}&\equiv&X_{j,(3,4,8)}^{(1,k,-1)}+X_{j-1,(3,4,8)}^{(1,k,-1)}X_{1,(2,5,8)}^{(-1,1,1)}\\
\nonumber
&=&X_{j,(3,4)}^{(1,k)}+X_{1,(2,5)}^{(-1,1)}X_{j-1,(2,3,4)}^{(-1,1,k)}
+(x_{1,2}-x_{1,8})X_{1,(2,5)}^{(-1,1)}X_{j-1,(2,3,4,8)}^{(-1,1,k,-1)}
\end{eqnarray}
and $X_{k,(6,7)}^{(k,1)}$ equals to
\begin{eqnarray}
\nonumber
x_{k,6}x_{1,7}&\equiv&X_{k,(3,4,8)}^{(1,k,-1)}X_{1,(2,5,8)}^{(-1,1,1)}\\
\nonumber
&\equiv&
X_{k+1,(3,4)}^{(1,k)}
+X_{1,(2,5)}^{(-1,1)}X_{k,(2,3,4)}^{(-1,1,k)}
+(x_{1,2}-x_{1,8})X_{1,(2,5)}^{(-1,1)}X_{k-1,(2,3,4,8)}^{(-1,1,k,-1)}.
\end{eqnarray}
Then, the matrix factorization (\ref{mf-r3-1-00-3}) is isomorphic to 
\begin{eqnarray}
\label{mf-r3-2-00-3}
&&K\left(
\left(
\begin{array}{c}
\vspace{0.1cm}\Lambda_{1,(1;6,7)}^{[k,1]}\\
\vspace{0.1cm}\vdots\\
\vspace{0.1cm}\Lambda_{k+1,(1;6,7)}^{[k,1]}\\
u_{2,(2,7,5,8)}^{[1,1]}-\sum_{j=2}^{k+1}X_{j-2,(2,3,4,8)}^{(-1,1,k,-1)}\Lambda_{j,(1;6,7)}^{[k,1]}
\end{array}
\right);
\left(
\begin{array}{c}
\vspace{0.1cm}x_{1,1}-X_{1,(2,3,4,5)}^{(-1,1,k,1)}\\
\vspace{0.1cm}\vdots\\
\vspace{0.1cm}x_{k+1,1}-X_{k+1,(2,3,4,5)}^{(-1,1,k,1)}\\
(x_{1,2}-x_{1,8})X_{1,(2,5)}^{(-1,1)}
\end{array}
\right)
\right)_{Q_5\left/\langle (x_{1,8}-x_{1,3})X_{k,(6,3)}^{(k,-1)}\rangle\right.}\{-k-1\}.
\end{eqnarray}
By Corollary \ref{induce-sq1}, there exist polynomials $C_1$, $C_2$, $\ldots$, $C_{k+1}$ $\in R_{(1,2,3,4,5)}^{(k+1,1,1,k,1)}$ and $C_0\in Q_5\left/\langle (x_{1,8}-x_{1,3})X_{k,(6,3)}^{(k,-1)}\rangle\right.$ satisfying $(x_{1,2}-x_{1,8})C_0\equiv C\in R_{(1,2,3,4,5)}^{(k+1,1,1,k,1)}$ $(\mod Q_5\left/\langle (x_{1,8}-x_{1,3})X_{k,(6,3)}^{(k,-1)}\rangle\right.)$ and we have an isomorphism to the factorization (\ref{mf-r3-2-00-3})
\begin{eqnarray}
\label{mf-r3-3-00-3}
\overline{S}_{(1,2;3,4,5)}^{[k+1,1;1,k,1]}\boxtimes
K\left(C_0;(x_{1,2}-x_{1,8})X_{1,(2,5)}^{(-1,1)}\right)_{Q_5\left/\langle (x_{1,8}-x_{1,3})X_{k,(6,3)}^{(k,-1)}\rangle\right.}\{-k-1\},
\end{eqnarray}
where
\begin{equation}
\label{common-mf-2}
\overline{S}_{(1,2;3,4,5)}^{[k+1,1;1,k,1]}
:=K\left(
\left(
\begin{array}{c}
\vspace{0.1cm}C_1\\
\vspace{0.1cm}\vdots\\
\vspace{0.1cm}C_{k+1}
\end{array}
\right);
\left(
\begin{array}{c}
x_{1,1}-X_{1,(2,3,4,5)}^{(-1,1,k,1)}\\
\vdots\\
x_{k+1,1}-X_{k+1,(2,3,4,5)}^{(-1,1,k,1)}
\end{array}
\right)
\right)_{R_{(1,2,3,4,5)}^{(k+1,1,1,k,1)}}.
\end{equation}
We choose isomorphisms of $Q_5\left/\langle (x_{1,8}-x_{1,3})X_{k,(6,3)}^{(k,-1)}\rangle\right.$ to be
\begin{eqnarray}
\nonumber
R_{13}&\simeq&R_{(1,2,3,4,5)}^{(k+1,1,1,k,1)}\oplus(x_{1,2}-x_{1,8})R_{(1,2,3,4,5)}^{(k+1,1,1,k,1)}\oplus \ldots \oplus x_{1,8}^{k-1}(x_{1,2}-x_{1,8})R_{(1,2,3,4,5)}^{(k+1,1,1,k,1)},\\
\nonumber
R_{14}&\simeq&R_{(1,2,3,4,5)}^{(k+1,1,1,k,1)}\oplus x_{1,8}R_{(1,2,3,4,5)}^{(k+1,1,1,k,1)}\oplus\ldots\oplus x_{1,8}^{k-1}R_{(1,2,3,4,5)}^{(k+1,1,1,k,1)}\oplus X_{k,(2,3,4,8)}^{(-1,1,k,-1)}R_{(1,2,3,4,5)}^{(k+1,1,1,k,1)}.
\end{eqnarray}
Then, the partial matrix factorization $K(C_0;(x_{1,2}-x_{1,8})X_{1,(2,5)}^{(-1,1)})_{Q_5\left/\langle (x_{1,8}-x_{1,3})X_{k,(6,3)}^{(k,-1)}\rangle\right.}$ is isomorphic to
\begin{eqnarray}
\nonumber
&&\xymatrix{R_{13}
\ar[rrrrr]_(.42){
\left(
\begin{array}{cc}
{}^{t}\mathfrak{0}_{k}&E_{k}(C)\\
\frac{C}{(x_{1,2}-x_{1,3})X_{k,(2,4)}^{(-1,k)}}&\mathfrak{0}_{k}
\end{array}
\right)
}
&&&&&R_{14}\{3-n\}
\ar[rrrrrrr]_(.55){
\left(
\begin{array}{cc}
\mathfrak{0}_{k}&(x_{1,2}-x_{1,3})X_{k,(2,4)}^{(-1,k)}X_{1,(2,5)}^{(-1,k)}\\
E_{k}(X_{1,(2,5)}^{(-1,k)})&{}^{t}\mathfrak{0}_{k}
\end{array}
\right)
}
&&&&&&&R_{13}}\\
\nonumber
&\simeq&K\left(
\frac{C}{(x_{1,2}-x_{1,3})X_{k,(2,4)}^{(-1,k)}};(x_{1,2}-x_{1,3})X_{k,(2,4)}^{(-1,k)}X_{1,(2,5)}^{(-1,k)}
\right)_{R_{(1,2,3,4,5)}^{(k+1,1,1,k,1)}}\\
\nonumber
&&\hspace{1cm}\oplus 
\bigoplus_{j=1}^{n}K\left(C;X_{k,(1,5)}^{(-1,k)}\right)_{R_{(1,2,3,4,5)}^{(k+1,1,1,k,1)}}\{2j\}.
\end{eqnarray}
Thus, the matrix factorization (\ref{mf-r3-3-00-3}) is decomposed into
\begin{eqnarray}
\label{mf-r3-fin-00-3}
&&\overline{S}_{(1,2;3,4,5)}^{[k+1,1;1,k,1]}\boxtimes
K\left(
\frac{C}{(x_{1,2}-x_{1,3})X_{k,(2,4)}^{(-1,k)}};(x_{1,2}-x_{1,3})X_{k,(2,4)}^{(-1,k)}X_{1,(2,5)}^{(-1,k)}
\right)_{R_{(1,2,3,4,5)}^{(k+1,1,1,k,1)}}\{-k-1\}\\
\nonumber
&&\hspace{1cm}\oplus 
\bigoplus_{j=1}^{k}\overline{S}_{(1,2;3,4,5)}^{[k+1,1;1,k,1]}\boxtimes
K\left(C;X_{1,(2,5)}^{(-1,k)}\right)_{R_{(1,2,3,4,5)}^{(k+1,1,1,k,1)}}\{-k-1+2j\}.
\end{eqnarray}
Secondly, we have
\begin{eqnarray}
\nonumber
\c\left(\input{figure/r3-lem-planar-k1-2-mf}\right)_n
&=&\overline{\Lambda}_{(1;6,7)}^{[k,1]}\boxtimes \overline{N}_{(8,6,4,3)}^{[1,k]} \boxtimes \overline{M}_{(2,7,5,8)}^{[1,1]}\\
\nonumber
&\simeq&
K\left(
\left(
\begin{array}{c}
\Lambda_{1,(1;6,7)}^{[k,1]}\\
\vspace{0.1cm}\vdots\\
\Lambda_{k+1,(1;6,7)}^{[k,1]}
\end{array}
\right);
\left(
\begin{array}{c}
x_{1,1}-X_{1,(6,7)}^{(k,1)}\\
\vspace{0.1cm}\vdots\\
x_{k+1,1}-X_{k+1,(6,7)}^{(k,1)}
\end{array}
\right)
\right)_{R_{(1,6,7)}^{(1,k,1)}}\\
&&
\nonumber
\hspace{0.5cm}\boxtimes
K\left(
\left(
\begin{array}{c}
A_{1,(8,6,4,3)}^{[1,k]}\\
\vdots\\
A_{k,(8,6,4,3)}^{[1,k]}\\
u_{k+1,(8,6,4,3)}^{[1,k]}(x_{1,8}-x_{1,3})
\end{array}
\right);
\left(
\begin{array}{c}
X_{1,(8,6)}^{(1,k)}-X_{1,(4,3)}^{(k,1)}\\
\vdots\\
X_{k,(8,6)}^{(1,k)}-X_{k,(4,3)}^{(k,1)}\\
X_{k,(6,3)}^{(k,-1)}
\end{array}
\right)
\right)_{R_{(8,6,4,3)}^{(1,k,k,1)}}\{-k+1\}\\
&&
\nonumber
\hspace{0.5cm}\boxtimes
K\left(
\left(
\begin{array}{c}
A_{1,(2,7,5,8)}^{[1,1]}\\
u_{2,(2,7,5,8)}^{[1,1]}
\end{array}
\right);
\left(
\begin{array}{c}
X_{1,(2,7)}^{(1,1)}-X_{1,(5,8)}^{(1,1)}\\
(x_{1,2}-x_{1,8})X_{1,(7,8)}^{(1,-1)}
\end{array}
\right)
\right)_{R_{(2,7,5,8)}^{(1,1,1,1)}}\{-1\}\\
&\simeq&
\label{mf-r3-1-10-3}
K\left(
\left(
\begin{array}{c}
\vspace{0.1cm}\Lambda_{1,(1;6,7)}^{[k,1]}\\
\vspace{0.1cm}\vdots\\
\vspace{0.1cm}\Lambda_{k+1,(1;6,7)}^{[k,1]}\\
u_{2,(2,7,5,8)}^{[1,1]}
\end{array}
\right);
\left(
\begin{array}{c}
\vspace{0.1cm}x_{1,1}-X_{1,(6,7)}^{(k,1)}\\
\vspace{0.1cm}\vdots\\
\vspace{0.1cm}x_{k+1,1}-X_{k+1,(6,7)}^{(k,1)}\\
(x_{1,2}-x_{1,8})X_{1,(7,8)}^{(1,-1)}
\end{array}
\right)
\right)_{Q_5\left/\langle X_{k,(6,3)}^{(k,-1)}\rangle\right.}\{-k\}.
\end{eqnarray}
The quotient $Q_5\left/\langle X_{k,(6,3)}^{(k,-1)}\rangle\right.$ has equations
\begin{eqnarray}
\nonumber
&&x_{j,6}=X_{j,(3,4,8)}^{(1,k,-1)} \hspace{1cm}(1 \leq j \leq k-1),\\
\nonumber
&&x_{k,6}=x_{1,3}X_{k-1,(4,8)}^{(k,-1)},\\
\nonumber
&&X_{k,(4,8)}^{(k,-1)}=0,\\
\nonumber
&&x_{1,7}=X_{1,(2,5,8)}^{(-1,1,1)}.
\end{eqnarray}
In the quotient $Q_5\left/\langle X_{k,(6,3)}^{(k,-1)}\rangle\right.$, $X_{j,(6,7)}^{(k,1)}$ $(1\leq j \leq k)$ equals to
\begin{eqnarray}
\nonumber
x_{j,6}+x_{j-1,6}x_{1,7}&\equiv&X_{j,(3,4,8)}^{(1,k,-1)}+X_{j-1,(3,4,8)}^{(1,k,-1)}X_{1,(2,5,8)}^{(-1,1,1)}\\
\nonumber
&=&X_{j,(3,4)}^{(1,k)}+X_{1,(2,5)}^{(-1,1)}X_{j-1,(2,3,4)}^{(-1,1,k)}
+(x_{1,2}-x_{1,8})X_{1,(2,5)}^{(-1,1)}X_{j-1,(2,3,4,8)}^{(-1,1,k,-1)},
\end{eqnarray}
and $X_{k+1,(6,7)}^{(k,1)}$ equals to
\begin{eqnarray}
\nonumber
x_{k,6}x_{1,7}&\equiv&x_{1,3}X_{k-1,(4,8)}^{(k,-1)}X_{1,(2,5,8)}^{(-1,1,1)}\\
\nonumber
&\equiv&X_{k+1,(3,4)}^{(1,k)}
+X_{1,(2,5)}^{(-1,1)}X_{k,(2,3,4)}^{(-1,1,k)}
+(x_{1,2}-x_{1,8})X_{1,(2,5)}^{(-1,1)}X_{k-1,(2,3,4,8)}^{(-1,1,k,-1)}.
\end{eqnarray}
Then, the matrix factorization (\ref{mf-r3-1-10-3}) is isomorphic to 
\begin{eqnarray}
\nonumber
&&K\left(
\left(
\begin{array}{c}
\vspace{0.1cm}\Lambda_{1,(1;6,7)}^{[k,1]}\\
\vspace{0.1cm}\vdots\\
\vspace{0.1cm}\Lambda_{k+1,(1;6,7)}^{[k,1]}\\
u_{2,(2,7,5,8)}^{[1,1]}-\sum_{j=2}^{k+1}X_{j-2,(2,3,4,8)}^{(-1,1,k,-1)}\Lambda_{j,(1;6,7)}^{[k,1]}
\end{array}
\right);
\left(
\begin{array}{c}
\vspace{0.1cm}x_{1,1}-X_{1,(2,3,4,5)}^{(-1,1,k,1)}\\
\vspace{0.1cm}\vdots\\
\vspace{0.1cm}x_{k+1,1}-X_{k+1,(2,3,4,5)}^{(-1,1,k,1)}\\
(x_{1,2}-x_{1,8})X_{1,(2,5)}^{(-1,1)}
\end{array}
\right)
\right)_{Q_5\left/\langle X_{k,(6,3)}^{(k,-1)}\rangle\right.}\{-k\}\\
\label{mf-r3-2-10-3}
&\simeq&
\overline{S}_{(1,2;3,4,5)}^{[k+1,1;1,k,1]}\boxtimes
K\left(C_0;(x_{1,2}-x_{1,8})X_{1,(2,5)}^{(-1,k)}\right)_{Q_5\left/\langle X_{k,(6,3)}^{(k,-1)}\rangle\right.}\{-k\}.
\end{eqnarray}
We choose isomorphisms of $Q_5\left/\langle X_{k,(6,3)}^{(k,-1)}\rangle\right.$ to be
\begin{eqnarray}
\nonumber
R_{14}&\simeq&R_{(1,2,3,4,5)}^{(k+1,1,1,k,1)}\oplus(x_{1,2}-x_{1,8})R_{(1,2,3,4,5)}^{(k+1,1,1,k,1)}\oplus \ldots \oplus x_{1,8}^{k-2}(x_{1,2}-x_{1,8})R_{(1,2,3,4,5)}^{(k+1,1,1,k,1)},\\
\nonumber
R_{15}&\simeq&R_{(1,2,3,4,5)}^{(k+1,1,1,k,1)}\oplus x_{1,8}R_{(1,2,3,4,5)}^{(k+1,1,1,k,1)}\oplus\ldots\oplus x_{1,8}^{k-2}R_{(1,2,3,4,5)}^{(k+1,1,1,k,1)}\oplus (-X_{k,(2,4,8)}^{(-1,k,-1)})R_{(1,2,3,4,5)}^{(k+1,1,1,k,1)}.
\end{eqnarray}
Then, the partial matrix factorization $K(C_0;(x_{1,2}-x_{1,8})X_{1,(2,5)}^{(-1,1)})_{Q_5\left/\langle X_{k,(6,3)}^{(k,-1)}\rangle\right.}$ is isomorphic to
\begin{eqnarray}
\nonumber
&&\xymatrix{R_{14}
\ar[rrrrr]_(.42){
\left(
\begin{array}{cc}
{}^{t}\mathfrak{0}_{k-1}&E_{k-1}(C)\\
\frac{C}{X_{k,(2,4)}^{(-1,k)}}&\mathfrak{0}_{k}
\end{array}
\right)
}
&&&&&R_{15}\{3-n\}
\ar[rrrrrrr]_(.55){
\left(
\begin{array}{cc}
\mathfrak{0}_{k-1}&X_{k,(2,4)}^{(-1,k)}X_{1,(2,5)}^{(-1,k)}\\
E_{k-1}(X_{1,(2,5)}^{(-1,k)})&{}^{t}\mathfrak{0}_{k-1}
\end{array}
\right)
}
&&&&&&&R_{14}}\\
\nonumber
&\simeq&K\left(
\frac{C}{X_{k,(2,4)}^{(-1,k)}};X_{k,(2,4)}^{(-1,k)}X_{1,(2,5)}^{(-1,k)}
\right)_{R_{(1,2,3,4,5)}^{(k+1,1,1,k,1)}}
\oplus 
\bigoplus_{j=1}^{k-1}K\left(C;X_{k,(1,5)}^{(-1,k)}\right)_{R_{(1,2,3,4,5)}^{(k+1,1,1,k,1)}}\{2j\}.
\end{eqnarray}
Thus, the matrix factorization (\ref{mf-r3-2-10-3}) is decomposed into
\begin{eqnarray}
\label{mf-r3-fin-10-3}
&&\overline{S}_{(1,2;3,4,5)}^{[k+1,1;1,k,1]}\boxtimes
K\left(
\frac{C}{X_{k,(2,4)}^{(-1,k)}};X_{k,(2,4)}^{(-1,k)}X_{1,(2,5)}^{(-1,k)}
\right)_{R_{(1,2,3,4,5)}^{(k+1,1,1,k,1)}}\{-k\}\\
\nonumber
&&\hspace{1cm}\oplus 
\bigoplus_{j=1}^{k-1}\overline{S}_{(1,2;3,4,5)}^{[k+1,1;1,k,1]}\boxtimes
K\left(C;X_{1,(2,5)}^{(-1,k)}\right)_{R_{(1,2,3,4,5)}^{(k+1,1,1,k,1)}}\{-k+2j\}.
\end{eqnarray}
Thirdly, we have
\begin{eqnarray}
\nonumber
\c\left(\input{figure/r3-lem-planar-k1-3-mf}\right)_n
&=&\overline{\Lambda}_{(1;6,7)}^{[k,1]}\boxtimes \overline{M}_{(8,6,4,3)}^{[1,k]} \boxtimes \overline{N}_{(2,7,5,8)}^{[1,1]}\\
\nonumber
&\simeq&
K\left(
\left(
\begin{array}{c}
\Lambda_{1,(1;6,7)}^{[k,1]}\\
\vspace{0.1cm}\vdots\\
\Lambda_{k+1,(1;6,7)}^{[k,1]}
\end{array}
\right);
\left(
\begin{array}{c}
x_{1,1}-X_{1,(6,7)}^{(k,1)}\\
\vspace{0.1cm}\vdots\\
x_{k+1,1}-X_{k+1,(6,7)}^{(k,1)}
\end{array}
\right)
\right)_{R_{(1,6,7)}^{(1,k,1)}}\\
&&
\nonumber
\hspace{0.5cm}\boxtimes
K\left(
\left(
\begin{array}{c}
A_{1,(8,6,4,3)}^{[1,k]}\\
\vdots\\
\vspace{0.1cm}A_{k,(8,6,4,3)}^{[1,k]}\\
u_{k+1,(8,6,4,3)}^{[1,k]}
\end{array}
\right);
\left(
\begin{array}{c}
X_{1,(8,6)}^{(1,k)}-X_{1,(4,3)}^{(k,1)}\\
\vdots\\
\vspace{0.1cm}X_{k,(8,6)}^{(1,k)}-X_{k,(4,3)}^{(k,1)}\\
(x_{1,8}-x_{1,3})X_{k,(6,3)}^{(k,-1)}
\end{array}
\right)
\right)_{R_{(8,6,4,3)}^{(1,k,k,1)}}\{-k\}\\
&&
\nonumber
\hspace{0.5cm}\boxtimes
K\left(
\left(
\begin{array}{c}
\vspace{0.1cm}A_{1,(2,7,5,8)}^{[1,1]}\\
u_{2,(2,7,5,8)}^{[1,1]}(x_{1,2}-x_{1,8})
\end{array}
\right);
\left(
\begin{array}{c}
\vspace{0.1cm}X_{1,(2,7)}^{(1,1)}-X_{1,(5,8)}^{(1,1)}\\
X_{1,(7,8)}^{(1,-1)}
\end{array}
\right)
\right)_{R_{(2,7,5,8)}^{(1,1,1,1)}}
\end{eqnarray}
\begin{eqnarray}
&\simeq&
\nonumber
K\left(
\left(
\begin{array}{c}
\vspace{0.1cm}\Lambda_{1,(1;6,7)}^{[k,1]}\\
\vspace{0.1cm}\vdots\\
\vspace{0.1cm}\Lambda_{k+1,(1;6,7)}^{[k,1]}\\
u_{2,(2,7,5,8)}^{[1,1]}(x_{1,2}-x_{1,8})
\end{array}
\right);
\left(
\begin{array}{c}
\vspace{0.1cm}x_{1,1}-X_{1,(6,7)}^{(k,1)}\\
\vspace{0.1cm}\vdots\\
\vspace{0.1cm}x_{k+1,1}-X_{k+1,(6,7)}^{(k,1)}\\
X_{1,(7,8)}^{(1,-1)}
\end{array}
\right)
\right)_{Q_5\left/\langle (x_{1,8}-x_{1,3})X_{k,(6,3)}^{(k,-1)}\rangle\right.}\{-k\}\\
\label{mf-r3-1-01-3}
&\simeq&
\overline{S}_{(1,2;3,4,5)}^{[k+1,1;1,k,1]}\boxtimes
K\left(C;X_{k,(2,5)}^{(-1,1)}\right)_{Q_5\left/\langle (x_{1,8}-x_{1,3})X_{k,(6,3)}^{(k,-1)}\rangle\right.}\{-k\}.
\end{eqnarray}
Since the partial matrix factorization $K\left(C;X_{1,(2,5)}^{(-1,k)}\right)_{Q_5\left/\langle (x_{1,8}-x_{1,3})X_{k,(6,3)}^{(k,-1)}\rangle\right.}$ is decomposed into
\begin{eqnarray}
\nonumber
&&\xymatrix{
R_{13}
\ar[rr]_(.4){E_{k+1}(C)}
&&R_{13}\{1-n\}
\ar[rrr]_(.55){E_{k+1}(X_{1,(2,5)}^{(-1,1)})}
&&&R_{13}}\\
\nonumber
&\simeq&
\bigoplus_{j=0}^{k}K\left(C;X_{1,(2,5)}^{(-1,1)}\right)_{R_{(1,2,3,4,5)}^{(k+1,1,1,k,1)}}\{2j\},
\end{eqnarray} 
then the matrix factorization (\ref{mf-r3-1-01-3}) is isomorphic to
\begin{eqnarray}
\label{mf-r3-fin-01-3}
&&
\bigoplus_{j=0}^{k}\overline{S}_{(1,2;3,4,5)}^{[k+1,1;1,k,1]}\boxtimes
K\left(C;X_{1,(2,5)}^{(-1,1)}\right)_{R_{(1,2,3,4,5)}^{(k+1,1,1,k,1)}}\{-k+2j\}.
\end{eqnarray}
Finally, we have
\begin{eqnarray}
\nonumber
\c\left(\input{figure/r3-lem-planar-k1-4-mf}\right)_n
&=&\overline{\Lambda}_{(1;6,7)}^{[k,1]}\boxtimes \overline{N}_{(8,6,4,3)}^{[1,k]} \boxtimes \overline{N}_{(2,7,5,8)}^{[1,1]}\\
\nonumber
&\simeq&
K\left(
\left(
\begin{array}{c}
\Lambda_{1,(1;6,7)}^{[k,1]}\\
\vspace{0.1cm}\vdots\\
\Lambda_{k+1,(1;6,7)}^{[k,1]}
\end{array}
\right);
\left(
\begin{array}{c}
x_{1,1}-X_{1,(6,7)}^{(k,1)}\\
\vspace{0.1cm}\vdots\\
x_{k+1,1}-X_{k+1,(6,7)}^{(k,1)}
\end{array}
\right)
\right)_{R_{(1,6,7)}^{(1,k,1)}}\\
&&
\nonumber
\hspace{0.5cm}\boxtimes
K\left(
\left(
\begin{array}{c}
A_{1,(8,6,4,3)}^{[1,k]}\\
\vdots\\
\vspace{0.1cm}A_{k,(8,6,4,3)}^{[1,k]}\\
u_{k+1,(8,6,4,3)}^{[1,k]}(x_{1,8}-x_{1,3})
\end{array}
\right);
\left(
\begin{array}{c}
X_{1,(8,6)}^{(1,k)}-X_{1,(4,3)}^{(k,1)}\\
\vdots\\
\vspace{0.1cm}X_{k,(8,6)}^{(1,k)}-X_{k,(4,3)}^{(k,1)}\\
X_{k,(6,3)}^{(k,-1)}
\end{array}
\right)
\right)_{R_{(8,6,4,3)}^{(1,k,k,1)}}\\
&&
\nonumber
\hspace{0.5cm}\boxtimes
K\left(
\left(
\begin{array}{c}
\vspace{0.1cm}A_{1,(2,7,5,8)}^{[1,1]}\\
u_{2,(2,7,5,8)}^{[1,1]}(x_{1,2}-x_{1,8})
\end{array}
\right);
\left(
\begin{array}{c}
\vspace{0.1cm}X_{1,(2,7)}^{(1,1)}-X_{1,(5,8)}^{(1,1)}\\
X_{1,(7,8)}^{(1,-1)}
\end{array}
\right)
\right)_{R_{(2,7,5,8)}^{(1,k,k,1)}}\{-k+1\}\\
&\simeq&
\nonumber
K\left(
\left(
\begin{array}{c}
\vspace{0.1cm}\Lambda_{1,(1;6,7)}^{[k,1]}\\
\vspace{0.1cm}\vdots\\
\vspace{0.1cm}\Lambda_{k+1,(1;6,7)}^{[k,1]}\\
u_{2,(2,7,5,8)}^{[1,1]}(x_{1,2}-x_{1,8})
\end{array}
\right);
\left(
\begin{array}{c}
\vspace{0.1cm}x_{1,1}-X_{1,(6,7)}^{(k,1)}\\
\vspace{0.1cm}\vdots\\
\vspace{0.1cm}x_{k+1,1}-X_{k+1,(6,7)}^{(k,1)}\\
X_{1,(7,8)}^{(1,-1)}
\end{array}
\right)
\right)_{Q_5\left/\langle X_{k,(6,3)}^{(k,-1)}\rangle\right.}\{-k+1\}\\
\label{mf-r3-1-11-3}
&\simeq&
\overline{S}_{(1,2;3,4,5)}^{[k+1,1;1,k,1]}\boxtimes
K\left(C;X_{1,(2,5)}^{(-1,1)}\right)_{Q_5\left/\langle X_{k,(6,3)}^{(k,-1)}\rangle\right.}\{-k+1\}.
\end{eqnarray}
Since the partial matrix factorization $K\left(C;X_{1,(2,5)}^{(-1,1)}\right)_{Q_5\left/\langle X_{k,(6,3)}^{(k,-1)}\rangle\right.}$ is decomposed into
\begin{eqnarray}
\nonumber
&&\xymatrix{
R_{15}
\ar[rr]_(.4){E_{k}(C)}
&&R_{15}\{1-n\}
\ar[rrr]_(.55){E_{k}(X_{1,(2,5)}^{(-1,1)})}
&&&R_{15}}\\
\nonumber
&\simeq&
\bigoplus_{j=0}^{k-1}K\left(C;X_{1,(2,5)}^{(-1,1)}\right)_{R_{(1,2,3,4,5)}^{(k+1,1,1,k,1)}}\{2j\},
\end{eqnarray} 
then the matrix factorization (\ref{mf-r3-1-11-3}) is isomorphic to
\begin{eqnarray}
\label{mf-r3-fin-11-3}
&&
\bigoplus_{j=0}^{k-1}\overline{S}_{(1,2;3,4,5)}^{[k+1,1;1,k,1]}\boxtimes
K\left(C;X_{1,(2,5)}^{(-1,1)}\right)_{R_{(1,2,3,4,5)}^{(k+1,1,1,k,1)}}\{-k+1+2j\}.
\end{eqnarray}
For these decompositions (\ref{mf-r3-fin-00-3}), (\ref{mf-r3-fin-10-3}), (\ref{mf-r3-fin-01-3}) and (\ref{mf-r3-fin-11-3}), the morphisms $\overline{\xi}_{+,1}$, $\overline{\xi}_{+,2}$, $\overline{\xi}_{+,3}$ and $\overline{\xi}_{+,4}$ of the complex (\ref{r3-lem-planar-comp-3}) transform into
\begin{eqnarray}
\nonumber
&&
\overline{\xi}_{+,1}\simeq
\left(
\begin{array}{ccc}
\id_{\overline{S}_{(1,2;3,4,5)}^{[k+1,1;1,k,1]}}\boxtimes(1,x_{1,2}-x_{1,3})
&\mathfrak{0}_{k-1}
&(-1)^{k-1}\id_{\overline{S}_{(1,2;3,4,5)}^{[k+1,1;1,k,1]}}\boxtimes(X_{k,(2,4)}^{(-1,k)},1)\\
&
&
(-1)^{k-2}X_{k-1,(2,4)}^{(-1,k)}\id\\
{}^{t}\mathfrak{0}_{k-1}
&E_{k-1}\left(\id\right)
&\vdots\\
&
&
X_{1,(2,4)}^{(-1,k)}\id
\end{array}
\right),\\
\nonumber
&&
\overline{\xi}_{+,2}\simeq
\left(
\begin{array}{cc}
\id_{\overline{S}_{(1,2;3,4,5)}^{[k+1,1;1,k,1]}}\boxtimes(1,(x_{1,2}-x_{1,3})X_{k+1,(2,4)}^{(-1,k)})&\mathfrak{0}_{k}\\
{}^{t}\mathfrak{0}_{k}&E_{k}\left(\id\right)
\end{array}
\right)
,\\
\nonumber
&&
\overline{\xi}_{+,3}\simeq
\left(
\begin{array}{cc}
\id_{\overline{S}_{(1,2;3,4,5)}^{[k+1,1;1,k,1]}}\boxtimes(1,-X_{k,(2,4)}^{(-1,k)})&\mathfrak{0}_{k-1}\\
{}^{t}\mathfrak{0}_{k-1}&E_{k-1}\left(\id\right)
\end{array}
\right)
,\\
\nonumber
&&
\overline{\xi}_{+,4}\simeq
-\left(
\begin{array}{cc}
&(-1)^{k-1}X_{k,(2,4)}^{(-1,k)}\id\\
E_{k}\left(\id\right)&\vdots\\
&X_{1,(2,4)}^{(-1,k)}\id
\end{array}
\right).
\end{eqnarray}
Then, the complex (\ref{r3-lem-planar-comp-3}) is isomorphic, in $\k^b(\HMF^{gr}_{R_{(1,2,3,4,5)}^{(k+1,1,1,k,1)},\omega_4})$ $(\omega_3=F_{k+1}(\mathbb{X}^{(k+1)}_{(1)})+F_{1}(\mathbb{X}^{(1)}_{(2)})-F_{1}(\mathbb{X}^{(1)}_{(3)})-F_{k}(\mathbb{X}^{(k)}_{(4)})-F_{1}(\mathbb{X}^{(1)}_{(5)}))$, to
\begin{equation}
\nonumber
\xymatrix{
-k-1\ar@{.}[d]&&&-k\ar@{.}[d]\\
{\overline{M}_3
\{(k+1)n\}
\left<k+1\right>
}
\ar[rrr]^{\id_{\overline{S}}\boxtimes(1,x_{1,2}-x_{1,3})}&&&
{\overline{M}_4
\{(k+1)n-1\}
\left<k+1\right>
}
,}
\end{equation}
where
\begin{eqnarray}
\label{mf-planar-1-2}
\overline{M}_3&=&
K\left(
\left(
\begin{array}{c}
\vspace{0.1cm}C_1\\
\vspace{0.1cm}\vdots\\
\vspace{0.1cm}C_{k+1}\\
\frac{C}{(x_{1,2}-x_{1,3})X_{k,(2,4)}^{(-1,k)}}
\end{array}
\right);
\left(
\begin{array}{c}
x_{1,1}-X_{1,(2,3,4,5)}^{(-1,1,k,1)}\\
\vdots\\
x_{k+1,1}-X_{k+1,(2,3,4,5)}^{(-1,1,k,1)}\\
(x_{1,2}-x_{1,3})X_{k,(2,4)}^{(-1,k)}X_{1,(2,5)}^{(-1,k)}
\end{array}
\right)
\right)_{R_{(1,2,3,4,5)}^{(k+1,1,1,k,1)}}\hspace{-1cm}\{-k-1\},\\
\label{mf-planar-2-2}
\overline{M}_4&=&
K\left(
\left(
\begin{array}{c}
\vspace{0.1cm}C_1\\
\vspace{0.1cm}\vdots\\
\vspace{0.1cm}C_{k+1}\\
\frac{C}{X_{k,(2,4)}^{(-1,k)}}
\end{array}
\right);
\left(
\begin{array}{c}
x_{1,1}-X_{1,(2,3,4,5)}^{(-1,1,k,1)}\\
\vdots\\
x_{k+1,1}-X_{k+1,(2,3,4,5)}^{(-1,1,k,1)}\\
X_{k,(2,4)}^{(-1,k)}X_{1,(2,5)}^{(-1,k)}
\end{array}
\right)
\right)_{R_{(1,2,3,4,5)}^{(k+1,1,1,k,1)}}\hspace{-1cm}\{-k\}.
\end{eqnarray}
By the way, we have
\begin{eqnarray}
\nonumber
&&\c\left(\input{figure/r3-lem-plus-k1-2-mf}\right)_n=\\
&&\xymatrix{
-k-1\ar@{.}[d]&&&-k\ar@{.}[d]\\
{
\c\left(\input{figure/r3-lem-planar-k1-5-mf}\right)_n
\begin{array}{c}
\{(k+1)n\}\\
\left<k+1\right>
\end{array}
}
\ar[rrr]^{\chi_{+,(2,1,6,3)}^{[1,k+1]}\boxtimes\id_{\overline{\Lambda}_{(6;4,5)}^{[k,1]}}}&&&
{
\c\left(\input{figure/r3-lem-planar-k1-6-mf}\right)_n
\begin{array}{c}
\{(k+1)n-1\}\\
\left<k+1\right>
\end{array}
}
},
\end{eqnarray}
\begin{eqnarray}
\nonumber
&&\c\left(\input{figure/r3-lem-planar-k1-5-mf}\right)_n\\
\label{mf-planar-3-2}
&&\simeq
K\left(
\left(
\begin{array}{c}
A_{1,(2,1,6,3)}^{[1,k+1]}\\
\vspace{0.1cm}\vdots\\
A_{k+1,(2,1,6,3)}^{[1,k+1]}\\
\vspace{0.2cm}u_{k+2,(2,1,6,3)}^{[1,k+1]}\\
\end{array}
\right);
\left(
\begin{array}{c}
X_{1,(1,2)}^{(k+1,1)}-X_{1,(3,4,5)}^{(1,k,1)}\\
\vspace{0.1cm}\vdots\\
X_{k+1,(1,2)}^{(k+1,1)}-X_{k+1,(3,4,5)}^{(1,k,1)}\\
(x_{1,2}-x_{1,3})X_{k+1,(1,3)}^{(k+1,-1)}
\end{array}
\right)
\right)_{Q_{10}}\{-k-1\},\\
\nonumber
&&\c\left(\input{figure/r3-lem-planar-k1-6-mf}\right)_n\\
\label{mf-planar-4-2}
&&\simeq
K\left(
\left(
\begin{array}{c}
A_{1,(2,1,6,3)}^{[1,k+1]}\\
\vspace{0.1cm}\vdots\\
A_{k+1,(2,1,6,3)}^{[1,k+1]}\\
u_{k+2,(2,1,6,3)}^{[1,k+1]}(x_{1,2}-x_{1,3})
\end{array}
\right);
\left(
\begin{array}{c}
X_{1,(1,2)}^{(k+1,1)}-X_{1,(3,4,5)}^{(1,k,1)}\\
\vspace{0.1cm}\vdots\\
X_{k+1,(1,2)}^{(k+1,1)}-X_{k+1,(3,4,5)}^{(1,k,1)}\\
X_{k+1,(1,3)}^{(k+1,-1)}
\end{array}
\right)
\right)_{Q_{10}}\{-k\},
\end{eqnarray}
where
\begin{eqnarray}
\nonumber	
Q_{10}&:=&R_{(1,2,3,4,5,6)}^{(k+1,1,1,k,1,k+1)}\left/\left<x_{1,6}-X_{1,(4,5)}^{(k,1)},\ldots,x_{k+1,6}-X_{k+1,(4,5)}^{(k,1)}\right>\right. \\
\nonumber
(&\simeq&R_{(1,2,3,4,5)}^{(k+1,1,1,k,1)}).
\end{eqnarray}
The right-hand side sequences 
$\left(
x_{1,1}-X_{1,(2,3,4,5)}^{(-1,1,k,1)},\ldots,x_{1,k+1}-X_{k+1,(2,3,4,5)}^{(-1,1,k,1)},X_{k,(2,4)}^{(-1,k)}X_{1,(2,5)}^{(-1,k)}
\right)$ of the matrix factorization $\overline{M}_4$ and 
$\left(
X_{1,(1,2)}^{(k+1,1)}-X_{1,(3,4,5)}^{(1,k,1)},\ldots,X_{k+1,(1,2)}^{(k+1,1)}-X_{k+1,(3,4,5)}^{(1,k,1)},X_{k+1,(1,3)}^{(k+1,-1)}
\right)$
of the matrix factorization (\ref{mf-planar-4-2}) transform to each other by a linear transformation over $R_{(1,2,3,4,5)}^{(k+1,1,1,k,1)}$.
Then, by Proposition \ref{equiv} and Theorem \ref{reg-eq}, we have
\begin{equation}
\label{isom3}
\overline{M}_4\simeq\c\left(\input{figure/r3-lem-planar-k1-6-mf}\right)_n.
\end{equation}
The sequences 
$\left(
x_{1,1}-X_{1,(2,3,4,5)}^{(-1,1,k,1)},\ldots,x_{1,k+1}-X_{k+1,(2,3,4,5)}^{(-1,1,k,1)},-X_{k+1,(2,3,4)}^{(-1,1,k)}X_{1,(2,5)}^{(-1,k)}
\right)$ of the matrix factorization $\overline{M}_3$ and 
$\left(
X_{1,(1,2)}^{(k+1,1)}-X_{1,(3,4,5)}^{(1,k,1)},\ldots,X_{k+1,(1,2)}^{(k+1,1)}-X_{k+1,(3,4,5)}^{(1,k,1)},(x_{1,2}-x_{1,3})X_{k+1,(1,3)}^{(k+1,-1)}
\right)$
of the matrix factorization (\ref{mf-planar-3-2}) also transform to each other.
We have
\begin{equation}
\label{isom4}
\overline{M}_3\simeq\c\left(\input{figure/r3-lem-planar-k1-5-mf}\right)_n.
\end{equation}
Thus, in $\k^b(\HMF^{gr}_{R_{(1,2,3,4,5)}^{(k+1,1,1,k,1)},\omega_4})$, we obtain
\begin{equation}
\nonumber
\c\left(\input{figure/r3-lem-plus-k1-1}\right)_n\simeq\c\left(\input{figure/r3-lem-plus-k1-2}\right)_n.
\end{equation}
\end{proof}
\begin{proof}[{\bf Proof of Proposition \ref{prop-r3} (4)}]
The complex for the diagram $\input{figure/r3-lem-minus-k1-1-text}$ is described as a complex of factorizations of $\k^b(\HMF^{gr}_{R_{(1,2,3,4,5)}^{(k+1,1,1,k,1)},\omega_4})$,
\begin{eqnarray}
\label{r3-lem-planar-comp-4}
&&\hspace{1.5cm}
\c\left(\input{figure/r3-lem-minus-k1-1-mf}\right)_n=\\
\nonumber
&&
\xymatrix{
k-1\ar@{.}[d]&k\ar@{.}[d]&k+1\ar@{.}[d]\\
{\c\left(\input{figure/r3-lem-planar-k1-4-mf}\right)_n
\hspace{-0.5cm}
\begin{array}{c}
{}_{\{-(k+1)n+2\}}\\
\left<k+1\right>
\end{array}}
\ar[r]^{\left(\hspace{-0.2cm}
\begin{array}{c}
\overline{\xi}_{-,1}\\
\overline{\xi}_{-,2}
\end{array}
\hspace{-0.2cm}\right)}&
{
\begin{array}{c}
\c\left(\input{figure/r3-lem-planar-k1-2-mf}\right)_n
\hspace{-0.5cm}
\begin{array}{c}
{}_{\{-(k+1)n+1\}}\\
\left<k+1\right>
\end{array}
\\
\bigoplus\\
\c\left(\input{figure/r3-lem-planar-k1-3-mf}\right)_n
\hspace{-0.5cm}
\begin{array}{c}
{}_{\{-(k+1)n+1\}}\\
\left<k+1\right>
\end{array}
\end{array}
}
\ar[r]^{\txt{$(\overline{\xi}_{-,3},\overline{\xi}_{-,4})$}}&
{\c\left(\input{figure/r3-lem-planar-k1-1-mf}\right)_n
\hspace{-0.5cm}
\begin{array}{c}
{}_{\{-(k+1)n\}}\\
\left<k+1\right>
\end{array}
}.
}
\end{eqnarray}
By the discussion of Proof of lemma \ref{prop-r3} (3), we also have
\begin{eqnarray}
\label{mf-r3-1-00-4}
\c\left(\input{figure/r3-lem-planar-k1-4-mf}\right)_n
&\simeq&\overline{S}_{(1,2;3,4,5)}^{[k+1,1;1,k,1]}\boxtimes
K\left(C;X_{1,(2,5)}^{(-1,1)}\right)_{Q_5\left/\langle X_{k,(6,3)}^{(k,-1)}\rangle\right.}\{-k+1\},\\
\label{mf-r3-1-10-4}
\c\left(\input{figure/r3-lem-planar-k1-2-mf}\right)_n
&\simeq&\overline{S}_{(1,2;3,4,5)}^{[k+1,1;1,k,1]}\boxtimes
K\left(C_0;(x_{1,2}-x_{1,8})X_{1,(2,5)}^{(-1,k)}\right)_{Q_5\left/\langle X_{k,(6,3)}^{(k,-1)}\rangle\right.}\{-k\},\\
\label{mf-r3-1-01-4}
\c\left(\input{figure/r3-lem-planar-k1-3-mf}\right)_n
&\simeq&\overline{S}_{(1,2;3,4,5)}^{[k+1,1;1,k,1]}\boxtimes
K\left(C;X_{k,(2,5)}^{(-1,1)}\right)_{Q_5\left/\langle (x_{1,8}-x_{1,3})X_{k,(6,3)}^{(k,-1)}\rangle\right.}\{-k\},\\
\label{mf-r3-1-11-4}
\c\left(\input{figure/r3-lem-planar-k1-1-mf}\right)_n
&\simeq&\overline{S}_{(1,2;3,4,5)}^{[k+1,1;1,k,1]}\boxtimes
K\left(C_0;(x_{1,2}-x_{1,8})X_{1,(2,5)}^{(-1,1)}\right)_{Q_5\left/\langle (x_{1,8}-x_{1,3})X_{k,(6,3)}^{(k,-1)}\rangle\right.}\{-k-1\},
\end{eqnarray}
The partial matrix factorization $K\left(C;X_{1,(2,5)}^{(-1,1)}\right)_{Q_5\left/\langle X_{k,(6,3)}^{(k,-1)}\rangle\right.}$ of (\ref{mf-r3-1-00-4}) is isomorphic to
\begin{eqnarray}
\nonumber
&&\xymatrix{R_{15}
\ar[rr]_(.42){E_{k}(C)}
&&R_{15}\{1-n\}
\ar[rr]_(.55){E_{k}(X_{(2,5)}^{(-1,1)})}
&&R_{15}}\\
\nonumber
&\simeq&
\bigoplus_{j=0}^{k-1}K\left(C;X_{1,(2,5)}^{(-1,1)}\right)_{R_{(1,2,3,4,5)}^{(k+1,1,1,k,1)}}\{2j\}.
\end{eqnarray}
The partial matrix factorization $K\left(C_0;(x_{1,2}-x_{1,8})X_{1,(2,5)}^{(-1,k)}\right)_{Q_5\left/\langle X_{k,(6,3)}^{(k,-1)}\rangle\right.}$ of (\ref{mf-r3-1-10-4}) is isomorphic to
\begin{eqnarray}
\nonumber
&&\xymatrix{R_{14}\ar[rrrr]_(.45){
\left(\hspace{-0.2cm}
\begin{array}{cc}
0&E_{k-1}(C)\\
\frac{C}{X_{k,(2,4)}^{(-1,k)}}&0
\end{array}\hspace{-0.2cm}
\right)
}
&&&&R_{15}\{3-n\}\ar[rrrrrr]_(.55){
\left(
\begin{array}{cc}
0&X_{k,(2,4)}^{(-1,k)}X_{1,(2,5)}^{(-1,1)}\\
E_{k-1}(X_{1,(2,5)}^{(-1,1)})&0
\end{array}
\right)
}
&&&&&&R_{14}}\\
\nonumber
&\simeq&
K\left(\frac{C}{X_{k,(2,4)}^{(-1,k)}};X_{k,(2,4)}^{(-1,k)}X_{1,(2,5)}^{(-1,1)}\right)_{R_{(1,2,3,4,5)}^{(k+1,1,1,k,1)}}
\oplus \bigoplus_{j=1}^{k-1}K\left(C;X_{1,(2,5)}^{(-1,1)}\right)_{R_{(1,2,3,4,5)}^{(k+1,1,1,k,1)}}\{2j\}.
\end{eqnarray} 
We consider an isomorphism of $Q_5\left/\langle (x_{1,8}-x_{1,3})X_{k,(6,3)}^{(k,-1)}\rangle\right.$ to be
\begin{eqnarray}
\nonumber
R_{16}
&:=&
R_{(1,2,3,4,5)}^{(k+1,1,1,k,1)}
\oplus(x_{1,8}-x_{1,3})R_{(1,2,3,4,5)}^{(k+1,1,1,k,1)}
\oplus\ldots
\oplus x_{1,8}^{k-2}(x_{1,8}-x_{1,3})R_{(1,2,3,4,5)}^{(k+1,1,1,k,1)}
\oplus X_{k,(2,3,4,8)}^{(-1,1,k,-1)}R_{(1,2,3,4,5)}^{(k+1,1,1,k,1)},\\
\nonumber
R_{17}
&:=&
R_{(1,2,3,4,5)}^{(k+1,1,1,k,1)}
\oplus(x_{1,2}-x_{1,8})R_{(1,2,3,4,5)}^{(k+1,1,1,k,1)}\\
\nonumber
&&\hspace{1.5cm}
\oplus(x_{1,8}-x_{1,3})(x_{1,2}-x_{1,8})R_{(1,2,3,4,5)}^{(k+1,1,1,k,1)}
\oplus\ldots
\oplus x_{1,8}^{k-2}(x_{1,8}-x_{1,3})(x_{1,2}-x_{1,8})R_{(1,2,3,4,5)}^{(k+1,1,1,k,1)}.
\end{eqnarray}
Then, the partial matrix factorization $K\left(C;X_{k,(2,5)}^{(-1,1)}\right)_{Q_5\left/\langle (x_{1,8}-x_{1,3})X_{k,(6,3)}^{(k,-1)}\rangle\right.}$ of (\ref{mf-r3-1-01-4}) is isomorphic to
\begin{eqnarray}
\nonumber
&&\xymatrix{R_{16}
\ar[rr]_(.42){E_{k+1}(C)}
&&R_{16}\{1-n\}
\ar[rr]_(.55){E_{k+1}(X_{(2,5)}^{(-1,1)})}
&&R_{16}}\\
\nonumber
&\simeq&
\bigoplus_{j=0}^{k}K\left(C;X_{1,(2,5)}^{(-1,1)}\right)_{R_{(1,2,3,4,5)}^{(k+1,1,1,k,1)}}\{2j\}.
\end{eqnarray}
The partial matrix factorization $K\left(C_0;(x_{1,2}-x_{1,8})X_{1,(2,5)}^{(-1,1)}\right)_{Q_5\left/\langle (x_{1,8}-x_{1,3})X_{k,(6,3)}^{(k,-1)}\rangle\right.}$ of (\ref{mf-r3-1-11-4}) is isomorphic to
\begin{eqnarray}
\nonumber
&&\xymatrix{R_{17}\ar[rrrr]_(.45){
\left(\hspace{-0.2cm}
\begin{array}{cc}
0&E_{k}(C)\\
\frac{C}{(x_{1,2}-x_{1,3})X_{k,(2,4)}^{(-1,k)}}&0
\end{array}\hspace{-0.2cm}
\right)
}
&&&&R_{16}\{3-n\}\ar[rrrrrr]_(.55){
\left(
\begin{array}{cc}
0&(x_{1,2}-x_{1,3})X_{k,(2,4)}^{(-1,k)}X_{1,(2,5)}^{(-1,1)}\\
E_{k}(X_{1,(2,5)}^{(-1,1)})&0
\end{array}
\right)
}
&&&&&&R_{17}}\\
\nonumber
&\simeq&
K\left(\frac{C}{(x_{1,2}-x_{1,3})X_{k,(2,4)}^{(-1,k)}};(x_{1,2}-x_{1,3})X_{k,(2,4)}^{(-1,k)}X_{1,(2,5)}^{(-1,1)}\right)_{R_{(1,2,3,4,5)}^{(k+1,1,1,k,1)}}
\oplus \bigoplus_{j=1}^{k}K\left(C;X_{1,(2,5)}^{(-1,1)}\right)_{R_{(1,2,3,4,5)}^{(k+1,1,1,k,1)}}\{2j\}.
\end{eqnarray} 
Then, the matrix factorizations (\ref{mf-r3-1-00-4}), (\ref{mf-r3-1-10-4}), (\ref{mf-r3-1-01-4}) and (\ref{mf-r3-1-11-4}) are decomposed as follows,
\begin{eqnarray}
\nonumber
\c\left(\input{figure/r3-lem-planar-k1-4-mf}\right)_n
&\simeq&
\bigoplus_{j=0}^{k-1}\overline{S}_{(1,2;3,4,5)}^{[k+1,1;1,k,1]}\boxtimes
K\left(
C;X_{1,(2,5)}^{(-1,1)}
\right)_{R_{(1,2,3,4,5)}^{(k+1,1,1,k,1)}}\{-k+1+2j\},
\\
\nonumber
\c\left(\input{figure/r3-lem-planar-k1-2-mf}\right)_n
&\simeq&
\overline{S}_{(1,2;3,4,5)}^{[k+1,1;1,k,1]}\boxtimes
K\left(
\frac{C}{X_{k,(2,4)}^{(-1,k)}};X_{k,(2,4)}^{(-1,k)}X_{1,(2,5)}^{(-1,1)}
\right)_{R_{(1,2,3,4,5)}^{(k+1,1,1,k,1)}}\{-k\}
\\
\nonumber
&&
\hspace{.5cm}\oplus
\bigoplus_{j=1}^{k-1}
\overline{S}_{(1,2;3,4,5)}^{[k+1,1;1,k,1]}\boxtimes
K\left(
C;X_{1,(2,5)}^{(-1,1)}\right
)_{R_{(1,2,3,4,5)}^{(k+1,1,1,k,1)}}\{-k+2j\},
\end{eqnarray}
\begin{eqnarray}
\nonumber
\c\left(\input{figure/r3-lem-planar-k1-4-mf}\right)_n
&\simeq&
\bigoplus_{j=0}^{k}
\overline{S}_{(1,2;3,4,5)}^{[k+1,1;1,k,1]}\boxtimes
K\left(
C;X_{1,(2,5)}^{(-1,1)}
\right)_{R_{(1,2,3,4,5)}^{(k+1,1,1,k,1)}}\{-k+2j\},
\\
\nonumber
\c\left(\input{figure/r3-lem-planar-k1-4-mf}\right)_n
&\simeq&
\overline{S}_{(1,2;3,4,5)}^{[k+1,1;1,k,1]}\boxtimes
K\left(
\frac{C}{(x_{1,2}-x_{1,3})X_{k,(2,4)}^{(-1,k)}};(x_{1,2}-x_{1,3})X_{k,(2,4)}^{(-1,k)}X_{1,(2,5)}^{(-1,1)}
\right)_{R_{(1,2,3,4,5)}^{(k+1,1,1,k,1)}}\hspace{-1cm}\{-k-1\}\\
\nonumber
&&
\hspace{.5cm}\oplus
\bigoplus_{j=1}^{k}
\overline{S}_{(1,2;3,4,5)}^{[k+1,1;1,k,1]}\boxtimes
K\left(
C;X_{1,(2,5)}^{(-1,1)}
\right)_{R_{(1,2,3,4,5)}^{(k+1,1,1,k,1)}}\{-k-1+2j\}.
\end{eqnarray}
For these decompositions, the morphisms $\overline{\xi}_{-,1}$, $\overline{\xi}_{-,2}$, $\overline{\xi}_{-,3}$ and $\overline{\xi}_{-,4}$ of the complex (\ref{r3-lem-planar-comp-4}) transform into
\begin{eqnarray}
\nonumber
\overline{\xi}_{-,1}&\simeq&
\left(
\begin{array}{cc}
\mathfrak{0}_{k-1}&\id_{\overline{S}_{(1,2;3,4,5)}^{[k+1,1;1,k,1]}}\boxtimes(X_{k,(2,4)}^{(-1,k)},1)\\
E_{k-1}(\id)&{}^{t}\mathfrak{0}_{k-1}
\end{array}
\right),\\
\nonumber
\overline{\xi}_{-,2}&\simeq&
\left(
\begin{array}{cc}
\mathfrak{0}_{k-1}&-X_{k,(2,4)}^{(-1,k)}\id\\
E_{k-1}(\id)&{}^{t}\mathfrak{0}_{k-1}\\
\mathfrak{0}_{k-1}&\id
\end{array}
\right),\\
\nonumber
\overline{\xi}_{-,3}&\simeq&
\left(
\begin{array}{cc}
\id_{\overline{S}_{(1,2;3,4,5)}^{[k+1,1;1,k,1]}}\boxtimes(x_{1,2}-x_{1,3},1)&\mathfrak{0}_{k}\\
\id_{\overline{S}_{(1,2;3,4,5)}^{[k+1,1;1,k,1]}}\boxtimes(-1,-X_{k,(2,4)}^{(-1,k)})&\mathfrak{0}_{k}\\
{}^{t}\mathfrak{0}_{k-1}&E_{k-1}(\id)
\end{array}
\right),\\
\nonumber
\overline{\xi}_{-,4}&\simeq&
-\left(
\begin{array}{cc}
\mathfrak{0}_{k}&\id_{\overline{S}_{(1,2;3,4,5)}^{[k+1,1;1,k,1]}}\boxtimes((x_{1,2}-x_{1,3})X_{k,(2,4)}^{(-1,k)},1)\\
E_{k}(\id)&{}^{t}\mathfrak{0}_{k}
\end{array}
\right).
\end{eqnarray}
Then, the complex (\ref{r3-lem-planar-comp-4}) is isomorphic, in $\k^b(\HMF^{gr}_{R_{(1,2,3,4,5)}^{(k+1,1,1,k,1)},\omega_4})$, to
\begin{equation}
\nonumber
\xymatrix{
k\ar@{.}[d]&&&k+1\ar@{.}[d]\\
{\overline{M}_4
\{1-(k+1)n\}
\left<k+1\right>
}
\ar[rrr]^{\id_{\overline{S}}\boxtimes(x_{1,2}-x_{1,3},1)}&&&
{\overline{M}_3
\{-(k+1)n\}
\left<k+1\right>
}
.}
\end{equation}
Since we have
\begin{equation}
\nonumber
\overline{M}_3\simeq\c\left(\input{figure/r3-lem-planar-k1-5-mf}\right)_n,
\overline{M}_4\simeq\c\left(\input{figure/r3-lem-planar-k1-6-mf}\right)_n,
\end{equation}
thus we obtain
\begin{equation}
\nonumber
\c\left(\input{figure/r3-lem-minus-k1-1}\right)_n\simeq\c\left(\input{figure/r3-lem-minus-k1-2}\right)_n.
\end{equation}
\end{proof}
\begin{proof}[{\bf Proof of Proposition \ref{prop-r3} (5)}]
The complex for the diagram $\input{figure/r3-lem-plus-k1-3-text}$ is described as a complex of factorizations of $\k^b(\HMF^{gr}_{R_{(1,2,3,4,5)}^{(k,1,1,1,k+1)},\omega_5})$,
\begin{eqnarray}
\label{r3-lem-planar-comp-5}
&&\hspace{1.5cm}
\c\left(\input{figure/r3-lem-plus-k1-3-mf}\right)_n=\\
\nonumber
&&
\xymatrix{
-k-1\ar@{.}[d]&-k\ar@{.}[d]&-k+1\ar@{.}[d]\\
{\c\left(\input{figure/r3-lem-planar-k1-7-mf}\right)_n
\hspace{-0.5cm}
\begin{array}{c}
{}_{\{(k+1)n\}}\\
\left<k+1\right>
\end{array}}
\ar[r]^{\left(\hspace{-0.2cm}
\begin{array}{c}
\overline{\sigma}_{+,1}\\
\overline{\sigma}_{+,2}
\end{array}
\hspace{-0.2cm}\right)}&
{
\begin{array}{c}
\c\left(\input{figure/r3-lem-planar-k1-8-mf}\right)_n
\hspace{-0.5cm}
\begin{array}{c}
{}_{\{(k+1)n-1\}}\\
\left<k+1\right>
\end{array}
\\
\bigoplus\\
\c\left(\input{figure/r3-lem-planar-k1-9-mf}\right)_n
\hspace{-0.5cm}
\begin{array}{c}
{}_{\{(k+1)n-1\}}\\
\left<k+1\right>
\end{array}
\end{array}
}
\ar[r]^{\txt{$(\overline{\sigma}_{+,3},\overline{\sigma}_{+,4})$}}&
{\c\left(\input{figure/r3-lem-planar-k1-10-mf}\right)_n
\hspace{-0.5cm}
\begin{array}{c}
{}_{\{(k+1)n-2\}}\\
\left<k+1\right>
\end{array}
}.
}
\end{eqnarray}
By Corollary \ref{cor2-11}, we have
\begin{eqnarray}
\label{mf-r3-1-00-5}
&&
\c\left(\input{figure/r3-lem-planar-k1-7-mf}\right)_n\simeq
K\left(
\left(
\begin{array}{c}
V_{1,(7,8;5)}^{[k,1]}\\
\vdots\\
V_{k+1,(7,8;5)}^{[k,1]}\\
u_{k+1,(6,1,7,4)}^{[1,k]}
\end{array}
\right)
;
\left(
\begin{array}{c}
X_{1,(7,8)}^{(k,1)}-x_{1,5}\\
\vdots\\
X_{k+1,(7,8)}^{(k,1)}-x_{k+1,5}\\
(x_{1,6}-x_{1,4})X_{k,(1,4)}^{(k,-1)}
\end{array}
\right)
\right)_{Q_{6}\left/\langle (x_{1,3}-x_{1,6})X_{1,(2,6)}^{(1,-1)}\rangle\right.}\{-2k-1\},
\end{eqnarray}
where $Q_{6}=R_{(1,2,3,4,5,6,7,8)}^{(k,1,1,1,k+1,1,k,1)}\left/\left<X_{1,(6,1)}^{(1,k)}-X_{1,(7,4)}^{(k,1)},\ldots,X_{k,(6,1)}^{(1,k)}-X_{k,(7,4)}^{(k,1)},X_{1,(3,2)}^{(1,1)}-X_{1,(8,6)}^{(1,1)}\right>\right.$. The quotient $Q_{6}\left/\langle (x_{1,3}-x_{1,6})X_{1,(2,6)}^{(1,-1)}\rangle\right.$ has equations
\begin{eqnarray}
\nonumber
&&x_{j,7}=X_{j,(1,4,6)}^{(k,-1,1)}\hspace{1cm}(1\leq j \leq k),\\
\nonumber
&&x_{1,8}=X_{1,(2,3,6)}^{(1,1,-1)},\\
\nonumber
&&(x_{1,3}-x_{1,6})X_{1,(2,6)}^{(1,-1)}=0.
\end{eqnarray}
Then, the matrix factorization (\ref{mf-r3-1-00-5}) is isomorphic to
\begin{equation}
\nonumber
K\left(
\left(
\begin{array}{c}
V_{1,(7,8;5)}^{[k,1]}\\
\vdots\\
V_{k+1,(7,8;5)}^{[k,1]}\\
u_{k+1,(6,1,7,4)}^{[1,k]}-V_{k+1,(7,8;5)}^{[k,1]}
\end{array}
\right)
;
\left(
\begin{array}{c}
X_{1,(1,2,3,4)}^{(k,1,1,-1)}-x_{1,5}\\
\vdots\\
X_{k+1,(1,2,3,4)}^{(k,1,1,-1)}-x_{k+1,5}\\
(x_{1,6}-x_{1,4})X_{k,(1,4)}^{(k,-1)}
\end{array}
\right)
\right)_{Q_{6}\left/\langle (x_{1,3}-x_{1,6})X_{1,(2,6)}^{(1,-1)}\rangle\right.}\{-2k-1\}.
\end{equation}
By Corollary \ref{induce-sq1}, there exist polynomials $D_1,\ldots,D_{k+1}$ $\in R_{(1,2,3,4,5)}^{(k,1,1,1,k+1)}$ and $D_0 \in Q_{6}\left/\langle (x_{1,3}-x_{1,6})X_{1,(2,6)}^{(1,-1)}\rangle\right.$ such that $D_0(x_{1,6}-x_{1,4})\equiv D \in R_{(1,2,3,4,5)}^{(k,1,1,1,k+1)}$ and we have an isomorphism to the above factorization
\begin{equation}
\label{mf-r3-2-00-5}
\overline{S}_{(1,2,3;4,5)}^{[k,1,1;1,k+1]}\boxtimes
K\left(D_0;(x_{1,6}-x_{1,4})X_{k,(1,4)}^{(k,-1)}\right)_{Q_{6}\left/\langle (x_{1,3}-x_{1,6})X_{1,(2,6)}^{(1,-1)}\rangle\right.}\{-2k-1\},
\end{equation}
where 
\begin{equation}
\nonumber
\overline{S}_{(1,2,3;4,5)}^{[k,1,1;1,k+1]}
:=K\left(
\left(
\begin{array}{c}
\vspace{.1cm}D_1\\
\vspace{.1cm}\vdots\\
\vspace{.2cm}D_{k+1}\\
D_0
\end{array}
\right)
;
\left(
\begin{array}{c}
X_{1,(1,2,3,4)}^{(k,1,1,-1)}-x_{1,5}\\
\vdots\\
X_{k+1,(1,2,3,4)}^{(k,1,1,-1)}-x_{k+1,5}\\
(x_{1,6}-x_{1,4})X_{k,(1,4)}^{(k,-1)}
\end{array}
\right)
\right)_{R_{(1,2,3,4,5)}^{(k,1,1,1,k+1)}}
\end{equation}
We consider isomorphisms of $Q_{6}\left/\langle (x_{1,3}-x_{1,6})X_{1,(2,6)}^{(1,-1)}\rangle\right.$ as $R_{(1,2,3,4,5)}^{(k,1,1,1,k+1)}$-module
\begin{eqnarray}
\nonumber
&&R_{18}:=R_{(1,2,3,4,5)}^{(k,1,1,1,k+1)}\oplus(x_{1,6}-x_{1,4})R_{(1,2,3,4,5)}^{(k,1,1,1,k+1)},\\
\nonumber
&&R_{19}:=R_{(1,2,3,4,5)}^{(k,1,1,1,k+1)}\oplus(x_{1,2}+x_{1,3}-x_{1,4}-x_{1,6})R_{(1,2,3,4,5)}^{(k,1,1,1,k+1)}.
\end{eqnarray}
The partial factorization $K\left(D_0;(x_{1,6}-x_{1,4})X_{k,(1,4)}^{(k,-1)}\right)_{Q_{6}\left/\langle (x_{1,3}-x_{1,6})X_{1,(2,6)}^{(1,-1)}\rangle\right.}$ is isomorphic to
\begin{eqnarray}
\nonumber
&&
\xymatrix{
R_{18}\ar[rrrrr]_(.45){
\left(
\begin{array}{cc}
0&D\\
\frac{D}{(x_{1,4}-x_{1,3})(x_{1,4}-x_{1,3})}&0
\end{array}
\right)}
&&&&&
R_{19}\{2k+1-n\}\ar[rrrrr]_(.55){
\left(
\begin{array}{cc}
0&(x_{1,4}-x_{1,3})X_{k,(1,2,4)}^{(k,1,-1)}\\
X_{k,(1,4)}^{(k,-1)}&0
\end{array}
\right)}
&&&&&
R_{18}.
}
\end{eqnarray}
Then, the matrix factorization (\ref{mf-r3-2-00-5}) is isomorphic to
\begin{eqnarray}
\label{mf-r3-fin-00-5}
&&\overline{S}_{(1,2,3;4,5)}^{[k,1,1;1,k+1]}\boxtimes
K\left(
\frac{D}{(x_{1,4}-x_{1,3})(x_{1,4}-x_{1,3})};(x_{1,4}-x_{1,3})X_{k,(1,2,4)}^{(k,1,-1)}
\right)_{R_{(1,2,3,4,5)}^{(k,1,1,1,k+1)}}\{-2k-1\}\\
\nonumber
&&\oplus
\overline{S}_{(1,2,3;4,5)}^{[k,1,1;1,k+1]}\boxtimes
K\left(
D;X_{k,(1,4)}^{(k,-1)}
\right)_{R_{(1,2,3,4,5)}^{(k,1,1,1,k+1)}}\{-2k+1\}.
\end{eqnarray}
By a similar discussion, we obtain
\begin{eqnarray}
\label{mf-r3-1-10-5}
&&
\c\left(\input{figure/r3-lem-planar-k1-8-mf}\right)_n\simeq
\overline{S}_{(1,2,3;4,5)}^{[k,1,1;1,k+1]}\boxtimes
K\left(D;X_{k,(1,4)}^{(k,-1)}\right)_{Q_{6}\left/\langle (x_{1,3}-x_{1,6})X_{1,(2,6)}^{(1,-1)}\rangle\right.}\{-2k\},
\\
\label{mf-r3-1-01-5}
&&
\c\left(\input{figure/r3-lem-planar-k1-9-mf}\right)_n\simeq
\overline{S}_{(1,2,3;4,5)}^{[k,1,1;1,k+1]}\boxtimes
K\left(D_0;(x_{1,6}-x_{1,4})X_{k,(1,4)}^{(k,-1)}\right)_{Q_{6}\left/\langle X_{1,(2,6)}^{(1,-1)}\rangle\right.}\{-2k\},
\\
\label{mf-r3-1-11-5}
&&
\c\left(\input{figure/r3-lem-planar-k1-10-mf}\right)_n\simeq
\overline{S}_{(1,2,3;4,5)}^{[k,1,1;1,k+1]}\boxtimes
K\left(D;X_{k,(1,4)}^{(k,-1)}\right)_{Q_{6}\left/\langle X_{1,(2,6)}^{(1,-1)}\rangle\right.}\{-2k+1\}.
\end{eqnarray}
\indent
The partial factorization $K\left(D;X_{k,(1,4)}^{(k,-1)}\right)_{Q_{6}\left/\langle (x_{1,3}-x_{1,6})X_{1,(2,6)}^{(1,-1)}\rangle\right.}$ of (\ref{mf-r3-1-10-5}) is isomorphic to
\begin{eqnarray}
\nonumber
\xymatrix{
R_{18}\ar[rr]_(.45){
\left(
\begin{array}{cc}
D&0\\
0&D
\end{array}
\right)}
&&
R_{18}\{2k-1-n\}\ar[rrrr]_(.55){
\left(
\begin{array}{cc}
X_{k,(1,4)}^{(k,-1)}&0\\
0&X_{k,(1,4)}^{(k,-1)}
\end{array}
\right)}
&&&&
R_{18}.
}
\end{eqnarray}
Then, the matrix factorization (\ref{mf-r3-1-10-5}) is decomposed into
\begin{equation}
\label{mf-r3-fin-10-5}
\overline{S}_{(1,2,3;4,5)}^{[k,1,1;1,k+1]}\boxtimes
K\left(D;X_{k,(1,4)}^{(k,-1)}\right)_{R_{(1,2,3,4,5)}^{(k,1,1,1,k+1)}}\{-2k\}\oplus
\overline{S}_{(1,2,3;4,5)}^{[k,1,1;1,k+1]}\boxtimes
K\left(D;X_{k,(1,4)}^{(k,-1)}\right)_{R_{(1,2,3,4,5)}^{(k,1,1,1,k+1)}}\{-2k+2\}.
\end{equation}
Since the quotient $Q_{6}\left/\langle X_{1,(2,6)}^{(1,-1)}\rangle\right.$ has equations
\begin{eqnarray}
\nonumber
&&x_{j,7}=X_{j,(1,2,4)}^{(k,1,-1)}\hspace{1cm}(1\leq j \leq k),\\
\nonumber
&&x_{1,8}=x_{1,3},\\
\nonumber
&&x_{1,6}=x_{1,2},
\end{eqnarray}
we have $Q_{6}\left/\langle X_{1,(2,6)}^{(1,-1)}\rangle\right.\simeq R_{(1,2,3,4,5)}^{(k,1,1,1,k+1)}$.
Then, the matrix factorization (\ref{mf-r3-1-01-5}) is isomorphic to
\begin{equation}
\label{mf-r3-fin-01-5}
\overline{S}_{(1,2,3;4,5)}^{[k,1,1;1,k+1]}\boxtimes
K\left(D;X_{k,(1,4)}^{(k,-1)}\right)_{R_{(1,2,3,4,5)}^{(k,1,1,1,k+1)}}\{-2k\}
\end{equation}
and the matrix factorization (\ref{mf-r3-1-11-5}) is isomorphic to
\begin{equation}
\label{mf-r3-fin-11-5}
\overline{S}_{(1,2,3;4,5)}^{[k,1,1;1,k+1]}\boxtimes
K\left(
\frac{D}{x_{1,2}-x_{1,4}};(x_{1,2}-x_{1,4})X_{k,(1,4)}^{(k,-1)}
\right)_{R_{(1,2,3,4,5)}^{(k,1,1,1,k+1)}}\{-2k+1\}.
\end{equation}
For these decompositions, the morphisms $\overline{\sigma}_{+,1}$, $\overline{\sigma}_{+,2}$, $\overline{\sigma}_{+,3}$ and $\overline{\sigma}_{+,4}$ of the complex (\ref{r3-lem-planar-comp-5}) transform into
\begin{eqnarray}
\nonumber
\overline{\sigma}_{+,1}&\simeq&
\left(
\begin{array}{cc}
\id_{\overline{S}_{(1,2,3;4,5)}^{[k,1,1;1,k+1]}}\boxtimes(1,(x_{1,2}-x_{1,4})(x_{1,3}-x_{1,4}))&0\\
0&\id
\end{array}
\right),\\
\nonumber
\overline{\sigma}_{+,2}&\simeq&
\left(
\id_{\overline{S}_{(1,2,3;4,5)}^{[k,1,1;1,k+1]}}\boxtimes(1,x_{1,3}-x_{1,4}),
\id_{\overline{S}_{(1,2,3;4,5)}^{[k,1,1;1,k+1]}}\boxtimes(x_{1,2}-x_{1,4},1)
\right),\\
\nonumber
\overline{\sigma}_{+,3}&\simeq&
\left(\id,\id_{\overline{S}_{(1,2,3;4,5)}^{[k,1,1;1,k+1]}}\boxtimes(x_{1,2}-x_{1,4},x_{1,2}-x_{1,4})\right),\\
\nonumber
\overline{\sigma}_{+,4}&\simeq&
-\id_{\overline{S}_{(1,2,3;4,5)}^{[k,1,1;1,k+1]}}\boxtimes(1,x_{1,2}-x_{1,4}).
\end{eqnarray}
Then, the complex (\ref{r3-lem-planar-comp-5}) is isomorphic, in $\k^b(\HMF^{gr}_{R_{(1,2,3,4,5)}^{(k,1,1,1,k+1)},\omega_5})$, to
\begin{equation}
\nonumber
\xymatrix{-k-1\ar@{.}[d]&&&-k\ar@{.}[d]
\\
\overline{M}_5\{(k+1)n\}
\ar[rrr]^{\id_{\overline{S}}\boxtimes(1,x_{1,3}-x_{1,4})}
&&&
\overline{M}_6\{(k+1)n-1\}
}
\end{equation}
where
\begin{eqnarray}
\nonumber
\overline{M}_5&\simeq&
\overline{S}_{(1,2,3;4,5)}^{[k,1,1;1,k+1]}\boxtimes
K\left(
\frac{D}{(x_{1,4}-x_{1,3})(x_{1,4}-x_{1,3})};(x_{1,4}-x_{1,3})X_{k,(1,2,4)}^{(k,1,-1)}
\right)_{R_{(1,2,3,4,5)}^{(k,1,1,1,k+1)}}\{-2k-1\},\\
\nonumber
\overline{M}_6&\simeq&
\overline{S}_{(1,2,3;4,5)}^{[k,1,1;1,k+1]}\boxtimes
K\left(D;X_{k,(1,4)}^{(k,-1)}\right)_{R_{(1,2,3,4,5)}^{(k,1,1,1,k+1)}}\{-2k\}.
\end{eqnarray}
We find
\begin{equation}
\nonumber
\overline{M}_5\simeq\left( \input{figure/r3-lem-planar-k1-11-mf}\right),\hspace{1cm}
\overline{M}_6\simeq\left( \input{figure/r3-lem-planar-k1-12-mf}\right).
\end{equation}
Thus, we obtain
\begin{equation}
\nonumber
\c\left(\input{figure/r3-lem-plus-k1-3}\right)_n\simeq\c\left(\input{figure/r3-lem-plus-k1-4}\right)_n.
\end{equation}
\end{proof}
\begin{proof}[{\bf Proof of Proposition \ref{prop-r3} (6)}]
The complex for the diagram $\input{figure/r3-lem-minus-k1-3-text}$ is described as a complex of factorizations of $\k^b(\HMF^{gr}_{R_{(1,2,3,4,5)}^{(k,1,1,1,k+1)},\omega_5})$,
\begin{eqnarray}
\label{r3-lem-planar-comp-6}
&&\hspace{1.5cm}
\c\left(\input{figure/r3-lem-minus-k1-3-mf}\right)_n=\\
\nonumber
&&
\xymatrix{
k-1\ar@{.}[d]&k\ar@{.}[d]&k+1\ar@{.}[d]\\
{\c\left(\input{figure/r3-lem-planar-k1-10-mf}\right)_n
\hspace{-0.5cm}
\begin{array}{c}
{}_{\{-(k+1)n+2\}}\\
\left<k+1\right>
\end{array}}
\ar[r]^{\left(\hspace{-0.2cm}
\begin{array}{c}
\overline{\sigma}_{-,1}\\
\overline{\sigma}_{-,2}
\end{array}
\hspace{-0.2cm}\right)}&
{
\begin{array}{c}
\c\left(\input{figure/r3-lem-planar-k1-8-mf}\right)_n
\hspace{-0.5cm}
\begin{array}{c}
{}_{\{-(k+1)n+1\}}\\
\left<k+1\right>
\end{array}
\\
\bigoplus\\
\c\left(\input{figure/r3-lem-planar-k1-9-mf}\right)_n
\hspace{-0.5cm}
\begin{array}{c}
{}_{\{-(k+1)n+1\}}\\
\left<k+1\right>
\end{array}
\end{array}
}
\ar[r]^{\txt{$(\overline{\sigma}_{-,3},\overline{\sigma}_{-,4})$}}&
{\c\left(\input{figure/r3-lem-planar-k1-7-mf}\right)_n
\hspace{-0.5cm}
\begin{array}{c}
{}_{\{-(k+1)n\}}\\
\left<k+1\right>
\end{array}
}.
}
\end{eqnarray}
By the discussion of Proof of Proposition \ref{prop-r3} (5), we have
\begin{eqnarray}
\label{mf-r3-1-11-6}
&&
\c\left(\input{figure/r3-lem-planar-k1-10-mf}\right)_n\simeq
\overline{S}_{(1,2,3;4,5)}^{[k,1,1;1,k+1]}\boxtimes
K\left(
\frac{D}{x_{1,2}-x_{1,4}};(x_{1,2}-x_{1,4})X_{k,(1,4)}^{(k,-1)}
\right)_{R_{(1,2,3,4,5)}^{(k,1,1,1,k+1)}}\{-2k+1\},\\
\label{mf-r3-1-10-6}
&&
\c\left(\input{figure/r3-lem-planar-k1-8-mf}\right)_n\simeq
\overline{S}_{(1,2,3;4,5)}^{[k,1,1;1,k+1]}\boxtimes
K\left(D;X_{k,(1,4)}^{(k,-1)}\right)_{Q_{6}\left/\langle (x_{1,3}-x_{1,6})X_{1,(2,6)}^{(1,-1)}\rangle\right.}\{-2k\},
\\
\label{mf-r3-1-01-6}
&&
\c\left(\input{figure/r3-lem-planar-k1-9-mf}\right)_n\simeq
\overline{S}_{(1,2,3;4,5)}^{[k,1,1;1,k+1]}\boxtimes
K\left(D;X_{k,(1,4)}^{(k,-1)}\right)_{R_{(1,2,3,4,5)}^{(k,1,1,1,k+1)}}\{-2k\},
\\
\label{mf-r3-1-00-6}
&&
\c\left(\input{figure/r3-lem-planar-k1-7-mf}\right)_n\simeq
\overline{S}_{(1,2,3;4,5)}^{[k,1,1;1,k+1]}\boxtimes
K\left(D_0;(x_{1,6}-x_{1,4})X_{k,(1,4)}^{(k,-1)}\right)_{Q_{6}\left/\langle (x_{1,3}-x_{1,6})X_{1,(2,6)}^{(1,-1)}\rangle\right.}\{-2k-1\},
\end{eqnarray}
where $Q_{6}=R_{(1,2,3,4,5,6,7,8)}^{(k,1,1,1,k+1,1,k,1)}\left/\left<X_{1,(6,1)}^{(1,k)}-X_{1,(7,4)}^{(k,1)},\ldots,X_{k,(6,1)}^{(1,k)}-X_{k,(7,4)}^{(k,1)},X_{1,(3,2)}^{(1,1)}-X_{1,(8,6)}^{(1,1)}\right>\right.$.
The partial factorization $K\left(D_0;(x_{1,6}-x_{1,4})X_{k,(1,4)}^{(k,-1)}\right)_{Q_{6}\left/\langle (x_{1,3}-x_{1,6})X_{1,(2,6)}^{(1,-1)}\rangle\right.}$ of (\ref{mf-r3-1-00-6}) is isomorphic to
\begin{eqnarray}
\nonumber
&&
\xymatrix{
R_{18}\ar[rrrrr]_(.45){
\left(
\begin{array}{cc}
0&D\\
\frac{D}{(x_{1,4}-x_{1,3})(x_{1,4}-x_{1,3})}&0
\end{array}
\right)}
&&&&&
R_{19}\{2k+1-n\}\ar[rrrrr]_(.55){
\left(
\begin{array}{cc}
0&(x_{1,4}-x_{1,3})X_{k,(1,2,4)}^{(k,1,-1)}\\
X_{k,(1,4)}^{(k,-1)}&0
\end{array}
\right)}
&&&&&
R_{18}.
}
\end{eqnarray}
Then, the matrix factorization (\ref{mf-r3-1-00-6}) is isomorphic to
\begin{eqnarray}
\label{mf-r3-fin-00-6}
&&\overline{S}_{(1,2,3;4,5)}^{[k,1,1;1,k+1]}\boxtimes
K\left(
\frac{D}{(x_{1,4}-x_{1,3})(x_{1,4}-x_{1,3})};(x_{1,4}-x_{1,3})X_{k,(1,2,4)}^{(k,1,-1)}
\right)_{R_{(1,2,3,4,5)}^{(k,1,1,1,k+1)}}\{-2k-1\}\\
\nonumber
&&\oplus
\overline{S}_{(1,2,3;4,5)}^{[k,1,1;1,k+1]}\boxtimes
K\left(
D;X_{k,(1,4)}^{(k,-1)}
\right)_{R_{(1,2,3,4,5)}^{(k,1,1,1,k+1)}}\{-2k+1\}.
\end{eqnarray}
The partial factorization $K\left(D;X_{k,(1,4)}^{(k,-1)}\right)_{Q_{6}\left/\langle (x_{1,3}-x_{1,6})X_{1,(2,6)}^{(1,-1)}\rangle\right.}$ of (\ref{mf-r3-1-10-6}) is isomorphic to
\begin{eqnarray}
\nonumber
\xymatrix{
R_{19}\ar[rr]_(.45){
\left(
\begin{array}{cc}
D&0\\
0&D
\end{array}
\right)}
&&
R_{19}\{2k-1-n\}\ar[rrrr]_(.55){
\left(
\begin{array}{cc}
X_{k,(1,4)}^{(k,-1)}&0\\
0&X_{k,(1,4)}^{(k,-1)}
\end{array}
\right)}
&&&&
R_{19}.
}
\end{eqnarray}
Then, the matrix factorization (\ref{mf-r3-1-10-6}) is decomposed into
\begin{equation}
\label{mf-r3-fin-10-6}
\overline{S}_{(1,2,3;4,5)}^{[k,1,1;1,k+1]}\boxtimes
K\left(D;X_{k,(1,4)}^{(k,-1)}\right)_{R_{(1,2,3,4,5)}^{(k,1,1,1,k+1)}}\{-2k\}\oplus
\overline{S}_{(1,2,3;4,5)}^{[k,1,1;1,k+1]}\boxtimes
K\left(D;X_{k,(1,4)}^{(k,-1)}\right)_{R_{(1,2,3,4,5)}^{(k,1,1,1,k+1)}}\{-2k+2\}.
\end{equation}
For these decompositions, the morphisms $\overline{\sigma}_{-,1}$, $\overline{\sigma}_{-,2}$, $\overline{\sigma}_{-,3}$ and $\overline{\sigma}_{-,4}$ of the complex (\ref{r3-lem-planar-comp-6}) transform into
\begin{eqnarray}
\nonumber
\overline{\sigma}_{-,1}&\simeq&
\left(
\begin{array}{c}
\id_{\overline{S}_{(1,2,3;4,5)}^{[k,1,1;1,k+1]}}\boxtimes(-x_{1,2}+x_{1,4},-x_{1,2}+x_{1,4})\\
\id
\end{array}
\right),\\
\nonumber
\overline{\sigma}_{-,2}&\simeq&
\id_{\overline{S}_{(1,2,3;4,5)}^{[k,1,1;1,k+1]}}\boxtimes(x_{1,2}-x_{1,4},1),\\
\nonumber
\overline{\sigma}_{-,3}&\simeq&
\left(
\begin{array}{cc}
0&\id_{\overline{S}_{(1,2,3;4,5)}^{[k,1,1;1,k+1]}}\boxtimes((x_{1,2}-x_{1,4})(x_{1,3}-x_{1,4}),1)\\
\id&0
\end{array}
\right),\\
\nonumber
\overline{\sigma}_{-,4}&\simeq&
-\left(
\begin{array}{c}
\id_{\overline{S}_{(1,2,3;4,5)}^{[k,1,1;1,k+1]}}\boxtimes(x_{1,3}-x_{1,4},1)\\
\id_{\overline{S}_{(1,2,3;4,5)}^{[k,1,1;1,k+1]}}\boxtimes(-1,-x_{1,2}+x_{1,4})
\end{array}
\right).
\end{eqnarray}
Then, the complex (\ref{r3-lem-planar-comp-6}) is isomorphic, in $\k^b(\HMF^{gr}_{R_{(1,2,3,4,5)}^{(k,1,1,1,k+1)},\omega_5})$, to
\begin{equation}
\nonumber
\xymatrix{k\ar@{.}[d]&&&k+1\ar@{.}[d]
\\
\overline{M}_6\{-(k+1)n+1\}
\ar[rrr]^{\id_{\overline{S}}\boxtimes(1,x_{1,3}-x_{1,4})}
&&&
\overline{M}_5\{-(k+1)n\}
}
\end{equation}
where
\begin{eqnarray}
\nonumber
\overline{M}_5&\simeq&
\overline{S}_{(1,2,3;4,5)}^{[k,1,1;1,k+1]}\boxtimes
K\left(
\frac{D}{(x_{1,4}-x_{1,3})(x_{1,4}-x_{1,3})};(x_{1,4}-x_{1,3})X_{k,(1,2,4)}^{(k,1,-1)}
\right)_{R_{(1,2,3,4,5)}^{(k,1,1,1,k+1)}}\{-2k-1\},\\
\nonumber
\overline{M}_6&\simeq&
\overline{S}_{(1,2,3;4,5)}^{[k,1,1;1,k+1]}\boxtimes
K\left(D;X_{k,(1,4)}^{(k,-1)}\right)_{R_{(1,2,3,4,5)}^{(k,1,1,1,k+1)}}\{-2k\}.
\end{eqnarray}
Thus, we obtain
\begin{equation}
\nonumber
\c\left(\input{figure/r3-lem-minus-k1-3}\right)_n\simeq\c\left(\input{figure/r3-lem-minus-k1-4}\right)_n.
\end{equation}
\end{proof}
\begin{proof}[{\bf Proof of Proposition \ref{prop-r3} (7)}]
The complex for the diagram $\input{figure/r3-lem-plus-1k-3-text}$ is described as a complex of factorizations of $\k^b(\HMF^{gr}_{R_{(1,2,3,4,5)}^{(k,1,1,1,k+1)},\omega_6})$,
\begin{eqnarray}
\label{r3-lem-planar-comp-7}
&&\hspace{1.5cm}
\c\left(\input{figure/r3-lem-plus-1k-3-mf}\right)_n=\\
\nonumber
&&
\xymatrix{
	-k-1\ar@{.}[d]&-k\ar@{.}[d]&-k+1\ar@{.}[d]\\
{\c\left(\input{figure/r3-lem-planar-1k-7-mf}\right)_n
\hspace{-0.5cm}
\begin{array}{c}
{}_{\{(k+1)n\}}\\
\left<k+1\right>
\end{array}}
\ar[r]^{\left(\hspace{-0.2cm}
\begin{array}{c}
\overline{\tau}_{+,1}\\
\overline{\tau}_{+,2}
\end{array}
\hspace{-0.2cm}\right)}&
{
\begin{array}{c}
\c\left(\input{figure/r3-lem-planar-1k-8-mf}\right)_n
\hspace{-0.5cm}
\begin{array}{c}
{}_{\{(k+1)n-1\}}\\
\left<k+1\right>
\end{array}
\\
\bigoplus\\
\c\left(\input{figure/r3-lem-planar-1k-9-mf}\right)_n
\hspace{-0.5cm}
\begin{array}{c}
{}_{\{(k+1)n-1\}}\\
\left<k+1\right>
\end{array}
\end{array}
}
\ar[r]^{\txt{$(\overline{\tau}_{+,3},\overline{\tau}_{+,4})$}}&
{\c\left(\input{figure/r3-lem-planar-1k-10-mf}\right)_n
\hspace{-0.5cm}
\begin{array}{c}
{}_{\{(k+1)n-2\}}\\
\left<k+1\right>
\end{array}
}.
}
\end{eqnarray}
By Corollary \ref{cor2-11}, we have
\begin{eqnarray}
\label{mf-r3-1-00-7}
&&
\c\left(\input{figure/r3-lem-planar-1k-7-mf}\right)_n\simeq
K\left(
\left(
\begin{array}{c}
V_{1,(7,8;5)}^{[1,k]}\\
\vdots\\
V_{k+1,(7,8;5)}^{[1,k]}\\
u_{k+1,(6,1,7,4)}^{[1,1]}
\end{array}
\right)
;
\left(
\begin{array}{c}
X_{1,(7,8)}^{(1,k)}-x_{1,5}\\
\vdots\\
X_{k+1,(7,8)}^{(1,k)}-x_{k+1,5}\\
(x_{1,6}-x_{1,4})X_{1,(1,4)}^{(1,-1)}
\end{array}
\right)
\right)_{Q_{7}\left/\langle (x_{1,3}-x_{1,6})X_{k,(2,6)}^{(k,-1)}\rangle\right.}\{-2k-1\},
\end{eqnarray}
where $Q_{7}=R_{(1,2,3,4,5,6,7,8)}^{(1,k,1,1,k+1,1,1,k)}\left/\left<X_{1,(6,1)}^{(1,1)}-X_{1,(7,4)}^{(1,1)},X_{1,(3,2)}^{(1,k)}-X_{1,(8,6)}^{(k,1)},\ldots,X_{k,(3,2)}^{(1,k)}-X_{k,(8,6)}^{(k,1)}\right>\right.$. This quotient $Q_{7}\left/\langle (x_{1,3}-x_{1,6})X_{k,(2,6)}^{(k,-1)}\rangle\right.$ has equations
\begin{eqnarray}
\nonumber
&&x_{1,7}=X_{1,(1,4,6)}^{(1,-1,1)},\\
\nonumber
&&x_{j,8}=X_{j,(2,3,6)}^{(k,1,-1)}\hspace{1cm}(1\leq j \leq k),\\
\nonumber
&&(x_{1,3}-x_{1,6})X_{k,(2,6)}^{(k,-1)}=0.
\end{eqnarray}
Then, the matrix factorization (\ref{mf-r3-1-00-7}) is isomorphic to
\begin{equation}
\nonumber
K\left(
\left(
\begin{array}{c}
V_{1,(7,8;5)}^{[1,k]}\\
\vdots\\
V_{k+1,(7,8;5)}^{[1,k]}\\
u_{k+1,(6,1,7,4)}^{[1,1]}-V_{k+1,(7,8;5)}^{[1,k]}
\end{array}
\right)
;
\left(
\begin{array}{c}
X_{1,(1,2,3,4)}^{(1,k,1,-1)}-x_{1,5}\\
\vdots\\
X_{k+1,(1,2,3,4)}^{(1,k,1,-1)}-x_{k+1,5}\\
(x_{1,6}-x_{1,4})X_{1,(1,4)}^{(1,-1)}
\end{array}
\right)
\right)_{Q_{7}\left/\langle (x_{1,3}-x_{1,6})X_{k,(2,6)}^{(k,-1)}\rangle\right.}\{-2k-1\}.
\end{equation}
By Corollary \ref{induce-sq1}, there exist polynomials $G_1,\ldots,G_{k+1}$ $\in R_{(1,2,3,4,5)}^{(1,k,1,1,k+1)}$ and $G_0 \in Q_{7}\left/\langle (x_{1,3}-x_{1,6})X_{k,(2,6)}^{(k,-1)}\rangle\right.$ such that $G_0(x_{1,6}-x_{1,4})\equiv G \in R_{(1,2,3,4,5)}^{(1,k,1,1,k+1)}$ and we have an isomorphism to the above factorization
\begin{equation}
\label{mf-r3-2-00-7}
\overline{S}_{(1,2,3;4,5)}^{[1,k,1;1,k+1]}\boxtimes
K\left(G_0;(x_{1,6}-x_{1,4})X_{1,(1,4)}^{(1,-1)}\right)_{Q_{7}\left/\langle (x_{1,3}-x_{1,6})X_{k,(2,6)}^{(k,-1)}\rangle\right.}\{-2k-1\},
\end{equation}
where 
\begin{equation}
\nonumber
\overline{S}_{(1,2,3;4,5)}^{[1,k,1;1,k+1]}
:=K\left(
\left(
\begin{array}{c}
\vspace{.1cm}G_1\\
\vspace{.1cm}\vdots\\
\vspace{.2cm}G_{k+1}\\
G_0
\end{array}
\right)
;
\left(
\begin{array}{c}
X_{1,(1,2,3,4)}^{(1,k,1,-1)}-x_{1,5}\\
\vdots\\
X_{k+1,(1,2,3,4)}^{(1,k,1,-1)}-x_{k+1,5}\\
(x_{1,6}-x_{1,4})X_{1,(1,4)}^{(1,-1)}
\end{array}
\right)
\right)_{R_{(1,2,3,4,5)}^{(1,k,1,1,k+1)}}
\end{equation}
We consider isomorphisms of $Q_{7}\left/\langle (x_{1,3}-x_{1,6})X_{k,(2,6)}^{(k,-1)}\rangle\right.$ as $R_{(1,2,3,4,5)}^{(1,k,1,1,k+1)}$-module
\begin{eqnarray}
\nonumber
&&R_{20}:=R_{(1,2,3,4,5)}^{(1,k,1,1,k+1)}
\oplus(x_{1,6}-x_{1,4})R_{(1,2,3,4,5)}^{(1,k,1,1,k+1)}
\oplus\ldots
\oplus x_{1,6}^{k-1}(x_{1,6}-x_{1,4})R_{(1,2,3,4,5)}^{(1,k,1,1,k+1)},\\
\nonumber
&&R_{21}:=R_{(1,2,3,4,5)}^{(1,k,1,1,k+1)}
\oplus x_{1,6}R_{(1,2,3,4,5)}^{(1,k,1,1,k+1)}
\oplus\ldots
\oplus x_{1,6}^{k-1}R_{(1,2,3,4,5)}^{(1,k,1,1,k+1)}
\oplus X_{k,(2,3,4,6)}^{(k,1,-1,-1)}R_{(1,2,3,4,5)}^{(1,k,1,1,k+1)}.
\end{eqnarray}
The partial factorization $K\left(G_0;(x_{1,6}-x_{1,4})X_{1,(1,4)}^{(1,-1)}\right)_{Q_{7}\left/\langle (x_{1,3}-x_{1,6})X_{k,(2,6)}^{(k,-1)}\rangle\right.}$ is isomorphic to
\begin{eqnarray}
\nonumber
&&
\xymatrix{
R_{20}\ar[rrrr]_(.45){
\left(
\begin{array}{cc}
{}^{t}\mathfrak{0}_k&E_{k}(G)\\
\frac{G}{X_{k+1,(2,3,4)}^{(k,1,-1)}}&\mathfrak{0}_k
\end{array}
\right)}
&&&&
R_{21}\{3-n\}\ar[rrrrrr]_(.55){
\left(
\begin{array}{cc}
\mathfrak{0}_k&(x_{1,3}-x_{1,4})X_{k+1,(1,2,4)}^{(1,k,-1)}\\
E_{k}(X_{1,(1,4)}^{(1,-1)})&{}^{t}\mathfrak{0}_k
\end{array}
\right)}
&&&&&&
R_{20}.
}
\end{eqnarray}
Then, the matrix factorization (\ref{mf-r3-2-00-7}) is isomorphic to
\begin{eqnarray}
\label{mf-r3-fin-00-7}
&&\overline{S}_{(1,2,3;4,5)}^{[1,k,1;1,k+1]}\boxtimes
K\left(
\frac{G}{X_{k+1,(2,3,4)}^{(k,1,-1)}};(x_{1,3}-x_{1,4})X_{k+1,(1,2,4)}^{(1,k,-1)}
\right)_{R_{(1,2,3,4,5)}^{(1,k,1,1,k+1)}}\{-2k-1\}\\
\nonumber
&&\oplus
\bigoplus_{j=1}^{k}\overline{S}_{(1,2,3;4,5)}^{[1,k,1;1,k+1]}\boxtimes
K\left(
G;X_{1,(1,4)}^{(1,-1)}
\right)_{R_{(1,2,3,4,5)}^{(1,k,1,1,k+1)}}\{-2k+1+2j\}.
\end{eqnarray}
By a similar discussion, we obtain
\begin{eqnarray}
\label{mf-r3-1-10-7}
&&
\c\left(\input{figure/r3-lem-planar-1k-8-mf}\right)_n\simeq
\overline{S}_{(1,2,3;4,5)}^{[1,k,1;1,k+1]}\boxtimes
K\left(G;X_{k,(1,4)}^{(1,-1)}\right)_{Q_{7}\left/\langle (x_{1,3}-x_{1,6})X_{k,(2,6)}^{(k,-1)}\rangle\right.}\{-2k\},
\\
\label{mf-r3-1-01-7}
&&
\c\left(\input{figure/r3-lem-planar-1k-9-mf}\right)_n\simeq
\overline{S}_{(1,2,3;4,5)}^{[1,k,1;1,k+1]}\boxtimes
K\left(G_0;(x_{1,6}-x_{1,4})X_{1,(1,4)}^{(1,-1)}\right)_{Q_{7}\left/\langle X_{k,(2,6)}^{(k,-1)}\rangle\right.}\{-2k\},
\\
\label{mf-r3-1-11-7}
&&
\c\left(\input{figure/r3-lem-planar-1k-10-mf}\right)_n\simeq
\overline{S}_{(1,2,3;4,5)}^{[1,k,1;1,k+1]}\boxtimes
K\left(G;X_{1,(1,4)}^{(1,-1)}\right)_{Q_{7}\left/\langle X_{k,(2,6)}^{(k,-1)}\rangle\right.}\{-2k+1\}.
\end{eqnarray}
\indent
The partial factorization $K\left(G;X_{1,(1,4)}^{(1,-1)}\right)_{Q_{7}\left/\langle (x_{1,3}-x_{1,6})X_{k,(2,6)}^{(k,-1)}\rangle\right.}$ of (\ref{mf-r3-1-10-7}) is isomorphic to
\begin{eqnarray}
\nonumber
\xymatrix{
R_{20}\ar[rr]_(.45){E_{k+1}(G)}
&&
R_{20}\{1-n\}\ar[rrr]_(.55){E_{k+1}(X_{1,(1,4)}^{(1,-1)})}
&&&
R_{20}.
}
\end{eqnarray}
Then, the matrix factorization (\ref{mf-r3-1-10-7}) is decomposed into
\begin{equation}
\label{mf-r3-fin-10-7}
\bigoplus_{j=0}^{k+1}\overline{S}_{(1,2,3;4,5)}^{[1,k,1;1,k+1]}\boxtimes
K\left(G;X_{1,(1,4)}^{(1,-1)}\right)_{R_{(1,2,3,4,5)}^{(1,k,1,1,k+1)}}\{-2k+2j\}.
\end{equation}
The quotient $Q_{7}\left/\langle X_{k,(2,6)}^{(k,-1)}\rangle\right.$ has equations
\begin{eqnarray}
\nonumber
&&x_{1,7}=X_{1,(1,2,4)}^{(1,1,-1)},\\
\nonumber
&&x_{j,8}=X_{j,(2,3,6)}^{(k,1,-1)}\hspace{1cm}(1\leq j \leq k),\\
\nonumber
&&X_{k,(2,6)}^{(k,-1)}=0.
\end{eqnarray}
We consider isomorphisms of $Q_{7}\left/\langle X_{k,(2,6)}^{(k,-1)}\rangle\right.$ as an $R_{(1,2,3,4,5)}^{(1,k,1,1,k+1)}$-module
\begin{eqnarray}
\nonumber
&&R_{22}:=R_{(1,2,3,4,5)}^{(1,k,1,1,k+1)}
\oplus(x_{1,6}-x_{1,4})R_{(1,2,3,4,5)}^{(1,k,1,1,k+1)}
\oplus\ldots
\oplus x_{1,6}^{k-2}(x_{1,6}-x_{1,4})R_{(1,2,3,4,5)}^{(1,k,1,1,k+1)},\\
\nonumber
&&R_{23}:=R_{(1,2,3,4,5)}^{(1,k,1,1,k+1)}
\oplus x_{1,6}R_{(1,2,3,4,5)}^{(1,k,1,1,k+1)}
\oplus\ldots
\oplus x_{1,6}^{k-2}R_{(1,2,3,4,5)}^{(1,k,1,1,k+1)}
\oplus X_{k-1,(2,4,6)}^{(k,-1,-1)}R_{(1,2,3,4,5)}^{(1,k,1,1,k+1)}.
\end{eqnarray}
The partial matrix factorization $K\left(G_0;(x_{1,6}-x_{1,4})X_{1,(1,4)}^{(1,-1)}\right)_{Q_{7}\left/\langle X_{k,(2,6)}^{(k,-1)}\rangle\right.}$ of (\ref{mf-r3-1-01-7}) is isomorphic to
\begin{eqnarray}
\nonumber
&&\xymatrix{
R_{22}
\ar[rrrr]_(.45){
\left(
\begin{array}{cc}
{}^{t}\mathfrak{0}_{k-1}&E_{k-1}(G)\\
\frac{G}{X_{k,(2,4)}^{(k,-1)}}&\mathfrak{0}_{k-1}
\end{array}
\right)}
&&&&
R_{23}\{3-n\}
\ar[rrrrrr]_(.55){
\left(
\begin{array}{cc}
\mathfrak{0}_{k-1}&X_{k+1,(1,2,4)}^{(1,k,-1)}\\
E_{k-1}(X_{1,(1,4)}^{(1,-1)})&{}^{t}\mathfrak{0}_{k-1}
\end{array}
\right)
}
&&&&&&
R_{22}
}\\
\nonumber
&\simeq&K\left(\frac{G}{X_{k,(2,4)}^{(k,-1)}};X_{k+1,(1,2,4)}^{(1,k,-1)}\right)_{R_{(1,2,3,4,5)}^{(1,k,1,1,k+1)}}\hspace{-1cm}\{-2k\}
\oplus\bigoplus_{j=1}^{k-1}
K\left(G;X_{1,(1,4)}^{(1,-1)}\right)_{R_{(1,2,3,4,5)}^{(1,k,1,1,k+1)}}\{-2k+2j\}
\end{eqnarray}
and the partial matrix factorization $K\left(G;X_{1,(1,4)}^{(1,-1)}\right)_{Q_{7}\left/\langle X_{k,(2,6)}^{(k,-1)}\rangle\right.}$ of (\ref{mf-r3-1-11-7}) is isomorphic to
\begin{eqnarray}
\nonumber
&&\xymatrix{
R_{22}
\ar[rr]_(.45){E_{k}(G)}
&&
R_{22}\{1-n\}
\ar[rr]_(.55){E_{k}(X_{1,(1,4)}^{(1,-1)})}
&&
R_{22}
}\\
\nonumber
&\simeq&\bigoplus_{j=0}^{k-1}
K\left(G;X_{1,(1,4)}^{(1,-1)}\right)_{R_{(1,2,3,4,5)}^{(1,k,1,1,k+1)}}\{-2k+2j\}
\end{eqnarray}
For these decompositions, the morphisms $\overline{\tau}_{+,1}$, $\overline{\tau}_{+,2}$, $\overline{\tau}_{+,3}$ and $\overline{\tau}_{+,4}$ of the complex (\ref{r3-lem-planar-comp-7}) transform into
\begin{eqnarray}
\nonumber
\overline{\tau}_{+,1}&\simeq&
\left(
\begin{array}{cc}
\id_{\overline{S}_{(1,2,3;4,5)}^{[1,k,1;1,k+1]}}\boxtimes(1,X_{k+1,(2,3,4)}^{(k,1,-1)})&\mathfrak{0}_{k}\\
{}^{t}\mathfrak{0}_{k}&E_{k}(\id)
\end{array}
\right),\\
\nonumber
\overline{\tau}_{+,2}&\simeq&
\left(
\begin{array}{ccc}
\id_{\overline{S}_{(1,2,3;4,5)}^{[1,k,1;1,k+1]}}\boxtimes(1,x_{1,3}-x_{1,4})&\mathfrak{0}_{k-1}&
(-1)^{k-1}\id_{\overline{S}_{(1,2,3;4,5)}^{[1,k,1;1,k+1]}}\boxtimes(X_{k,(2,4)}^{(k,-1)},1)\\
&&(-1)^{k-2}\id_{\overline{S}_{(1,2,3;4,5)}^{[1,k,1;1,k+1]}}\boxtimes(X_{k-1,(2,4)}^{(k,-1)},X_{k-1,(2,4)}^{(k,-1)})\\
{}^{t}\mathfrak{0}_{k-1}&E_{k-1}\left(\id\right)&\vdots\\
&&\id_{\overline{S}_{(1,2,3;4,5)}^{[1,k,1;1,k+1]}}\boxtimes(X_{1,(2,4)}^{(k,-1)},X_{1,(2,4)}^{(k,-1)})
\end{array}
\right),\\
\nonumber
\overline{\tau}_{+,3}&\simeq&
\left(
\begin{array}{cc}
&(-1)^{k-1}\id_{\overline{S}_{(1,2,3;4,5)}^{[1,k,1;1,k+1]}}\boxtimes(X_{k,(2,4)}^{(k,-1)},X_{k,(2,4)}^{(k,-1)})\\
E_{k}(\id)&\vdots\\
&\id_{\overline{S}_{(1,2,3;4,5)}^{[1,k,1;1,k+1]}}\boxtimes(X_{1,(2,4)}^{(k,-1)},X_{1,(2,4)}^{(k,-1)})
\end{array}
\right),\\
\nonumber
\overline{\tau}_{+,4}&\simeq&
-\left(
\begin{array}{cc}
\id_{\overline{S}_{(1,2,3;4,5)}^{[1,k,1;1,k+1]}}\boxtimes(1,X_{k,(2,4)}^{(k,-1)})&\mathfrak{0}_{k-1}\\
{}^{t}\mathfrak{0}_{k-1}&E_{k-1}(\id)
\end{array}
\right).
\end{eqnarray}
Then, the complex (\ref{r3-lem-planar-comp-7}) is isomorphic, in $\k^b(\HMF^{gr}_{R_{(1,2,3,4,5)}^{(1,k,1,1,k+1)},\omega_6})$, to
\begin{equation}
\nonumber
\xymatrix{-k-1\ar@{.}[d]&&&-k\ar@{.}[d]
\\
\overline{M}_7\{(k+1)n\}
\ar[rrr]^{\id_{\overline{S}}\boxtimes(1,x_{1,3}-x_{1,4})}
&&&
\overline{M}_8\{(k+1)n-1\}
}
\end{equation}
where
\begin{eqnarray}
\nonumber
\overline{M}_7&\simeq&
\overline{S}_{(1,2,3;4,5)}^{[1,k,1;1,k+1]}\boxtimes
K\left(
\frac{G}{X_{k+1,(2,3,4)}^{(k,1,-1)}};(x_{1,3}-x_{1,4})X_{k+1,(1,2,4)}^{(1,k,-1)}
\right)_{R_{(1,2,3,4,5)}^{(1,k,1,1,k+1)}}\{-2k-1\},\\
\nonumber
\overline{M}_8&\simeq&
\overline{S}_{(1,2,3;4,5)}^{[1,k,1;1,k+1]}\boxtimes
K\left(\frac{G}{X_{k,(2,4)}^{(k,-1)}};X_{k+1,(1,2,4)}^{(1,k,-1)}
\right)_{R_{(1,2,3,4,5)}^{(1,k,1,1,k+1)}}\{-2k\}.
\end{eqnarray}
We have
\begin{equation}
\nonumber
\overline{M}_7\simeq\left( \input{figure/r3-lem-planar-1k-11-mf}\right),\hspace{1cm}
\overline{M}_8\simeq\left( \input{figure/r3-lem-planar-1k-12-mf}\right).
\end{equation}
Thus, we obtain
\begin{equation}
\nonumber
\c\left(\input{figure/r3-lem-plus-1k-3}\right)_n\simeq\c\left(\input{figure/r3-lem-plus-1k-4}\right)_n.
\end{equation}
\end{proof}
\begin{proof}[{\bf Proof of Proposition \ref{prop-r3} (8)}]
The complex for the diagram $\input{figure/r3-lem-minus-1k-3-text}$ is described as a complex of factorizations of $\k^b(\HMF^{gr}_{R_{(1,2,3,4,5)}^{(1,k,1,1,k+1)},\omega_6})$,
\begin{eqnarray}
\label{r3-lem-planar-comp-8}
&&\hspace{1.5cm}
\c\left(\input{figure/r3-lem-minus-1k-3-mf}\right)_n=\\
\nonumber
&&
\xymatrix{
k-1\ar@{.}[d]&k\ar@{.}[d]&k+1\ar@{.}[d]\\
{\c\left(\input{figure/r3-lem-planar-1k-10-mf}\right)_n
\hspace{-0.5cm}
\begin{array}{c}
{}_{\{-(k+1)n+2\}}\\
\left<k+1\right>
\end{array}}
\ar[r]^{\left(\hspace{-0.2cm}
\begin{array}{c}
\overline{\tau}_{-,1}\\
\overline{\tau}_{-,2}
\end{array}
\hspace{-0.2cm}\right)}&
{
\begin{array}{c}
\c\left(\input{figure/r3-lem-planar-1k-9-mf}\right)_n
\hspace{-0.5cm}
\begin{array}{c}
{}_{\{-(k+1)n+1\}}\\
\left<k+1\right>
\end{array}
\\
\bigoplus\\
\c\left(\input{figure/r3-lem-planar-1k-8-mf}\right)_n
\hspace{-0.5cm}
\begin{array}{c}
{}_{\{-(k+1)n+1\}}\\
\left<k+1\right>
\end{array}
\end{array}
}
\ar[r]^{\txt{$(\overline{\tau}_{-,3},\overline{\tau}_{-,4})$}}&
{\c\left(\input{figure/r3-lem-planar-1k-7-mf}\right)_n
\hspace{-0.5cm}
\begin{array}{c}
{}_{\{-(k+1)n\}}\\
\left<k+1\right>
\end{array}
}.
}
\end{eqnarray}
By the discussion of Proof of Proposition \ref{prop-r3} (7), we have
\begin{eqnarray}
\label{mf-r3-1-11-8}
&&
\c\left(\input{figure/r3-lem-planar-1k-10-mf}\right)_n\simeq
\overline{S}_{(1,2,3;4,5)}^{[1,k,1;1,k+1]}\boxtimes
K\left(G;X_{1,(1,4)}^{(1,-1)}\right)_{Q_{7}\left/\langle X_{k,(2,6)}^{(k,-1)}\rangle\right.}\{-2k+1\},\\
\label{mf-r3-1-01-8}
&&
\c\left(\input{figure/r3-lem-planar-1k-9-mf}\right)_n\simeq
\overline{S}_{(1,2,3;4,5)}^{[1,k,1;1,k+1]}\boxtimes
K\left(G;X_{1,(1,4)}^{(1,-1)}\right)_{Q_{7}\left/\langle (x_{1,3}-x_{1,6})X_{k,(2,6)}^{(k,-1)}\rangle\right.}\{-2k\},
\\
\label{mf-r3-1-10-8}
&&
\c\left(\input{figure/r3-lem-planar-1k-8-mf}\right)_n\simeq
\overline{S}_{(1,2,3;4,5)}^{[1,k,1;1,k+1]}\boxtimes
K\left(G_0;(x_{1,6}-x_{1,4})X_{1,(1,4)}^{(1,-1)}\right)_{Q_{7}\left/\langle X_{k,(2,6)}^{(k,-1)}\rangle\right.}\{-2k\},
\\
\label{mf-r3-1-00-8}
&&
\c\left(\input{figure/r3-lem-planar-1k-7-mf}\right)_n\simeq
\overline{S}_{(1,2,3;4,5)}^{[1,k,1;1,k+1]}\boxtimes
K\left(G_0;(x_{1,6}-x_{1,4})X_{1,(1,4)}^{(1,-1)}\right)_{Q_{7}\left/\langle (x_{1,3}-x_{1,6})X_{k,(2,6)}^{(k,-1)}\rangle\right.}\{-2k-1\},
\end{eqnarray}
where $Q_{7}=R_{(1,2,3,4,5,6,7,8)}^{(1,k,1,1,k+1,1,1,k)}\left/\left<X_{1,(6,1)}^{(1,1)}-X_{1,(7,4)}^{(1,1)},X_{1,(3,2)}^{(1,k)}-X_{1,(8,6)}^{(k,1)},\ldots,X_{k,(3,2)}^{(1,k)}-X_{k,(8,6)}^{(k,1)}\right>\right.$.\\
\indent
The partial factorization $K\left(G;X_{1,(1,4)}^{(1,-1)}\right)_{Q_{7}\left/\langle X_{k,(2,6)}^{(k,-1)}\rangle\right.}$ of (\ref{mf-r3-1-11-8}) is isomorphic to
\begin{eqnarray}
\nonumber
&&
\xymatrix{
R_{23}\ar[rr]_(.45){E_{k}(G)}
&&
R_{23}\{1-n\}\ar[rr]_(.55){E_{k}(X_{1,(1,4)}^{(1,-1)})}
&&
R_{23}
}\\
\nonumber
&\simeq&
\bigoplus_{j=0}^{k-1}K\left(G;X_{1,(1,4)}^{(1,-1)}\right)_{R_{(1,2,3,4,5)}^{(1,k,1,1,k+1)}}\{2j\}.
\end{eqnarray}
Then, the matrix factorization (\ref{mf-r3-1-11-8}) is isomorphic to
\begin{eqnarray}
\label{mf-r3-fin-11-8}
\bigoplus_{j=0}^{k-1}\overline{S}_{(1,2,3;4,5)}^{[1,k,1;1,k+1]}\boxtimes
K\left(G;X_{1,(1,4)}^{(1,-1)}\right)_{R_{(1,2,3,4,5)}^{(1,k,1,1,k+1)}}\{-2k+1+2j\}.
\end{eqnarray}
The partial matrix factorization $K\left(G_0;(x_{1,6}-x_{1,4})X_{1,(1,4)}^{(1,-1)}\right)_{Q_{7}\left/\langle X_{k,(2,6)}^{(k,-1)}\rangle\right.}$ of (\ref{mf-r3-1-01-8}) is isomorphic to
\begin{eqnarray}
\nonumber
&&\xymatrix{
R_{22}
\ar[rrrr]_(.45){
\left(
\begin{array}{cc}
{}^{t}\mathfrak{0}_{k-1}&E_{k-1}(G)\\
\frac{G}{X_{k,(2,4)}^{(k,-1)}}&\mathfrak{0}_{k-1}
\end{array}
\right)}
&&&&
R_{23}\{3-n\}
\ar[rrrrrr]_(.55){
\left(
\begin{array}{cc}
\mathfrak{0}_{k-1}&X_{k+1,(1,2,4)}^{(1,k,-1)}\\
E_{k-1}(X_{1,(1,4)}^{(1,-1)})&{}^{t}\mathfrak{0}_{k-1}
\end{array}
\right)
}
&&&&&&
R_{22}
}\\
\nonumber
&\simeq&K\left(\frac{G}{X_{k,(2,4)}^{(k,-1)}};X_{k+1,(1,2,4)}^{(1,k,-1)}\right)_{R_{(1,2,3,4,5)}^{(1,k,1,1,k+1)}}
\oplus\bigoplus_{j=1}^{k-1}
K\left(G;X_{1,(1,4)}^{(1,-1)}\right)_{R_{(1,2,3,4,5)}^{(1,k,1,1,k+1)}}\{2j\}
\end{eqnarray}
Then, the matrix factorization (\ref{mf-r3-1-01-8}) is isomorphic to
\begin{eqnarray}
\label{mf-r3-fin-01-8}
&&
\hspace{-1cm}
\overline{S}_{(1,2,3;4,5)}^{[1,k,1;1,k+1]}\boxtimes K\left(\frac{G}{X_{k,(2,4)}^{(k,-1)}};X_{k+1,(1,2,4)}^{(1,k,-1)}\right)_{R_{(1,2,3,4,5)}^{(1,k,1,1,k+1)}}\hspace{-1cm}\{-2k\}
\\
\nonumber
&&
\oplus\bigoplus_{j=1}^{k-1}
\overline{S}_{(1,2,3;4,5)}^{[1,k,1;1,k+1]}\boxtimes K\left(G;X_{1,(1,4)}^{(1,-1)}\right)_{R_{(1,2,3,4,5)}^{(1,k,1,1,k+1)}}\{-2k+2j\}.
\end{eqnarray}
We consider isomorphisms of $Q_{7}\left/\langle (x_{1,3}-x_{1,6})X_{k,(2,6)}^{(k,-1)}\rangle\right.$ as an $R_{(1,2,3,4,5)}^{(1,k,1,1,k+1)}$-module
\begin{eqnarray}
\nonumber
&&R_{24}:=R_{(1,2,3,4,5)}^{(1,k,1,1,k+1)}
\oplus(x_{1,6}-x_{1,4})R_{(1,2,3,4,5)}^{(1,k,1,1,k+1)}\\
\nonumber
&&\hspace{1.5cm}\oplus(x_{1,3}-x_{1,6})(x_{1,6}-x_{1,4})R_{(1,2,3,4,5)}^{(1,k,1,1,k+1)}
\oplus\ldots
\oplus x_{1,6}^{k-2}(x_{1,3}-x_{1,6})(x_{1,6}-x_{1,4})R_{(1,2,3,4,5)}^{(1,k,1,1,k+1)},\\
\nonumber
&&R_{25}:=R_{(1,2,3,4,5)}^{(1,k,1,1,k+1)}
\oplus(x_{1,3}-x_{1,6})R_{(1,2,3,4,5)}^{(1,k,1,1,k+1)}
\oplus\ldots
\oplus x_{1,6}^{k-2}(x_{1,3}-x_{1,6})R_{(1,2,3,4,5)}^{(1,k,1,1,k+1)}
\oplus X_{k,(2,3,4,6)}^{(k,1,-1,-1)}R_{(1,2,3,4,5)}^{(1,k,1,1,k+1)}.
\end{eqnarray}
The partial factorization $K\left(G;X_{k,(1,4)}^{(k,-1)}\right)_{Q_{7}\left/\langle (x_{1,3}-x_{1,6})X_{k,(2,6)}^{(k,-1)}\rangle\right.}$ of (\ref{mf-r3-1-10-8}) is isomorphic to
\begin{eqnarray}
\nonumber
&&
\xymatrix{
R_{25}\ar[rr]_(.4){E_{k+1}(G)}
&&
R_{25}\{1-n\}\ar[rr]_(.6){E_{k+1}(X_{1,(1,4)}^{(1,-1)})}
&&
R_{25}
}\\
\nonumber
&\simeq&\bigoplus_{j=0}^{k}K\left(G;X_{1,(1,4)}^{(1,-1)}\right)_{R_{(1,2,3,4,5)}^{(1,k,1,1,k+1)}}\{2j\}
\end{eqnarray}
Then, the matrix factorization (\ref{mf-r3-1-10-8}) is decomposed into
\begin{equation}
\label{mf-r3-fin-10-8}
\bigoplus_{j=0}^{k}\overline{S}_{(1,2,3;4,5)}^{[1,k,1;1,k+1]}\boxtimes
K\left(G;X_{1,(1,4)}^{(1,-1)}\right)_{R_{(1,2,3,4,5)}^{(1,k,1,1,k+1)}}\{-2k+2j\}.
\end{equation}
The partial matrix factorization $K\left(G_0;(x_{1,6}-x_{1,4})X_{1,(1,4)}^{(1,-1)}\right)_{Q_{7}\left/\langle (x_{1,3}-x_{1,6})X_{k,(2,6)}^{(k,-1)}\rangle\right.}$ of (\ref{mf-r3-1-00-8}) is isomorphic to
\begin{eqnarray}
\nonumber
&&\xymatrix{
R_{24}
\ar[rrrr]_(.45){
\left(
\begin{array}{cc}
{}^{t}\mathfrak{0}_{k}&E_{k}(G)\\
\frac{G}{X_{k+1,(2,3,4)}^{(k,1,-1)}}&\mathfrak{0}_{k}
\end{array}
\right)}
&&&&
R_{25}\{3-n\}
\ar[rrrrrr]_(.55){
\left(
\begin{array}{cc}
\mathfrak{0}_{k}&(x_{1,3}-x_{1,4})X_{k+1,(1,2,4)}^{(1,k,-1)}\\
E_{k}(X_{1,(1,4)}^{(1,-1)})&{}^{t}\mathfrak{0}_{k}
\end{array}
\right)
}
&&&&&&
R_{24}
}\\
\nonumber
&\simeq&
K\left(
\frac{G}{X_{k+1,(2,3,4)}^{(k,1,-1)}};(x_{1,3}-x_{1,4})X_{k+1,(1,2,4)}^{(1,k,-1)}
\right)_{R_{(1,2,3,4,5)}^{(1,k,1,1,k+1)}}
\oplus\bigoplus_{j=1}^{k}
K\left(G;X_{1,(1,4)}^{(1,-1)}\right)_{R_{(1,2,3,4,5)}^{(1,k,1,1,k+1)}}\{2j\}.
\end{eqnarray}
Then, the matrix factorization (\ref{mf-r3-1-00-8}) is isomorphic to
\begin{eqnarray}
\label{mf-r3-fin-00-8}
&&
\hspace{-1cm}
\overline{S}_{(1,2,3;4,5)}^{[1,k,1;1,k+1]}\boxtimes
K\left(
\frac{G}{X_{k+1,(2,3,4)}^{(k,1,-1)}};(x_{1,3}-x_{1,4})X_{k+1,(1,2,4)}^{(1,k,-1)}
\right)_{R_{(1,2,3,4,5)}^{(1,k,1,1,k+1)}}\{-2k-1\}\\
\nonumber
&&
\bigoplus_{j=1}^{k}
\overline{S}_{(1,2,3;4,5)}^{[1,k,1;1,k+1]}\boxtimes
K\left(G;X_{1,(1,4)}^{(1,-1)}\right)_{R_{(1,2,3,4,5)}^{(1,k,1,1,k+1)}}\{-2k-1+2j\}.
\end{eqnarray}
For decompositions (\ref{mf-r3-fin-11-8}), (\ref{mf-r3-fin-01-8}), (\ref{mf-r3-fin-10-8}) and (\ref{mf-r3-fin-00-8}), the morphisms $\overline{\tau}_{-,1}$, $\overline{\tau}_{-,2}$, $\overline{\tau}_{-,3}$ and $\overline{\tau}_{-,4}$ of the complex (\ref{r3-lem-planar-comp-8}) transform into
\begin{eqnarray}
\nonumber
\overline{\tau}_{-,1}&\simeq&
\left(
\begin{array}{cc}
\frak{0}_{k-1}&\id_{\overline{S}_{(1,2,3;4,5)}^{[1,k,1;1,k+1]}}\boxtimes(X_{k,(2,4)}^{(k,-1)},1)\\
E_{k-1}(\id)&{}^{t}\frak{0}_{k-1}
\end{array}
\right)
,\\
\nonumber
\overline{\tau}_{-,2}&\simeq&
\left(
\begin{array}{cc}
\frak{0}_{k-1}&-\id_{\overline{S}_{(1,2,3;4,5)}^{[1,k,1;1,k+1]}}\boxtimes(X_{k,(2,4)}^{(k,-1)},X_{k,(2,4)}^{(k,-1)})\\
E_{k-1}(\id)&{}^{t}\frak{0}_{k-1}\\
\frak{0}_{k-1}&\id
\end{array}
\right),\\
\nonumber
\overline{\tau}_{-,3}&\simeq&
\left(
\begin{array}{cc}
\id_{\overline{S}_{(1,2,3;4,5)}^{[1,k,1;1,k+1]}}\boxtimes(x_{1,3}-x_{1,4},1)&\frak{0}_{k-1}\\
-\id_{\overline{S}_{(1,2,3;4,5)}^{[1,k,1;1,k+1]}}\boxtimes(1,X_{k,(2,4)}^{(k,-1)})&\frak{0}_{k-1}\\
{}^{t}\frak{0}_{k-1}&E_{k-1}(\id)
\end{array}
\right)
,\\
\nonumber
\overline{\tau}_{-,4}&\simeq&
-\left(
\begin{array}{cc}
\frak{0}_{k}&\id_{\overline{S}_{(1,2,3;4,5)}^{[1,k,1;1,k+1]}}\boxtimes(X_{k+1,(2,3,4)}^{(k,1,-1)},1)\\
E_{k}(\id)&{}^{t}\frak{0}_{k}
\end{array}
\right)
.
\end{eqnarray}
Then, the complex (\ref{r3-lem-planar-comp-8}) is isomorphic, in $\k^b(\HMF^{gr}_{R_{(1,2,3,4,5)}^{(1,k,1,1,k+1)},\omega_6})$, to
\begin{equation}
\nonumber
\xymatrix{k\ar@{.}[d]&&&k+1\ar@{.}[d]
\\
\overline{M}_8\{-(k+1)n+1\}
\ar[rrr]^{\id_{\overline{S}}\boxtimes(x_{1,3}-x_{1,4},1)}
&&&
\overline{M}_7\{-(k+1)n\}.
}
\end{equation}
Thus, we obtain
\begin{equation}
\nonumber
\c\left(\input{figure/r3-lem-minus-1k-3}\right)_n\simeq\c\left(\input{figure/r3-lem-minus-1k-4}\right)_n.
\end{equation}
\end{proof}
%
%
%
%
\subsection{Proof of Proposition \ref{prop-twist}}\label{prop-twist-proof}
%
%
\begin{proof}[{\bf Proof of Proposition \ref{prop-twist} (1)}]
The complex of matrix factorization for the diagram $\input{figure/fig-plus-1k-tri1-text}$ is the following object of $\k^b(\HMF^{gr}_{R_{(1,2,3)}^{(k+1,1,k)},\omega_7})$, 
$\omega_7=F_{k+1}(\mathbb{X}^{(k+1)}_{(1)})-F_{1}(\mathbb{X}^{(1)}_{(2)})-F_{k}(\mathbb{X}^{(k)}_{(3)})$,
\begin{equation}
\label{com-twist1}
\xymatrix{
&&-k\ar@{.}[d]&&&-k+1\ar@{.}[d]\\
\c\left(\input{figure/fig-plus-1k-tri1-mf}\right)_n
\hspace{-1cm}&\hspace{-1cm}=\hspace{-1cm}
&\hspace{-.7cm}\c\left(\input{figure/fig-resolution-1k-1-mf}\right)_n\{kn\}\left<k\right>
\ar[rrr]^{\id_{\overline{\Lambda}_{(1;4,5)}^{[k,1]}}\boxtimes\chi_{+,(5,4,3,2)}^{[1,k]}}
&&&\c\left(\input{figure/fig-resolution-1k-2-mf}\right)_n\{kn-1\}\left<k\right>\\
&\hspace{-1cm}=\hspace{-1cm}
&\hspace{-1cm}\overline{\Lambda}_{(1;4,5)}^{[k,1]}\boxtimes\overline{M}_{(5,4,3,2)}^{[1,k]}\{kn\}\left<k\right>
\ar[rrr]^{\id_{\overline{\Lambda}_{(1;4,5)}^{[k,1]}}\boxtimes
\left(\id_{\overline{S}}\boxtimes(1,x_{1,5}-x_{1,2})\right)}
&&&\overline{\Lambda}_{(1;4,5)}^{[k,1]}\boxtimes\overline{N}_{(5,4,3,2)}^{[1,k]}\{kn-1\}\left<k\right>.
}
\end{equation}
We have
\begin{eqnarray}
\nonumber
\overline{\Lambda}_{(1;4,5)}^{[k,1]}\boxtimes\overline{M}_{(5,4,3,2)}^{[1,k]}&=&
K\left(
\left(
\begin{array}{c}
\Lambda_{1,(1;4,5)}^{[k,1]}\\
\vdots\\
\Lambda_{k+1,(1;4,5)}^{[k,1]}
\end{array}
\right)
;
\left(
\begin{array}{c}
x_{1,1}-X_{1,(4,5)}^{(k,1)}\\
\vdots\\
x_{k+1,1}-X_{k+1,(4,5)}^{(k,1)}
\end{array}
\right)
\right)_{R_{(1,4,5)}^{(k+1,k,1)}}\\
\nonumber
&&\boxtimes
K\left(
\left(
\begin{array}{c}
A_{1,(5,4,3,2)}^{[1,k]}\\
\vdots\\
A_{k,(5,4,3,2)}^{[1,k]}\\
v_{k+1,(5,4,3,2)}^{[1,k]}
\end{array}
\right)
;
\left(
\begin{array}{c}
X_{1,(4,5)}^{(k,1)}-X_{1,(2,3)}^{(1,k)}\\
\vdots\\
X_{k,(4,5)}^{(k,1)}-X_{k,(2,3)}^{(1,k)}\\
(x_{1,5}-x_{1,2})X_{k,(2,4)}^{(-1,k)}
\end{array}
\right)
\right)_{R_{(2,3,4,5)}^{(1,k,k,1)}}\{-k\}\\
\nonumber
&\simeq&K\left(
\left(
\begin{array}{c}
\Lambda_{1,(1;2,3)}^{[1,k]}\\
\vdots\\
\Lambda_{k+1,(1;2,3)}^{[1,k]}
\end{array}
\right)
;
\left(
\begin{array}{c}
x_{1,1}-X_{1,(2,3)}^{(1,k)}\\
\vdots\\
x_{k+1,1}-X_{k+1,(2,3)}^{(1,k)}
\end{array}
\right)
\right)_{R_{(1,2,3,5)}^{(k+1,1,k,1)}\left/\left<X_{k+1,(2,3,5)}^{(1,k,-1)}\right>\right.}
\hspace{-1cm}\{-k\}\\
\label{mf-twist1}
&\simeq&\overline{\Lambda}_{(1;2,3)}^{[1,k]}
\boxtimes\left(
\xymatrix{
R_{(1,2,3,5)}^{(k+1,1,k,1)}\left/\left<X_{k+1,(2,3,5)}^{(1,k,-1)}\right>\right.
\ar[r]
&0\
\ar[r]
&R_{(1,2,3,5)}^{(k+1,1,k,1)}\left/\left<X_{k+1,(2,3,5)}^{(1,k,-1)}\right>\right.}
\right)\{-k\},
\end{eqnarray}
\begin{eqnarray}
\nonumber
\overline{\Lambda}_{(1;4,5)}^{[k,1]}\boxtimes\overline{N}_{(5,4,3,2)}^{[1,k]}&=&
K\left(
\left(
\begin{array}{c}
\Lambda_{1,(1;4,5)}^{[k,1]}\\
\vdots\\
\Lambda_{k+1,(1;4,5)}^{[k,1]}
\end{array}
\right)
;
\left(
\begin{array}{c}
x_{1,1}-X_{1,(4,5)}^{(k,1)}\\
\vdots\\
x_{k+1,1}-X_{k+1,(4,5)}^{(k,1)}
\end{array}
\right)
\right)_{R_{(1,4,5)}^{(k+1,k,1)}}\\
\nonumber
&&\boxtimes
K\left(
\left(
\begin{array}{c}
A_{1,(5,4,3,2)}^{[1,k]}\\
\vdots\\
A_{k,(5,4,3,2)}^{[1,k]}\\
v_{k+1,(5,4,3,2)}^{[1,k]}(x_{1,5}-x_{1,2})
\end{array}
\right)
;
\left(
\begin{array}{c}
X_{1,(4,5)}^{(k,1)}-X_{1,(2,3)}^{(1,k)}\\
\vdots\\
X_{k,(4,5)}^{(k,1)}-X_{k,(2,3)}^{(1,k)}\\
X_{k,(2,4)}^{(-1,k)}
\end{array}
\right)
\right)_{R_{(2,3,4,5)}^{(1,k,k,1)}}\{-k+1\}\\
\nonumber
&\simeq&K\left(
\left(
\begin{array}{c}
\Lambda_{1,(1;2,3)}^{[1,k]}\\
\vdots\\
\Lambda_{k+1,(1;2,3)}^{[1,k]}
\end{array}
\right)
;
\left(
\begin{array}{c}
x_{1,1}-X_{1,(2,3)}^{(1,k)}\\
\vdots\\
x_{k+1,1}-X_{k+1,(2,3)}^{(1,k)}
\end{array}
\right)
\right)_{R_{(1,2,3,5)}^{(k+1,1,k,1)}\left/\left<X_{k,(3,5)}^{(k,-1)}\right>\right.}\hspace{-1cm}\{-k+1\}\\
\label{mf-twist2}
&\simeq&
\overline{\Lambda}_{(1;2,3)}^{[1,k]}
\boxtimes\left(
\xymatrix{
R_{(1,2,3,5)}^{(k+1,1,k,1)}\left/\left<X_{k,(3,5)}^{(k,-1)}\right>\right.
\ar[r]
&0\
\ar[r]
&R_{(1,2,3,5)}^{(k+1,1,k,1)}\left/\left<X_{k,(3,5)}^{(k,-1)}\right>\right.}
\right)\{-k+1\}.
\end{eqnarray}
We consider isomorphisms as an $R_{(1,2,3)}^{(k+1,1,k)}$-module
\begin{eqnarray}
\nonumber
R_{(1,2,3,5)}^{(k+1,1,k,1)}\left/\left<X_{k+1,(2,3,5)}^{(1,k,-1)}\right>\right.
&\simeq&R_{(1,2,3)}^{(k+1,1,k)}\oplus x_{1,5}R_{(1,2,3)}^{(k+1,1,k)}\oplus\ldots
\oplus x_{1,5}^{k-1}R_{(1,2,3)}^{(k+1,1,k)}\oplus X_{k,(3,5)}^{(k,-1)}R_{(1,2,3)}^{(k+1,1,k)},\\
\nonumber
R_{(1,2,3,5)}^{(k+1,1,k,1)}\left/\left<X_{k,(3,5)}^{(k,-1)}\right>\right.
&\simeq&R_{(1,2,3)}^{(k+1,1,k)}\oplus x_{1,5}R_{(1,2,3)}^{(k+1,1,k)}\oplus\ldots
\oplus x_{1,5}^{k-1}R_{(1,2,3)}^{(k+1,1,k)}.
\end{eqnarray}
The matrix factorization (\ref{mf-twist1}) is decomposed into
\begin{equation}
\nonumber
\bigoplus_{j=0}^{k}\overline{\Lambda}_{(1;2,3)}^{[1,k]}\{-k+2j\}.
\end{equation}
The matrix factorization (\ref{mf-twist2}) is decomposed into
\begin{equation}
\nonumber
\bigoplus_{j=0}^{k-1}\overline{\Lambda}_{(1;2,3)}^{[1,k]}\{-k+1+2j\}.
\end{equation}
These decompositions change the morphism $\id_{\overline{\Lambda}_{(1;4,5)}^{[k,1]}}\boxtimes\left(\id_{\overline{S}}\boxtimes(1,x_{1,5}-x_{1,2})\right)$ into 
\begin{equation}
\nonumber
\left(
\begin{array}{cc}
E_{k}(\id_{\overline{\Lambda}_{(1;2,3)}^{[1,k]}})&{}^{t}\frak{0}_{k}
\end{array}
\right).
\end{equation}
Thus, the complex (\ref{com-twist1}) is homotopic to
\begin{eqnarray}
\nonumber
&&\xymatrix{
-k\ar@{.}[d]&&-k+1\ar@{.}[d]\\
\c\left(\input{figure/fig-planar-1k-tri1-mf}\right)_n\{kn+k\}\left<k\right>
\ar[rr]&&0
}
\\
\nonumber
&\simeq&
\c\left(\input{figure/fig-planar-1k-tri1-mf}\right)_n\{kn+k\}\left<k\right>\left[-k\right].
\end{eqnarray}
It is obvious that we have the other isomorphism of Proposition \ref{prop-twist} (1) by the symmetry.
\end{proof}
%
%
\begin{proof}[{\bf Proof of Proposition \ref{prop-twist} (2)}]
The complex of matrix factorization for the diagram $\input{figure/fig-minus-1k-tri1-text}$ is the following object of $\k^b(\HMF^{gr}_{R_{(1,2,3)}^{(k+1,1,k)},\omega_7})$
\begin{equation}
\label{com-twist2}
\xymatrix{
&&k-1\ar@{.}[d]&&&k\ar@{.}[d]\\
\c\left(\input{figure/fig-minus-1k-tri1-mf}\right)_n
\hspace{-1cm}&\hspace{-1cm}=\hspace{-1cm}
&\hspace{-.7cm}\c\left(\input{figure/fig-resolution-1k-2-mf}\right)_n\{-kn+1\}\left<k\right>
\ar[rrr]^{\id_{\overline{\Lambda}_{(1;4,5)}^{[k,1]}}\boxtimes\chi_{-,(5,4,3,2)}^{[1,k]}}
&&&\c\left(\input{figure/fig-resolution-1k-1-mf}\right)_n\{-kn\}\left<k\right>\\
&\hspace{-1cm}=\hspace{-1cm}
&\hspace{-1cm}\overline{\Lambda}_{(1;4,5)}^{[k,1]}\boxtimes\overline{N}_{(5,4,3,2)}^{[1,k]}\{-kn+1\}\left<k\right>
\ar[rrr]^{\id_{\overline{\Lambda}_{(1;4,5)}^{[k,1]}}\boxtimes
\left(\id_{\overline{S}}\boxtimes(x_{1,5}-x_{1,2},1)\right)}
&&&\overline{\Lambda}_{(1;4,5)}^{[k,1]}\boxtimes\overline{M}_{(5,4,3,2)}^{[1,k]}\{-kn\}\left<k\right>.
}
\end{equation}
By the discussion of Proof of Proposition \ref{prop-twist} (1), we have
\begin{eqnarray}
\nonumber
&&\overline{\Lambda}_{(1;4,5)}^{[k,1]}\boxtimes\overline{N}_{(5,4,3,2)}^{[1,k]}
\simeq
\overline{\Lambda}_{(1;2,3)}^{[1,k]}
\boxtimes\left(
\xymatrix{
R_{(1,2,3,5)}^{(k+1,1,k,1)}\left/\left<X_{k,(3,5)}^{(k,-1)}\right>\right.
\ar[r]
&0\
\ar[r]
&R_{(1,2,3,5)}^{(k+1,1,k,1)}\left/\left<X_{k,(3,5)}^{(k,-1)}\right>\right.}
\right)\{-k+1\},
\\
\nonumber
&&\overline{\Lambda}_{(1;4,5)}^{[k,1]}\boxtimes\overline{M}_{(5,4,3,2)}^{[1,k]}
\simeq
\overline{\Lambda}_{(1;2,3)}^{[1,k]}
\boxtimes\left(
\xymatrix{
R_{(1,2,3,5)}^{(k+1,1,k,1)}\left/\left<X_{k+1,(2,3,5)}^{(1,k,-1)}\right>\right.
\ar[r]
&0\
\ar[r]
&R_{(1,2,3,5)}^{(k+1,1,k,1)}\left/\left<X_{k+1,(2,3,5)}^{(1,k,-1)}\right>\right.}
\right)\{-k\}.
\end{eqnarray}
We consider isomorphisms as an $R_{(1,2,3)}^{(k+1,1,k)}$-module
\begin{eqnarray}
\nonumber
R_{(1,2,3,5)}^{(k+1,1,k,1)}\left/\left<X_{k,(3,5)}^{(k,-1)}\right>\right.
&\simeq&R_{(1,2,3)}^{(k+1,1,k)}\oplus x_{1,5}R_{(1,2,3)}^{(k+1,1,k)}\oplus\ldots
\oplus x_{1,5}^{k-1}R_{(1,2,3)}^{(k+1,1,k)},\\
\nonumber
R_{(1,2,3,5)}^{(k+1,1,k,1)}\left/\left<X_{k+1,(2,3,5)}^{(1,k,-1)}\right>\right.
&\simeq&R_{(1,2,3)}^{(k+1,1,k)}\oplus (x_{1,5}-x_{1,2})R_{(1,2,3)}^{(k+1,1,k)}\oplus\ldots
\oplus x_{1,5}^{k-1}(x_{1,5}-x_{1,2})R_{(1,2,3)}^{(k+1,1,k)}.
\end{eqnarray}
Then, $\overline{\Lambda}_{(1;4,5)}^{[k,1]}\boxtimes\overline{N}_{(5,4,3,2)}^{[1,k]}$ is decomposed into
\begin{equation}
\nonumber
\bigoplus_{j=0}^{k-1}\overline{\Lambda}_{(1;2,3)}^{[1,k]}\{-k+1+2j\}.
\end{equation}
$\overline{\Lambda}_{(1;4,5)}^{[k,1]}\boxtimes\overline{N}_{(5,4,3,2)}^{[1,k]}$ is decomposed into
\begin{equation}
\nonumber
\bigoplus_{j=0}^{k}\overline{\Lambda}_{(1;2,3)}^{[1,k]}\{-k+2j\}.
\end{equation}
These decompositions change the morphism $\id_{\overline{\Lambda}_{(1;4,5)}^{[k,1]}}\boxtimes\left(\id_{\overline{S}}\boxtimes(x_{1,5}-x_{1,2},1)\right)$ into
\begin{equation}
\nonumber
\left(
\begin{array}{c}
\frak{0}_{k}\\
E_{k}(\id)
\end{array}
\right).
\end{equation}
Thus, the complex (\ref{com-twist2}) is homotopic to
\begin{eqnarray}
\nonumber
&&\xymatrix{
k-1\ar@{.}[d]&&&k\ar@{.}[d]\\
0
\ar[rrr]&&&\c\left(\input{figure/fig-planar-1k-tri1-mf}\right)_n\{-kn-k\}\left<k\right>
}\\
\nonumber
&\simeq&\c\left(\input{figure/fig-planar-1k-tri1-mf}\right)_n\{-kn-k\}\left<k\right>\left[k\right].
\end{eqnarray}
It is obvious that we have the other isomorphism of Proposition \ref{prop-twist} (2) by the symmetry.
\end{proof}
%
%
\begin{proof}[{\bf Proof of Proposition \ref{prop-twist} (3)}]
The complex of matrix factorization for the diagram $\input{figure/fig-plus-1k-tri2-text}$ is the following object of $\k^b(\HMF^{gr}_{R_{(1,2,3)}^{(k+1,1,k)},-\omega_7})$
\begin{equation}
\label{com-twist3}
\xymatrix{
&&-k\ar@{.}[d]&&&-k+1\ar@{.}[d]\\
\c\left(\input{figure/fig-plus-1k-tri2-mf}\right)_n
\hspace{-1cm}&\hspace{-1cm}=\hspace{-1cm}
&\hspace{-.7cm}\c\left(\input{figure/fig-resolution-1k-3-mf}\right)_n\{kn\}\left<k\right>
\ar[rrr]^(.45){\chi_{+,(2,3,4,5)}^{[1,k]}\boxtimes\id_{\overline{V}_{(4,5;1)}^{[k,1]}}}
&&&\c\left(\input{figure/fig-resolution-1k-4-mf}\right)_n\{kn-1\}\left<k\right>\\
&\hspace{-1cm}=\hspace{-1cm}&
\hspace{-1cm}\overline{M}_{(2,3,4,5)}^{[1,k]}\boxtimes\overline{V}_{(4,5;1)}^{[k,1]}\{kn\}\left<k\right>
\ar[rrr]^(.4){\left(\id_{\overline{S}}\boxtimes(1,x_{1,2}-x_{1,5})\right)\boxtimes\id_{\overline{V}_{(4,5;1)}^{[k,1]}}}
&&&\overline{N}_{(2,3,4,5)}^{[1,k]}\boxtimes\overline{V}_{(4,5;1)}^{[k,1]}\{kn-1\}\left<k\right>.
}
\end{equation}
We have
\begin{eqnarray}
\nonumber
\overline{M}_{(2,3,4,5)}^{[1,k]}\boxtimes\overline{V}_{(4,5;1)}^{[k,1]}&=&
K\left(
\left(
\begin{array}{c}
A_{1,(2,3,4,5)}^{[1,k]}\\
\vdots\\
A_{k,(2,3,4,5)}^{[1,k]}\\
v_{k+1,(2,3,4,5)}^{[1,k]}
\end{array}
\right)
;
\left(
\begin{array}{c}
X_{1,(2,3)}^{(1,k)}-X_{1,(4,5)}^{(k,1)}\\
\vdots\\
X_{k,(2,3)}^{(1,k)}-X_{k,(4,5)}^{(k,1)}\\
(x_{1,2}-x_{1,5})X_{k,(3,5)}^{(k,-1)}
\end{array}
\right)
\right)_{R_{(2,3,4,5)}^{(1,k,k,1)}}\{-k\}\\
\nonumber
&&\boxtimes
K\left(
\left(
\begin{array}{c}
V_{1,(4,5;1)}^{[k,1]}\\
\vdots\\
V_{k+1(4,5;1)}^{[k,1]}
\end{array}
\right)
;
\left(
\begin{array}{c}
X_{1,(4,5)}^{(k,1)}-x_{1,1}\\
\vdots\\
X_{k+1,(4,5)}^{(k,1)}-x_{k+1,1}
\end{array}
\right)
\right)_{R_{(1,4,5)}^{(k+1,k,1)}}\{-k\}
\\
\nonumber
&\simeq&K\left(
\left(
\begin{array}{c}
V_{1,(2,3;1)}^{[1,k]}\\
\vdots\\
V_{k+1(2,3;1)}^{[1,k]}
\end{array}
\right)
;
\left(
\begin{array}{c}
X_{1,(2,3)}^{(1,k)}-x_{1,1}\\
\vdots\\
X_{k+1,(2,3)}^{(1,k)}-x_{k+1,1}
\end{array}
\right)
\right)_{R_{(1,2,3,5)}^{(k+1,1,k,1)}\left/\left<X_{k+1,(2,3,5)}^{(1,k,-1)}\right>\right.}
\hspace{-1cm}\{-2k\}\\
\label{mf-twist3}
&\simeq&\overline{V}_{(2,3;1)}^{[1,k]}
\boxtimes\left(
\xymatrix{
R_{(1,2,3,5)}^{(k+1,1,k,1)}\left/\left<X_{k+1,(2,3,5)}^{(1,k,-1)}\right>\right.
\ar[r]
&0\
\ar[r]
&R_{(1,2,3,5)}^{(k+1,1,k,1)}\left/\left<X_{k+1,(2,3,5)}^{(1,k,-1)}\right>\right.}
\right)\{-k\},
\end{eqnarray}
\begin{eqnarray}
\nonumber
\overline{N}_{(2,3,4,5)}^{[1,k]}\boxtimes\overline{V}_{(4,5;1)}^{[k,1]}&=&
K\left(
\left(
\begin{array}{c}
A_{1,(2,3,4,5)}^{[1,k]}\\
\vdots\\
A_{k,(2,3,4,5)}^{[1,k]}\\
v_{k+1,(2,3,4,5)}^{[1,k]}(x_{1,2}-x_{1,5})
\end{array}
\right)
;
\left(
\begin{array}{c}
X_{1,(2,3)}^{(1,k)}-X_{1,(4,5)}^{(k,1)}\\
\vdots\\
X_{k,(2,3)}^{(1,k)}-X_{k,(4,5)}^{(k,1)}\\
X_{k,(3,5)}^{(k,-1)}
\end{array}
\right)
\right)_{R_{(2,3,4,5)}^{(1,k,k,1)}}\{-k+1\}\\
\nonumber
&&
\boxtimes
K\left(
\left(
\begin{array}{c}
V_{1,(4,5;1)}^{[k,1]}\\
\vdots\\
V_{k+1(4,5;1)}^{[k,1]}
\end{array}
\right)
;
\left(
\begin{array}{c}
X_{1,(4,5)}^{(k,1)}-x_{1,1}\\
\vdots\\
X_{k+1,(4,5)}^{(k,1)}-x_{k+1,1}
\end{array}
\right)
\right)_{R_{(1,4,5)}^{(k+1,k,1)}}\{-k\}
\\
\nonumber
&\simeq&K\left(
\left(
\begin{array}{c}
V_{1,(2,3;1)}^{[1,k]}\\
\vdots\\
V_{k+1(2,3;1)}^{[1,k]}
\end{array}
\right)
;
\left(
\begin{array}{c}
X_{1,(2,3)}^{(1,k)}-x_{1,1}\\
\vdots\\
X_{k+1,(2,3)}^{(1,k)}-x_{k+1,1}
\end{array}
\right)
\right)_{R_{(1,2,3,5)}^{(k+1,1,k,1)}\left/\left<X_{k,(3,5)}^{(k,-1)}\right>\right.}
\hspace{-1cm}\{-2k+1\}\\
\label{mf-twist4}
&\simeq&
\overline{V}_{(2,3;1)}^{[1,k]}
\boxtimes\left(
\xymatrix{
R_{(1,2,3,5)}^{(k+1,1,k,1)}\left/\left<X_{k,(3,5)}^{(k,-1)}\right>\right.
\ar[r]
&0\
\ar[r]
&R_{(1,2,3,5)}^{(k+1,1,k,1)}\left/\left<X_{k,(3,5)}^{(k,-1)}\right>\right.}
\right)\{-k+1\}.
\end{eqnarray}
We consider isomorphisms as an $R_{(1,2,3)}^{(k+1,1,k)}$-module
\begin{eqnarray}
\nonumber
R_{(1,2,3,5)}^{(k+1,1,k,1)}\left/\left<X_{k+1,(2,3,5)}^{(1,k,-1)}\right>\right.
&\simeq&R_{(1,2,3)}^{(k+1,1,k)}\oplus x_{1,5}R_{(1,2,3)}^{(k+1,1,k)}\oplus\ldots
\oplus x_{1,5}^{k-1}R_{(1,2,3)}^{(k+1,1,k)}\oplus X_{k,(3,5)}^{(k,-1)}R_{(1,2,3)}^{(k+1,1,k)},\\
\nonumber
R_{(1,2,3,5)}^{(k+1,1,k,1)}\left/\left<X_{k,(3,5)}^{(k,-1)}\right>\right.
&\simeq&R_{(1,2,3)}^{(k+1,1,k)}\oplus x_{1,5}R_{(1,2,3)}^{(k+1,1,k)}\oplus\ldots
\oplus x_{1,5}^{k-1}R_{(1,2,3)}^{(k+1,1,k)}.
\end{eqnarray}
The matrix factorization (\ref{mf-twist3}) is decomposed into
\begin{equation}
\nonumber
\bigoplus_{j=0}^{k}\overline{V}_{(2,3;1)}^{[1,k]}\{-k+2j\}.
\end{equation}
The matrix factorization (\ref{mf-twist4}) is decomposed into
\begin{equation}
\nonumber
\bigoplus_{j=0}^{k-1}\overline{V}_{(2,3;1)}^{[1,k]}\{-k+1+2j\}.
\end{equation}
These decompositions change the morphism $\left(\id_{\overline{S}}\boxtimes(1,x_{1,2}-x_{1,5})\right)\boxtimes\id_{\overline{V}_{(4,5;1)}^{[k,1]}}$ into 
\begin{equation}
\nonumber
\left(
\begin{array}{cc}
E_{k}(\id_{\overline{V}_{(2,3;1)}^{[1,k]}})&{}^{t}\frak{0}_{k}
\end{array}
\right).
\end{equation}
Thus, the complex (\ref{com-twist3}) is homotopic to
\begin{eqnarray}
\nonumber
&&\xymatrix{
-k\ar@{.}[d]&&&-k+1\ar@{.}[d]\\
\c\left(\input{figure/fig-planar-1k-tri2-mf}\right)_n\{kn+k\}\left<k\right>
\ar[rrr]&&&0
}\\
\nonumber
&\simeq&\c\left(\input{figure/fig-planar-1k-tri2-mf}\right)_n\{kn+k\}\left<k\right>\left[-k\right].
\end{eqnarray}
It is obvious that we have the other isomorphism of Proposition \ref{prop-twist} (3) by the symmetry.
\end{proof}
%
%
\begin{proof}[{\bf Proof of Proposition \ref{prop-twist} (4)}]
The complex of matrix factorization for the diagram $\input{figure/fig-minus-1k-tri2-text}$ is the following object of $\k^b(\HMF^{gr}_{R_{(1,2,3)}^{(k+1,1,k)},-\omega_7})$
\begin{equation}
\label{com-twist4}
\xymatrix{
&&k-1\ar@{.}[d]&&&k\ar@{.}[d]\\ 
\c\left(\input{figure/fig-minus-1k-tri2-mf}\right)_n
\hspace{-1cm}&\hspace{-1cm}=\hspace{-1cm}
&\hspace{-.7cm}\c\left(\input{figure/fig-resolution-1k-4-mf}\right)_n\{-kn+1\}\left<k\right>
\ar[rrr]^{\chi_{-,(2,3,4,5)}^{[1,k]}\boxtimes\id_{\overline{V}_{(4,5;1)}^{[k,1]}}}
&&&\c\left(\input{figure/fig-resolution-1k-3-mf}\right)_n\{-kn\}\left<k\right>\\
&\hspace{-1cm}=\hspace{-1cm}
&\hspace{-1cm}\overline{N}_{(2,3,4,5)}^{[1,k]}\boxtimes\overline{V}_{(4,5;1)}^{[k,1]}\{-kn+1\}\left<k\right>
\ar[rrr]^(.4){\left(\id_{\overline{S}}\boxtimes(x_{1,2}-x_{1,5},1)\right)\boxtimes\id_{\overline{V}_{(4,5;1)}^{[k,1]}}}
&&&\overline{M}_{(2,3,4,5)}^{[1,k]}\boxtimes\overline{V}_{(4,5;1)}^{[k,1]}\{-kn\}\left<k\right>.
}
\end{equation}
By the discussion of Proof of Proposition \ref{prop-twist} (3), we have
\begin{eqnarray}
\nonumber
&&\overline{N}_{(2,3,4,5)}^{[1,k]}\boxtimes\overline{V}_{(4,5;1)}^{[k,1]}
\simeq
\overline{V}_{(2,3;1)}^{[1,k]}
\boxtimes\left(
\xymatrix{
R_{(1,2,3,5)}^{(k+1,1,k,1)}\left/\left<X_{k,(3,5)}^{(k,-1)}\right>\right.
\ar[r]
&0\
\ar[r]
&R_{(1,2,3,5)}^{(k+1,1,k,1)}\left/\left<X_{k,(3,5)}^{(k,-1)}\right>\right.}
\right)\{-k+1\},
\\
\nonumber
&&\overline{M}_{(2,3,4,5)}^{[1,k]}\boxtimes\overline{V}_{(4,5;1)}^{[k,1]}
\simeq
\overline{V}_{(2,3;1)}^{[1,k]}
\boxtimes\left(
\xymatrix{
R_{(1,2,3,5)}^{(k+1,1,k,1)}\left/\left<X_{k+1,(2,3,5)}^{(1,k,-1)}\right>\right.
\ar[r]
&0\
\ar[r]
&R_{(1,2,3,5)}^{(k+1,1,k,1)}\left/\left<X_{k+1,(2,3,5)}^{(1,k,-1)}\right>\right.}
\right)\{-k\}.
\end{eqnarray}
We consider isomorphisms as an $R_{(1,2,3)}^{(k+1,1,k)}$-module
\begin{eqnarray}
\nonumber
R_{(1,2,3,5)}^{(k+1,1,k,1)}\left/\left<X_{k,(3,5)}^{(k,-1)}\right>\right.
&\simeq&R_{(1,2,3)}^{(k+1,1,k)}\oplus x_{1,5}R_{(1,2,3)}^{(k+1,1,k)}\oplus\ldots
\oplus x_{1,5}^{k-1}R_{(1,2,3)}^{(k+1,1,k)},\\
\nonumber
R_{(1,2,3,5)}^{(k+1,1,k,1)}\left/\left<X_{k+1,(2,3,5)}^{(1,k,-1)}\right>\right.
&\simeq&R_{(1,2,3)}^{(k+1,1,k)}\oplus (x_{1,2}-x_{1,5})R_{(1,2,3)}^{(k+1,1,k)}\oplus\ldots
\oplus x_{1,5}^{k-1}(x_{1,2}-x_{1,5})R_{(1,2,3)}^{(k+1,1,k)}.
\end{eqnarray}
These isomorphisms change the morphism $\left(\id_{\overline{S}}\boxtimes(x_{1,2}-x_{1,5},1)\right)\boxtimes\id_{\overline{V}_{(4,5;1)}^{[k,1]}}$ into
\begin{equation}
\nonumber
\left(
\begin{array}{c}
\frak{0}_{k}\\
E_{k}(\id)
\end{array}
\right).
\end{equation}
Thus, the complex (\ref{com-twist4}) is homotopic to
\begin{eqnarray}
\nonumber
&&\xymatrix{
k-1\ar@{.}[d]&&&k\ar@{.}[d]\\
0
\ar[rrr]&&&\c\left(\input{figure/fig-planar-1k-tri2-mf}\right)_n\{-kn-k\}\left<k\right>
}\\
\nonumber
&\simeq&\c\left(\input{figure/fig-planar-1k-tri2-mf}\right)_n\{-kn-k\}\left<k\right>\left[k\right].
\end{eqnarray}
It is obvious that we have another isomorphism of Proposition \ref{prop-twist} (4) by the symmetry.
\end{proof}
%
%
%
%
\newpage
\appendix
\section{Virtual link case}\label{app1}
The virtual link theory is given by Kauffman \cite{Kau}.
For a given link diagram, it is represented in Gauss word (see definition \cite{Kau}) and the Gauss word recovers the given link diagram by an operation.
However, there is a Gauss word which has no link diagram obtained by the operation.
Kauffman introduced a virtual crossing, then defined a new topological class called a virtual link diagram.
An object of this class has invariance under isotopy and local moves called the virtual Reidemeister moves.
Roughly speaking, they consist of the Reidemeister moves for crossings and virtual crossings.
Moreover, we can naturally generalize the virtual link diagram into the colored a virtual link diagram.\\
\begin{figure}[htb]
\input{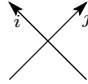}
\caption{$[i,j]$-colored virtual crossing}
\end{figure}\\
\indent
We consider the following virtual $[i,j]$-colored crossing assigned formal indexes.
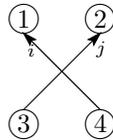
\begin{center}
\input{figure/figcolor-virtual-mf}\\[2em]
\end{center}
We define a matrix factorization for a virtual $[i,j]$-colored crossing to be
\begin{equation}
\nonumber
\c\left(\input{figure/figcolor-virtual-mf}\right)_n:=\overline{L}_{(1,4)}^{[i]}\boxtimes\overline{L}_{(2,3)}^{[j]}.
\end{equation}
\begin{theorem}
By this definition of a matrix factorization for colored virtual crossing, we obtain isomorphisms in $\k^b(\HMF^{gr}_{R,\omega})$ corresponding to the colored virtual Reidemeister moves.
\end{theorem}
\begin{proof}
We naturally obtain this claim by the structure of matrix factorizations.
\end{proof}
Khovanov homology for a virtual link is introduced by V. O. Manturov\cite{Man}.\\
\begin{problem}
How does a relationship between virtual Khovanov homology and virtual Khovanov-Rozansky homology in the case $n=2$ exist?
\end{problem}
%
%
%
%
\newpage
\section{Normalized MOY bracket}\label{NMOY}
Hitoshi Murakami, Tomotada Ohtsuki and Shuji Yamada gave
a polynomial-valued regular link invariant\footnote{The regular link invariant is invariant under the Reidemeister moves II and III. 
It is well-known that we obtain the link invariant from a regular link invariant by taking multiplication of a suitable power of $q$.} with a bracket form associated with $U_q(\mathfrak{sl}_n)$ and $i$ anti-symmetric tensor product of the vector representation , called the MOY bracket. 
It is defined on the set of a colored oriented link diagrams whose component has a coloring from $1$ to $n-1$ \cite{MOY}.
Slightly speaking, this regular link invariant is associated with the quantum group $U_q(\mathfrak{sl}_n)$ 
and its fundamental representations $\land^i V_n$ $(i=1,\ldots,n-1)$, where $V_n$ is the $n$ dimensional vector representation of $U_q(\mathfrak{sl}_n)$.
It is well-known that we obtain a link invariant by normalizing a regular link invariant.\\
\indent The normalized MOY bracket $\left< \cdot \right>_n$ is defined as follows. It locally expands $\pm$-crossings with coloring from $1$ to $n-1$ into a linear combination of planar diagrams with coloring from $1$ to $n$ as follows,
\begin{eqnarray}
\label{figplus-color1}
\left< \input{figure/figplus-color}\right>_n &=& 
\sum_{k=0}^{j} (-1)^{-k+j-i}q^{k+in-i^2+(i-j)^2+2(i-j)}
\left< \input{figure/figmoysmooth2}\right>_n \hspace{1cm} {\rm for}\,\, i \geq j,\\[1em]
\label{figplus-color2}
\left< \input{figure/figplus-color}\right>_n &=& 
\sum_{k=0}^{i} (-1)^{-k+i-j}q^{k+jn-j^2+(j-i)^2+2(j-i)}
\left< \input{figure/figmoysmooth1}\right>_n \hspace{1cm} {\rm for}\,\, i < j,\\[1em]
\label{figminus-color1}
\left< \input{figure/figminus-color}\right>_n &=& 
\sum_{k=0}^{i} (-1)^{k+j-i}q^{-k-jn+j^2-(j-i)^2-2(j-i)}
\left< \input{figure/figmoysmooth1}\right>_n \hspace{1cm} {\rm for}\,\, i\leq j,\\[1em]
\label{figminus-color2}
\left< \input{figure/figminus-color}\right>_n &=& 
\sum_{k=0}^{j} (-1)^{k+i-j}q^{-k-in+i^2-(i-j)^2-2(i-j)}
\left< \input{figure/figmoysmooth2}\right>_n \hspace{1cm} {\rm for}\,\, i > j,
\end{eqnarray}
where the edge colored $0$ vanishes and the bracket $\left<\cdot\right>_n$ assigns $0$ to a diagram having an edge with the coloring number which is greater than $n$.\\
\begin{remark}\label{remark-moy}
\begin{itemize}
\item[\bf(1)]These expansions do not change the outside diagram of the local crossing.
\item[\bf(2)]We consider that the trivial representation of $U_q(\mathfrak{sl}_n)$ is running on the edge with $n$-coloring as the quantum link invariant. 
\item[\bf(3)]
For example, in the case $j>n-i$, the diagram of the $(n-i)$-term in Equation $(\ref{figplus-color1})$ has the edge with coloring $i+j$. Then, this term equals $0$, because $i+j$ is greater than $n$.
\item[\bf(4)]A {\bf colored planar diagram} is built on some trivalent diagrams combinatorially glued by the above expansion.
\begin{figure}[hbt]
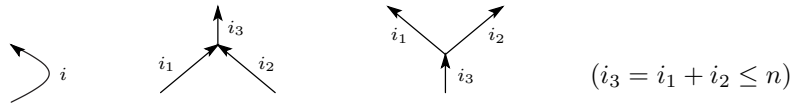

\input{figure/figmoy3-color-i} \hspace{1cm} \input{figure/figgluing-in-3valent-color-i} \hspace{1cm}   \input{figure/figgluing-out-3valent-color-i}\hspace{1cm} {\rm(}$i_3=i_1+i_2\leq n${\rm )}
\caption{$i$-colored line $\&$ ($i_1$,$i_2$,$i_3$)-colored trivalent diagrams}
\end{figure}
\item[\bf(5)]The HOMFLY-PT polynomial is the same with MOY bracket in the case $i=j=1$.
\end{itemize}
\end{remark}
The following relations is the same as ones defined by Murakami, Ohtsuki and Yamada in \cite{MOY}.\\
\indent
For a colored planar diagram $D$ which consists of the disjoint union of colored planar diagrams $D_1$ and $D_2$,
the bracket $\left< D \right>_n$ is defined by the product of $\left< D_1 \right>_n$ and $\left< D_2 \right>_n$,
\begin{eqnarray}
\left< D \right>_n &=& \left< D_1 \right>_n \left< D_2 \right>_n .
\end{eqnarray}
\indent
A closed loop with coloring $i$ ($1 \leq i \leq n$) evaluates $\displaystyle\left[ n \atop i \right]_q$,
\begin{eqnarray}
\left< 
\unitlength 0.1in
\begin{picture}(3.0,1.50)(2.2500,-4.0)
%
\special{pn 8}%
\special{ar 350 350 120 120  0.0000000 6.2831853}%
\put(5.25000,-3.5000){\makebox(0,0){${}_{i}$}}%
\end{picture}%
\right>_n &=& \left[ n \atop i \right]_q \hspace{1cm} (i=1,\ldots , n),
\end{eqnarray}
where $[n]_q=\frac{q^n -q^{-n}}{q-q^{-1}}$, $[n]_q!=[n]_q[n-1]_q\ldots [2]_q[1]_q$ and $\left[n \atop m \right]_q = \frac{[n]_q!}{[m]_q![n-m]_q!}$.\\
\indent
Normalized MOY bracket has the following relations between values of the bracket for some planar diagrams:
\begin{eqnarray}
\label{ass-diag}
\left< \input{figure/fig-ass-dia1}\right>_n &=& \left< \input{figure/fig-ass-dia2}\right>_n  ,\\[-0.1em]
\left< \input{figure/fig-coass-dia1}\right>_n &=& \left< \input{figure/fig-coass-dia2}\right>_n ,
\end{eqnarray}
where $1 \leq i_1,i_2,i_3 \leq n-2$, $i_5 = i_1 + i_2\leq n-1$, $i_6 = i_2 + i_3\leq n-1$ and $i_4 = i_1 + i_2 + i_3 \leq n$. And
\begin{eqnarray}
\label{planar-relation2}
\left< \input{figure/fig-bubble-color}\right>_n &=& 
\left[ i_3 \atop i_1 \right]_q  \left< \input{figure/fig-line-color}\right>_n ,
\\[-0.1em]
\label{planar-relation1}
\left< \input{figure/fig-bubble-color1}\right>_n &=&
\left[ n - i_1 \atop i_2 \right]_q  \left< \input{figure/fig-line-color1}\right>_n ,
\end{eqnarray}
where $1 \leq i_1,i_2 \leq n-1$ and $2 \leq i_3 = i_1 + i_2 \leq n$.\\
\indent
Moreover, we have
\begin{eqnarray*}
\left< \input{figure/figsquare1j--k--k+1j-k--k--1j}\right>_n = 
\left[ i_1-1 \atop i_2 \right]_q  \left< \input{figure/figsquare1i-1}\right>_n +
\left[ i_1-1 \atop i_2-1 \right]_q  \left< \input{figure/figsquare1i-1--i-1+1--1i-1}\right>_n,
\end{eqnarray*}\\[2em]
\begin{equation*}
\left< \input{figure/figsquare1j--j+1--j1--j+1--1j-rev}\right>_n=
\left< \input{figure/figsquare1j-rev}\right>_n +
\left[ n-j-1 \right]_q  \left< \input{figure/figsquare1j--j-1--1j-rev}\right>_n.
\end{equation*}\\[2em]
\indent
Some more relations exist between the values of the bracket for other colored planar diagrams. 
But we leave out the relations because we will not discuss them in following section. See \cite{MOY} about other relations. 
\begin{theorem}\label{MOY}
The bracket $\left< \cdot \right>_n$ is invariant under the Reidemeister moves $I$, $II$ and $III$.
\end{theorem}
\begin{proof}
The proof of invariance under the Reidemeister moves $II$ and $III$ is the same with the proof in \cite{MOY}.\\ 
Therefore, it suffices to show invariance under the Reidemeister move I.
When $\pm$-crossings have the colorings $i$ and $j$ such that $i=j$, $\pm$-curls appear.
We consider $+$-curl. \\
\indent
First, we consider the case $i<n-i$. By the equations (\ref{figplus-color1}), (\ref{planar-relation1}) and (\ref{planar-relation2}), we have
\begin{eqnarray*}
\left< \input{figure/fig-r1-p} \right>_n &=& \sum_{k=0}^{i} (-1)^{-k}q^{k+in-i^2}\left< \input{figure/figmoy-r1-smooth1}\right>_n\\[-0.1em]
&=& \sum_{k=0}^{i} (-1)^{-k}q^{k+in-i^2}
\left[ n-k \atop i \right]_q
\left< \input{figure/figmoy-r1-smooth2}\right>_n\\[-0.1em]
 &=& \sum_{k=0}^{i} (-1)^{-k}q^{k+in-i^2}
\left[ n-k \atop i \right]_q
\left[ i \atop k \right]_q
\left< \input{figure/figmoy-r1-smooth3}\right>_n.
\end{eqnarray*}
We have the following lemma.
\begin{lemma}\label{q-relation}
$$
A_{n,i}:=\sum_{k=0}^{i} (-1)^{-k}q^{k+in-i^2}
\left[ n-k \atop i \right]_q
\left[ i \atop k \right]_q =1
$$
\end{lemma}
\begin{proof}[{\bf Proof of Lemma \ref{q-relation}}]
We prove the lemma by induction to $i$. If $i=1$, then we have
$$
A_{n,1}=q^{1-n} [n]_q - q^{-n} [n-1]_q = q^{1-n} (q^{n-1}+ \ldots +q^{1-n}) - q^{-n} (q^{n-2}+ \ldots +q^{2-n}) = 1.
$$
We show that $A_{n,i}=A_{n-1,i-1}$.
Let $A_{n,i}^{(k)}$ be the $k$-th term of $A_{n,i}$,
$$
A_{n,i}^{(k)}=\left\{ 
\begin{array}{ll}(-1)^{-k}q^{k+in-i^2}
\displaystyle \left[ \displaystyle  n-k \atop \displaystyle i \right]_q
\displaystyle \left[ \displaystyle  i \atop \displaystyle k \right]_q  & if \hspace{0.5cm}0 \leq k \leq i \\
0 & otherwise,
\end{array}
\right.
$$
and we set
$$
T_k = \left\{ 
\begin{array}{ll}
(-1)^{-k} q^{-i^2+in+i}\displaystyle \left[ \displaystyle  n-1-k \atop \displaystyle i \right]_q
\displaystyle \left[ \displaystyle  i-1 \atop \displaystyle k \right]_q & if \hspace{0.5cm}0 \leq k \leq i-1 \\
0 & otherwise.
\end{array}
\right.
$$
Then, we have
\begin{eqnarray*}
A_{n,i}&=&\sum_{k=0}^{ i } A_{n,i}^{(k)}\\[-0.1em]
&=&\sum_{k=0}^{ i } -T_{k-1} + A_{n-1,i-1}^{(k)} +T_{k}\\[-0.1em]
&=&\sum_{k=0}^{ i-1 }A_{n-1,i-1}^{(k)} =A_{n-1,i-1}. 
\end{eqnarray*}
Hence, we obtain $A_{n,i}=A_{n-i+1,1}=1$ by induction of $i$.
\end{proof}
\indent
Next, we consider the case of $i>n-i$. By Remark \ref{remark-moy} (2), we have
\begin{eqnarray*}
\left< \input{figure/fig-r1-p} \right>_n 
&=& \sum_{k=0}^{n-i} (-1)^{-k}q^{k+in-i^2}
\left[ n-k \atop i \right]_q
\left[ i \atop k \right]_q
\left< \input{figure/figmoy-r1-smooth3}\right>_n\\
&=& \sum_{k=0}^{n-i} (-1)^{-k}q^{k+(n-i)n-(n-i)^2}
\left[ n-k \atop n-i \right]_q
\left[ n-i \atop k \right]_q
\left< \input{figure/figmoy-r1-smooth3}\right>_n\\
&=&A_{n,n-i}\left< \input{figure/figmoy-r1-smooth3}\right>_n.
\end{eqnarray*}
By Lemma \ref{q-relation}, we have $A_{n,n-i}=A_{i+1,1}=1$. Hence, we have 
$$
\left< \input{figure/fig-r1-p} \right>_n = \left< \input{figure/figmoy-r1-smooth3}\right>_n.
$$
We similarly have that
$$
\left< \input{figure/fig-r1-m} \right>_n = \left< \input{figure/figmoy-r1-smooth3}\right>_n.
$$
\end{proof}
%
%
%
%

\end{document}